\documentclass{article}

\usepackage[letterpaper, margin=1in]{geometry}
\usepackage{amsmath,amssymb}
\usepackage{graphicx}
\usepackage[amsmath,thmmarks,framed]{ntheorem}
\usepackage{authblk}
\usepackage{booktabs}
\usepackage{multirow}
\usepackage{setspace}
\usepackage{xcolor}
\usepackage{comment}
\definecolor{darkgreen}{rgb}{0.0,0.5,0.0}
\usepackage[colorlinks,linkcolor=red,citecolor=darkgreen,urlcolor=blue]{hyperref}
\hypersetup{
  pdftitle={Stochastic Gradient Leimkuhler--Matthews Integrator with Second Order Discretization Error},
}
\usepackage[numbers,square,comma,sort&compress]{natbib}

\newtheorem{theorem}{Theorem}
\newtheorem{lemma}{Lemma}
\newtheorem{corollary}{Corollary}
\newtheorem{remark}{Remark}
\newtheorem{assumption}{Assumption}
\newenvironment{proof}[1][Proof]{\par\medskip\noindent\textit{#1.}\ }{\hfill$\square$\par\medskip}

\newcommand{\Le}{\leqslant}
\newcommand{\Ge}{\geqslant}
\newcommand{\D}{\mathrm{d}}
\newcommand{\MSE}{\mathrm{MSE}}
\newcommand{\norm}[1]{\lVert #1 \rVert}
\newcommand{\pnorm}[1]{{\vert\kern-0.25ex\vert\kern-0.25ex\vert} #1 {\vert\kern-0.25ex\vert\kern-0.25ex\vert}}
\newcommand{\Var}{\mathrm{Var}}

\title{Stochastic Gradient Leimkuhler--Matthews Integrator \\
	with Second Order Discretization Error}
\author[1]{Xuda Ye}
\author[2]{Yuliang Wang}
\author[1]{Guang Lin}
\author[1]{Cindy Xiangrui Kong}
\affil[1]{Department of Mathematics, Purdue University, West Lafayette, IN 47907}
\affil[2]{Department of Mathematics, Duke University, Durham, NC 27708}
\date{\today}

\numberwithin{equation}{section}

\begin{document}
\maketitle

\begin{abstract}
Stochastic gradient Langevin dynamics (SGLD) has a discretization error of first order in the step size $h$, even with the full gradient. We study the stochastic gradient Leimkuhler--Matthews (SGLM) integrator, which has the same cost as SGLD, one stochastic gradient per step, and reuses each Gaussian vector in consecutive steps. With the full gradient, SGLM is the Leimkuhler--Matthews (LM) integrator, whose discretization error for the invariant distribution is of second order. For a potential strongly convex outside a ball, we prove the bound $\mathcal O(\frac{1}{Nh} + \Sigma^2 h^2 + h^4)$ on the mean square error (MSE) of time averages of SGLM over $N$ steps, where $\Sigma$ bounds the variance of the gradient noise. The proof uses a modified class that absorbs the non-Markovian part of SGLM. The discretization error enters this bound as $h^4$ and the MSE bound of SGLD as $h^2$, while a larger batch reduces only the term $\Sigma^2 h^2$ in both bounds. In Bayesian logistic regression on MNIST, SGLM has a smaller estimated MSE than SGLD. The coefficient of $h$ in the numerical bias of SGLM is given explicitly by the covariance of the gradient noise and the continuous Poisson equation. The coefficient for SGLD contains in addition the Hessian of the potential, which comes from the discretization error. Through Riesz representors, the coefficients of the numerical bias determine the coefficients of the Wasserstein-1 errors of SGLM and SGLD.


\end{abstract}

\section{Introduction}
\label{sec: intro}

Sampling from the Boltzmann distribution $\pi(x) \propto e^{-U(x)}$ on $\mathbb R^d$ is a central task in statistical physics \citep{lelievre2016partial, leimkuhler2015molecular} and in Bayesian inference \citep{welling2011bayesian, dalalyan2017user, durmus2017nonasymptotic}. The simplest stochastic differential equation that leaves $\pi$ invariant is the overdamped Langevin dynamics
\begin{equation}
    \D X_t = -\nabla U(X_t)\,\D t + \sqrt2\,\D B_t , \qquad X_0 = x \in \mathbb R^d ,
    \label{eq: langevin}
\end{equation}
where $(B_t)_{t \Ge 0}$ is a standard Brownian motion in $\mathbb R^d$ and $U$ a confining potential; time averages along the solution $(X_t)_{t \Ge 0}$ approximate expectations under $\pi$ in the long-time limit \citep{roberts1996exponential, lamberton2002recursive}. A numerical integrator of \eqref{eq: langevin} with step size $h > 0$ produces samples $(Z_k)_{k \Ge 0}$ that approximate the invariant distribution $\pi$; such integrators form the basis of Langevin Monte Carlo \citep{dalalyan2017user, durmus2017nonasymptotic, cheng2018sharp, durmus2019high, vempala2019rapid, majka2020nonasymptotic}. Characterizing the accuracy of these samples is a basic question of numerical analysis.

The mean square error (MSE) of time averages is a natural criterion for this accuracy. For a test function $f$, with $\pi(f) := \int_{\mathbb R^d} f\,\D\pi$, the MSE of the time average over $N$ steps is
\begin{equation}
    \MSE(N, h) := \mathbb E \Biggl[ \biggl| \frac1N \sum_{k=0}^{N-1} f(Z_k) - \pi(f) \biggr|^2 \Biggr] .
    \label{eq: mse}
\end{equation}
The Wasserstein distance, the Kullback--Leibler (KL) divergence and the total variation distance compare the law of the samples with $\pi$, and are hard to evaluate in high dimension. The MSE instead concerns a single numerical trajectory, and accounts for both the fluctuation of the time average and the numerical bias, the deviation of its long-time limit from $\pi(f)$ \citep{mattingly2010convergence, chen2015convergence, teh2016consistency, vollmer2016exploration}.

In Bayesian inference on large data sets, $\nabla U$ is replaced by an unbiased stochastic estimator $\nabla \hat U(x; \xi)$ computed from a mini-batch, with $\xi$ a random variable independent of the Brownian motion \citep{welling2011bayesian, chen2015convergence, teh2016consistency, vollmer2016exploration}. The standard sampler of this kind is stochastic gradient Langevin dynamics (SGLD) \citep{welling2011bayesian}, the Euler--Maruyama (EM) integrator \citep{kloeden1992numerical, milstein2004stochastic} of \eqref{eq: langevin} with the stochastic gradient,
\begin{equation}
    Z_{k+1} = Z_k - h\, \nabla \hat U(Z_k; \xi_k) + \sqrt{2h}\, \chi_{k+1} ,
    \label{eq: sgld}
    \tag{SGLD}
\end{equation}
where $(\chi_k)_{k \Ge 0}$ is a sequence of independent standard Gaussian vectors in $\mathbb R^d$ and $\xi_k$ is drawn afresh at every step. The discretization error of \eqref{eq: sgld} is of first order in $h$: its numerical bias is of order $h$ even with the full gradient (Section~\ref{subsec: related sgld}). The existing stochastic gradient integrators that preserve a second order discretization error are built on stochastic gradient Hamiltonian Monte Carlo (SGHMC) \citep{chen2014stochastic} with a symmetric splitting integrator \citep{chen2015convergence} or on the underdamped Langevin dynamics \citep{gouraud2025hmc}, whose splitting integrators are of second order with one gradient evaluation per step \citep{monmarche2021high}; the underdamped dynamics with stochastic gradients is studied in \citep{leimkuhler2016adaptive, matthews2018langevin, paulin2025sampling, lu2025mse}.

We study the \emph{stochastic gradient Leimkuhler--Matthews} (SGLM) integrator, which has the cost of \eqref{eq: sgld}, one stochastic gradient per step, and differs from it in only one respect: the Gaussian vector $\chi_{k+1}$ is replaced by the average of $\chi_k$ and $\chi_{k+1}$, so that each Gaussian vector is used twice, in consecutive steps,
\begin{equation}
    Z_{k+1} = Z_k - h\, \nabla \hat U(Z_k; \xi_k) + \sqrt{\frac h2}\, \bigl( \chi_k + \chi_{k+1} \bigr) .
    \label{eq: sglm}
    \tag{SGLM}
\end{equation}
With the full gradient, SGLM is the Leimkuhler--Matthews (LM) integrator \citep{leimkuhler2013rational, vilmart2015postprocessed},
\begin{equation}
    Z_{k+1} = Z_k - h\, \nabla U(Z_k) + \sqrt{\frac h2}\, \bigl( \chi_k + \chi_{k+1} \bigr) .
    \label{eq: lm}
    \tag{LM}
\end{equation}
The reuse of the Gaussian vector gives LM a second order discretization error for the invariant distribution $\pi$ \citep{leimkuhler2013rational, leimkuhler2014long, vilmart2015postprocessed}, although its samples $(Z_k)_{k \Ge 0}$ do not form a Markov chain and its weak error on finite time intervals is of first order. LM evaluates $\nabla U$ once per step, while almost all other integrators of \eqref{eq: langevin} with a second order discretization error are of predictor--corrector type, with at least two sequential evaluations of $\nabla U$ per step (Section~\ref{subsec: related lm}). With a stochastic gradient, each of these evaluations carries its own noise. This raises the question whether SGLM preserves the second order discretization error of LM.

We provide a positive answer to this question. Our contributions are summarized as follows:

\subsubsection*{1. SGLM separates the stochastic gradient error from the discretization error.} Under regularity conditions on $U$, $f$ and the stochastic gradient that require strong convexity of $U$ only outside a ball, we prove in Theorem~\ref{thm: sglm mse} the MSE bound for SGLM,
\begin{flalign}
	\makebox[0.33\linewidth][l]{(SGLM, variance $\Sigma$)} &
	\MSE(N, h) \Le C \Bigl( \frac{1}{Nh} + \Sigma^2 h^2 + h^4 \Bigr) , &&
	\label{eq: intro sg mse}
\end{flalign}
where $\Sigma$ bounds the variance of the stochastic gradient, and $C$ depends on the regularity constants of $U$ and $f$, the dimension $d$ and the initial value $Z_0$, but not on $N$, $h$ or $\Sigma$. In \eqref{eq: intro sg mse}, $\frac{1}{Nh}$ comes from the fluctuation of the time average, $\Sigma^2 h^2$ from the square of the stochastic gradient error, and $h^4$ from the square of the discretization error. The same separation holds for the UBU integrator with stochastic gradients (SG-UBU) \citep{lu2025mse}, and was previously known only for integrators of the underdamped dynamics \citep{lu2025mse, paulin2026theoretical}. The \eqref{eq: sgld} integrator does not have it, since its stochastic gradient error and its discretization error are both of first order in $h$. With the full gradient, $\Sigma = 0$ and \eqref{eq: intro sg mse} yields the MSE bound $\frac{1}{Nh} + h^4$ of LM. To our knowledge, this is also the first MSE bound for the LM integrator.

The separation in \eqref{eq: intro sg mse} is advantageous in the mini-batch regime: if each step averages $P$ independent stochastic gradients, the variance bound $\Sigma$ is replaced by $\Sigma/P$, so the contribution of the gradient noise to the MSE becomes $\Sigma^2 h^2 / P^2$, with $C$ in \eqref{eq: intro sg mse} independent of $P$. Increasing the batch size $P$ at a fixed step size $h$ therefore reduces the stochastic gradient error. By contrast, the MSE of \eqref{eq: sgld} keeps the term $h^2$ from the square of its discretization error at every batch size, and only a smaller $h$ reduces it. Since the gradients of a batch are evaluated in parallel while the steps are sequential, SGLM with a large batch size $P$ has a smaller error bound than SGLD under the same cost.

The proof of \eqref{eq: intro sg mse} runs on the Markovian form of \eqref{eq: sglm}, a Markov chain $(Y_k)_{k \Ge 0}$ with $Z_k = Y_k + \sqrt{h/2}\,\chi_k$. It builds on the discrete Poisson equation of our earlier work \citep{lu2025mse, ye2026mse}, combined with a modified class of test function, potential, generator and discrete Poisson solution that absorbs the non-Markovian part of the integrator, the Gaussian shift $\sqrt{h/2}\,\chi_k$.

\subsubsection*{2. The coefficients of the numerical bias are explicit.} The time averages of $f$ along \eqref{eq: sglm} converge almost surely to $\pi_h(f)$, for a numerical invariant distribution $\pi_h$ independent of the initial value, and the numerical bias $\pi_h(f) - \pi(f)$ is of order $h$ with an explicit coefficient given in Corollary~\ref{cor: bias coefficient},
\begin{flalign}
	\makebox[0.33\linewidth][l]{(SGLM, numerical bias)} &
	\lim_{h \to 0} \frac{\pi_h(f) - \pi(f)}{h} = \frac12\, \pi\bigl( \nabla^2 \phi_0(x) : \mathcal C(x) \bigr) , &&
	\label{eq: intro bias coefficient}
\end{flalign}
where $\mathcal C(x)$ is the covariance of the gradient noise $\nabla \hat U(x;\xi) - \nabla U(x)$, $\phi_0$ solves the continuous Poisson equation \eqref{eq: continuous poisson}, $-\mathcal L \phi_0 = f - \pi(f)$ with $\mathcal L$ the generator of \eqref{eq: langevin}, and $:$ is the full contraction of two matrices. The coefficient comes entirely from the gradient noise: with the full gradient, $\mathcal C(x) = 0$ and the numerical bias of LM has no linear part. The gradient noise enters only through the local weak error of a single step, through its covariance $\mathcal C(x)$, as for SG-UBU in \citep[Section~3.4]{lu2025mse}.

The discrete Poisson equation also characterizes the numerical bias of \eqref{eq: sgld}. Its time averages converge to $\bar\pi_h(f)$, with $\bar\pi_h$ the invariant distribution of SGLD, and by Corollary~\ref{cor: sgld coefficient}
\begin{flalign}
	\makebox[0.28\linewidth][l]{(SGLD, numerical bias)} &
	\lim_{h \to 0} \frac{\bar\pi_h(f) - \pi(f)}{h} = \frac12\, \pi\Bigl( \nabla^2 \phi_0(x) : \bigl( \nabla^2 U(x) + \mathcal C(x) \bigr) \Bigr) , &&
	\label{eq: intro sgld coefficient}
\end{flalign}
a coefficient studied formally for $d = 1$ in \citep{vollmer2016exploration}. Besides the covariance $\mathcal C(x)$, the coefficient of SGLD carries the Hessian $\nabla^2 U(x)$, which comes from the discretization error of EM. With the full gradient, $\mathcal C(x) = 0$ and \eqref{eq: intro sgld coefficient} gives the coefficient of EM; with a batch of $P$ independent stochastic gradients, the coefficient of SGLM is divided by $P$, while that of SGLD tends to the coefficient of EM as $P$ grows. In \eqref{eq: intro bias coefficient} and \eqref{eq: intro sgld coefficient}, the gradient noise and the discretization error both enter linearly through $\nabla^2 \phi_0$, so the solution $\phi_0$ of the continuous Poisson equation is an essential structure of the numerical bias of the two stochastic gradient integrators, not only a tool for upper bounds on the MSE.

\subsubsection*{3. The coefficients of the Wasserstein-1 error are explicit.} The suprema of the absolute values of the coefficients \eqref{eq: intro bias coefficient} and \eqref{eq: intro sgld coefficient} over the test functions with $\norm{\nabla f}_\infty \Le 1$ are the limits of the normalized Wasserstein-1 errors $\mathcal W_1(\pi_h, \pi)/h$ and $\mathcal W_1(\bar\pi_h, \pi)/h$ as $h \to 0$ (Corollary~\ref{cor: w1 limits}):
\begin{flalign}
	\makebox[0.35\linewidth][l]{(SGLM, $\mathcal W_1$ coefficient)} &
	\lim_{h \to 0} \frac{\mathcal W_1(\pi_h, \pi)}{h} = \frac12\, \norm{\psi_{\mathcal C}\, \pi}_{\mathrm{KR}} , &&
	\label{eq: intro w1 sglm} \\
	\makebox[0.35\linewidth][l]{(SGLD, $\mathcal W_1$ coefficient)} &
	\lim_{h \to 0} \frac{\mathcal W_1(\bar\pi_h, \pi)}{h} = \frac12\, \norm{( \psi_{\nabla^2 U} + \psi_{\mathcal C} )\, \pi}_{\mathrm{KR}} , &&
	\label{eq: intro w1 sgld}
\end{flalign}
where $\psi_{\mathcal C}$ and $\psi_{\nabla^2 U}$ are the Riesz representors of $\mathcal C(x)$ and $\nabla^2 U(x)$ defined in \eqref{eq: riesz}, and $\norm{\cdot}_{\mathrm{KR}}$ is the Kantorovich--Rubinstein norm in Section~\ref{subsec: notation}. The normalized biases $(\pi_h - \pi)/h$ and $(\bar\pi_h - \pi)/h$ converge to $\frac12 \psi_{\mathcal C}\, \pi$ and $\frac12 (\psi_{\nabla^2 U} + \psi_{\mathcal C})\, \pi$ in this norm at the rate $h^{\frac12} (1 + \log\frac1h)$ (Theorem~\ref{thm: kr expansion}). SGLM does not always have a smaller coefficient of the Wasserstein-1 error than SGLD: for a gradient noise of large amplitude supported away from the minimizer of $U$, the order is reversed (Lemma~\ref{lem: w1 counterexample}). With a batch of $P$ independent stochastic gradients, $\psi_{\mathcal C}/P$ replaces $\psi_{\mathcal C}$ in both limits, with $\mathcal C(x)$ the covariance of a single stochastic gradient. For a strongly convex potential, the coefficient of SGLM is at most that of SGLD once the batch size reaches a threshold linear in $\sup_x \norm{\mathcal C(x)}$ (Theorem~\ref{thm: batch threshold}).

\subsubsection*{Summary of the results}

Three features distinguish these results. First, SGLM separates the stochastic gradient error from the discretization error at the cost of SGLD, so that the batch size and the step size control the two errors independently; for an integrator of the overdamped dynamics, this separation is new.

Second, the results on the numerical bias are equalities, not only upper bounds. The coefficients \eqref{eq: intro bias coefficient} and \eqref{eq: intro sgld coefficient} identify the numerical biases of SGLM and SGLD, and \eqref{eq: intro w1 sglm} and \eqref{eq: intro w1 sgld} their Wasserstein-1 errors, to first order in $h$, uniformly over the Lipschitz test functions for the Wasserstein-1 errors. The mainstream coupling methods bound these errors from above without identifying their coefficients; here the coefficients are the Kantorovich--Rubinstein norms of the signed measures given by the Riesz representors, and the shift of the invariant distribution caused by the gradient noise and by the discretization is characterized exactly.

Third, the exact coefficients make the comparison of SGLM with SGLD precise. The coefficient of SGLD keeps the Hessian of the potential at every batch size, while that of SGLM is divided by the batch size. SGLM can be less accurate than SGLD at small batch sizes (Lemma~\ref{lem: w1 counterexample}), and, for a strongly convex potential, Theorem~\ref{thm: batch threshold} gives an explicit threshold of the batch size beyond which SGLM is more accurate than SGLD in the Wasserstein-1 distance, to first order in $h$. This threshold does not depend explicitly on the dimension $d$; it depends only on the eigenvalues of $\nabla^2 U(x)$ and $\mathcal C(x)$.

As applications of the analysis of the coefficients, Section~\ref{sec: representors} describes three approaches to bias correction. The numerical bias of SGLM is $\frac h2\, \pi(f \psi_{\mathcal C})$ to first order in $h$, so it is removed once $\pi(f \psi_{\mathcal C})$ is computed, with the explicit representor $\psi_{\mathcal C}$ for the covariances of Table~\ref{tab: explicit representors} or by Richardson--Romberg extrapolation in the step size. The third approach, the shifting correction, replaces a sample $Z \sim \pi_h$ of SGLM by $Z + h \nabla \varphi_{\mathcal C}(Z)$, where $\varphi_{\mathcal C}$ solves \eqref{eq: shift potential}, and formally removes the first order bias coefficient for all test functions at once. The function $\varphi_{\mathcal C}$ is approximated by minimizing the energy \eqref{eq: extrapolated energy}, which requires only the samples of SGLM at the step sizes $h$ and $2h$, so the shifting correction is completely data-driven.

In the numerical experiments of Section~\ref{sec: numerics}, the evaluated coefficients of the numerical bias of SGLM and SGLD in a 4D Gaussian mixture distribution approach the exact coefficients \eqref{eq: intro bias coefficient} and \eqref{eq: intro sgld coefficient} as $h$ decreases, and the evaluated coefficients of the Wasserstein-1 errors in a radial distribution approach those of \eqref{eq: intro w1 sglm} and \eqref{eq: intro w1 sgld}. In a 16D product of ring distributions, the approaches to bias correction of Section~\ref{sec: representors} reduce the numerical bias of SGLM. In Bayesian logistic regression on MNIST, SGLM has a smaller estimated MSE than SGLD with the same step size and batch size.

Section~\ref{sec: related} reviews the related works. Section~\ref{sec: preliminaries} states the assumptions, the Markovian forms \eqref{eq: lm markov} and \eqref{eq: sglm markov}, and the estimates of the discrete Poisson equation. Section~\ref{sec: modified} constructs the modified class, Section~\ref{sec: local} estimates the local errors of the LM and SGLM integrators against the modified dynamics, and Section~\ref{sec: mse} states the main results. Section~\ref{sec: representors} describes the approaches to bias correction, and Section~\ref{sec: numerics} presents the numerical experiments. Section~\ref{sec: conclusion} concludes, and the proofs are in the appendices.

\section{Related works}
\label{sec: related}

This section reviews the work related to \eqref{eq: lm} and \eqref{eq: sglm}. Section~\ref{subsec: related lm} reviews the integrators of \eqref{eq: langevin} with a second order discretization error and the works on the LM integrator. Section~\ref{subsec: related sgld} reviews the sampling error of \eqref{eq: sgld}, and Section~\ref{subsec: related underdamped} the stochastic gradient integrators of the underdamped Langevin dynamics.

\subsection{Second order integrators in the overdamped case}
\label{subsec: related lm}

Among the integrators of \eqref{eq: langevin} that evaluate $\nabla U$ alone, only a few have a discretization error of order $h^2$ for the invariant distribution. They include LM \citep{leimkuhler2013rational, vilmart2015postprocessed} and stochastic Heun \citep{milstein2004stochastic, debrabant2010runge, mannella2025numerical}. The large family of stochastic Runge--Kutta (SRK) integrators \citep{rossler2010runge, debrabant2010runge, li2019stochastic, yang2026accelerating} contains further examples, such as SRK-I/II \citep{yang2026accelerating, ye2026mse}. Beyond LM, almost all these integrators are of predictor--corrector type: the second evaluation of $\nabla U$ in a step is taken at a state computed from the first. With a stochastic gradient, each of these evaluations carries its own noise, while LM evaluates $\nabla U$ once per step and extends to stochastic gradients by replacing that single evaluation; this is the reason we study LM and its stochastic gradient form SGLM. A detailed comparison of these integrators is given in \citep{ye2026mse}. Integrators of high order for the invariant distribution that also use derivatives of $\nabla U$ have been constructed \citep{talay1990second, abdulle2014high}.

The LM integrator was introduced in \citep{leimkuhler2013rational} as the high friction limit of the BAOAB splitting of the underdamped Langevin dynamics. Its second order discretization error for the invariant distribution was proved by Talay--Tubaro expansions \citep{talay1990expansion}. These expansions were derived for the weak error, whose first order coefficient decays exponentially in time \citep{leimkuhler2014long}, and for the invariant distribution of splitting integrators on the torus \citep{leimkuhler2015computation}. The second order discretization error of LM was also proved by writing the integrator as the Markovian form \eqref{eq: lm markov} followed by the Gaussian shift $\sqrt{h/2}\,\chi_k$ as a postprocessor \citep{vilmart2015postprocessed}. The postprocessor argument uses the sufficient conditions for high order on the invariant distribution of \citep{abdulle2014high}, which rest on the Talay--Tubaro expansion \citep{talay1990expansion} and on the weak backward error analysis of \citep{debussche2012weak, kopec2013weak}. In high dimension and under weak or sparse interactions among the variables, the bias of the unadjusted underdamped Langevin algorithm on marginals of a few coordinates depends on the dimension only through logarithmic factors \citep{chen2026delocalization}, extending the delocalization of bias of the unadjusted overdamped Langevin algorithm \citep{chen2026convergence}. This delocalization holds for all large friction parameters and therefore carries over to \eqref{eq: lm}.
These works control the numerical bias of LM; none of them bounds the MSE \eqref{eq: mse} or treats stochastic gradients.

The modified class of this paper serves a purpose different from that of the modified equations of \citep{abdulle2014high}, which are used to construct integrators. Here the integrator is fixed, and the modified potential of Section~\ref{sec: modified} defines the modified dynamics \eqref{eq: modified langevin}, against which the Markovian form \eqref{eq: lm markov} has a local weak error of order $h^3$. The modified class also provides the modified discrete Poisson solution, to which the MSE analysis of \citep{lu2025mse, ye2026mse} is applied. The order $h^3$ of this local weak error appears in the Talay--Tubaro expansions \citep{leimkuhler2014long, leimkuhler2015computation}, as the vanishing at equilibrium of its coefficient of order $h$, and in the postprocessor argument \citep{vilmart2015postprocessed}. The expansions in \citep{leimkuhler2014long} rely on the exponential decay in time of the derivatives of the Kolmogorov solution of Section~\ref{subsec: poisson}, for which they refer to \citep{talay1990second}. Theorem~\ref{thm: u decay} proves this decay on $\mathbb R^d$ for a potential satisfying Assumption~\ref{ass: U}, extending \citep[Theorem~1]{ye2026mse} to derivatives of order $4$.

\subsection{Sampling error analysis for SGLD}
\label{subsec: related sgld}

The \eqref{eq: sgld} integrator converges weakly to the solution of \eqref{eq: langevin} as its step size tends to zero \citep{sato2014approximation}. Its sampling error has been analyzed along the following lines, and Table~\ref{tab: sgld} collects representative bounds; see \citep{nemeth2021stochastic} for a review.
\begin{itemize}
	\setlength{\itemsep}{0pt}
	\item Recurrence relations of the integrator, combined with strong or weak error estimates, give nonasymptotic bounds in Wasserstein distance \citep{dalalyan2017user, nagapetyan2017true, brosse2018promises, barkhagen2021stochastic, chau2021stochastic, zhang2023nonasymptotic, shaw2025random}. Weighted transportation cost inequalities give bounds in $\mathcal W_2$ for nonconvex potentials \citep{raginsky2017nonconvex}. Under log-concavity at infinity the variance of the stochastic gradient does not change the rate of weak convergence \citep{majka2020nonasymptotic}, and a conductance analysis gives bounds in total variation for targets that are not log-concave \citep{zou2021faster}. For targets log-concave outside a compact set, the law of the $k$-th iterate of SGLD converges geometrically in $\mathcal W_1$ to the numerical invariant distribution $\bar\pi_h$ of SGLD \citep{li2024geometric}.
	\item The continuous Poisson equation \eqref{eq: continuous poisson} controls the error of the time averages of numerical schemes \citep{mattingly2010convergence, brosse2019tamed} and of SGLD and SGHMC, including their MSE \citep{chen2015convergence, teh2016consistency, vollmer2016exploration, xu2018global, li2023improving}. These results require its solution $\phi_0$ to have bounded or polynomially growing derivatives.
	\item The KL divergence between the laws of SGLD and of \eqref{eq: langevin} is bounded uniformly in time \citep{li2022sharp}. Bounds in KL divergence also follow from convex optimization on the Wasserstein space \citep{durmus2019analysis} and, under a logarithmic Sobolev inequality, from the nearly Gaussian structure of the gradient noise \citep{das2023utilising}.
	\item Variance-reduced and control-variate estimators of the gradient improve these guarantees \citep{dubey2016variance, chatterji2018theory, baker2019control, zou2019sampling, kinoshita2022improved}. The numerical bias of SGLD is reduced by Richardson--Romberg extrapolation \citep{durmus2016stochastic}, and its cost at a prescribed accuracy is reduced by multilevel Monte Carlo \citep{giles2016multilevel, giles2020multilevel}. The extrapolation, however, introduces negative weights and does not increase the accuracy of the trajectory itself. The cost of the gradient noise is quantified by oracle lower bounds \citep{chatterji2022oracle}, by the dependence of the MSE on the mini-batch size under a fixed budget \citep{chen2019convergence}, and by a comparison of SGLD with EM at a prescribed mean square accuracy \citep{majka2026computational}.
\end{itemize}

\begin{table}[htb]
	\centering
	\setlength{\tabcolsep}{4pt}
	\begingroup
	\setstretch{1.2}
	\begin{tabular}{@{}p{48mm}p{20mm}p{80mm}@{}}
		\toprule
		reference & assumption & sampling error \\
		\midrule
		Theorem~5 \citep{dalalyan2017user} & global & $\mathcal W_2(\mu_k, \pi) \lesssim e^{-\lambda kh} + h + \Sigma^{\frac12} h^{\frac12}$ \\
		Theorem~5, Corollary~12 \citep{nagapetyan2017true} & global & $\mathbb E |X_{kh} - Z_k|^2 \lesssim e^{-\lambda kh} + h$ \\
		Corollary~18 \citep{durmus2019analysis} & global & $\mathrm{KL}(\bar\mu_N \| \pi) \lesssim (Nh)^{-1} + h$ \\
		Theorem~3.8 \citep{shaw2025random} & global & $\mathcal W_2(\mu_k, \pi) \lesssim e^{-\lambda kh} + h^{\frac12}$ \\
		Proposition~1, Theorem~1 \citep{paulin2026sgld} & global & $\mathcal W_2(\mu_k, \pi) \lesssim e^{-\lambda kh} + h$ \\
		Proposition~10 \citep{raginsky2017nonconvex} & non-global & $\mathcal W_2(\mu_k, \pi) \lesssim (\delta^{\frac14} + h^{\frac14})\, kh + e^{-\lambda kh}$ \\
		Theorem~4.5 \citep{zou2021faster} & non-global & $\norm{\mu_k - \pi}_{\mathrm{TV}} \lesssim e^{-\lambda kh} + h^{\frac12} + \epsilon$ \\
		Theorem~9 \citep{vollmer2016exploration} & non-global & $\MSE(N,h) \lesssim (Nh)^{-1} + h^2$ \\
		Theorem~3 \citep{chen2015convergence} & non-global & $\MSE(N,h) \lesssim (Nh)^{-1} + h^2$ \\
		Corollaries~3.1 and~5.1 \citep{li2022sharp} & log-Sobolev & $\mathrm{KL}(\mu_k \| \nu_k) \lesssim h^2$, \ $\mathcal W_2(\mu_k, \pi) \lesssim e^{-\lambda kh} + h$ \\
		Theorem~5 \citep{das2023utilising} & log-Sobolev & $\mathrm{KL}(\mu_k \| \pi) \lesssim e^{-\lambda kh} + h$ \\
		\midrule
		SGLM, Theorem~\ref{thm: sglm mse} & non-global & $\MSE(N,h) \lesssim (Nh)^{-1} + \Sigma^2 h^2 + h^4$ \\
		SGLM, Corollary~\ref{cor: w1 limits} & non-global & $\lim_{h \to 0} \frac{\mathcal W_1(\pi_h, \pi)}{h} = \frac12 \norm{\psi_{\mathcal C}\, \pi}_{\mathrm{KR}}$ \\
		SGLD, Corollary~\ref{cor: w1 limits} & non-global & $\lim_{h \to 0} \frac{\mathcal W_1(\bar\pi_h, \pi)}{h} = \frac12 \norm{( \psi_{\nabla^2 U} + \psi_{\mathcal C} )\, \pi}_{\mathrm{KR}}$ \\
		\bottomrule
	\end{tabular}
	\endgroup
	\caption{Sampling error results for \eqref{eq: sgld} and for \eqref{eq: sglm}. In the assumption column, global marks a strongly convex or convex potential, non-global a dissipative potential, one convex outside a ball, or, in \citep{vollmer2016exploration, chen2015convergence}, bounds on the solution of the Poisson equation with moment conditions in place of convexity, and log-Sobolev a target satisfying a logarithmic Sobolev inequality; the references impose in addition conditions on the regularity, the initial distribution and the moments of the gradient noise. Here $\mu_k$ denotes the law of the $k$-th iterate of SGLD, $\bar\pi_h$ and $\pi_h$ the numerical invariant distributions of SGLD and of SGLM, $\bar\mu_N$ the average of $\mu_1, \dots, \mu_N$, and $\nu_k$ the law of the solution of \eqref{eq: langevin} at time $kh$. Further, $X_{kh}$ is the solution of \eqref{eq: langevin} coupled with the $k$-th iterate $Z_k$ of SGLD, $\delta$ the relative variance of the gradient noise in \citep{raginsky2017nonconvex}, and $\epsilon \in (0,1)$ the truncation parameter in \citep{zou2021faster}, which fixes the radius of the truncated target on which the constants of that bound depend. Finally, $\Sigma$ is the bound on the variance of the gradient noise in Assumption~\ref{ass: sg}, $\psi_{\mathcal C}$ and $\psi_{\nabla^2 U}$ the Riesz representors \eqref{eq: riesz} of its covariance $\mathcal C(x)$ in \eqref{eq: covariance} and of $\nabla^2 U(x)$, $\norm{\cdot}_{\mathrm{KR}}$ the Kantorovich--Rubinstein norm, $\mathcal W_p$ the Wasserstein distance of order $p$, and $\lambda > 0$ a rate of convergence. Other constants, depending on the dimension, the potential and the gradient noise, are suppressed.}
	\label{tab: sgld}
\end{table}

The numerical bias of \eqref{eq: sgld} is of first order in $h$. The bias of the time average over $N$ steps is $\mathcal O\bigl( (Nh)^{-1} + h \bigr)$ \citep[Theorem~9]{vollmer2016exploration}. The uniform-in-time KL bound of \citep[Corollary~3.1]{li2022sharp}, of order $h^2$, combined with the geometric ergodicity of SGLD for a potential strongly convex outside a ball \citep{li2024geometric}, yields $\mathcal W_1(\bar\pi_h, \pi) \Le Ch$ \citep[Section~5.1]{li2022sharp}. Corollary~\ref{cor: w1 limits} gives the exact coefficients of $h$ in $\mathcal W_1(\bar\pi_h, \pi)$ and $\mathcal W_1(\pi_h, \pi)$, determined through Riesz representors by the bias coefficients of Corollaries~\ref{cor: bias coefficient} and \ref{cor: sgld coefficient}. For a strongly convex potential with a Lipschitz gradient and gradient noise that is Lipschitz in mean square, $\mathcal W_2(\bar\pi_h, \pi) \Le Ch$ when the gradient noise has a finite moment of order $4$, while under a finite moment of order $2$ only, the bound is of order $h^{\frac12}$ \citep[Theorem~1]{paulin2026sgld}. A family of spiked gradient noises shows that the order $h^{\frac12}$ cannot be improved uniformly over gradient noises of a fixed variance \citep[Lemma~1]{paulin2026sgld}. Each noise of this family has finite moments of all orders, while their moments of order $12$ are unbounded over the family; Assumption~\ref{ass: sg}, whose bound on the moment of order $12$ of the stochastic gradient has a fixed constant $M$, excludes the family as a whole. The first order term has two sources: the gradient noise and the discretization error of EM. The discretization error of EM gives a numerical bias of first order even with the full gradient \citep{talay1990expansion}; see also \citep[Section~4.1]{vollmer2016exploration}.

A different line of work removes the numerical bias. The proximal sampler \citep{lee2021structured, chen2022improved}, with its restricted Gaussian oracle implemented by first order rejection sampling \citep{chen2026high}, reaches a prescribed accuracy in total variation with a number of gradient queries polylogarithmic in the inverse accuracy \citep{altschuler2024faster, chen2026high}, and stochastic gradients with subexponential tails preserve this rate for targets satisfying a log-Sobolev or a Poincar\'e inequality \citep{chen2026stochastic}. These samplers need a warm start, an approximate proximal step and a rejection step with a random number of gradient queries in each iteration, and light-tailed gradient noise: with stochastic gradients of bounded variance only, a number of queries proportional to the inverse accuracy is necessary and sufficient \citep{chen2026stochastic}. SGLD and SGLM use one stochastic gradient per step and no warm start, and their numerical bias, quantified in the present paper, can be made smaller than the fluctuation of the time average by the choice of the step size and the batch size.

\subsection{Stochastic gradient integrators for underdamped Langevin dynamics}
\label{subsec: related underdamped}

The underdamped Langevin dynamics augments the position with a velocity and is discretized by splitting integrators such as BAOAB \citep{leimkuhler2013rational, leimkuhler2015molecular}, OBABO \citep{bussi2007accurate} and UBU \citep{sanz2021wasserstein, paulin2024correction}. For convex and smooth potentials, the unadjusted chain of a splitting integrator of second order with one gradient evaluation per step contracts in Wasserstein distance at a rate independent of the dimension \citep{monmarche2021high}.

SGHMC adds a friction to compensate for the gradient noise \citep{chen2014stochastic}, within a general recipe for such samplers \citep{ma2015complete}. Adaptive thermostats correct for the gradient noise \citep{ding2014bayesian, leimkuhler2016adaptive, shang2015covariance}, with exponential convergence of the law \citep{leimkuhler2020hypocoercivity} and a quantified bias from mini-batching \citep{sekkat2023minibatching}. Variance reduction accelerates SGHMC \citep{li2019variance, zou2018stochastic, zou2019stochastic, chatterji2018theory, zou2021convergence}, and multilevel coupling gives unbiased estimators \citep{chada2026unbiased}. For nonconvex optimization, the complexity bounds of SGHMC can be tighter than those of \eqref{eq: sgld} \citep{gao2022global}. Wasserstein bounds are proved without log-concavity \citep{akyildiz2024nonasymptotic, chak2023reflection}, through contraction rates of stochastic gradient integrators \citep{leimkuhler2023contraction}, and with dimension-free rates for log-concave targets \citep{gouraud2025hmc}.

Consider a stochastic gradient sampler such as SGHMC, with unbiased gradients and a Markov integrator whose local weak error is of order $h^{K+1}$. If the derivatives of the solution of the continuous Poisson equation \eqref{eq: continuous poisson} are bounded, the bias of its time average is of order $(Nh)^{-1} + h^K$ and its MSE of order $(Nh)^{-1} + h^{2K}$ \citep[Theorems~2 and~3]{chen2015convergence}. Applied to the Markovian form \eqref{eq: lm markov}, which has weak order $1$ \citep{ye2026mse}, this framework gives a discretization term of order $h^2$ in the MSE, whereas the modified class of Section~\ref{sec: modified} gives $h^4$. The derivative bounds that this framework requires follow from the decay \eqref{eq: u bound} of Theorem~\ref{thm: u decay}.

For SG-UBU and a globally convex potential, the discrete Poisson equation bounds the MSE of the time average \citep{lu2025mse}. The numerical bias is of first order in $h$, with a coefficient proportional to the variance of the stochastic gradient. For a finite-sum potential whose components are each strongly convex, UBU combined with a symmetric forward and backward sweep over the mini-batches reaches a bias of order $h^2$ with a single mini-batch per step \citep{paulin2025sampling}. The constant of this bound grows with the number of mini-batches, and the order rests on the symmetry of the whole sweep, which lasts a time equal to $2h$ times the number of mini-batches. For SG-UBU and strongly log-concave targets, Gaussian convolution inequalities bound the Wasserstein distance between the numerical invariant distribution and $\pi$ \citep{paulin2026theoretical}. The bound is the discretization error of the UBU integrator plus a stochastic gradient error of first order in $h$, which vanishes with the full gradient.

Like SG-UBU, \eqref{eq: sglm} evaluates one gradient per step and has a discretization error of second order, and the separation in \eqref{eq: intro sg mse} holds without the global convexity assumed for SG-UBU in \citep{lu2025mse}. Unlike the sweep of \citep{paulin2025sampling}, SGLM draws $\xi_k$ independently, so its numerical bias keeps the first order term \eqref{eq: intro bias coefficient} from the gradient noise. As for SG-UBU, this term comes from the gradient noise only and vanishes with the full gradient, whereas the coefficient \eqref{eq: intro sgld coefficient} of \eqref{eq: sgld} also carries the Hessian $\nabla^2 U(x)$ from the discretization error of EM. UBU and LM are analyzed differently: UBU has strong error of second order \citep{sanz2021wasserstein, paulin2024correction}, whereas the Markovian form \eqref{eq: lm markov} of LM has strong error of first order \citep{ye2026mse}. The second order discretization error of LM is recovered through the modified class of Section~\ref{sec: modified}.

\section{Preliminaries}
\label{sec: preliminaries}

\subsection{Notation and assumptions}
\label{subsec: notation}

Throughout, $\langle \cdot , \cdot \rangle$ is the Euclidean inner product, $\norm{\,\cdot\,}$ the operator norm of a tensor, which for a vector is the Euclidean norm $|\cdot|$, $I_d$ is the $d \times d$ identity matrix, $A \succeq B$ means that $A - B$ is positive semidefinite, $v^{\otimes j}$ is the $j$-fold tensor power, $A : B$ the full contraction of two tensors of the same order, $\mathrm{Tr}[A]$ the trace of a square matrix $A$, and $\mathrm{Sym}\, T$ the average of a tensor $T$ over the permutations of its indices. Derivatives of $U$ contract with tensors as in \citep{ye2026mse}: for instance, for a matrix $Q$ and a vector $v$, $\nabla^3 U(x)[Q] \in \mathbb R^d$ and $\nabla^3 U(x)[v] \in \mathbb R^{d\times d}$ have entries
\begin{equation*}
    \bigl( \nabla^3 U(x)[Q] \bigr)_i = \sum_{j,k=1}^d \partial_i\partial_j\partial_k U(x)\, Q_{jk} , \qquad
    \bigl( \nabla^3 U(x)[v] \bigr)_{ij} = \sum_{k=1}^d \partial_i\partial_j\partial_k U(x)\, v_k .
\end{equation*}
We write $\Delta U = \sum_i \partial_i^2 U$ and $\nabla \Delta U$ for its gradient, $\pi(\psi) = \int \psi\,\D\pi$, $\mathbb E^x$ for the expectation of the Langevin dynamics \eqref{eq: langevin} started at $X_0 = x$, and $\mathbb E_\xi$ and $\mathbb P_\xi$ for the expectation and the probability over the random variable $\xi$ of the stochastic gradient at fixed $x$. For a tensor-valued function $g$ and $j \Ge 0$, $\norm{\nabla^j g}_\infty := \sup_x \norm{\nabla^j g(x)}$, with $|\cdot|$ for $j = 1$; for a Lipschitz function $f$, $\norm{\nabla f}_\infty$ is its Lipschitz constant. For a random tensor $X$ and $p \Ge 1$, $\pnorm{X}_p := ( \mathbb E \norm{X}^p )^{\frac1p}$, with the Euclidean norm $|X|$ in place of $\norm{X}$ when $X$ is a vector. Most proofs do not track the dependence of the constants on the dimension $d$, so the norms used in the paper are equivalent, up to constants depending only on $d$.

For a signed measure $\nu$ with $\nu(\mathbb R^d) = 0$ and $\int |x|\, |\nu|(\D x) < \infty$, the Kantorovich--Rubinstein norm is
\begin{equation*}
    \norm{\nu}_{\mathrm{KR}} := \sup_{\norm{\nabla f}_\infty \Le 1} \int_{\mathbb R^d} f\, \D\nu ,
\end{equation*}
with the supremum over the Lipschitz functions $f$ with $\norm{\nabla f}_\infty \Le 1$, so that $\mathcal W_1(\mu, \rho) = \norm{\mu - \rho}_{\mathrm{KR}}$ for probability measures $\mu$ and $\rho$ with finite first moments by the Kantorovich--Rubinstein duality. For a function $\psi$ with $\pi(\psi) = 0$ and $\pi(|x|\,|\psi|) < \infty$, $\psi\pi$ is the signed measure $\psi(x)\, \pi(\D x)$.

Three conventions fix the notation of the paper. The solution of the discrete Poisson equation of \eqref{eq: langevin} is written $\phi_h$, the subscript $h$ marking its dependence on the step size. Every member of the modified class of Section~\ref{sec: modified} carries a tilde: the modified test function $\tilde f_h$, the modified potential $\tilde U_h$ with its Boltzmann distribution, the modified target $\tilde\pi_h$, the modified generator $\tilde{\mathcal L}_h$, and the modified discrete Poisson solution $\tilde\phi_h$. So do the modified dynamics $\tilde X_t$, defined together with $\tilde{\mathcal L}_h$, and its solution $\tilde Y_k(t)$ started at $Y_k$. The hat marks the stochastic gradient: $\nabla \hat U(x;\xi)$ is the estimator and $U(x)$ the exact potential, and $\hat{\mathcal F}_k$ and $\hat{\mathbb E}_k$ in \eqref{eq: sg filtration} are the filtration and the conditional expectation of \eqref{eq: sglm markov}. Table~\ref{tab: notation} lists the objects of \eqref{eq: langevin} next to their modified counterparts.

\begin{table}[htbp]
\centering
\begin{tabular}{lllll}
\toprule
 & \multicolumn{2}{c}{original} & \multicolumn{2}{c}{modified} \\
\cmidrule(lr){2-3} \cmidrule(lr){4-5}
object & symbol & definition & symbol & definition \\
\midrule
test function & $f$ & Assumption~\ref{ass: f} & $\tilde f_h$ & \eqref{eq: heat shift} \\
potential & $U$ & Assumption~\ref{ass: U} & $\tilde U_h$ & \eqref{eq: modified potential} \\
Langevin dynamics & $X_t$ & \eqref{eq: langevin} & $\tilde X_t$, $\tilde Y_k(t)$ & \eqref{eq: modified langevin} \\
generator & $\mathcal L$ & \eqref{eq: generator} & $\tilde{\mathcal L}_h$ & \eqref{eq: modified generator} \\
invariant distribution & $\pi$ & $\pi \propto e^{-U}$ & $\tilde\pi_h$ & $\tilde\pi_h \propto e^{-\tilde U_h}$, \eqref{eq: modified potential} \\
discrete Poisson solution & $\phi_h$ & \eqref{eq: discrete poisson solution} & $\tilde\phi_h$ & \eqref{eq: modified poisson} \\
\bottomrule
\end{tabular}
\caption{Notation of the original objects and of their modified counterparts.}
\label{tab: notation}
\end{table}

For the full gradient \eqref{eq: lm} integrator, we make the following two assumptions on the potential $U$ and the test function $f$. Strong convexity of $U$ is required only outside a bounded ball.

\begin{assumption}
	\label{ass: U}
	$U \in C^7(\mathbb R^d)$ satisfies
	\begin{equation*}
		|\nabla U(0)| \Le M , \qquad \max_{2 \Le i \Le 7}\, \norm{\nabla^i U(x)} \Le M , \qquad \forall x \in \mathbb R^d ,
	\end{equation*}
	and there are constants $m > 0$ and $R \Ge 0$ such that $\nabla^2 U(x) \succeq m I_d$ for $|x| \Ge R$.
\end{assumption}

\begin{assumption}
	\label{ass: f}
	$f \in C^4(\mathbb R^d)$ satisfies
	\begin{equation*}
		\max_{1 \Le j \Le 4}\, \norm{\nabla^j f(x)} \Le L , \qquad \forall x \in \mathbb R^d .
	\end{equation*}
\end{assumption}

For SGLM, additional boundedness conditions are required for the stochastic gradient $\nabla \hat U(x;\xi)$.

\begin{assumption}
	\label{ass: sg}
	The estimator is unbiased and its $L^{12}$ norm with respect to $\xi$ grows linearly,
	\begin{equation*}
		\mathbb E_\xi\big[\nabla \hat U(x;\xi)\big] = \nabla U(x), \qquad
		\Big( \mathbb E_\xi \big[ |\nabla \hat U(x;\xi)|^{12} \big] \Big)^{\frac{1}{12}} \Le M \bigl( 1 + |x| \bigr) , \qquad
		\forall x \in \mathbb R^d ,
	\end{equation*}
	$\hat U(\cdot;\xi) \in C^2(\mathbb R^d)$ for every $\xi$ with $(x,\xi) \mapsto \nabla^2 \hat U(x;\xi)$ measurable, and its Hessian has a second moment bound, uniform on unit balls,
	\begin{equation*}
		\mathbb E_\xi \bigg[ \sup_{|x'-x| \Le 1} \norm{\nabla^2 \hat U(x';\xi)}^2 \bigg] \Le M^2 , \qquad
		\forall x \in \mathbb R^d ,
	\end{equation*}
	with $M$ the constant of Assumption~\ref{ass: U}, and its variance is bounded by a constant $\Sigma \Le M^2$,
	\begin{equation*}
		\mathbb E_\xi \big[ |\nabla \hat U(x;\xi) - \nabla U(x)|^2 \big] \Le \Sigma , \qquad
		\forall x \in \mathbb R^d .
	\end{equation*}
\end{assumption}

The constraint $\Sigma \Le M^2$ in Assumption~\ref{ass: sg} is imposed for the convenience of analysis. The bound on the Hessian of $\hat U$ in Assumption~\ref{ass: sg} ensures the uniqueness of the invariant distribution of \eqref{eq: sglm markov}, as required by Harris's ergodic theorem \citep{hairer2008yet}. A stronger condition is the almost sure uniform bound $\norm{\nabla^2 \hat U(x;\xi)} \Le M$ for all $x \in \mathbb R^d$ and almost every $\xi$. Although Assumption~\ref{ass: sg} requires the potential form, $\nabla \hat U(x;\xi)$ does not need to be a gradient: the proofs use only the boundedness of its derivative $\nabla^2 \hat U(x;\xi)$.

From here on, the letter $C$ denotes a positive constant depending only on $m$, $M$, $R$ and $d$, but not on $h$, $N$, $P$, $\Sigma$ or the initial state, and it changes from line to line. The letter $h_0$ denotes an upper bound in $(0,1]$ on the step size $h$, depending only on $m$, $M$, $R$ and $d$. In each statement, $h_0$ is the smallest of the bounds its proof uses, and every statement below holds for $h \in (0,h_0]$.

Although this paper focuses on SGLM, the proofs are carried out also for the full gradient LM, as the stochastic gradient case follows directly in most places. When we speak of the LM integrators in the plural, the statement covers both the full gradient and the stochastic gradient cases.

\subsection{Markovian form of the LM integrators}

Although the \eqref{eq: lm} integrator is non-Markovian, it can be written in a Markovian form for the convenience of analysis. Let $(B_t)_{t \Ge 0}$ be a standard Brownian motion in $\mathbb R^d$, the one driving \eqref{eq: langevin}, and realize the Gaussian vectors as its increments over the steps,
\begin{equation*}
	\chi_k = \frac1{\sqrt{h}}\big(B_{(k+1)h} - B_{kh}\big), \qquad k \Ge 0 .
\end{equation*}
Let $\mathcal F_k$ be the sigma algebra generated by the Brownian motion up to the time $kh$, the past before the step $k$, and $\mathbb E_k$ the conditional expectation given that past,
\begin{equation}
	\mathcal F_k := \sigma\bigl( B_s : 0 \Le s \Le kh \bigr) , \qquad
	\mathbb E_k[\,\cdot\,] := \mathbb E\bigl[\,\cdot\, \big| \mathcal F_{k}\bigr] ,
	\label{eq: filtration}
\end{equation}
so that $\chi_k$ is $\mathcal F_{k+1}$-measurable and independent of $\mathcal F_{k}$. We construct a Markov chain $(Y_k)_{k\Ge0}$ by
\begin{equation}
	Y_{k+1} = Y_k - h\, \nabla U\Big(Y_k + \sqrt{\frac h2}\, \chi_k\Big) + \sqrt{2h}\, \chi_k .
	\label{eq: lm markov}
	\tag{LM'}
\end{equation}
Then $Z_k = Y_k + \sqrt{h/2}\, \chi_k$ satisfies \eqref{eq: lm}, so the LM integrator is equivalent to \eqref{eq: lm markov}, the form in which it is analyzed as a postprocessed integrator in \citep{vilmart2015postprocessed}; here $Y_k$ is $\mathcal F_{k}$-measurable and $Z_k$ is $\mathcal F_{k+1}$-measurable. Similarly, the \eqref{eq: sglm} integrator is equivalent to
\begin{equation}
	Y_{k+1} = Y_k - h\, \nabla \hat U\Big(Y_k + \sqrt{\frac h2}\, \chi_k; \xi_k\Big) + \sqrt{2h}\, \chi_k ,
	\label{eq: sglm markov}
	\tag{SGLM'}
\end{equation}
where the random variables $\xi_k$ are independent and identically distributed and independent of $(B_t)_{t \Ge 0}$.
For \eqref{eq: sglm markov} the past before the step $k$ also contains $\xi_0, \dots, \xi_{k-1}$, and we define 
\begin{equation}
	\hat{\mathcal F}_k := \sigma\bigl( B_s ,\ \xi_j : 0 \Le s \Le kh ,\ 0 \Le j \Le k-1 \bigr) , \qquad
	\hat{\mathbb E}_k[\,\cdot\,] := \mathbb E\bigl[\,\cdot\, \big| \hat{\mathcal F}_{k}\bigr] ,
	\label{eq: sg filtration}
\end{equation}
so that $\chi_k$ and $\xi_k$ are independent of $\hat{\mathcal F}_k$, and $Y_k$ is $\hat{\mathcal F}_{k}$-measurable and $Z_k$ is $\hat{\mathcal F}_{k+1}$-measurable.

Since $\chi_k$ is determined by the Brownian motion on $[kh,(k+1)h]$, we call \eqref{eq: lm markov} and \eqref{eq: sglm markov} the Markovian forms of LM and SGLM. The analysis through the discrete Poisson equation is carried out on the Markovian form $(Y_k)_{k \Ge 0}$, not on the original form $(Z_k)_{k \Ge 0}$, and their connection is made by the modified class of Section~\ref{sec: modified}. The Markovian form \eqref{eq: lm markov} has strong and weak order $1$ \citep{ye2026mse}, hence the second order discretization error of LM comes entirely from the Gaussian shift $\sqrt{h/2}\, \chi_k$ between $Y_k$ and $Z_k$.

\subsection{Elliptic estimates of the discrete Poisson equation}
\label{subsec: poisson}

The generator of the overdamped Langevin dynamics \eqref{eq: langevin} is
\begin{equation}
    \mathcal L \psi = \Delta \psi - \nabla U \cdot \nabla \psi ,
    \label{eq: generator}
\end{equation}
and $e^{t\mathcal L}$, $t \Ge 0$, is its semigroup, $e^{t\mathcal L} f(x) = \mathbb E^x[f(X_t)]$, where $X_t$ solves \eqref{eq: langevin} with $X_0 = x$. The Dirichlet form of $\mathcal L$ \citep{bakry2014analysis} gives the integration by parts
\begin{equation}
    \pi(g_1\, \mathcal L g_2) = -\pi(\nabla g_1 \cdot \nabla g_2)
    \label{eq: generator ibp}
\end{equation}
for $g_1, g_2 \in C^2(\mathbb R^d)$ whose derivatives of orders $0$ to $2$ grow at most polynomially.

For the test function $f$ of Assumption~\ref{ass: f}, the Kolmogorov solution $u(x,t) := e^{t\mathcal L} f(x) - \pi(f)$ solves the parabolic equation $\partial_t u = \mathcal L u$ with $u(\cdot,0) = f - \pi(f)$. By the exponential convergence of $u(x,t)$ established in Theorem~\ref{thm: u decay} below, the discrete Poisson solution \citep{lu2025mse, ye2026mse}
\begin{equation}
    \phi_h(x) := h \sum_{k=0}^{\infty} u(x,kh) = h\, \bigl( \mathcal I - e^{h\mathcal L} \bigr)^{-1} \bigl( f(x) - \pi(f) \bigr) ,
    \label{eq: discrete poisson solution}
\end{equation}
with $\mathcal I$ the identity operator, is well defined and, by telescoping the series with $e^{h\mathcal L} u(\cdot, kh) = u(\cdot, (k+1)h)$ and the exponential decay \eqref{eq: u bound} of $u$, satisfies the discrete Poisson equation
\begin{equation}
    \frac{\phi_h(x) - (e^{h\mathcal L} \phi_h)(x)}{h} = f(x) - \pi(f) , \qquad x \in \mathbb R^d .
    \label{eq: discrete poisson}
\end{equation}
As $h \to 0$, \eqref{eq: discrete poisson} formally becomes the continuous Poisson equation
\begin{equation}
    -\mathcal L \phi_0(x) = f(x) - \pi(f) , \qquad x \in \mathbb R^d ,
    \label{eq: continuous poisson}
\end{equation}
studied in \citep{mattingly2010convergence, vollmer2016exploration, brosse2019tamed}. By \eqref{eq: u bound} of Theorem~\ref{thm: u decay} below, $\phi_0 = \int_0^\infty u(\cdot, t)\, \D t$ solves \eqref{eq: continuous poisson}, with derivatives of orders $1$ to $4$ bounded by $C L$. Among the functions with a bounded gradient, \eqref{eq: continuous poisson} determines $\phi_0$ up to an additive constant, as the proof of Corollary~\ref{cor: bias coefficient} shows.

\begin{theorem}
\label{thm: u decay}
    Let Assumptions~\ref{ass: U} and \ref{ass: f} hold. There are $\lambda > 0$ and $C$, depending only on $m$, $M$, $R$ and $d$, such that for all $x \in \mathbb R^d$ and $t \Ge 0$,
    \begin{equation}
        |u(x,t)| \Le C L\, (1 + |x|)\, e^{-\lambda t} , \qquad
        \max_{1 \Le j \Le 4}\, \norm{\nabla^j u(x,t)} \Le C L\, e^{-\lambda t} .
        \label{eq: u bound}
    \end{equation}
    Consequently, for every $h \in (0,\lambda^{-1}]$ the discrete Poisson solution \eqref{eq: discrete poisson solution} satisfies
    \begin{equation}
        \max_{1 \Le j \Le 4}\, \norm{\nabla^j \phi_h(x)} \Le C L , \qquad \forall x \in \mathbb R^d .
        \label{eq: poisson bound}
    \end{equation}
    \emph{(Proof in Appendix~\ref{app: elliptic}.)}
\end{theorem}

Theorem~\ref{thm: u decay} requires one order of derivative of $u$ more than Theorem 1 of \citep{ye2026mse}, and its proof is almost the same; the details are given in Appendix~\ref{app: elliptic}. The proof uses the derivatives of $U$ up to order $5$ only, although Assumption~\ref{ass: U} provides bounds up to order $7$; the $2$ extra orders are consumed by the modified potential of Section~\ref{subsec: modified}, whose construction loses $2$ derivatives. Corollary~\ref{cor: modified poisson} applies it to the modified potential $\tilde U_h$ and the modified test function $\tilde f_h$, whose discrete Poisson solution is written $\tilde\phi_h$.

\section{Modified class for the LM integrator}
\label{sec: modified}

The MSE \eqref{eq: mse} averages the test function along the samples $(Z_k)_{k \Ge 0}$ of the \eqref{eq: lm} integrator, while the analysis runs on the Markovian form \eqref{eq: lm markov}. The Gaussian shift between the two is absorbed by the \emph{modified class}: a modified test function $\tilde f_h$, a modified potential $\tilde U_h$ with its modified generator $\tilde{\mathcal L}_h$, and a modified discrete Poisson solution $\tilde\phi_h$. The modified class is local: its construction does not involve the long-time behavior. Section~\ref{subsec: modified test function} defines the modified test function $\tilde f_h$ as the heat flow of $f$, which splits the MSE into the time average error of the Markovian form and a bias against a reference measure. Section~\ref{subsec: modified} defines the modified potential $\tilde U_h$ by the correction that makes the bias of order $h^2$, cut off far from the origin so that the bounds of Assumption~\ref{ass: U} up to order $5$ hold uniformly in $h$, with $2M$ and $\frac m2$ in place of $M$ and $m$. Section~\ref{subsec: modified bias} proves that the bias against the modified target $\tilde\pi_h$ is of order $h^2$. Section~\ref{subsec: modified poisson} defines the modified generator $\tilde{\mathcal L}_h$ and the modified discrete Poisson solution $\tilde\phi_h$.

\subsection{Modified test function $\tilde f_h$}
\label{subsec: modified test function}

The positions $Z_k$ and $Y_k$ in \eqref{eq: lm} and \eqref{eq: lm markov} differ by the Gaussian shift $\sqrt{h/2}\, \chi_k \sim \mathcal N(0, \frac h2 I_d)$, which is independent of $\mathcal F_{k}$. The conditional mean of $f(Z_k)$ given $\mathcal F_{k}$ is therefore the heat semigroup applied to $f$ at time $\frac h4$, evaluated at $Y_k$. This motivates the definition of the \emph{modified test function}
\begin{equation}
	\tilde f_h(y) := (e^{\frac h4\Delta }f)(y) = 
	\mathbb E\biggl[ f\Bigl( y + \sqrt{\frac h2}\, \chi \Bigr) \biggr] , \quad \chi \sim \mathcal N(0, I_d) ,
	\label{eq: heat shift}
\end{equation}
for which $\mathbb E_k\bigl[ f(Z_k) \bigr] = \tilde f_h(Y_k)$. Since $\tilde f_h$ is an average of translates of $f$, every derivative of $\tilde f_h$ is the same average of the corresponding derivative of $f$, and $\tilde f_h$ satisfies Assumption~\ref{ass: f} with the same constant $L$ for every $h$. Expanding $f$ about $y$ to the order $4$ along $\sqrt{h/2}\, \chi$,
\begin{equation*}
	f\Bigl( y + \sqrt{\frac h2}\, \chi \Bigr) = f(y) + \sqrt{\frac h2}\, \nabla f(y) \cdot \chi + \frac h4\, \nabla^2 f(y)[\chi, \chi] + \frac16 \Bigl( \frac h2 \Bigr)^{\frac32} \nabla^3 f(y)\bigl[ \chi^{\otimes 3} \bigr] + \varrho_4^{\chi}(y) ,
\end{equation*}
where the remainder satisfies
\begin{equation*}
	|\varrho_4^{\chi}(y)| \Le \frac{L}{24} \Bigl( \frac h2 \Bigr)^2 |\chi|^4 .
\end{equation*}
Under the expectation the terms of odd order in $\chi$ vanish, hence
\begin{equation}
    \Bigl| \tilde f_h(y) - f(y) - \frac h4\, \Delta f(y) \Bigr| \Le \frac{L}{24}\, \Bigl( \frac h2 \Bigr)^2\, \mathbb E |\chi|^4 \Le C L h^2 , \qquad \forall y \in \mathbb R^d,
    \label{eq: shift expansion}
\end{equation}
where $\Delta f = \mathrm{Tr}[\nabla^2 f] = \sum_{i=1}^d \partial_i^2 f$.
The modified test function $\tilde f_h$ absorbs the Gaussian shift, so that the analysis runs on the positions $Y_k$ of the Markovian form \eqref{eq: lm markov}.

The modified test function $\tilde f_h(Y_k)$ is an unbiased estimate of $f(Z_k)$ under $\mathcal F_{k}$. Define the difference $\mathfrak M_k := f(Z_k) - \tilde f_h(Y_k)$; then $\mathbb E_k[\mathfrak M_k] = 0$ by \eqref{eq: heat shift}, and since $f$ is $L$-Lipschitz,
\begin{equation}
    \mathbb E_k\bigl[ \mathfrak M_k^2 \bigr] \Le \mathbb E_k\Bigl[ \bigl| f(Z_k) - f(Y_k) \bigr|^2 \Bigr] \Le L^2\, \frac h2\, \mathbb E |\chi|^2 = \frac{d h}{2}\, L^2 .
    \label{eq: martingale variance}
\end{equation}
The $\mathfrak M_k$ are $\mathcal F_{k+1}$-measurable martingale differences, so their sum has additive variance. For every probability measure $\rho$ on $\mathbb R^d$ with a finite first moment, inserting $f(Z_k) = \tilde f_h(Y_k) + \mathfrak M_k$ into the MSE \eqref{eq: mse}, applying Cauchy--Schwarz inequality to the three resulting terms and then the additive variance with \eqref{eq: martingale variance} gives
\begin{equation}
\begin{aligned}
    \MSE(N,h)
    &= \mathbb E \Biggl[ \biggl| \frac1N \sum_{k=0}^{N-1} \bigl( \tilde f_h(Y_k) - \rho(\tilde f_h) \bigr) + \bigl( \rho(\tilde f_h) - \pi(f) \bigr) + \frac1N \sum_{k=0}^{N-1} \mathfrak M_k \biggr|^2 \Biggr] \\
    &\Le 3\, \mathbb E \Biggl[ \biggl| \underbrace{\frac1N \sum_{k=0}^{N-1} \bigl( \tilde f_h(Y_k) - \rho(\tilde f_h) \bigr)}_{\text{time average error}} \biggr|^2 \Biggr] + 3\, \bigl( \underbrace{\rho(\tilde f_h) - \pi(f)}_{\text{modified bias}} \bigr)^2 + \frac{3 d L^2}{2}\, \frac{h}{N} .
\end{aligned}
    \label{eq: mse split}
\end{equation}
The term $\frac{3 d L^2}{2}\, \frac hN$ of \eqref{eq: mse split} is absorbed into $\frac{1}{Nh}$, so it remains to estimate the time average error and the modified bias. The measure $\rho$ is to be chosen so that the modified bias and the discretization part of the time average error are of order $h^2$.

The measure $\rho$ cannot simply be $\pi$. With $\rho = \pi$, the expansion \eqref{eq: shift expansion} gives
\begin{equation*}
	\pi(\tilde f_h) - \pi(f) = \frac h4\, \pi(\Delta f) + \mathcal O(L h^2) ,
\end{equation*}
whose leading coefficient $\pi(\Delta f)$ is in general nonzero. The choice of $\rho$ has to compensate for the Laplacian term, so that the modified bias $\rho(\tilde f_h) - \pi(f)$ is of order $h^2$.

\subsection{Modified potential $\tilde U_h$}
\label{subsec: modified}

We aim to find a distribution $\rho$ such that the modified bias $\rho(\tilde f_h) - \pi(f)$ is of order $h^2$. Since $\tilde f_h = f + \frac h4 \Delta f + \mathcal O(h^2)$ by \eqref{eq: shift expansion}, we consider the ansatz
\begin{equation}
	\rho \propto \exp\Bigl( - U - \frac h4\, V \Bigr) ,
	\label{eq: ansatz}
\end{equation}
where $V$ is a function to be determined and the coefficient $\frac h4$ matches that of $\Delta f$ in \eqref{eq: shift expansion}. We derive $V = |\nabla U|^2 - \Delta U$ by a perturbation analysis. The perturbed potential $U + \frac h4 V$, once made rigorous by a cutoff far from the origin, is the \emph{modified potential} $\tilde U_h$, and its Boltzmann distribution $\tilde\pi_h \propto \exp(-\tilde U_h)$ is the \emph{modified target}. Finally, we prove the modified bias estimate $\tilde\pi_h(\tilde f_h) - \pi(f) = \mathcal O(h^2)$.

\subsubsection*{Derivation of $V$ by perturbation}

For a test function $\psi$, the ansatz \eqref{eq: ansatz} gives
\begin{equation*}
\rho(\psi) = \frac{\pi( \psi e^{-\frac h4 V} ) }{ \pi( e^{-\frac h4 V})} = \frac{\int_{\mathbb R^d} \psi e^{-U-\frac h4V}\D x}{\int_{\mathbb R^d} e^{-U-\frac h4V}\D x}.
\end{equation*}
Expanding $e^{-\frac h4 V} = 1 - \frac h4 V + \mathcal O(h^2)$ in the numerator and in the denominator of $\rho(\psi)$, we obtain
\begin{equation}
    \rho(\psi) = \pi(\psi) - \frac h4 \bigl( \pi(\psi V) - \pi(\psi)\, \pi(V) \bigr) + \mathcal O(h^2) ,
    \label{eq: rho expansion}
\end{equation}
the term of order $h$ being $-\frac h4$ times the covariance of $\psi$ and $V$ under $\pi$. Taking $\psi = \tilde f_h$ in \eqref{eq: rho expansion} and inserting the expansion \eqref{eq: shift expansion} of the modified test function,
\begin{equation}
    \rho(\tilde f_h) - \pi(f) = \frac h4 \Bigl( \pi(\Delta f) - \pi(f V) + \pi(f)\, \pi(V) \Bigr) + \mathcal O( h^2) .
    \label{eq: bias order h}
\end{equation}
The modified bias is therefore of order $h^2$ exactly when the term of order $h$ in \eqref{eq: bias order h} vanishes,
\begin{equation}
    \pi(\Delta f) = \pi(f V) - \pi(f)\, \pi(V) .
    \label{eq: V condition}
\end{equation}
We look for a single $V$ satisfying \eqref{eq: V condition} for every $f$. The integration by parts \eqref{eq: generator ibp}, first with $g_1 = 1$ and $g_2 = f$ and then with $g_1 = f$ and $g_2 = U$, gives
\begin{equation}
    \pi(\Delta f) = \pi\bigl( \nabla f \cdot \nabla U \bigr) = \pi\bigl( f\, |\nabla U|^2 - f\, \Delta U \bigr) .
    \label{eq: integration by parts}
\end{equation}
Hence we select the potential function
\begin{equation}
    V := |\nabla U|^2 - \Delta U = -\mathcal L U ,
    \label{eq: V}
\end{equation}
for which $\pi(V) = -\pi(\mathcal L U) = 0$ by the invariance of $\pi$.
This computation is formal, and the cutoff introduced below changes the measure \eqref{eq: ansatz}; Lemma~\ref{lem: shifted bias} proves the bias estimate for the modified target.

\subsubsection*{Cutoff in the modified potential}

The naive construction $U + \frac h4 V$ of the modified potential does not satisfy Assumption~\ref{ass: U}, because the derivatives of $|\nabla U|^2$ grow with $|x|$. Since Theorem~\ref{thm: u decay} must be applied to the modified potential with constants independent of $h$, the modification is cut off far from the origin. Fix a smooth cutoff function $\theta \in C^\infty(\mathbb R^d)$ with $0 \Le \theta \Le 1$, $\theta = 1$ on $\{|x| \Le 1\}$ and $\theta = 0$ on $\{|x| \Ge 2\}$. With $\delta \in (0,1]$ fixed by \eqref{eq: delta bound} in the proof of Lemma~\ref{lem: modified class} below, we define the modified potential by
\begin{equation}
    \tilde U_h := U + \frac h4\, V\, \theta_h , \qquad \theta_h(x) := \theta\Bigl( \frac{h x}{\delta} \Bigr) ,
    \label{eq: modified potential}
\end{equation}
so that $\tilde U_h = U + \frac h4 V$ on the ball $\{|x| \Le \delta/h\}$ and $\tilde U_h = U$ on $\{|x| \Ge 2\delta/h\}$. In particular, since $0 \Le \theta_h \Le 1$ and $|V(y)| \Le C(1 + |y|)^2$ by Assumption~\ref{ass: U},
\begin{equation}
    |\tilde U_h(x) - U(x)| \Le \frac h4 \max_{|y| \Le 2\delta/h} |V(y)| \Le \frac{C}{h} , \qquad x \in \mathbb R^d .
    \label{eq: modified potential difference}
\end{equation}
The modified target $\tilde\pi_h \propto e^{-\tilde U_h}$ is therefore a perturbation of $\pi$ by a bounded potential, although the bound \eqref{eq: modified potential difference} is not uniform in $h$. The derivatives of $\tilde U_h - U$ of orders $1$ to $5$, however, are bounded uniformly in $x$ and $h$: by \eqref{eq: modified class} in Lemma~\ref{lem: modified class} below they are at most $C h (1 + |x|)$, and they vanish for $|x| > 2\delta/h$.

\begin{lemma}
\label{lem: modified class}
    Let Assumption~\ref{ass: U} hold. There are $\delta \in (0,1]$ and $h_0 \in (0, \delta]$ such that for every $h \in (0,h_0]$, the modified potential $\tilde U_h \in C^5(\mathbb R^d)$ satisfies
    \begin{equation}
        |\nabla \tilde U_h(0)| \Le 2M , \qquad \max_{2 \Le i \Le 5}\, \norm{\nabla^i \tilde U_h(x)} \Le 2M , \qquad \forall x \in \mathbb R^d ,
        \label{eq: modified class bounds}
    \end{equation}
    the convexity condition $\nabla^2 \tilde U_h(x) \succeq \frac m2 I_d$ for $|x| \Ge R$, and
    \begin{equation}
        \max_{1 \Le i \Le 5}\, \norm{\nabla^i \tilde U_h(x) - \nabla^i U(x)} \Le C h \bigl( 1 + |x| \bigr) , \qquad \forall x \in \mathbb R^d .
        \label{eq: modified class}
    \end{equation}
    \emph{(Proof in Appendix~\ref{sec: modified class proof}.)}
\end{lemma}

Theorem~\ref{thm: u decay} therefore applies to $\tilde U_h$, since its proof uses the bounds of Assumption~\ref{ass: U} up to order $5$ only (Appendix~\ref{app: elliptic}). The constants, in particular the rate $\lambda$, depend only on $m$, $M$, $R$ and $d$, not on $h$. In the rest of the paper, the number $\delta$ is fixed by Lemma~\ref{lem: modified class} and absorbed into the constant $C$.

\subsection{Estimation of the modified bias $\tilde\pi_h(\tilde f_h) - \pi(f)$}
\label{subsec: modified bias}

The cancellation \eqref{eq: V condition} between the heat flow of the test function and the modification of the potential holds for the modified target $\tilde\pi_h$ as it does for the ansatz \eqref{eq: ansatz}, since the cutoff acts only where $\pi$ has exponentially small mass. Lemma~\ref{lem: shifted bias} carries the perturbation argument that produced $V$ over to the modified target.

\begin{lemma}
\label{lem: shifted bias}
    Let Assumption~\ref{ass: U} hold. For every $h \in (0,h_0]$ and every $f \in C^1(\mathbb R^d)$ with $\norm{\nabla f}_\infty < \infty$,
    \begin{equation}
        \bigl| \tilde\pi_h(\tilde f_h) - \pi(f) \bigr| \Le C h^2 \norm{\nabla f}_\infty .
        \label{eq: shifted bias}
    \end{equation}
    \emph{(Proof in Appendix~\ref{sec: shifted bias proof}.)}
\end{lemma}

\begin{remark}
The Gaussian potential $U = \frac12 |x|^2$, taken without the cutoff, illustrates the modification. Here $V = |x|^2 - d$, so the modified potential equals $(1 + \frac h2) \frac{|x|^2}{2}$ up to an additive constant and $\tilde\pi_h$ has variance $(1 + \frac h2)^{-1}$ in each coordinate. The chain becomes $Y_{k+1} = (1-h) Y_k + \sqrt{2h}\,(1 - \frac h2)\, \chi_k$, whose stationary variance $1 - \frac h2$ agrees with that of $\tilde\pi_h$ up to $\mathcal O(h^2)$. The shifted variable $Z_k = Y_k + \sqrt{h/2}\, \chi_k$ has stationary variance $1$, in agreement with the exactness of \eqref{eq: lm} for Gaussian targets established in \citep[Remark~4.1]{vilmart2015postprocessed}.
\end{remark}

\subsection{Modified discrete Poisson equation $\tilde\phi_h$}
\label{subsec: modified poisson}

Let the \emph{modified dynamics} be the overdamped Langevin dynamics driven by the potential $\tilde U_h$,
\begin{equation}
    \D \tilde X_t = -\nabla \tilde U_h(\tilde X_t)\,\D t + \sqrt2\,\D B_t ,
    \label{eq: modified langevin}
\end{equation}
with the \emph{modified generator}
\begin{equation}
    \tilde{\mathcal L}_h \psi = \Delta \psi - \nabla \tilde U_h \cdot \nabla \psi ,
    \label{eq: modified generator}
\end{equation}
and its semigroup $e^{t\tilde{\mathcal L}_h}$. The invariant distribution of \eqref{eq: modified langevin} is the modified target $\tilde\pi_h$. The time average error in \eqref{eq: mse split} with $\rho = \tilde\pi_h$ is analyzed through the discrete Poisson equation of the modified generator \eqref{eq: modified generator}. Let $\tilde\phi_h$ be the discrete Poisson solution \eqref{eq: discrete poisson solution} with $\tilde{\mathcal L}_h$, $\tilde\pi_h$ and $\tilde f_h$ in place of $\mathcal L$, $\pi$ and $f$, so that
\begin{equation}
    \frac{\tilde\phi_h(x) - (e^{h\tilde{\mathcal L}_h} \tilde\phi_h)(x)}{h} = \tilde f_h(x) - \tilde\pi_h(\tilde f_h) , \qquad x \in \mathbb R^d .
    \label{eq: modified poisson}
\end{equation}

\begin{corollary}
\label{cor: modified poisson}
    Let Assumptions~\ref{ass: U} and \ref{ass: f} hold. For every $h \in (0,h_0]$, $\tilde\phi_h$ is well defined and
    \begin{equation}
        \max_{1 \Le j \Le 4}\, \norm{\nabla^j \tilde\phi_h(x)} \Le C L , \qquad \forall x \in \mathbb R^d .
        \label{eq: modified poisson bound}
    \end{equation}
\end{corollary}

Corollary~\ref{cor: modified poisson} is Theorem~\ref{thm: u decay} applied to the modified class. By Lemma~\ref{lem: modified class}, $\tilde U_h$ satisfies the bounds that the proof of Theorem~\ref{thm: u decay} uses, and $\tilde f_h$ satisfies Assumption~\ref{ass: f} with the constant $L$. The constants in \eqref{eq: modified poisson bound} therefore depend only on $m$, $M$, $R$ and $d$, and the threshold $h_0$ includes $\lambda^{-1}$.

\section{Local error of the LM integrators}
\label{sec: local}

This section compares one step of the Markovian form \eqref{eq: lm markov} with one step of the modified dynamics \eqref{eq: modified langevin}, and derives from this comparison the local weak error used in the MSE analysis. Section~\ref{subsec: bias cutoff} bounds the moments of the chain uniformly in time, so that the conditional one-step estimates can be averaged over the chain. It then bounds the sizes of a step and matches the moments of orders $1$ to $3$ of the chain and of the modified dynamics, with a mismatch of order $h^3$. Section~\ref{subsec: weak expansion} turns the matching into the local weak error of the LM integrator, and of the SGLM integrator, whose gradient noise adds an explicit term of order $h^2$ and a remainder of order $h^3$.

Taking $\rho = \tilde\pi_h$ in \eqref{eq: mse split}, Lemma~\ref{lem: shifted bias} bounds the modified bias by $C L h^2$, and the MSE satisfies
\begin{equation}
    \MSE(N,h) \Le 3\, \mathbb E \Biggl[ \biggl| \underbrace{\frac1N \sum_{k=0}^{N-1} \bigl( \tilde f_h(Y_k) - \tilde\pi_h(\tilde f_h) \bigr)}_{\text{time average error}} \biggr|^2 \Biggr] + C L^2 \Big(h^4 + \frac{dh}{N}\Big).
    \label{eq: mse split modified}
\end{equation}
The time average error is centered at $\tilde\pi_h(\tilde f_h)$, so it is controlled through the discrete Poisson equation of a dynamics that preserves $\tilde\pi_h$. This is the modified dynamics \eqref{eq: modified langevin}, whose discrete Poisson solution $\tilde\phi_h$ in \eqref{eq: modified poisson} has the right side $\tilde f_h - \tilde\pi_h(\tilde f_h)$. The discrete Poisson solution transfers local errors into the time average error \citep{lu2025mse, ye2026mse}. The local error studied in this section is therefore the error of one step of the Markovian form \eqref{eq: lm markov} against one step of the modified dynamics from the same state.

Since the chains \eqref{eq: lm markov} and \eqref{eq: sglm markov} are time-homogeneous, statements about a single step are made at the index $0$ and apply at any $k$; Lemmas~\ref{lem: moments}, \ref{lem: weak} and \ref{lem: sg weak} are stated at the index $k$, in the form the MSE analysis uses. One step of the \eqref{eq: lm markov} chain is
\begin{equation}
    Y_1 - Y_0 = \sqrt{2h}\, \chi_0 - h\, \nabla U(Z_0) , \qquad Z_0 = Y_0 + \sqrt{\frac h2}\, \chi_0 ,
    \label{eq: chain step}
\end{equation}
where $\chi_0 \sim \mathcal N(0, I_d)$ is independent of $\mathcal F_0$, and one step of the \eqref{eq: sglm markov} chain is
\begin{equation}
    Y_1 - Y_0 = \sqrt{2h}\, \chi_0 - h\, \nabla \hat U(Z_0;\xi_0) ,
    \label{eq: sg chain step}
\end{equation}
where $Z_0$ is given by \eqref{eq: chain step}, and $\chi_0$ and $\xi_0$ are independent of $\hat{\mathcal F}_0$. For each $k \Ge 0$, let $(\tilde Y_k(t))_{t \Ge 0}$ solve the modified dynamics \eqref{eq: modified langevin} from $\tilde Y_k(0) = Y_k$, driven by $(B_{kh+t} - B_{kh})_{t \Ge 0}$, whose increment over $[0,h]$ is $\sqrt h\, \chi_k$. Then $\mathbb E_k[\psi(\tilde Y_k(h))] = (e^{h\tilde{\mathcal L}_h} \psi)(Y_k)$. From the same state $Y_0$, the chains \eqref{eq: lm markov} and \eqref{eq: sglm markov} have the same $Z_0$ and the same $(\tilde Y_0(t))_{t \Ge 0}$, and differ only in $Y_1$. For both chains, the local error of one step against the modified dynamics is $Y_1 - \tilde Y_0(h)$.

\subsection{Moment bounds and moment matching}
\label{subsec: bias cutoff}

The MSE analysis squares the polynomial weight of the local weak errors of Section~\ref{subsec: weak expansion}, and therefore needs moments of order $12$ of the chain, uniformly in time. The following lemma provides them. Its proof is standard and follows the Lyapunov function arguments in \citep{ye2026mse}, with Lemma~\ref{lem: modified class} for the modified dynamics.

\begin{lemma}
\label{lem: moments}
    Let Assumption~\ref{ass: U} hold. There are $h_0 > 0$ and $C$ such that for the chain \eqref{eq: lm markov}, every $h \in (0,h_0]$ and every $k \Ge 0$,
    \begin{equation}
        \mathbb E |Y_k|^{12} \Le C \bigl( 1 + |Y_0| \bigr)^{12} , \qquad
        \max\bigl\{ \mathbb E |\tilde Y_k(h) - Y_k|^{12} ,\; \mathbb E |Y_{k+1} - Y_k|^{12} \bigr\} \Le C h^{6} \bigl( 1 + |Y_0| \bigr)^{12} .
        \label{eq: moment bounds}
    \end{equation}
    If Assumption~\ref{ass: sg} additionally holds, the same moment bounds hold for the chain \eqref{eq: sglm markov}, with the same $h_0$ and $C$. \emph{(Proof in Appendix~\ref{sec: moments proof}.)}
\end{lemma}

The next lemma bounds the size of one step of the chain and of the modified dynamics, and the size of the local error between them; these sizes control the remainders of all the expansions below. Lemmas~\ref{lem: sizes} and \ref{lem: moment matching} are stated for \eqref{eq: lm markov} only; the local weak error estimate \eqref{eq: sg weak} of \eqref{eq: sglm markov} relies on them only through the full gradient step \eqref{eq: chain step}.

\begin{lemma}
\label{lem: sizes}
    Let Assumption~\ref{ass: U} hold. Then, for the chain \eqref{eq: lm markov} and every $s \in [0,h]$,
    \begin{equation}
        \pnorm{Y_1 - Y_0}_4 \Le C \sqrt h \bigl( 1 + |Y_0| \bigr) , \quad
        \pnorm{\tilde Y_0(s) - Y_0}_4 \Le C \sqrt h \bigl( 1 + |Y_0| \bigr) , \quad
        \pnorm{Y_1 - \tilde Y_0(h)}_4 \Le C h^{\frac32} \bigl( 1 + |Y_0| \bigr) ,
        \label{eq: sizes}
    \end{equation}
    where $\pnorm{\cdot}_4$ is taken under $\mathbb E_0$. \emph{(Proof in Appendix~\ref{sec: sizes proof}.)}
\end{lemma}

The moments of one step of the modified dynamics and of the chain are expanded in Lemmas~\ref{lem: exact moments} and \ref{lem: step moments} of Appendix~\ref{app: local}. Compared with the exact dynamics \eqref{eq: langevin}, the chain matches the moments of orders $2$ and $3$ to $\mathcal O(h^3)$ but misses the mean by $\frac{h^2}{4}\nabla V$, with $V$ as in \eqref{eq: V}; the drift correction of the modified dynamics removes exactly this mismatch. The cutoff keeps the cancellation on the ball $\{|Y_0| \Le \delta/h\}$, and outside the ball the remaining mismatch is again of order $h^3$, because $h|Y_0|/\delta \Ge 1$ there.

\begin{lemma}
\label{lem: moment matching}
    Let Assumption~\ref{ass: U} hold. Then, for the chain \eqref{eq: lm markov},
    \begin{equation}
        \Bigl\| \mathbb E_0\bigl[ (Y_1 - Y_0)^{\otimes j} \bigr] - \mathbb E_0\bigl[ (\tilde Y_0(h) - Y_0)^{\otimes j} \bigr] \Bigr\| \Le C h^3 \bigl( 1 + |Y_0| \bigr)^{6} , \qquad j = 1, 2, 3 .
        \label{eq: moment matching}
    \end{equation}
    \emph{(Proof in Appendix~\ref{sec: matching proof}.)}
\end{lemma}

\subsection{Local weak error expansion}
\label{subsec: weak expansion}

The local weak error of the LM integrator on a test function $G$ is $\mathbb E_0\bigl[ G(Y_1) - G(\tilde Y_0(h)) \bigr]$, and the MSE analysis takes $G = \tilde\phi_h$. Because the local error is of size $h^{\frac32}$, against $h^2$ for the SRK integrators of \citep{ye2026mse}, $G$ is expanded along a step to the order $3$: the terms of orders $1$ to $3$ are matched by Lemma~\ref{lem: moment matching}, and the remainder is controlled by Lemma~\ref{lem: sizes}. This is why Assumption~\ref{ass: f} asks for $4$ derivatives of $f$.

\begin{lemma}
\label{lem: weak}
    Let Assumption~\ref{ass: U} hold and let $G \in C^4(\mathbb R^d)$ have bounded derivatives of orders $1$ to $4$. Then, for the chain \eqref{eq: lm markov} and every $k \Ge 0$,
    \begin{equation}
        \Bigl| \mathbb E_k \bigl[ G(Y_{k+1}) - G(\tilde Y_k(h)) \bigr] \Bigr|
        \Le C h^3 \sum_{j=1}^{4} \norm{\nabla^j G}_\infty\, \bigl( 1 + |Y_k| \bigr)^{6} .
        \label{eq: weak}
    \end{equation}
    \emph{(Proof in Appendix~\ref{sec: weak proof}.)}
\end{lemma}

With the stochastic gradient, the gradient noise at $x$ has the covariance
\begin{equation}
    \mathcal C(x) := \mathbb E_\xi\Bigl[ \bigl( \nabla \hat U(x;\xi) - \nabla U(x) \bigr) \otimes \bigl( \nabla \hat U(x;\xi) - \nabla U(x) \bigr) \Bigr] , \qquad \mathrm{Tr}\bigl[ \mathcal C(x) \bigr] \Le \Sigma ,
    \label{eq: covariance}
\end{equation}
where the bound on the trace is Assumption~\ref{ass: sg}. The step \eqref{eq: sg chain step} adds to the step \eqref{eq: chain step} the gradient noise at $Z_0$, multiplied by $-h$. Given $\hat{\mathcal F}_0$ and $\chi_0$ this noise has mean $0$, so the mean of the step is unchanged, while its moment of order $2$ gains $h^2\, \hat{\mathbb E}_0[\mathcal C(Z_0)]$, as for \eqref{eq: sgld}.

Lemma~\ref{lem: sg weak} keeps this gain as an explicit term of order $h^2$, bounded through \eqref{eq: trace bound} by $\frac{\Sigma h^2}{2}\, \norm{\nabla^2 G}_\infty$, against the order $h^3$ of \eqref{eq: weak}; the rest of the gradient noise also contributes only a remainder of order $h^3$.

\begin{lemma}
\label{lem: sg weak}
    Let Assumptions~\ref{ass: U} and \ref{ass: sg} hold and let $G \in C^4(\mathbb R^d)$ have bounded derivatives of orders $1$ to $4$. Then, for the chain \eqref{eq: sglm markov} and every $k \Ge 0$,
    \begin{equation}
        \qquad \biggl| \hat{\mathbb E}_k \bigl[ G(Y_{k+1}) - G(\tilde Y_k(h)) \bigr] - \frac{h^2}{2}\, \nabla^2 G(Y_k) : \hat{\mathbb E}_k\bigl[ \mathcal C(Z_k) \bigr] \biggr| \\
        \Le C h^3 \sum_{j=1}^{4} \norm{\nabla^j G}_\infty\, \bigl( 1 + |Y_k| \bigr)^{6} .
    \label{eq: sg weak}
    \end{equation}
    \emph{(Proof in Appendix~\ref{sec: sg weak proof}.)}
\end{lemma}

\section{Main results on the LM integrators}
\label{sec: mse}

This section states the main results: the MSE bounds for the samples of the \eqref{eq: lm} and \eqref{eq: sglm} integrators in Section~\ref{subsec: mse analysis}, the numerical bias of SGLM in Section~\ref{subsec: numerical bias}, its comparison with the numerical bias of \eqref{eq: sgld} in Section~\ref{subsec: sgld}, and the explicit coefficients of $h$ in $\mathcal W_1(\pi_h, \pi)$ and $\mathcal W_1(\bar\pi_h, \pi)$, given by the Riesz representors, with their comparison over the batch size, in Section~\ref{subsec: w1}. Their proofs are in Appendices~\ref{app: mse}, \ref{app: bias}, \ref{app: sgld} and \ref{app: w1}.

\subsection{MSE analysis of SGLM}
\label{subsec: mse analysis}

The samples of the \eqref{eq: lm} and \eqref{eq: sglm} integrators start from a given $Z_0 \in \mathbb R^d$. The analysis below runs on the Markovian form from a given $Y_0$, with $Z_k = Y_k + \sqrt{h/2}\,\chi_k$ for every $k \Ge 0$; Appendix~\ref{app: mse} applies it to the samples from $Z_0$ by starting the Markovian form at $Y_1 = Z_0 - h\nabla U(Z_0) + \sqrt{h/2}\,\chi_0$, with $\nabla \hat U(Z_0;\xi_0)$ in place of $\nabla U(Z_0)$ for SGLM, so that $Z_k = Y_k + \sqrt{h/2}\,\chi_k$ for every $k \Ge 1$.

With $\rho = \tilde\pi_h$, the MSE satisfies \eqref{eq: mse split modified}, and it remains to bound its time average error. Telescoping $\tilde\phi_h(Y_{k+1}) - \tilde\phi_h(Y_k)$ along the chain of the Markovian form, with the solution $\tilde\phi_h$ of \eqref{eq: modified poisson}, gives the time-averaged error formula of \citep{lu2025mse, ye2026mse},
\begin{equation}
    \frac1N \sum_{k=0}^{N-1} \bigl( \tilde f_h(Y_k) - \tilde\pi_h(\tilde f_h) \bigr)
    = \frac{\tilde\phi_h(Y_0) - \tilde\phi_h(Y_N)}{Nh} + \frac1N \sum_{k=0}^{N-1} \Upsilon_k + \frac1N \sum_{k=0}^{N-1} \Psi_k ,
    \label{eq: time average formula}
\end{equation}
with the discretization term and the martingale difference
\begin{equation}
    \Upsilon_k := \frac{\tilde\phi_h(Y_{k+1}) - \tilde\phi_h(\tilde Y_k(h))}{h} , \qquad
    \Psi_k := \frac{\tilde\phi_h(\tilde Y_k(h)) - \tilde\phi_h(Y_k)}{h} + \tilde f_h(Y_k) - \tilde\pi_h(\tilde f_h) .
    \label{eq: Upsilon Psi}
\end{equation}
The conditional mean of $\Upsilon_k$, under $\mathbb E_k$ for \eqref{eq: lm} and $\hat{\mathbb E}_k$ for \eqref{eq: sglm}, is the local weak error of Lemmas~\ref{lem: weak} and \ref{lem: sg weak} with $G = \tilde\phi_h$, divided by $h$; the lemmas apply because \eqref{eq: modified poisson bound} bounds the derivatives of $\tilde\phi_h$ of orders $1$ to $4$ by $C L$. The term $\Psi_k$ has conditional mean $0$ by \eqref{eq: modified poisson}.

\begin{theorem}
\label{thm: lm mse}
    Let Assumptions~\ref{ass: U} and \ref{ass: f} hold, and let $(Z_k)_{k \Ge 0}$ be the samples of \eqref{eq: lm} from a given $Z_0 \in \mathbb R^d$. Then, for every $h \in (0,h_0]$ and every $N \Ge 1$,
    \begin{equation}
        \MSE(N,h) \Le C L^2 \biggl[ \bigl( 1 + |Z_0| \bigr)^2 \frac{1}{Nh} + \bigl( 1 + |Z_0| \bigr)^{12} h^4 \biggr] .
        \label{eq: lm mse}
    \end{equation}
    \emph{(Proof in Appendix~\ref{sec: lm mse proof}.)}
\end{theorem}

The term $\frac{1}{Nh}$ is the order of the variance of a time average over a trajectory of length $Nh$, and $h^4$ is the square of the bias of order $h^2$ of \eqref{eq: lm}, which comes from Lemma~\ref{lem: shifted bias} and \eqref{eq: weak}. Balancing them gives, for $N \Ge h_0^{-5}$, the step $h = N^{-\frac15} \Le h_0$ and $\MSE(N,h) \Le C L^2 (1+|Z_0|)^{12} N^{-\frac45}$.

For the stochastic gradient, the explicit term of Lemma~\ref{lem: sg weak}, divided by $h$, adds to the conditional mean of $\Upsilon_k$ a term of order $\Sigma h$, whose square gives the order $\Sigma^2 h^2$ in \eqref{eq: sglm mse}.

\begin{theorem}
\label{thm: sglm mse}
    Let Assumptions~\ref{ass: U}, \ref{ass: f} and \ref{ass: sg} hold, and let $(Z_k)_{k \Ge 0}$ be the samples of \eqref{eq: sglm} from a given $Z_0 \in \mathbb R^d$. Then, for every $h \in (0,h_0]$ and every $N \Ge 1$,
    \begin{equation}
        \MSE(N,h) \Le C L^2 \biggl[ \bigl( 1 + |Z_0| \bigr)^2 \Bigl( \frac{1}{Nh} + \Sigma^2 h^2 \Bigr) + \bigl( 1 + |Z_0| \bigr)^{12} h^4 \biggr] .
        \label{eq: sglm mse}
    \end{equation}
    \emph{(Proof in Appendix~\ref{sec: sglm mse proof}.)}
\end{theorem}

With a batch of $P$ independent stochastic gradients, the chain \eqref{eq: sglm} reads
\begin{equation*}
    Z_{k+1} = Z_k - \frac hP \sum_{p=1}^{P} \nabla \hat U\bigl(Z_k; \xi_k^{(p)}\bigr) + \sqrt{\frac h2}\, \bigl( \chi_k + \chi_{k+1} \bigr) ,
\end{equation*}
with $\xi_k^{(1)}, \dots, \xi_k^{(P)}$ independent and drawn afresh at every step. By Minkowski's inequality, the averaged estimator keeps the bound $M(1 + |x|)$ on its $L^{12}$ norm and the bound on its Hessian, uniformly in $P$, and its variance is at most $\Sigma/P$, so it satisfies Assumption~\ref{ass: sg} with $\Sigma/P$ in place of $\Sigma$; the proof of Theorem~\ref{thm: sglm mse} is unchanged, and \eqref{eq: sglm mse} holds with $\Sigma^2 h^2 / P^2$ in place of $\Sigma^2 h^2$. 

\subsection{Numerical bias estimate of SGLM}
\label{subsec: numerical bias}

The MSE bounds control the time average along one trajectory. This subsection identifies the long-time limit of that average for \eqref{eq: sglm}, and its distance to $\pi$, the numerical bias. The Markovian form \eqref{eq: sglm markov} is ergodic by Harris's ergodic theorem \citep{hairer2008yet}, applied with the drift condition behind Lemma~\ref{lem: moments} and a minorization that uses the Hessian bound of Assumption~\ref{ass: sg}.

\begin{lemma}
\label{lem: invariant}
    Let Assumptions~\ref{ass: U} and \ref{ass: sg} hold. For every $h \in (0,h_0]$, \eqref{eq: sglm markov} has a unique invariant distribution $\mu_h$, with $\mu_h(|y|^{12}) \Le C$. Let $\pi_h := \mu_h * \mathcal N(0, \frac h2 I_d)$. Then, for every $f$ satisfying Assumption~\ref{ass: f} and every given $Z_0 \in \mathbb R^d$, the samples of \eqref{eq: sglm} satisfy
    \begin{equation}
        \lim_{N \to \infty} \frac1N \sum_{k=0}^{N-1} f(Z_k) = \pi_h(f) \qquad \text{almost surely} .
        \label{eq: ergodic limit}
    \end{equation}
    \emph{(Proof in Appendix~\ref{sec: invariant proof}.)}
\end{lemma}

The measure $\pi_h$ is determined by the chain \eqref{eq: sglm markov} alone, through $\mu_h$, and does not depend on the initial value $Z_0$. Since \eqref{eq: sglm} is not Markovian, it has no invariant distribution in the usual sense; by \eqref{eq: ergodic limit}, the measure $\pi_h$ can be viewed as the numerical invariant distribution of SGLM.

The next theorem expands the numerical bias $\pi_h(f) - \pi(f)$ to the order $h$, with a coefficient given by $\tilde\phi_h$ and the covariance of the gradient noise, defined in \eqref{eq: covariance} as
\begin{equation*}
    \mathcal C(x) = \mathbb E_\xi\Bigl[ \bigl( \nabla \hat U(x;\xi) - \nabla U(x) \bigr) \otimes \bigl( \nabla \hat U(x;\xi) - \nabla U(x) \bigr) \Bigr] , \qquad x \in \mathbb R^d .
\end{equation*}
The bias estimate follows the same strategy as that of SG-UBU in \citep[Section~3.4]{lu2025mse}. Since $\mathcal C(x)$ is positive semidefinite, Assumption~\ref{ass: sg} gives $\norm{\mathcal C(x)} \Le \mathrm{Tr}\bigl[ \mathcal C(x) \bigr] \Le \Sigma$.

\begin{theorem}
\label{thm: numerical bias}
    Let Assumptions~\ref{ass: U}, \ref{ass: f} and \ref{ass: sg} hold. For every $h \in (0,h_0]$, with $Y_0 \sim \mu_h$, $\chi_0 \sim \mathcal N(0,I_d)$ independent of $Y_0$, and $Z_0 = Y_0 + \sqrt{h/2}\,\chi_0$, so that $Z_0 \sim \pi_h$, the numerical invariant distribution of \eqref{eq: sglm},
    \begin{equation}
        \Bigl| \pi_h(f) - \pi(f) - \frac h2\, \mathbb E\bigl[ \nabla^2 \tilde\phi_h(Y_0) : \mathcal C(Z_0) \bigr] \Bigr| \Le C L h^2 .
        \label{eq: numerical bias}
    \end{equation}
    In particular, $|\pi_h(f) - \pi(f)| \Le C L ( \Sigma h + h^2 )$. \emph{(Proof in Appendix~\ref{sec: numerical bias proof}.)}
\end{theorem}

The bound $C L (\Sigma h + h^2)$ of Theorem~\ref{thm: numerical bias} follows from the uniform bound \eqref{eq: modified poisson bound} on the modified discrete Poisson solution $\tilde\phi_h$, together with $\norm{\mathcal C(x)} \Le \Sigma$. As $h \to 0$, the coefficient $\frac12\, \mathbb E\bigl[ \nabla^2 \tilde\phi_h(Y_0) : \mathcal C(Z_0) \bigr]$ of $h$ in \eqref{eq: numerical bias} converges to an expression in the solution $\phi_0$ of the continuous Poisson equation \eqref{eq: continuous poisson}.

\begin{corollary}
\label{cor: bias coefficient}
    Let Assumptions~\ref{ass: U}, \ref{ass: f} and \ref{ass: sg} hold. Then, for \eqref{eq: sglm},
    \begin{equation}
        \lim_{h \to 0} \frac{\pi_h(f) - \pi(f)}{h} = \frac12\, \pi\bigl( \nabla^2 \phi_0(x) : \mathcal C(x) \bigr) .
        \label{eq: bias coefficient}
    \end{equation}
    \emph{(Proof in Appendix~\ref{sec: bias coefficient proof}.)}
\end{corollary}

Corollary~\ref{cor: bias coefficient} gives a rigorous characterization of the linear part of the numerical bias of \eqref{eq: sglm}: its coefficient of order $h$ comes entirely from the covariance $\mathcal C(x)$ of the stochastic gradient, and vanishes with the full gradient. This is intrinsically different from \eqref{eq: sgld}, whose numerical bias carries a first order discretization error even with the full gradient, so that its coefficient of order $h$ remains nonzero for some test functions as the variance of the stochastic gradient tends to zero (Corollary~\ref{cor: w1 limits}). For SG-UBU, the corresponding limit is identified only formally in \citep{lu2025mse}; for SGLM, Corollary~\ref{cor: bias coefficient} makes the coefficient explicit and proves the limit.

Linear response theory \citep{kubo1957statistical, hairer2010simple} also explains the coefficient \eqref{eq: bias coefficient}. To first order in $h$, the gradient noise of \eqref{eq: sglm} adds the diffusion $\frac h2\, \mathcal C : \nabla^2$ to the generator \eqref{eq: generator}, and the invariant distribution responds to this perturbation by $\frac h2\, \pi(\nabla^2 \phi_0 : \mathcal C)$ for the test function $f$. The invariant distributions of Langevin integrators are viewed in the same way, as perturbations of $\pi$ depending on the step size, in \citep{leimkuhler2015computation}. With the solution of the continuous Poisson equation \eqref{eq: continuous poisson} written as $\phi_0 = \int_0^\infty e^{t\mathcal L} \bigl( f - \pi(f) \bigr) \D t$, two integrations by parts give the Green--Kubo form \citep{green1954markoff, kubo1957statistical} of the coefficient,
\begin{equation*}
    \frac12\, \pi\bigl( \nabla^2 \phi_0 : \mathcal C \bigr) = \frac12 \int_0^\infty \pi\Bigl( \pi^{-1}
     \nabla^2 : (\mathcal C \pi)\; e^{t\mathcal L} \bigl( f - \pi(f) \bigr) \Bigr) \D t , \qquad \nabla^2 : (\mathcal C \pi) = \sum_{i,j=1}^d \partial_i \partial_j \bigl( \mathcal C_{ij}\, \pi \bigr) ,
\end{equation*}
that is, the time integral of the correlation between $\pi^{-1} \nabla^2 : (\mathcal C \pi)$ at the initial value and $f$ at the time $t$, along \eqref{eq: langevin} started from $\pi$.

\subsection{Numerical bias comparison with SGLD}
\label{subsec: sgld}

The MSE analysis through the discrete Poisson equation applies to \eqref{eq: sgld} without extra effort: the samples of SGLD form a Markov chain, so no modified class is needed, and the time-averaged error formula \eqref{eq: time average formula} holds with $Z_k$, $f$, $\pi(f)$, the discrete Poisson solution $\phi_h$ of \eqref{eq: discrete poisson solution} and the solution $X_h$ of \eqref{eq: langevin} from $X_0 = Z_k$ in place of $Y_k$, $\tilde f_h$, $\tilde\pi_h(\tilde f_h)$, $\tilde\phi_h$ and $\tilde Y_k(h)$. For the samples of SGLD from a given $Z_0 \in \mathbb R^d$, it gives
\begin{equation}
    \MSE(N,h) = \mathcal O\Bigl( \frac{1}{Nh} + \Sigma^2 h^2 + h^2 \Bigr) ,
    \label{eq: sgld mse}
\end{equation}
with the constant depending also on $L$ and $Z_0$. The term $\Sigma^2 h^2$ comes from the gradient noise, as in \eqref{eq: sglm mse}, and the last term $h^2$ from the discretization error of EM, where \eqref{eq: sglm mse} has $h^4$. Since $\Sigma \Le M^2$, \eqref{eq: sgld mse} reduces to $\MSE(N,h) = \mathcal O\bigl( (Nh)^{-1} + h^2 \bigr)$, the bound already proved through the continuous Poisson equation \eqref{eq: continuous poisson} \citep[Theorem~9]{vollmer2016exploration}, \citep[Theorem~3]{chen2015convergence}, and we do not restate its proof.

The arguments of Lemma~\ref{lem: invariant}, based on Harris's ergodic theorem, show that \eqref{eq: sgld} has a unique invariant distribution $\bar\pi_h$, with $\bar\pi_h(|x|^{12}) \Le C$, and that the time averages of $f$ along SGLD converge almost surely to $\bar\pi_h(f)$ (Lemma~\ref{lem: sgld invariant} in Appendix~\ref{app: sgld}). The invariant distribution $\bar\pi_h$ differs from the numerical invariant distribution $\pi_h$ of \eqref{eq: sglm} in the coefficient of order $h$ of the numerical bias. In the next theorem, this coefficient carries, besides the covariance $\mathcal C(x)$ of the gradient noise, the one-step error of EM against \eqref{eq: langevin}, through $\nabla^2 U\, \nabla U - \nabla \Delta U$ in the mean of the step and the Hessian $\nabla^2 U$ in its moment of order $2$.

\begin{theorem}
\label{thm: sgld bias}
    Let Assumptions~\ref{ass: U}, \ref{ass: f} and \ref{ass: sg} hold. For every $h \in (0,h_0]$, with $Z_0 \sim \bar\pi_h$, the invariant distribution of \eqref{eq: sgld},
    \begin{equation}
    \begin{aligned}
        & \biggl| \bar\pi_h(f) - \pi(f) - \frac h2\, \mathbb E\Bigl[ \nabla^2 \phi_h(Z_0) : \bigl( 2\nabla^2 U(Z_0) + \mathcal C(Z_0) \bigr) \Bigr] \\
        & \quad + \frac h2\, \mathbb E\Bigl[ \nabla \phi_h(Z_0) \cdot \bigl( \nabla^2 U(Z_0)\, \nabla U(Z_0) - \nabla \Delta U(Z_0) \bigr) \Bigr] \biggr| \Le C L h^2 .
    \end{aligned}
        \label{eq: sgld bias}
    \end{equation}
    In particular, $|\bar\pi_h(f) - \pi(f)| \Le C L h$. \emph{(Proof in Appendix~\ref{sec: sgld bias proof}.)}
\end{theorem}

As $h \to 0$, the coefficient of $h$ in \eqref{eq: sgld bias} converges to an expression in the solution $\phi_0$ of \eqref{eq: continuous poisson}, in which only the Hessian $\nabla^2 \phi_0$ remains.

\begin{corollary}
\label{cor: sgld coefficient}
    Let Assumptions~\ref{ass: U}, \ref{ass: f} and \ref{ass: sg} hold. Then, for \eqref{eq: sgld},
    \begin{equation}
        \lim_{h \to 0} \frac{\bar\pi_h(f) - \pi(f)}{h} = \frac12\, \pi\Bigl( \nabla^2 \phi_0(x) : \bigl( \nabla^2 U(x) + \mathcal C(x) \bigr) \Bigr) .
        \label{eq: sgld coefficient}
    \end{equation}
    \emph{(Proof in Appendix~\ref{sec: sgld coefficient proof}.)}
\end{corollary}

Corollary~\ref{cor: sgld coefficient} is not entirely new: the coefficient of order $h$ of the numerical bias of \eqref{eq: sgld} has been studied formally in \citep[Section~4.1]{vollmer2016exploration}, for $d = 1$ and under the assumption that the chain is ergodic, where the gradient noise adds to the coefficient of EM a term in the variance of the stochastic gradient, and the modified SGLD of \citep[Section~4.2]{vollmer2016exploration} removes this term. Corollary~\ref{cor: sgld coefficient} proves the limit in $\mathbb R^d$ and writes the coefficient in terms of the solution $\phi_0$ of \eqref{eq: continuous poisson}. The \eqref{eq: sglm} integrator treats the discretization bias in a different and clean way: it preserves the second order discretization error of the LM integrator, so its coefficient \eqref{eq: bias coefficient} carries the gradient noise only.

\begin{remark}
Corollary~\ref{cor: sgld coefficient} also requires the stochastic gradient to have a finite moment of high order, as in Assumption~\ref{ass: sg}. For \eqref{eq: sgld} with a strongly convex potential, \citep[Theorem~1]{paulin2026sgld} bounds the numerical bias in $\mathcal W_2$ by a term of order $h$ if the moment of order $4$ of the gradient noise under $\pi$ is finite; if only the moment of order $2$ is finite, the bound could deteriorate to the order $h^{\frac12}$.
\end{remark}

The coefficients \eqref{eq: bias coefficient} of \eqref{eq: sglm} and \eqref{eq: sgld coefficient} of \eqref{eq: sgld} differ by $\frac12\, \pi\bigl( \nabla^2 \phi_0(x) : \nabla^2 U(x) \bigr)$, the contribution of the discretization error of EM. This term does not depend on the gradient noise: it is the coefficient of SGLD with the full gradient, and it remains at every batch size. The coefficient of SGLM, by contrast, vanishes with $\mathcal C(x)$, and with a batch of $P$ independent stochastic gradients, whose average has the covariance $\mathcal C(x)/P$, it is divided by $P$. For a given test function, the two coefficients are not ordered in absolute value, since $\nabla^2 \phi_0$ is indefinite in general and the terms in $\nabla^2 U(x)$ and in $\mathcal C(x)$ can have opposite signs. Section~\ref{subsec: w1} compares the suprema of their absolute values over the test functions with $\norm{\nabla f}_\infty \Le 1$, and Lemma~\ref{lem: w1 counterexample} gives a gradient noise for which even the supremum of the absolute value of the coefficient of SGLD over these test functions is smaller than that of SGLM.

The expectation in \eqref{eq: sgld bias} contains the gradient $\nabla \phi_h$, through the error of EM in the mean of the step, while the limit \eqref{eq: sgld coefficient} of Corollary~\ref{cor: sgld coefficient} contains the Hessian $\nabla^2 \phi_0$ only. The two forms are connected by an integration by parts that holds under $\pi$ and is not available under $\bar\pi_h$. By \eqref{eq: generator}, for each component $i$,
\begin{equation*}
    \mathcal L(\partial_i U) = \Delta \partial_i U - \nabla U \cdot \nabla \partial_i U = \bigl( \nabla \Delta U - \nabla^2 U\, \nabla U \bigr)_i ,
\end{equation*}
so the vector multiplying $\nabla \phi_h$ in \eqref{eq: sgld bias} is, componentwise, minus the generator applied to $\nabla U$,
\begin{equation}
    \nabla^2 U\, \nabla U - \nabla \Delta U = -\mathcal L \nabla U .
    \label{eq: mean error generator}
\end{equation}
For every $\psi \in C^3(\mathbb R^d)$ with bounded derivatives of orders $1$ to $3$, \eqref{eq: generator ibp} with $g_1 = \partial_i \psi$ and $g_2 = \partial_i U$, summed over $i$, and \eqref{eq: mean error generator} give
\begin{equation}
\begin{aligned}
    & \pi\Bigl( \nabla \psi(x) \cdot \bigl( \nabla^2 U(x)\, \nabla U(x) - \nabla \Delta U(x) \bigr) \Bigr) = -\sum_{i=1}^{d} \pi\bigl( \partial_i \psi(x)\, \mathcal L \partial_i U(x) \bigr) \\
    & \quad = \sum_{i=1}^{d} \pi\bigl( \nabla \partial_i \psi(x) \cdot \nabla \partial_i U(x) \bigr) = \pi\bigl( \nabla^2 \psi(x) : \nabla^2 U(x) \bigr) .
\end{aligned}
    \label{eq: sgld integration by parts}
\end{equation}

At a fixed step size, the expectation in \eqref{eq: sgld bias} is taken under $\bar\pi_h$, while \eqref{eq: sgld integration by parts} with $\psi = \phi_h$ holds under $\pi$. Moving the terms of order $h$ from $\bar\pi_h$ to $\pi$ within the remainder $C L h^2$ requires Theorem~\ref{thm: sgld bias} for a test function built from $\nabla \phi_h$ and $\nabla^2 \phi_h$, whose derivatives of orders $1$ to $4$ involve the derivatives of $\phi_h$ up to order $6$, beyond the order $4$ of Theorem~\ref{thm: u decay}. The gradient $\nabla \phi_h$ therefore stays in \eqref{eq: sgld bias}.

As $h \to 0$, $\bar\pi_h$ converges weakly to $\pi$ with moments bounded uniformly in $h$, and $\nabla \phi_h$ and $\nabla^2 \phi_h$ converge to $\nabla \phi_0$ and $\nabla^2 \phi_0$ uniformly on compact sets, so the expectation in \eqref{eq: sgld bias} converges to the same expression under $\pi$. There \eqref{eq: sgld integration by parts} with $\psi = \phi_0$ applies, and the term in $\nabla \phi_0$ cancels half of the term in $2 \nabla^2 U$. Theorem~\ref{thm: sgld bias} thus holds at every step size $h \in (0,h_0]$ with a remainder of order $h^2$ but keeps both $\nabla \phi_h$ and $\nabla^2 \phi_h$, while the Hessian form of the coefficient in Corollary~\ref{cor: sgld coefficient} is reached only in the limit $h \to 0$, without a rate.

The gap between Theorem~\ref{thm: sgld bias} and Corollary~\ref{cor: sgld coefficient} also reveals a feature of the discrete Poisson equation. The solution $\phi_h$ of \eqref{eq: discrete poisson} compares one step of \eqref{eq: sgld} with one step of \eqref{eq: langevin}, so only the differences of their moments enter, and the derivatives of $\phi_h$ up to order $4$ suffice for the nonasymptotic bound \eqref{eq: sgld bias}. The solution $\phi_0$ of \eqref{eq: continuous poisson} compares one step of SGLD with the generator $\mathcal L$ instead: the coefficient of order $h$ then involves the derivatives of $\phi_0$ up to order $4$, and a remainder of order $h^2$ requires those up to order $6$, beyond Theorem~\ref{thm: u decay}. With $\phi_0$ alone, the terms that do not cancel at a fixed step size are therefore hard to characterize under the present assumptions.

\subsection{Explicit coefficient in Wasserstein-1 error}
\label{subsec: w1}

The coefficients \eqref{eq: bias coefficient} and \eqref{eq: sgld coefficient} depend on the test function $f$ through $\nabla^2 \phi_0$. The suprema of their absolute values over the test functions with $\norm{\nabla f}_\infty \Le 1$ are the limits of $\mathcal W_1(\pi_h, \pi)/h$ and $\mathcal W_1(\bar\pi_h, \pi)/h$ as $h \to 0$ (Corollary~\ref{cor: w1 limits}). To evaluate these suprema, both coefficients are written as $\pi(f \psi)$, with a function $\psi$ that does not depend on $f$. The construction of $\psi$ uses the Poincar\'e inequality of $\pi$.

\begin{lemma}
\label{lem: poincare}
    Let Assumption~\ref{ass: U} hold. There is a constant $C_P$, depending only on $m$, $M$, $R$ and $d$, such that
    \begin{equation}
        \pi\bigl( (g - \pi(g))^2 \bigr) \Le C_P\, \pi\bigl( |\nabla g|^2 \bigr)
        \label{eq: poincare}
    \end{equation}
    for every $g \in C^1(\mathbb R^d)$ with a bounded gradient.\\[4pt]
    \emph{(Directly derived from \citep[Corollary~1.6 and Theorem~1.4]{bakry2008simple}; see Appendix~\ref{sec: riesz well-posedness} for further discussion.)}
\end{lemma}

For a bounded measurable function $\mathcal A : \mathbb R^d \to \mathbb R^{d \times d}$ such that $\mathcal A(x)$ is a symmetric matrix for every $x \in \mathbb R^d$, define the Riesz representor of $\mathcal A$ as the function $\psi_{\mathcal A} \in L^2(\pi)$ with $\pi(\psi_{\mathcal A}) = 0$ and
\begin{equation}
    \pi\bigl( \psi_{\mathcal A}(x)\, (-\mathcal L \omega)(x) \bigr) = \pi\bigl( \nabla^2 \omega(x) : \mathcal A(x) \bigr)
    \label{eq: riesz}
\end{equation}
for every $\omega \in C^4(\mathbb R^d)$ with bounded derivatives of orders $1$ to $4$. The next theorem shows that the Riesz representor is well defined.

\begin{theorem}
\label{thm: riesz}
    Let Assumption~\ref{ass: U} hold, and let $\mathcal A : \mathbb R^d \to \mathbb R^{d \times d}$ be bounded and measurable, with $\mathcal A(x)$ a symmetric matrix for every $x \in \mathbb R^d$. There is a unique $\psi_{\mathcal A} \in L^2(\pi)$ with $\pi(\psi_{\mathcal A}) = 0$ that satisfies \eqref{eq: riesz} for every $\omega \in C^4(\mathbb R^d)$ with bounded derivatives of orders $1$ to $4$, the map $\mathcal A \mapsto \psi_{\mathcal A}$ is linear, and
    \begin{equation}
        \pi\bigl( \psi_{\mathcal A}^2 \bigr) \Le d\, (1 + M C_P) \norm{\mathcal A}_\infty^2 ,
        \label{eq: riesz bound}
    \end{equation}
    with $C_P$ the constant of Lemma~\ref{lem: poincare}.
    \emph{(Proof in Appendix~\ref{sec: riesz well-posedness}.)}
\end{theorem}

For $f$ satisfying Assumption~\ref{ass: f}, \eqref{eq: riesz} with $\omega = \phi_0$ gives $\pi(\nabla^2 \phi_0 : \mathcal A) = \pi(f \psi_{\mathcal A})$, so the coefficients \eqref{eq: bias coefficient} and \eqref{eq: sgld coefficient} read
\begin{equation}
    \lim_{h \to 0} \frac{\pi_h(f) - \pi(f)}{h} = \frac12\, \pi( f \psi_{\mathcal C} ) , \qquad
    \lim_{h \to 0} \frac{\bar\pi_h(f) - \pi(f)}{h} = \frac12\, \pi\bigl( f\, ( \psi_{\nabla^2 U} + \psi_{\mathcal C} ) \bigr) .
    \label{eq: coefficients riesz}
\end{equation}

The representor of the Hessian $\nabla^2 U$ is explicit. The vector $\nabla^2 U\, \nabla U - \nabla \Delta U$ in \eqref{eq: mean error generator} is the gradient of $\frac12 |\nabla U|^2 - \Delta U$, so \eqref{eq: sgld integration by parts} and \eqref{eq: generator ibp} give
\begin{equation*}
    \pi\bigl( \nabla^2 \omega : \nabla^2 U \bigr) = \pi\Bigl( \nabla \omega \cdot \nabla \Bigl( \frac12 |\nabla U|^2 - \Delta U \Bigr) \Bigr) = \pi\Bigl( \Bigl( \frac12 |\nabla U|^2 - \Delta U \Bigr) (-\mathcal L \omega) \Bigr) .
\end{equation*}
Since $\pi(\mathcal L \omega) = 0$, the mean of $\frac12 |\nabla U|^2 - \Delta U$ may be subtracted; it equals $-\frac12 \pi(\Delta U)$, since $\pi(\mathcal L U) = 0$ gives $\pi(|\nabla U|^2) = \pi(\Delta U)$. Hence \eqref{eq: riesz} holds for $\mathcal A = \nabla^2 U$ with
\begin{equation}
    \psi_{\nabla^2 U}(x) = \frac12\, |\nabla U(x)|^2 - \Delta U(x) + \frac12\, \pi(\Delta U) .
    \label{eq: psi hessian}
\end{equation}

Through the representors, the normalized biases $(\pi_h - \pi)/h$ and $(\bar\pi_h - \pi)/h$ converge to the signed measures $\frac12 \psi_{\mathcal C}\, \pi$ and $\frac12 (\psi_{\nabla^2 U} + \psi_{\mathcal C})\, \pi$ in the Kantorovich--Rubinstein norm, at the rate $h^{\frac12} (1 + \log\frac1h)$.

\begin{theorem}
\label{thm: kr expansion}
    Let Assumptions~\ref{ass: U} and \ref{ass: sg} hold. Then, for every $h \in (0,h_0]$, the numerical invariant distribution $\pi_h$ of \eqref{eq: sglm} and the invariant distribution $\bar\pi_h$ of \eqref{eq: sgld} satisfy
    \begin{align}
        & \Bigl\| \pi_h - \pi - \frac h2\, \psi_{\mathcal C}\, \pi \Bigr\|_{\mathrm{KR}} \Le C h^{\frac32} \Bigl( 1 + \log\frac1h \Bigr) ,
        \label{eq: kr sglm} \\
        & \Bigl\| \bar\pi_h - \pi - \frac h2\, ( \psi_{\nabla^2 U} + \psi_{\mathcal C} )\, \pi \Bigr\|_{\mathrm{KR}} \Le C h^{\frac32} \Bigl( 1 + \log\frac1h \Bigr) .
        \label{eq: kr sgld}
    \end{align}
    \emph{(Proof in Appendix~\ref{sec: kr expansion proof}.)}
\end{theorem}

Since $\mathcal W_1(\pi_h, \pi) = \norm{\pi_h - \pi}_{\mathrm{KR}}$ and $\mathcal W_1(\bar\pi_h, \pi) = \norm{\bar\pi_h - \pi}_{\mathrm{KR}}$, the triangle inequality for $\norm{\cdot}_{\mathrm{KR}}$ and Theorem~\ref{thm: kr expansion} give, for every $h \in (0,h_0]$,
\begin{equation}
    \biggl| \mathcal W_1(\pi_h, \pi) - \frac h2\, \norm{\psi_{\mathcal C}\, \pi}_{\mathrm{KR}} \biggr|
    + \biggl| \mathcal W_1(\bar\pi_h, \pi) - \frac h2\, \norm{( \psi_{\nabla^2 U} + \psi_{\mathcal C} )\, \pi}_{\mathrm{KR}} \biggr|
    \Le C h^{\frac32} \Bigl( 1 + \log\frac1h \Bigr) .
    \label{eq: w1 limits}
\end{equation}
The bound \eqref{eq: w1 limits} determines the Wasserstein-1 errors to first order in $h$, but it is not guaranteed to be sharp in the discretization error: its remainder is only of order $h^{\frac32}$, up to the factor $1 + \log\frac1h$. Enhancing this order to $h^2$ is left for future studies; Remark~\ref{rem: h2} discusses the globally strongly convex case.

\begin{remark}
\label{rem: h2}
When $U$ is globally strongly convex, that is, $R = 0$ in Assumption~\ref{ass: U}, and the covariance $\mathcal C(x)$ is smooth, we expect the remainder of \eqref{eq: kr sglm}, \eqref{eq: kr sgld} and \eqref{eq: w1 limits} to improve to $C h^2$. The factor $h^{\frac12} (1 + \log\frac1h)$ is lost where a step of the chain acts on a Lipschitz test function and where $\psi_{\mathcal C}$ is only known to lie in $L^2(\pi)$. Under global strong convexity, the derivatives of the solution of \eqref{eq: langevin} with respect to its initial point decay at the rate $e^{-mt}$, so $\psi_{\mathcal C}$, like the explicit $\psi_{\nabla^2 U}$ of \eqref{eq: psi hessian}, should be smooth with polynomially growing derivatives. The steps of the chain can then be applied to the smooth measures $\psi_{\mathcal C}\, \pi$ and $\psi_{\nabla^2 U}\, \pi$ instead of the test function, and each of these losses would become a factor $h$.
\end{remark}

Since the remainder in \eqref{eq: w1 limits} is $o(h)$, dividing \eqref{eq: w1 limits} by $h$ and letting $h \to 0$ gives the explicit limits of $\mathcal W_1(\pi_h, \pi)/h$ and $\mathcal W_1(\bar\pi_h, \pi)/h$.

\begin{corollary}
\label{cor: w1 limits}
    Let Assumptions~\ref{ass: U} and \ref{ass: sg} hold. Then the numerical invariant distribution $\pi_h$ of \eqref{eq: sglm} and the invariant distribution $\bar\pi_h$ of \eqref{eq: sgld} satisfy
    \begin{equation*}
        \lim_{h \to 0} \frac{\mathcal W_1(\pi_h, \pi)}{h} = \frac12\, \norm{\psi_{\mathcal C}\, \pi}_{\mathrm{KR}} , \qquad
        \lim_{h \to 0} \frac{\mathcal W_1(\bar\pi_h, \pi)}{h} = \frac12\, \norm{( \psi_{\nabla^2 U} + \psi_{\mathcal C} )\, \pi}_{\mathrm{KR}} .
    \end{equation*}
    Moreover, $\norm{\psi_{\nabla^2 U}\, \pi}_{\mathrm{KR}} > 0$, and $\norm{\psi_{\mathcal C}\, \pi}_{\mathrm{KR}} > 0$ if $\mathcal C(x)$ does not vanish identically. \emph{(Proof in Appendix~\ref{sec: w1 limits proof}.)}
\end{corollary}

Corollary~\ref{cor: w1 limits} is nontrivial: it gives an accurate characterization of the coefficient of $h$ in the Wasserstein-1 error, as an equality. This makes it distinct from the mainstream coupling methods of numerical analysis, which bound the Wasserstein-1 error from above without identifying its coefficient. Although the representor $\psi_{\mathcal C}$ is difficult to evaluate in practice, Corollary~\ref{cor: w1 limits} and \eqref{eq: w1 limits} are the decisive answer to how much SGLM and SGLD shift the invariant distribution: up to the remainder $C h^{\frac32} (1 + \log\frac1h)$ of \eqref{eq: w1 limits}, the shift in the Wasserstein-1 distance is $\frac h2\, \norm{\psi_{\mathcal C}\, \pi}_{\mathrm{KR}}$ for SGLM and $\frac h2\, \norm{( \psi_{\nabla^2 U} + \psi_{\mathcal C} )\, \pi}_{\mathrm{KR}}$ for SGLD.

Corollary~\ref{cor: w1 limits} thus compares SGLM with SGLD through the representors: the limit of $\mathcal W_1(\pi_h, \pi)/h$ is at most that of $\mathcal W_1(\bar\pi_h, \pi)/h$ exactly when
\begin{equation}
    \norm{( \psi_{\nabla^2 U} + \psi_{\mathcal C} )\, \pi}_{\mathrm{KR}} \Ge \norm{\psi_{\mathcal C}\, \pi}_{\mathrm{KR}} .
    \label{eq: w1 comparison}
\end{equation}
This inequality fails in general: Lemma~\ref{lem: w1 counterexample} gives a gradient noise of large amplitude, supported away from the minimizer of $U$, that violates it.

\begin{lemma}
\label{lem: w1 counterexample}
    Let $d = 1$, $U(x) = \frac12 x^2$ and $\nabla \hat U(x;\xi) = x + 10\, \xi\, a(x)$, where $\xi$ takes the values $1$ and $-1$ with probability $\frac12$ each and
    \begin{equation}
        a(x) := \exp\Bigl( 1 - \frac{1}{1 - 4 (|x| - 2)^2} \Bigr) \quad \text{if } \bigl| |x| - 2 \bigr| < \frac12 , \qquad a(x) := 0 \quad \text{otherwise} .
        \label{eq: shell noise}
    \end{equation}
    Then Assumptions~\ref{ass: U} and \ref{ass: sg} hold, and
    \begin{equation*}
        \frac12\, \norm{\psi_{\mathcal C}\, \pi}_{\mathrm{KR}} \approx 1.230 > 1.115 \approx \frac12\, \norm{( \psi_{\nabla^2 U} + \psi_{\mathcal C} )\, \pi}_{\mathrm{KR}} ,
    \end{equation*}
    so \eqref{eq: w1 comparison} fails. \emph{(Proof in Appendix~\ref{sec: w1 counterexample proof}.)}
\end{lemma}

With a batch of $P$ independent stochastic gradients, $\mathcal C(x)/P$ replaces the covariance $\mathcal C(x)$ of a single stochastic gradient. As $P$ grows, the limit of $\mathcal W_1(\pi_h, \pi)/h$ tends to $0$, while that of $\mathcal W_1(\bar\pi_h, \pi)/h$ tends to $\frac12 \norm{\psi_{\nabla^2 U}\, \pi}_{\mathrm{KR}} > 0$; for a strongly convex potential, \eqref{eq: w1 comparison} holds once $P$ reaches an explicit threshold.

\begin{theorem}
\label{thm: batch threshold}
    Let Assumptions~\ref{ass: U} with $R = 0$ and \ref{ass: sg} hold, so that $\nabla^2 U(x) \succeq m I_d$ for every $x \in \mathbb R^d$, and let $\mathcal C(x)$ be the covariance of a single stochastic gradient. Then, for every batch size $P$ satisfying
    \begin{equation}
        P \Ge \frac{9\, \norm{\nabla^2 U}_\infty\, \norm{\mathcal C}_\infty}{m^2} = \frac{9}{m^2} \sup_{x\in\mathbb R^d} \norm{\nabla^2U(x)} \sup_{x\in\mathbb R^d} \norm{\mathcal C(x)},
        \label{eq: batch threshold}
    \end{equation}
    the limit of $\mathcal W_1(\pi_h, \pi)/h$ for SGLM with a batch of size $P$ is at most that of $\mathcal W_1(\bar\pi_h, \pi)/h$ for SGLD with the same batch, which by Corollary~\ref{cor: w1 limits}, with $\psi_{\mathcal C}/P$ in place of $\psi_{\mathcal C}$, reads
    \begin{equation}
        \Big\|\Big( \psi_{\nabla^2 U} + \frac1P\psi_{\mathcal C} \Big) \pi\Big\|_{\mathrm{KR}} \Ge \frac1P\, \norm{\psi_{\mathcal C}\, \pi}_{\mathrm{KR}} .
        \label{eq: w1 comparison batch}
    \end{equation}
    \emph{(Proof in Appendix~\ref{sec: batch threshold proof}.)}
\end{theorem}

The norm $\norm{\cdot}$ in \eqref{eq: batch threshold} is the operator norm of Section~\ref{subsec: notation}, so the lower bound of $P$ does not depend explicitly on the dimension $d$; it depends only on the eigenvalues of $\nabla^2 U(x)$ and $\mathcal C(x)$.

\section{Bias correction via the Riesz representors}
\label{sec: representors}

By \eqref{eq: coefficients riesz}, the numerical bias $\pi_h(f) - \pi(f)$ of SGLM equals $\frac h2\, \pi(f \psi_{\mathcal C})$ to first order in $h$, so it is removed once $\pi(f \psi_{\mathcal C})$ is computed. Subtracting $\frac h2\, \pi(f \psi_{\mathcal C})$ from $\pi_h(f)$ requires the value of $\pi(f \psi_{\mathcal C})$ for the given test function $f$, not the representor $\psi_{\mathcal C}$ itself. This section describes three approaches to bias correction. For the covariances of Section~\ref{subsec: explicit representors}, $\pi(f \psi_{\mathcal C})$ is computed with the explicit representor $\psi_{\mathcal C}$. The Richardson--Romberg extrapolation of Section~\ref{subsec: extrapolation} requires only the time averages at the step sizes $h$ and $2h$. The shifting correction of Section~\ref{subsec: shifting} improves the samples themselves by solving the function $\varphi_{\mathcal C}$ from a variational form. Section~\ref{subsec: kr computation} describes the computation of the Kantorovich--Rubinstein norms that give the coefficients of the Wasserstein-1 errors. The three approaches also apply to SGLD, with $\bar\pi_h$ and $\psi_{\nabla^2 U} + \psi_{\mathcal C}$ in place of $\pi_h$ and $\psi_{\mathcal C}$ by \eqref{eq: coefficients riesz}, where $\psi_{\nabla^2 U}$ is explicit by \eqref{eq: psi hessian}; we focus on SGLM for simplicity.

\subsection{Covariances with an explicit representor}
\label{subsec: explicit representors}

Under Assumptions~\ref{ass: U} and \ref{ass: sg}, define the vector function $b_{\mathcal C}(x)$ by
\begin{equation}
    b_{\mathcal C}(x) := \mathcal C(x)\, \nabla U(x) - \operatorname{div} \mathcal C(x) ,
    \label{eq: representor field}
\end{equation}
with $(\operatorname{div} \mathcal C)_j = \sum_{i=1}^d \partial_i \mathcal C_{ij}$. An integration by parts against $\pi \propto e^{-U}$ gives, for $\omega$ as in \eqref{eq: riesz},
\begin{equation}
    \pi\bigl( \nabla^2 \omega : \mathcal C \bigr) = \pi\bigl( \nabla \omega \cdot b_{\mathcal C} \bigr) .
    \label{eq: representor ibp}
\end{equation}
When $b_{\mathcal C} = \nabla g$ for a function $g$ whose derivatives grow at most polynomially, \eqref{eq: generator ibp} and the uniqueness in Theorem~\ref{thm: riesz} give $\psi_{\mathcal C} = g - \pi(g)$. Table~\ref{tab: explicit representors} lists covariances and potentials for which $b_{\mathcal C}$ is a gradient.

\begin{table}[htb]
\centering
\small
\begin{tabular}{lll}
\toprule
covariance $\mathcal C(x)$ & potential $U$ & representor $\psi_{\mathcal C}(x)$, up to a constant \\
\midrule
$a\bigl( U(x) \bigr)\, I_d$ & any & $\displaystyle \int_0^{U(x)} a(s)\, \D s - a\bigl( U(x) \bigr)$ \\[2.5ex]
$\mathcal C(x) \in \mathbb R$, $d = 1$ & any & $\displaystyle \int_0^{x} \mathcal C(y)\, U'(y)\, \D y - \mathcal C(x)$ \\[2.5ex]
$\mathrm{diag}\bigl( c_1(x_1), \dots, c_d(x_d) \bigr)$ & $\displaystyle U(x) = \sum_{i=1}^d U_i(x_i)$ & $\displaystyle \sum_{i=1}^d \Bigl( \int_0^{x_i} c_i(y)\, U_i'(y)\, \D y - c_i(x_i) \Bigr)$ \\[2.5ex]
$a_1(|x|)\, I_d + a_2(|x|)\, \frac{x \otimes x}{|x|^2}$ & $U(x) = u(|x|)$ & $\displaystyle \int_0^{|x|} \Bigl( (a_1 + a_2)(r)\, u'(r) - \frac{(d-1)\, a_2(r)}{r} \Bigr) \D r - (a_1 + a_2)(|x|)$ \\
\bottomrule
\end{tabular}
\caption{Covariances whose Riesz representor $\psi_{\mathcal C}$ is explicit, with $\psi_{\mathcal C}$ given up to the additive constant fixed by $\pi(\psi_{\mathcal C}) = 0$. The functions $a$, $a_1$, $a_2$ and $c_1, \dots, c_d$ are of class $C^2$, and $a_2(r) = O(r^2)$ as $r \to 0$.}
\label{tab: explicit representors}
\end{table}

With a constant $a = \sigma^2$, the row $a(U(x))\, I_d$ of Table~\ref{tab: explicit representors} gives $\psi_{\mathcal C} = \sigma^2 \bigl( U - \pi(U) \bigr)$. For a general covariance in dimension $d \Ge 2$, the field $b_{\mathcal C}$ is not a gradient and $\psi_{\mathcal C}$ has no closed form; the extrapolation of Section~\ref{subsec: extrapolation} and the energy \eqref{eq: extrapolated energy} of Section~\ref{subsec: shifting} do not require $\psi_{\mathcal C}$ directly.

\subsection{Richardson--Romberg extrapolation}
\label{subsec: extrapolation}

Without the representor, the term of order $h$ of the numerical bias of a given test function is removed by extrapolation in the step size, as for the invariant distribution of an ergodic diffusion \citep{lemaire2015invariant} and for SGLD \citep{durmus2016stochastic}. For $h \in (0, h_0/2]$ and a function $f$ with $\norm{\nabla f}_\infty < \infty$, \eqref{eq: kr sglm} at the step sizes $h$ and $2h$ gives
\begin{equation}
    \biggl| \frac{\pi_{2h}(f) - \pi_h(f)}{h} - \frac12\, \pi(f \psi_{\mathcal C}) \biggr| \Le C \norm{\nabla f}_\infty\, h^{\frac12} \Bigl( 1 + \log\frac1h \Bigr) ,
    \label{eq: extrapolation sglm}
\end{equation}
so the coefficient $\frac12\, \pi(f \psi_{\mathcal C})$ is approximated by $(\pi_{2h}(f) - \pi_h(f))/h$, computed from the time averages of $f$ along runs of \eqref{eq: sglm} with the step sizes $h$ and $2h$; the error is of order $h$ when the remainder of \eqref{eq: kr sglm} is of order $h^2$, as expected in Remark~\ref{rem: h2}. Subtracting $\pi_{2h}(f) - \pi_h(f)$, which is $h$ times the approximation of $\frac12\, \pi(f \psi_{\mathcal C})$ in \eqref{eq: extrapolation sglm}, from $\pi_h(f)$ gives the approximation $2\pi_h(f) - \pi_{2h}(f)$ of $\pi(f)$, and \eqref{eq: kr sglm} gives
\begin{equation}
    \bigl| 2\pi_h(f) - \pi_{2h}(f) - \pi(f) \bigr| \Le C \norm{\nabla f}_\infty\, h^{\frac32} \Bigl( 1 + \log\frac1h \Bigr) ,
    \label{eq: extrapolation estimate}
\end{equation}
and the error is of order $h^2$ when the remainder of \eqref{eq: kr sglm} is of order $h^2$.

Extrapolation corrects the expectations of given test functions but not the samples, whose laws remain $\pi_h$ and $\pi_{2h}$; it does not raise the order of the Wasserstein-1 error in \eqref{eq: w1 limits}. Moreover, \eqref{eq: extrapolation sglm} holds for each given test function and does not mean that the empirical samples of $\pi_h$ and $\pi_{2h}$ approximate $\psi_{\mathcal C}$ itself: for the gradient noise $\sigma \xi$ with $\xi \sim \mathcal N(0, I_d)$, $\psi_{\mathcal C} = \sigma^2 \bigl( U - \pi(U) \bigr)$ by Table~\ref{tab: explicit representors}, which is not read off the samples.

\subsection{Postprocessed shifting correction}
\label{subsec: shifting}

The extrapolation of Section~\ref{subsec: extrapolation} improves the expectations $\pi_h(f)$, but not the samples of SGLM themselves. We instead define a vector field that shifts the samples of SGLM to more accurate ones. The vector field is $h\, \nabla \varphi_{\mathcal C}$, with $\varphi_{\mathcal C}$ the solution of
\begin{equation}
    \mathcal L \varphi_{\mathcal C} = \frac12\, \psi_{\mathcal C} , \qquad \pi(\varphi_{\mathcal C}) = 0 .
    \label{eq: shift potential}
\end{equation}
By \eqref{eq: generator ibp}, $\varphi_{\mathcal C}$ is the minimizer of the functional
\begin{equation}
    \mathcal E(\varphi) := \frac12\, \pi\bigl( |\nabla \varphi|^2 \bigr) + \frac12\, \pi\bigl( \psi_{\mathcal C}\, \varphi \bigr) , \qquad \pi(\varphi) = 0 ,
    \label{eq: shift energy}
\end{equation}
and since $\psi_{\mathcal C} \in L^2(\pi)$ and $\pi(\psi_{\mathcal C}) = 0$, the minimizer over the functions $\varphi$ with $\varphi, |\nabla \varphi| \in L^2(\pi)$ exists and is unique by the Poincar\'e inequality \eqref{eq: poincare} and the Riesz representation theorem. The shifting correction replaces a sample $Z \sim \pi_h$ of SGLM by the postprocessed sample
\begin{equation}
    Z + h\, \nabla \varphi_{\mathcal C}(Z) ,
    \label{eq: shifting correction}
\end{equation}
where the map $x \mapsto x + h\, \nabla \varphi_{\mathcal C}(x)$ is a postprocessor in the sense of \citep{vilmart2015postprocessed}, that is, a change of variables applied to the samples to raise the order for the invariant distribution. Formally, let $f$ be a smooth test function, and assume that the remainder of \eqref{eq: kr sglm} is of order $h^2$, as expected in Remark~\ref{rem: h2}. A Taylor expansion of $f$ at $x$ gives
\begin{equation*}
    f\bigl( x + h \nabla \varphi_{\mathcal C}(x) \bigr) = f(x) + h\, \nabla f(x) \cdot \nabla \varphi_{\mathcal C}(x) + \mathcal O(h^2) ,
\end{equation*}
and the expectation under $\pi_h$ gives
\begin{equation*}
    \pi_h\bigl( f(x + h \nabla \varphi_{\mathcal C}(x)) \bigr) = \pi_h(f) + h\, \pi_h(\nabla f \cdot \nabla \varphi_{\mathcal C}) + \mathcal O(h^2) .
\end{equation*}
By \eqref{eq: kr sglm} applied to $f$ and to $\nabla f \cdot \nabla \varphi_{\mathcal C}$,
\begin{equation*}
    \pi_h(f) = \pi(f) + \frac h2\, \pi(f \psi_{\mathcal C}) + \mathcal O(h^2) , \qquad
    \pi_h(\nabla f \cdot \nabla \varphi_{\mathcal C}) = \pi(\nabla f \cdot \nabla \varphi_{\mathcal C}) + \mathcal O(h) ,
\end{equation*}
so that
\begin{equation*}
    \pi_h\bigl( f(x + h \nabla \varphi_{\mathcal C}(x)) \bigr) = \pi(f) + \frac h2\, \pi(f \psi_{\mathcal C}) + h\, \pi(\nabla f \cdot \nabla \varphi_{\mathcal C}) + \mathcal O(h^2) .
\end{equation*}
By \eqref{eq: generator ibp} with $g_1 = f$ and $g_2 = \varphi_{\mathcal C}$, and by \eqref{eq: shift potential},
\begin{equation*}
    \pi(\nabla f \cdot \nabla \varphi_{\mathcal C}) = -\pi(f\, \mathcal L \varphi_{\mathcal C}) = -\frac12\, \pi(f \psi_{\mathcal C}) ,
\end{equation*}
so the terms $\frac h2\, \pi(f \psi_{\mathcal C})$ and $h\, \pi(\nabla f \cdot \nabla \varphi_{\mathcal C})$ cancel, and formally
\begin{equation*}
    \pi_h\bigl( f(x + h \nabla \varphi_{\mathcal C}(x)) \bigr) = \pi(f) + \mathcal O(h^2) .
\end{equation*}
The error of the postprocessed samples is thus of order $h^2$ for every test function at once. The shifting correction formally raises the order of the Wasserstein-1 error in \eqref{eq: w1 limits} to $h^2$. The rigorous study of the shifting correction is left for future work.

When $\psi_{\mathcal C}$ is explicit, as in Section~\ref{subsec: explicit representors}, $\varphi_{\mathcal C}$ is approximated by minimizing \eqref{eq: shift energy} over a class of functions, with $\pi$ replaced by the average over the empirical samples of $\pi_h$. In general, the energy \eqref{eq: shift energy} uses $\psi_{\mathcal C}$ only through $\frac12\, \pi(\psi_{\mathcal C}\, \varphi)$, the coefficient \eqref{eq: coefficients riesz} of the test function $\varphi$, which \eqref{eq: extrapolation sglm} approximates. Replacing $\frac12\, \pi(\psi_{\mathcal C}\, \varphi)$ by $(\pi_{2h}(\varphi) - \pi_h(\varphi))/h$, and $\pi$ by $\pi_h$ in $\frac12\, \pi(|\nabla \varphi|^2)$, gives
\begin{equation}
    \mathcal E_h(\varphi) := \frac12\, \pi_h\bigl( |\nabla \varphi|^2 \bigr) + \frac{\pi_{2h}(\varphi) - \pi_h(\varphi)}{h} .
    \label{eq: extrapolated energy}
\end{equation}
The energy $\mathcal E_h$ requires only the empirical samples of $\pi_h$ and $\pi_{2h}$, from runs of \eqref{eq: sglm} with the step sizes $h$ and $2h$, so it is completely data-driven: it involves neither $U$, $\nabla U$ nor $\mathcal C(x)$, and one minimization gives the shift. It does not change when a constant is added to $\varphi$, so its minimization needs no constraint on $\varphi$. A minimizer $\varphi$ of $\mathcal E_h$ satisfies $h\, \pi_h(\nabla f \cdot \nabla \varphi) = \pi_h(f) - \pi_{2h}(f)$ for every $f$, so the shift $x \mapsto x + h \nabla \varphi(x)$ formally carries $\pi_h(f)$ to the approximation $2\pi_h(f) - \pi_{2h}(f)$ of \eqref{eq: extrapolation estimate}, up to terms of order $h^2$, for every $f$ at once. The step size $h$ cannot be too small for the minimization to be stable: the fluctuation of the empirical averages that replace $\pi_{2h}(\varphi)$ and $\pi_h(\varphi)$ is divided by $h$ in \eqref{eq: extrapolated energy}, and the minimizer of $\mathcal E_h$ with these empirical averages fits this fluctuation.

\begin{remark}
The loss function $\mathcal E_h$ of \eqref{eq: extrapolated energy} is closely related to the loss in action matching \citep{neklyudov2023action}, which learns, from samples of a family of distributions $(q_t)_{t \in [0,1]}$, a vector field $\nabla s_t$ that transports them, by minimizing
\begin{equation*}
    q_0(s_0) - q_1(s_1) + \int_0^1 q_t\Bigl( \frac12\, |\nabla s_t|^2 + \partial_t s_t \Bigr)\, \D t
\end{equation*}
over functions $s_t(x)$. With $q_0 = \pi_{2h}$, $q_1 = \pi_h$, the function $s_t = h \varphi$ independent of $t$, and $q_t$ replaced by $\pi_h$ in the integral, this loss equals $h^2\, \mathcal E_h(\varphi)$. A loss of the same form, whose $2$ terms are expectations estimated by sampling, is also written in \citep{cao2022learning}. In these works, the vector field transports the samples between the given distributions, while the shifting correction applies it beyond them.

Since $\pi_{2h} - \pi_h$ is of order $h$, a minimizer $\varphi$ of \eqref{eq: extrapolated energy} formally gives $\pi_{2h}\bigl( f(x + h \nabla \varphi(x)) \bigr) = \pi_h(f) + \mathcal O(h^2)$, so the map $x \mapsto x + h \nabla \varphi(x)$ carries $\pi_{2h}$ to $\pi_h$, and the shifting correction applies the same map once more, from $\pi_h$ toward $\pi$. Formally, up to terms of order $h^2$, the shifting correction is thus the extrapolation at the time $2$ of the Wasserstein geodesic from $\pi_{2h}$ at the time $0$ to $\pi_h$ at the time $1$, the operation used in \citep{gallouet2024geodesic, gallouet2025metric} for time discretizations of Wasserstein gradient flows.
\end{remark}

\subsection{Calculation of the Kantorovich--Rubinstein norm}
\label{subsec: kr computation}

The coefficients of the Wasserstein-1 errors in Corollary~\ref{cor: w1 limits} are $\frac12 \norm{\psi\,\pi}_{\mathrm{KR}}$, with $\psi = \psi_{\mathcal C}$ for SGLM and $\psi = \psi_{\nabla^2 U} + \psi_{\mathcal C}$ for SGLD, so they require the norm itself and not only the representor. On a convex compact domain, with no flux through its boundary, the supremum that defines the norm in Section~\ref{subsec: notation} equals the least cost of a flux $F$ that carries the negative part of $\psi\pi$ onto its positive part \citep[(2)]{li2018parallel}:
\begin{equation}
    \norm{\psi\pi}_{\mathrm{KR}} = \min\Bigl\{ \int_{\mathbb R^d} |F(x)|\, \D x : \operatorname{div} F = -\psi\,\pi \Bigr\} .
    \label{eq: kr flux}
\end{equation}

For the computation, $\psi\,\pi$ is restricted to a box centered at the origin that carries all but a negligible part of $\pi$, with $n$ cells per side: the function $f$ of the supremum is placed at the centers of the cells, with forward differences of Euclidean norm at most $1$; the flux $F$ is placed on the faces; and $\psi\,\pi$ is replaced by its cell averages. The primal--dual iteration of \citet{li2018parallel} solves the two discrete problems together, and refining the grid approximates \eqref{eq: kr flux}.
The grid has $n^d$ cells, so the cost grows exponentially with the dimension: in $d=4$, $n = 100$ already gives $10^8$ cells. Computing Kantorovich--Rubinstein norms in high dimensions is intrinsically difficult, and it is possibly unnecessary in the current setting, because the Wasserstein-1 error is not a \emph{quantity of interest} in real computations: samples are used to estimate the expectations $\pi(f)$ of given test functions, and the bias of each of them is given to first order by the coefficient \eqref{eq: coefficients riesz}, which \eqref{eq: extrapolation sglm} approximates without the norm or the representor.

\section{Numerical experiments}
\label{sec: numerics}

The numerical experiments study the convergence of SGLM and SGLD and compare the two integrators in various tasks. Sections~\ref{subsec: mixture} and~\ref{subsec: radial} verify the leading terms of order $h$ of the numerical bias and of the Wasserstein-1 error in Section~\ref{sec: mse}, Section~\ref{subsec: rings} applies the bias corrections of Section~\ref{sec: representors} to SGLM, and Section~\ref{sec:mnist} compares the MSE of SGLM and SGLD in a Bayesian logistic regression on MNIST.

\subsection{Evaluating the coefficients in a 4D Gaussian mixture distribution}
\label{subsec: mixture}

\begin{figure}[htb]
\centering
\includegraphics[width=0.96\textwidth]{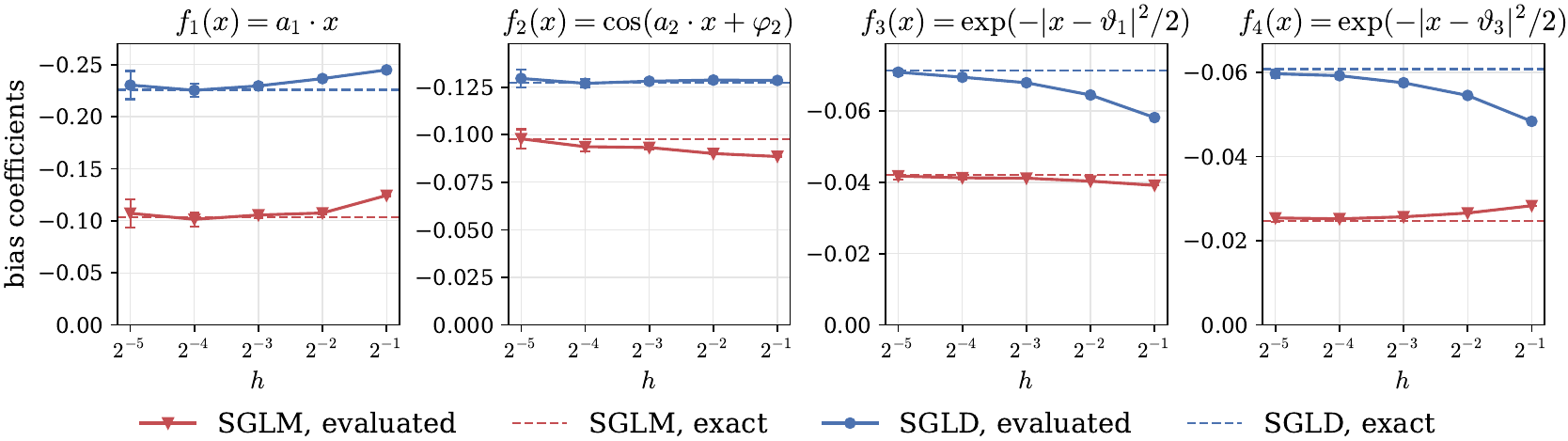}
\caption{Gaussian noise: the evaluated coefficients $(\pi_h(f) - \pi(f))/h$ of SGLM and $(\bar\pi_h(f) - \pi(f))/h$ of SGLD, with error bars of $2$ standard errors, and the exact coefficients \eqref{eq: coefficients riesz}, dashed. The vertical axis starts at $0$ and is reversed for negative exact coefficients.}
\label{fig: bias gaussian}
\end{figure}

The target $\pi$ is a mixture of $5$ Gaussian distributions on $\mathbb R^4$ with unequal weights, means $\vartheta_1, \dots, \vartheta_5$ near the vertices of a regular simplex of edge $4.5$, and covariance matrices of random orientation, so that $U = -\log \pi$ is not convex. The stochastic gradient $\nabla \hat U(x;\xi)$ adds to $\nabla U(x)$ either a Gaussian noise $\xi \sim \mathcal N(0, I_4)$, with $\mathcal C(x) = I_4$, or a non-Gaussian noise, a vector of independent centered exponential random variables of unit variance, rotated and scaled depending on $x$, with $\mathcal C(x) = \bigl( 2 - \pi(x)/\pi(\vartheta_1) \bigr)\, I_4$. Both covariances belong to the first row of Table~\ref{tab: explicit representors}, so the exact coefficients of the numerical bias follow from \eqref{eq: coefficients riesz} by quadrature. The test functions $f_1, \dots, f_4$, given in the panel titles of Figure~\ref{fig: bias gaussian} with vectors $a_1, a_2$ of lengths $1$, $0.6$ and the phase $\varphi_2$, have explicit expectations under $\pi$; they are taken from a pool of $12$, among the functions whose exact coefficients exceed $0.024$ in absolute value for both schemes and both noises.

\begin{figure}[htb]
\centering
\includegraphics[width=0.96\textwidth]{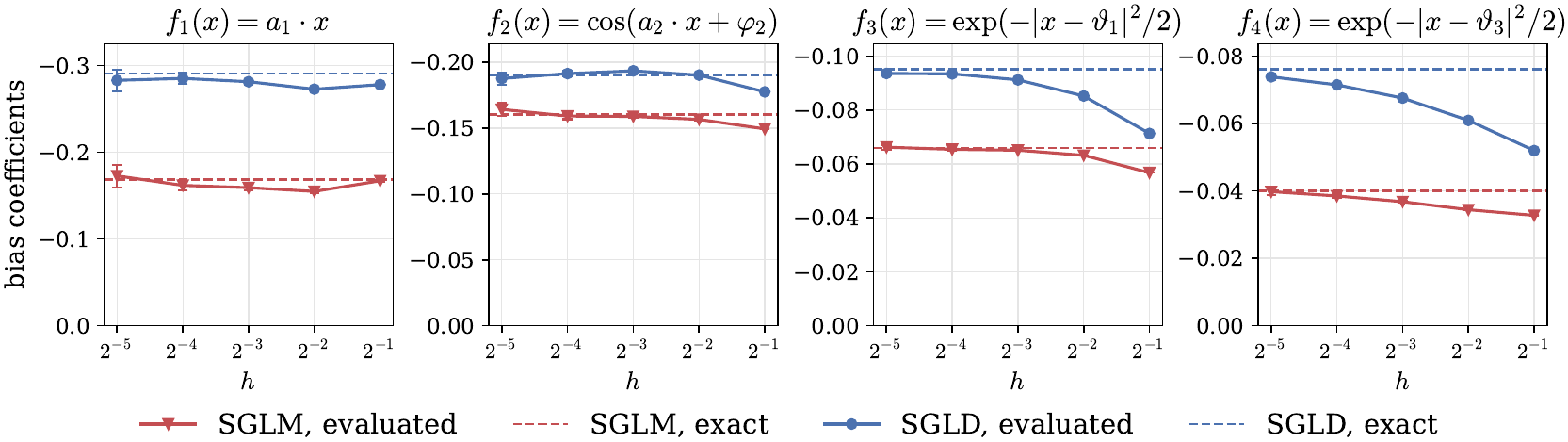}
\caption{Non-Gaussian noise: the evaluated coefficients $(\pi_h(f) - \pi(f))/h$ of SGLM and $(\bar\pi_h(f) - \pi(f))/h$ of SGLD, with error bars of $2$ standard errors, and the exact coefficients \eqref{eq: coefficients riesz}, dashed. The vertical axis starts at $0$ and is reversed for negative exact coefficients.}
\label{fig: bias nongaussian}
\end{figure}

The evaluated coefficients $(\pi_h(f) - \pi(f))/h$ of SGLM and $(\bar\pi_h(f) - \pi(f))/h$ of SGLD are obtained from the time averages of $16$ chains, each over the time $5 \times 10^7$, at $h = 2^{-1}, \dots, 2^{-5}$. In Figures~\ref{fig: bias gaussian} and~\ref{fig: bias nongaussian}, they approach the exact coefficients as $h$ decreases, for every test function, both schemes and both noises; at $h = 2^{-5}$, the relative difference is at most $3.8\%$. The agreement under the non-Gaussian noise, whose third moments do not vanish, is in line with \eqref{eq: bias coefficient} and \eqref{eq: sgld coefficient}, in which the gradient noise enters only through its covariance $\mathcal C(x)$. For each of these test functions and both noises, the exact coefficient of SGLM is smaller in absolute value than that of SGLD, between $41\%$ and $84\%$ of it.

\begin{table}[htb]
\centering
\small
\begin{tabular}{llccccc}
\toprule
noise & integrator & $\frac12 \norm{\psi\,\pi}_{\mathrm{KR}}$ & $f_1$ & $f_2$ & $f_3$ & $f_4$ \\
\midrule
\multirow{2}{*}{Gaussian} & SGLM & $0.5643$ & $18.3\%$ & $28.9\%$ & $12.3\%$ & $7.2\%$ \\
 & SGLD & $0.9824$ & $23.0\%$ & $21.6\%$ & $12.0\%$ & $10.2\%$ \\
\midrule
\multirow{2}{*}{non-Gaussian} & SGLM & $0.9701$ & $17.4\%$ & $27.5\%$ & $11.2\%$ & $6.8\%$ \\
 & SGLD & $1.3636$ & $21.3\%$ & $23.2\%$ & $11.5\%$ & $9.2\%$ \\
\bottomrule
\end{tabular}
\caption{Exact coefficients $\frac12 \norm{\psi\,\pi}_{\mathrm{KR}}$ of the Wasserstein-1 errors in Corollary~\ref{cor: w1 limits}, with $\psi = \psi_{\mathcal C}$ for SGLM and $\psi = \psi_{\nabla^2 U} + \psi_{\mathcal C}$ for SGLD. For each test function $f$ of Figure~\ref{fig: bias gaussian}, the percentage in the table shows $|\pi(f\psi)| / \bigl( \norm{\nabla f}_\infty \norm{\psi\, \pi}_{\mathrm{KR}} \bigr) \Le 1$, which is the ratio of its exact coefficient of the numerical bias, normalized by $\norm{\nabla f}_\infty$, to the exact coefficient of the Wasserstein-1 error.}
\label{tab: w1 coefficients}
\end{table}

Table~\ref{tab: w1 coefficients} lists the exact coefficients of the Wasserstein-1 errors in Corollary~\ref{cor: w1 limits}, calculated as in Section~\ref{subsec: kr computation} and extrapolated from grids with $n = 25, 50, 100$. The exact coefficient of SGLM is $57\%$ of that of SGLD under the Gaussian noise and $71\%$ under the non-Gaussian noise, so \eqref{eq: w1 comparison} holds for both noises. The ratios in Table~\ref{tab: w1 coefficients} are at most $0.29$, so the Wasserstein-1 error, a supremum over all test functions $f$ with $\norm{\nabla f}_\infty \Le 1$, exceeds the numerical bias of each of these test functions, normalized by $\norm{\nabla f}_\infty$, by a factor larger than $3$.

\subsection{Convergence of Wasserstein-1 error in a radial distribution}
\label{subsec: radial}

For a radial potential $U(x) = u(|x|)$ and the gradient noise $\sigma \xi$ with $\xi \sim \mathcal N(0, I_d)$, the covariance is $\mathcal C(x) = \sigma^2 I_d$, so both integrators commute with rotations. Hence the distributions $\pi_h$, $\bar\pi_h$ and $\pi$ are rotation invariant, and the Wasserstein-1 distance between two of them is determined by their radial densities. The radial parts of \eqref{eq: sgld} and of the Markovian form \eqref{eq: sglm markov} are Markov chains on $[0, \infty)$, whose invariant distributions are computed by discretizing their stationarity equations, and the representors $\psi_{\mathcal C} = \sigma^2 \bigl( U - \pi(U) \bigr)$ and $\psi_{\nabla^2 U}$ of \eqref{eq: psi hessian} are radial, with $\norm{\psi\,\pi}_{\mathrm{KR}} = \int_0^\infty \bigl| \int_{|x| \Le r} \psi \, \D\pi \bigr| \, \D r$. The evaluated coefficients $\mathcal W_1(\pi_h, \pi)/h$ and $\mathcal W_1(\bar\pi_h, \pi)/h$ at a given $h$, and the exact coefficients $\frac12 \norm{\psi\,\pi}_{\mathrm{KR}}$ of Corollary~\ref{cor: w1 limits}, are thus computed without sampling in every dimension. The potential is
\begin{equation}
    u(r) = (r^2 - r)\, \frac{r^2}{r^2 + 1/4} ,
    \label{eq: radial potential}
\end{equation}
which grows quadratically and is not globally convex, with a ring of minima at $|x| \approx 0.64$. Since $u(r) = -4r^3 + \mathcal O(r^4)$ near $r = 0$, $U$ is not three times differentiable at the origin and does not satisfy Assumption~\ref{ass: U}; the example mainly serves to observe the convergence of the Wasserstein-1 error numerically.

\begin{figure}[htb]
\centering
\includegraphics[width=0.46\textwidth]{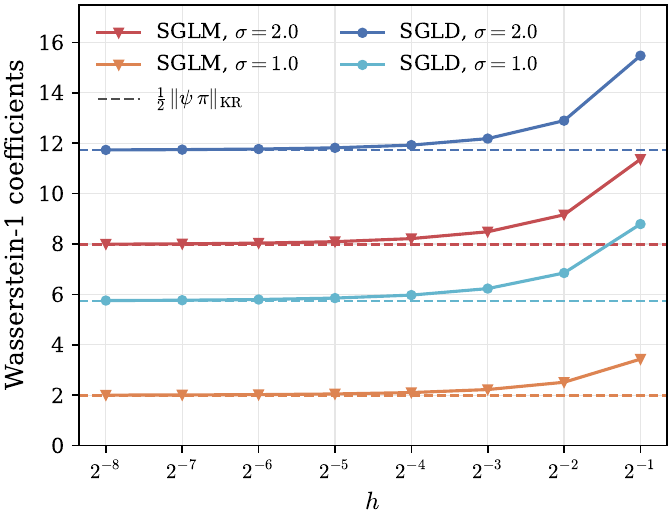}
\includegraphics[width=0.46\textwidth]{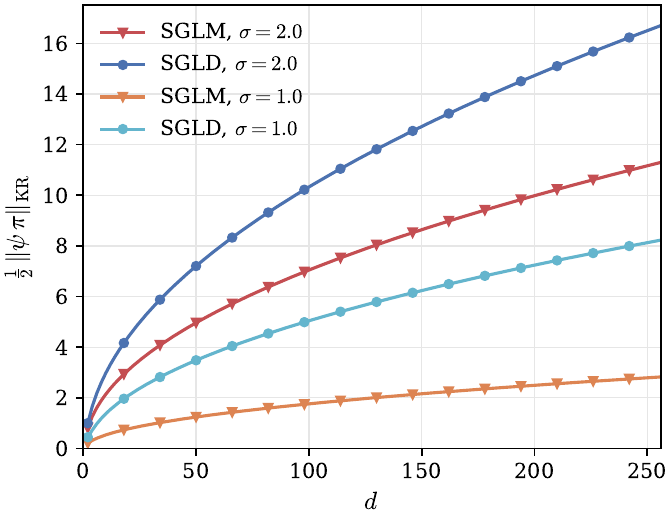}
\caption{Radial potential \eqref{eq: radial potential} with the gradient noise $\sigma \xi$, $\xi \sim \mathcal N(0, I_d)$. Left: the evaluated coefficients $\mathcal W_1(\pi_h, \pi)/h$ of SGLM and $\mathcal W_1(\bar\pi_h, \pi)/h$ of SGLD in dimension $d = 128$, with the exact coefficients $\frac12 \norm{\psi\,\pi}_{\mathrm{KR}}$ of Corollary~\ref{cor: w1 limits}, dashed. Right: the exact coefficients against $d = 2, 4, \dots, 256$.}
\label{fig: radial}
\end{figure}

In Figure~\ref{fig: radial}, left, the evaluated coefficients converge to the exact coefficients for $\sigma = 1.0, 2.0$, and their distance to the exact coefficients halves with $h$; at $h = 2^{-8}$ it is between $0.09\%$ and $0.32\%$ of the exact coefficients. In this example the remainder of \eqref{eq: w1 limits} is therefore of order $h^2$, the order expected in Remark~\ref{rem: h2} under global strong convexity, although $U$ is not convex.

In Figure~\ref{fig: radial}, right, the exact coefficients grow like $\sqrt d$ for $d = 2, 4, \dots, 256$: from $d = 128$ to $d = 256$ with $\sigma = 1.0$, they are multiplied by $1.416$ for SGLM and by $1.434$ for SGLD, close to $\sqrt 2$. This is the growth in the Gaussian case $u(r) = r^2$, where both exact coefficients are proportional to $\mathbb E|X| \sim \sqrt d$ with $X \sim \mathcal N(0, I_d)$.

The exact coefficient of SGLM is below that of SGLD for every $d$ and both $\sigma$, so \eqref{eq: w1 comparison} holds throughout, and the advantage of SGLM is larger when the variance of the gradient noise is smaller: the ratio of the exact coefficient of SGLM to that of SGLD is $0.46$ at $d = 2$ and $0.34$ at $d = 256$ for $\sigma = 1.0$, against $0.81$ and $0.68$ for $\sigma = 2.0$. The exact coefficient of SGLM is proportional to $\sigma^2$, while that of SGLD contains the representor $\psi_{\nabla^2 U}$ of the discretization error of EM, which does not depend on $\sigma$.

\subsection{Bias correction of SGLM in a 16D product of ring distributions}
\label{subsec: rings}

The target $\pi$ on $\mathbb R^{16}$ is a product of $8$ ring distributions in rotated coordinates,
\begin{equation}
    U(x) = \sum_{i=1}^{8} u_i\Bigl( \sqrt{y_{2i-1}^2 + y_{2i}^2} \Bigr) , \qquad y = Q x , \qquad u_i(r) = \frac{r^2}{2} + \beta_i\, e^{-r^2/2} ,
    \label{eq: ring potential}
\end{equation}
with a random orthogonal matrix $Q$ and $\beta_1, \dots, \beta_8$ drawn uniformly from $[2, 4]$. Each $u_i$ has a ring of minima at $r = \sqrt{2 \log \beta_i}$, so $U$ is not convex. The gradient noise added to $\nabla U(x)$ is either a Gaussian one $\xi \sim \mathcal N(0, I_{16})$, or a non-Gaussian one: on each pair $(x_{2i-1}, x_{2i})$, a vector of $2$ independent centered exponential random variables of unit variance, rotated by the angle $x_{2i-1} + x_{2i}$. Both noises have the covariance $\mathcal C(x) = I_{16}$, so $\psi_{\mathcal C} = U - \pi(U)$ by the first row of Table~\ref{tab: explicit representors}. The third moments of the non-Gaussian noise do not vanish and depend on $x$.

SGLM is run at the step sizes $h = 2^{-4}$ and $2h = 2^{-3}$, with $32$ chains, each over the time $2 \times 10^7$, at each step size. The expectations $\pi_h(f)$ and $\pi_{2h}(f)$ are replaced by the time averages of $f$ over the chains. The means $\pi(f)$ and the exact coefficients $\frac12\, \pi(f \psi_{\mathcal C}) = \frac12 \bigl( \pi(fU) - \pi(f)\, \pi(U) \bigr)$ are each computed with $1.024 \times 10^{10}$ exact samples of $\pi$. The test functions for evaluating the numerical bias are chosen as
\begin{align*}
    f_1(x) &= x_1 , & f_2(x) &= \frac14 \sum_{j=1}^{16} x_j , \\
    f_3(x) &= \cos x_1 , & f_4(x) &= \sin(x_1 + x_2) , \\
    f_5(x) &= \sqrt{1 + x_1^2} , & f_6(x) &= e^{-(x_1^2 + x_2^2)/2} , \\
    f_7(x) &= \log\bigl( 1 + e^{x_1 - x_2} \bigr) , & f_8(x) &= \sqrt{1 + |x|^2} , \\
    f_9(x) &= \cos(x_1 + x_2) , & f_{10}(x) &= \sqrt{1 + x_1^2 + x_2^2} , \\
    f_{11}(x) &= \frac{x_1 x_2}{\sqrt{1 + x_1^2 + x_2^2}} , & f_{12}(x) &= \log\bigl( 1 + x_1^2 + x_2^2 + x_3^2 + x_4^2 \bigr) .
\end{align*}

For each test function, Tables~\ref{tab: bias correction gaussian} and~\ref{tab: bias correction nongaussian} list the approximations of $\pi(f)$ in Table~\ref{tab: approximations} at $h = 2^{-4}$, each with $\pi(f)$ subtracted, where $\varphi_{\mathcal C}$ of \eqref{eq: shift potential} is approximated by minimizing \eqref{eq: extrapolated energy}.

\begin{table}[htb]
\centering
\small
\begin{tabular}{ll}
\toprule
method & approximation of $\pi(f)$ \\
\midrule
SGLM & $\pi_h(f)$ \\
explicit coefficient & $\pi_h(f) - \frac h2\, \pi(f \psi_{\mathcal C})$ \\
extrapolation & $2\pi_h(f) - \pi_{2h}(f)$ \\
shifting correction & $\pi_h\bigl( f(x + h \nabla \varphi_{\mathcal C}(x)) \bigr)$ \\
\bottomrule
\end{tabular}
\caption{The approximations of $\pi(f)$ by SGLM and by the approaches of Section~\ref{sec: representors}, listed in the columns SGLM, explicit, extrapolation and shift of Tables~\ref{tab: bias correction gaussian} and~\ref{tab: bias correction nongaussian}.}
\label{tab: approximations}
\end{table}

The entry of SGLM is the numerical bias $\pi_h(f) - \pi(f)$, and the entries of explicit, extrapolation and shift are the errors of the approaches of Section~\ref{sec: representors}. The column explicit uses the exact coefficient, as in Section~\ref{subsec: explicit representors}, and the column extrapolation is the error bounded in \eqref{eq: extrapolation estimate}.

The column shift uses the postprocessed samples \eqref{eq: shifting correction}, with $\varphi_{\mathcal C}$ approximated by a neural network with $3$ hidden layers of width $64$, trained on $1.68 \times 10^7$ states of SGLM at each of the step sizes $h$ and $2h$, stored once per unit of time along $32$ chains and thus covering the total time $1.68 \times 10^7$. The change $\pi_h\bigl( f(x + h \nabla \varphi_{\mathcal C}(x)) \bigr) - \pi_h(f)$ is computed on states of runs at $h$ independent of the training states and added to the time average of $f$ at $h$.

\begin{table}[htb]
\centering
\small
\begin{tabular}{lcccc}
\toprule
 & SGLM & explicit & extrapolation & shift \\
\midrule
$f_{1}$ & $-8.64\times10^{-5}$ & $\boldsymbol{-8.57\times10^{-5}}$ & $-3.00\times10^{-4}$ & $-5.36\times10^{-4}$ \\
$f_{2}$ & $\boldsymbol{-1.83\times10^{-4}}$ & $-1.85\times10^{-4}$ & $-5.42\times10^{-4}$ & $-8.86\times10^{-4}$ \\
$f_{3}$ & $-4.86\times10^{-3}$ & $\boldsymbol{-6.73\times10^{-5}}$ & $1.49\times10^{-4}$ & $1.86\times10^{-4}$ \\
$f_{4}$ & $\boldsymbol{-2.63\times10^{-5}}$ & $-2.68\times10^{-5}$ & $-9.72\times10^{-5}$ & $-2.61\times10^{-4}$ \\
$f_{5}$ & $5.70\times10^{-3}$ & $\boldsymbol{1.08\times10^{-4}}$ & $-2.63\times10^{-4}$ & $-2.29\times10^{-4}$ \\
$f_{6}$ & $-2.99\times10^{-3}$ & $\boldsymbol{-1.51\times10^{-5}}$ & $4.67\times10^{-5}$ & $1.50\times10^{-4}$ \\
$f_{7}$ & $3.86\times10^{-3}$ & $6.38\times10^{-5}$ & $-2.60\times10^{-4}$ & $\boldsymbol{-1.71\times10^{-6}}$ \\
$f_{8}$ & $3.61\times10^{-2}$ & $\boldsymbol{5.07\times10^{-4}}$ & $-1.48\times10^{-3}$ & $-1.99\times10^{-3}$ \\
$f_{9}$ & $-3.97\times10^{-3}$ & $\boldsymbol{-2.48\times10^{-7}}$ & $2.29\times10^{-5}$ & $1.11\times10^{-4}$ \\
$f_{10}$ & $1.00\times10^{-2}$ & $\boldsymbol{1.59\times10^{-4}}$ & $-3.98\times10^{-4}$ & $-5.42\times10^{-4}$ \\
$f_{11}$ & $9.84\times10^{-5}$ & $\boldsymbol{-3.18\times10^{-5}}$ & $-1.08\times10^{-4}$ & $1.35\times10^{-4}$ \\
$f_{12}$ & $1.11\times10^{-2}$ & $\boldsymbol{1.19\times10^{-4}}$ & $-3.47\times10^{-4}$ & $-3.33\times10^{-4}$ \\
\bottomrule
\end{tabular}
\caption{Gaussian noise: the numerical bias $\pi_h(f) - \pi(f)$ of SGLM at $h = 2^{-4}$ and the errors of the explicit coefficient, the extrapolation and the shifting correction of Section~\ref{sec: representors}. The entry of smallest absolute value in each row is in bold.}
\label{tab: bias correction gaussian}
\end{table}

\begin{table}[htb]
\centering
\small
\begin{tabular}{lcccc}
\toprule
 & SGLM & explicit & extrapolation & shift \\
\midrule
$f_{1}$ & $\boldsymbol{1.99\times10^{-4}}$ & $2.00\times10^{-4}$ & $6.42\times10^{-4}$ & $2.44\times10^{-3}$ \\
$f_{2}$ & $5.56\times10^{-5}$ & $5.38\times10^{-5}$ & $\boldsymbol{-3.11\times10^{-5}}$ & $1.32\times10^{-3}$ \\
$f_{3}$ & $-4.76\times10^{-3}$ & $3.22\times10^{-5}$ & $\boldsymbol{4.79\times10^{-6}}$ & $3.09\times10^{-4}$ \\
$f_{4}$ & $\boldsymbol{4.94\times10^{-4}}$ & $\boldsymbol{4.94\times10^{-4}}$ & $-9.99\times10^{-4}$ & $5.61\times10^{-4}$ \\
$f_{5}$ & $5.63\times10^{-3}$ & $\boldsymbol{3.87\times10^{-5}}$ & $-1.91\times10^{-4}$ & $-4.42\times10^{-4}$ \\
$f_{6}$ & $-2.98\times10^{-3}$ & $\boldsymbol{-5.15\times10^{-6}}$ & $5.24\times10^{-5}$ & $7.70\times10^{-5}$ \\
$f_{7}$ & $4.02\times10^{-3}$ & $\boldsymbol{2.27\times10^{-4}}$ & $2.61\times10^{-4}$ & $1.13\times10^{-3}$ \\
$f_{8}$ & $3.61\times10^{-2}$ & $\boldsymbol{5.43\times10^{-4}}$ & $-1.44\times10^{-3}$ & $-1.76\times10^{-3}$ \\
$f_{9}$ & $-4.00\times10^{-3}$ & $\boldsymbol{-3.50\times10^{-5}}$ & $-8.72\times10^{-5}$ & $-1.78\times10^{-4}$ \\
$f_{10}$ & $1.00\times10^{-2}$ & $\boldsymbol{1.47\times10^{-4}}$ & $-4.08\times10^{-4}$ & $-3.72\times10^{-4}$ \\
$f_{11}$ & $2.07\times10^{-4}$ & $\boldsymbol{7.66\times10^{-5}}$ & $1.37\times10^{-4}$ & $3.47\times10^{-4}$ \\
$f_{12}$ & $1.11\times10^{-2}$ & $\boldsymbol{7.75\times10^{-5}}$ & $-4.24\times10^{-4}$ & $-6.31\times10^{-4}$ \\
\bottomrule
\end{tabular}
\caption{Non-Gaussian noise: the same quantities as in Table~\ref{tab: bias correction gaussian}, with the entry of smallest absolute value in each row in bold.}
\label{tab: bias correction nongaussian}
\end{table}

In general, the explicit coefficient has the smallest errors: it removes the term of order $h$ of the bias with the exact coefficient, and what remains is the remainder of \eqref{eq: kr sglm} and the fluctuation of the time averages. The extrapolation and the shifting correction have larger errors for some test functions, mainly from the additional error of the run at $2h = 2^{-3}$. Formally, the term of order $h^2$ of $\pi_{2h}(f) - \pi(f)$ is $4$ times that of $\pi_h(f) - \pi(f)$, so $2\pi_h(f) - \pi_{2h}(f) - \pi(f)$ carries the term of order $h^2$ of $\pi_h(f) - \pi(f)$ with the factor $-2$, against $1$ for the explicit coefficient; for $f_4$ under the non-Gaussian noise, Table~\ref{tab: bias correction nongaussian} gives $-9.99\times10^{-4}$ for the extrapolation against $4.94\times10^{-4}$ for the explicit coefficient. The shifting correction inherits the additional error of the run at $2h$, since the term of order $h$ of $\pi_h\bigl( f(x + h \nabla \varphi_{\mathcal C}(x)) \bigr) - \pi_h(f)$ is formally $h\, \pi_h(\nabla f \cdot \nabla \varphi_{\mathcal C}) = \pi_h(f) - \pi_{2h}(f)$ by the minimization of \eqref{eq: extrapolated energy}; the entries of shift for $f_8$, $f_{10}$ and $f_{12}$ have, in both tables, the sign and size of those of extrapolation. The largest error of the shifting correction, $2.44\times10^{-3}$ for $f_1$ under the non-Gaussian noise, comes from the fluctuation of the stored states on which $\varphi_{\mathcal C}$ is trained, which cover the total time $1.68 \times 10^7$, against $6.4 \times 10^8$ for the time averages: on the stored states, $(\pi_{2h}(f_1) - \pi_h(f_1))/h$ is $-3.70\times10^{-2}$, with the standard error $1.69\times10^{-2}$, against $-7.08\times10^{-3}$ from the time averages, and the change caused by the shift, $2.24\times10^{-3}$, is close to $-h$ times $-3.70\times10^{-2}$. Nevertheless, the largest absolute error over $f_1, \dots, f_{12}$ is $3.61\times10^{-2}$ for SGLM under both noises, against at most $5.43\times10^{-4}$, $1.48\times10^{-3}$ and $2.44\times10^{-3}$ for the explicit coefficient, the extrapolation and the shifting correction; the corrections reduce the largest absolute error of SGLM by factors between $14.8$ and $71.1$, so the bias correction works in this example.

\subsection{Bayesian logistic regression on MNIST}
\label{sec:mnist}

\begin{figure}[htb]
  \centering
  \includegraphics[width=0.72\textwidth]{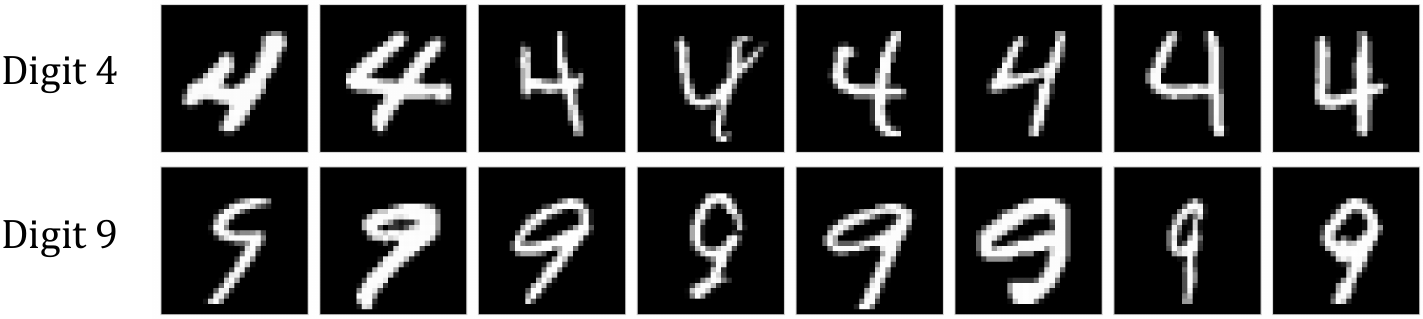}
  \caption{Examples of the digits $4$ and $9$ of the MNIST data set.}
  \label{fig:mnist-data}
\end{figure}

Bayesian logistic regression is applied to the handwritten digits $4$ and $9$ of the MNIST data set~\citep{lecun1998gradient}, shown in Figure~\ref{fig:mnist-data}, with the training set restricted to these digits. Each of the $n = 10{,}791$ training images, $5{,}342$ of the digit $4$ and $5{,}449$ of the digit $9$, is flattened, its pixel intensities are divided by $255$, and a constant $1$ is appended for the intercept, giving the vector $a_i \in \mathbb R^{785}$, with the label $y_i = 0$ for the digit $4$ and $y_i = 1$ for the digit $9$. With the prior $\mathcal N(0, \sigma_0^2 I_{785})$, $\sigma_0^2 = 10^{-3}$, the target is $\pi \propto e^{-U}$ on $\mathbb R^{785}$ with
\begin{equation}
    U(x) = \frac{|x|^2}{2\sigma_0^2} + \sum_{i=1}^{n} \Bigl( \log\bigl( 1 + e^{a_i \cdot x} \bigr) - y_i\, a_i \cdot x \Bigr) .
    \label{eq:mnist-potential}
\end{equation}
The stochastic gradient averages $P$ indices $\xi = (\xi_1, \dots, \xi_P)$, drawn independently and uniformly from $\{1, \dots, n\}$ at every step, independently of the Gaussian variables of the integrator:
\begin{equation}
    \nabla \hat U(x;\xi) = \frac{x}{\sigma_0^2} + \frac{n}{P} \sum_{p=1}^{P} \Bigl( \frac{1}{1 + e^{-a_{\xi_p} \cdot x}} - y_{\xi_p} \Bigr)\, a_{\xi_p} .
    \label{eq:mnist-minibatch-gradient}
\end{equation}

The test functions measure the posterior in directions of different curvature. Let $x_*$ be the minimizer of $U$, computed by L-BFGS with the full gradient, and let $\lambda_1 \Ge \dots \Ge \lambda_8$ be the $8$ largest eigenvalues of $\nabla^2 U(x_*)$, with orthonormal eigenvectors $v_1, \dots, v_8$. In the standardized coordinates
\begin{equation}
    z_j(x) = \sqrt{\lambda_j}\, v_j \cdot (x - x_*) , \qquad j = 1, \dots, 8 ,
    \label{eq:mnist-local-coordinates}
\end{equation}
the $16$ bounded test functions are
\begin{align}
    f_{C_j}(x) &= \cos z_j(x) , && j = 1, \dots, 8 ,
    \label{eq:mnist-cosine-observables} \\
    f_{G_j}(x) &= e^{-z_j(x)^2/2} , && j = 1, \dots, 4 ,
    \label{eq:mnist-gaussian-1d-observables} \\
    f_{G_{j,k}}(x) &= e^{-(z_j(x)^2 + z_k(x)^2)/2} , && (j,k) \in \{(1,2), (3,4), (5,6), (7,8)\} .
    \label{eq:mnist-gaussian-2d-observables}
\end{align}
The cosines $f_{C_j}$ probe single directions, and the Gaussians $f_{G_j}$ and $f_{G_{j,k}}$ the concentration near $x_*$ in $1$ or $2$ coordinates.

The step size follows \citet[Section~4]{brosse2018promises} for logistic regression with the Gaussian prior, $h = (\sigma_0^{-2} + \delta/4)^{-1} \approx 8.5099 \times 10^{-6}$, where $\delta$ is the largest eigenvalue of $\sum_{i=1}^n a_i \otimes a_i$. For each batch size $P = 1024, 2048, 4096$, SGLD and SGLM are run in $20$ pairs of trajectories, the $2$ trajectories of a pair using the same indices $\xi$ and the same Gaussian variables. Every trajectory starts at $x_*$, discards its first $2^{14}$ steps and retains the next $2^{17}$ states. The reference values $\widehat\pi_{\mathrm{ref}}(f)$ are the pooled averages of $4$ independent chains of NUTS~\citep{hoffman2014no} with the full gradient, which leaves $\pi$ invariant, each with $2^{19}$ retained states after $2^{14}$ adaptation steps. The chains use the implementation of NUTS in Pyro~\citep{bingham2019pyro}, in the coordinates whitened by the Cholesky factor of $\nabla^2 U(x_*)$, with the target acceptance rate $0.9$, the maximal tree depth $10$, and the step size and a diagonal mass matrix adapted during the adaptation steps. The $4$ chains agree closely for all $16$ test functions, their Monte Carlo errors are small compared with the errors below, and these Monte Carlo errors are included in the resampling intervals below.

For each integrator, batch size and test function $f$, the estimates use the first $N$ retained states, for every multiple $N$ of $8192$ up to $2^{17}$. With $Z_{r,k}$ the $k$th retained state of the repetition $r = 1, \dots, 20$, let
\[
  \widehat\mu_r=\frac{1}{N}\sum_{k=0}^{N-1}f(Z_{r,k}),
  \qquad
  \overline\mu=\frac{1}{20}\sum_{r=1}^{20}\widehat\mu_r
\]
be the time average and its mean over the repetitions. The estimated finite-time bias and the variance of the time averages are
\[
  \widehat b=\overline\mu-\widehat\pi_{\mathrm{ref}}(f),
  \qquad
  \widehat V=\frac{1}{20}\sum_{r=1}^{20}
    (\widehat\mu_r-\overline\mu)^2.
\]
The reported estimate of the MSE \eqref{eq: mse} is
\begin{equation}
  \widehat{\MSE}_f(N,h)
  =\frac{1}{20}\sum_{r=1}^{20}
   \bigl(\widehat\mu_r-\widehat\pi_{\mathrm{ref}}(f)\bigr)^2
  =\widehat V+\widehat b^{\,2}.
  \label{eq:mnist-mse}
\end{equation}
The variance is that of the time averages over the repetitions, not of single states, and the bias is that of time averages of finite length, not the numerical bias of the invariant distribution.

\begin{figure}[htb]
  \centering
  \includegraphics[width=\textwidth]{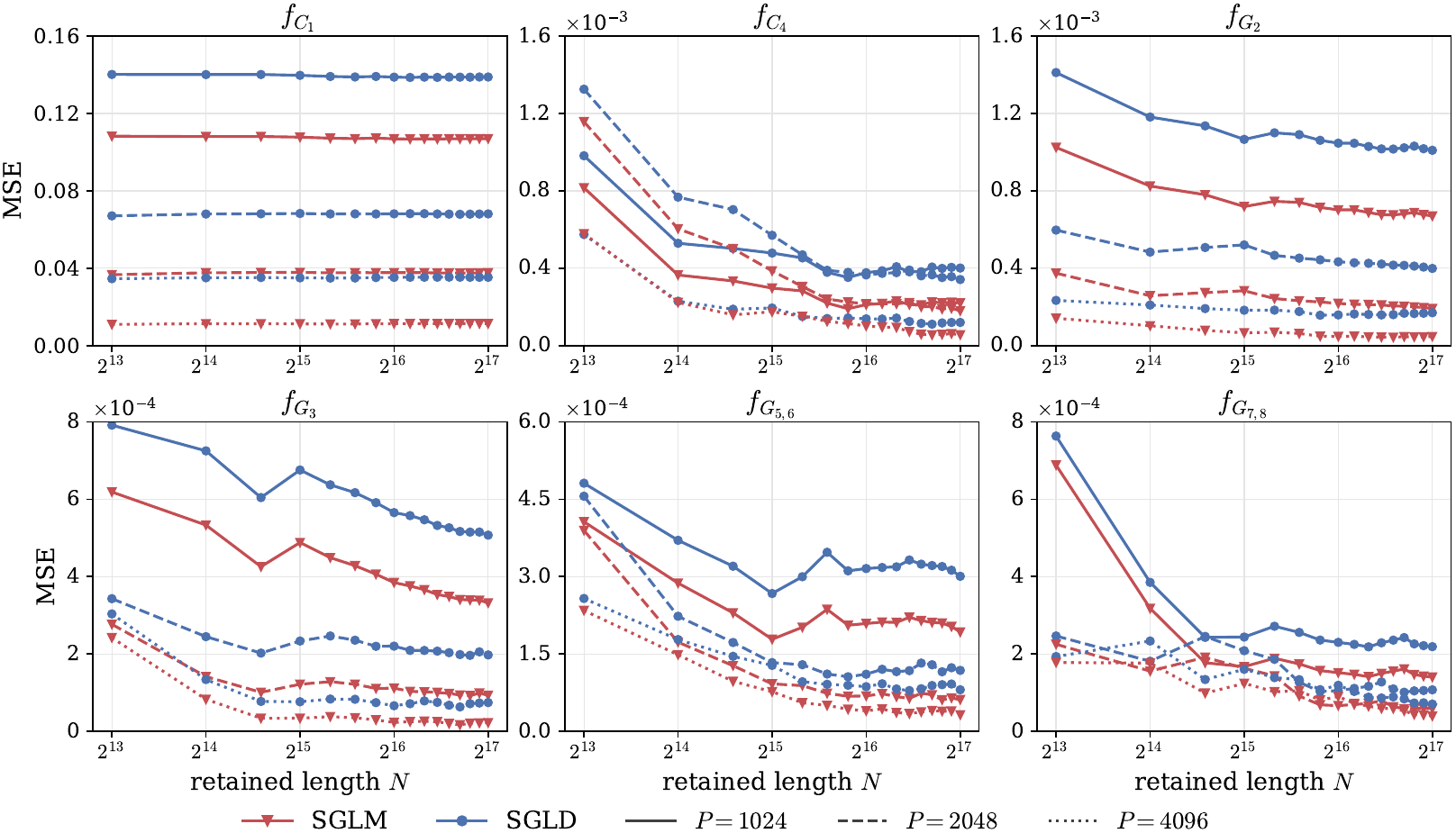}
  \caption{Estimated MSE \eqref{eq:mnist-mse} of the time averages of SGLD and SGLM for $f_{C_1}, f_{C_4}, f_{G_2}, f_{G_3}, f_{G_{5,6}}$ and $f_{G_{7,8}}$ against the retained length $N$, with $20$ pairs of trajectories for each batch size $P$.}
  \label{fig:mnist-mse}
\end{figure}

Figure~\ref{fig:mnist-mse} shows $2$ test functions of each form, $f_{C_1}, f_{C_4}, f_{G_2}, f_{G_3}, f_{G_{5,6}}$ and $f_{G_{7,8}}$, which together involve all $8$ coordinates \eqref{eq:mnist-local-coordinates}. For $f_{C_1}$ the estimated MSE is nearly constant in $N$, so longer averages bring little improvement in the tested range. For $f_{C_4}$ and $f_{G_{7,8}}$, several curves decrease markedly at small $N$ and then level off. For $f_{G_2}$, $f_{G_3}$ and $f_{G_{5,6}}$, longer averages are particularly beneficial for SGLM at $P = 4096$, whose errors fall below those of SGLD. The leveling off of several curves indicates a part of the error from the numerical bias, which longer averages at fixed $h$ and $P$ cannot remove.

At the largest retained length $N = 2^{17}$, SGLM has the smaller estimated MSE for all $16$ test functions at each of the $3$ batch sizes, that is, in all $48$ comparisons, and in each of them the $95\%$ interval of $\widehat{\MSE}_{\mathrm{SGLD}} - \widehat{\MSE}_{\mathrm{SGLM}}$ obtained by resampling the $20$ pairs with replacement lies above $0$. For each test function and batch size, the comparison at $N = 2^{17}$ uses the quantities
\begin{gather*}
    \text{MSE Reduction} = 100 \biggl( 1 - \frac{\widehat{\MSE}_{\mathrm{SGLM}}}{\widehat{\MSE}_{\mathrm{SGLD}}} \biggr) \% , \\
    \text{Bias Ratio} = \frac{|\widehat b_{\mathrm{SGLM}}|}{|\widehat b_{\mathrm{SGLD}}|} , \qquad
    \text{Variance Ratio} = \frac{\widehat V_{\mathrm{SGLM}}}{\widehat V_{\mathrm{SGLD}}} .
\end{gather*}
A positive MSE reduction and ratios below $1$ favor SGLM. Table~\ref{tab:mnist-endpoint} gives the median of each quantity over the $16$ test functions.

\begin{table}[htb]
	\centering
	\small
	\begin{tabular}{cccc}
		\toprule
		$P$ & MSE Reduction & Bias Ratio & Variance Ratio \\
		\midrule
		1024 & 34.6\% & 0.784 & 1.020 \\
		2048 & 46.5\% & 0.677 & 1.018 \\
		4096 & 63.3\% & 0.516 & 1.032 \\
		\bottomrule
	\end{tabular}
	\caption{Comparison of SGLM with SGLD at $N = 2^{17}$: the median over the $16$ test functions of each quantity, computed for each test function.}
	\label{tab:mnist-endpoint}
\end{table}

The median MSE reduction grows with the batch size, from $34.6\%$ to $63.3\%$, while the median variance ratio stays close to $1$ and the median bias ratio decreases from $0.784$ to $0.516$: the gain of SGLM comes from smaller finite-time biases, not from smaller fluctuations. This trend in $P$ is consistent with the separation of the errors in Theorem~\ref{thm: sglm mse}: a batch of $P$ independent stochastic gradients reduces the contribution of the gradient noise to $\Sigma^2 h^2 / P^2$, while the discretization error contributes $h^4$ for SGLM and $h^2$ for SGLD in \eqref{eq: sgld mse}. In the bias, the coefficient \eqref{eq: bias coefficient} of SGLM is divided by $P$, while the coefficient \eqref{eq: sgld coefficient} of SGLD keeps the contribution of the discretization error of EM at every batch size.

\section{Conclusion}
\label{sec: conclusion}

For a potential strongly convex outside a ball, this paper proves MSE bounds for the SGLM integrator and its full gradient form LM, and gives the explicit coefficients of $h$ in the numerical bias and in the Wasserstein-1 error of SGLM and SGLD. The samples of \eqref{eq: sglm} and \eqref{eq: lm} do not form a Markov chain; the modified class of Section~\ref{sec: modified} moves the analysis to the Markovian forms \eqref{eq: sglm markov} and \eqref{eq: lm markov}, and the discrete Poisson equation of \citep{lu2025mse, ye2026mse} converts their local weak errors against the modified dynamics into the bound $\frac{1}{Nh} + \Sigma^2 h^2 + h^4$ of Theorem~\ref{thm: sglm mse}, with the bound $\frac{1}{Nh} + h^4$ of Theorem~\ref{thm: lm mse} for the full gradient. The stochastic gradient error and the discretization error thus separate for an integrator of the overdamped dynamics, as for SG-UBU \citep{lu2025mse} and unlike \eqref{eq: sgld}.

The numerical bias of SGLM is of first order in $h$, with the coefficient \eqref{eq: bias coefficient} given by the covariance $\mathcal C(x)$ of the gradient noise. For \eqref{eq: sgld}, the same discrete Poisson equation gives the coefficient \eqref{eq: sgld coefficient}, which carries the Hessian $\nabla^2 U(x)$ besides the covariance $\mathcal C(x)$. This contribution of the discretization error of EM remains at every batch size, while the coefficient of SGLM is divided by the batch size. For a given test function, the two coefficients are not ordered in absolute value, since $\nabla^2 \phi_0$ is indefinite in general. In both coefficients, the gradient noise and the discretization error enter linearly through $\nabla^2 \phi_0$: the Poisson equation describes the structure of the numerical bias of both stochastic gradient integrators, beyond its role in the upper bounds on the MSE. The gap between Theorem~\ref{thm: sgld bias}, which holds at every step size with $\phi_h$, and Corollary~\ref{cor: sgld coefficient}, which reaches the Hessian form with $\phi_0$ only in the limit $h \to 0$, also shows that the discrete Poisson equation requires fewer derivatives than the continuous one for nonasymptotic bounds on the numerical bias.

Through the Riesz representors of $\mathcal C(x)$ and $\nabla^2 U(x)$, the coefficients of $h$ in the Wasserstein-1 errors are explicit: $\mathcal W_1(\pi_h, \pi)/h$ and $\mathcal W_1(\bar\pi_h, \pi)/h$ converge to $\frac12 \norm{\psi_{\mathcal C}\, \pi}_{\mathrm{KR}}$ and $\frac12 \norm{(\psi_{\nabla^2 U} + \psi_{\mathcal C})\, \pi}_{\mathrm{KR}}$ (Corollary~\ref{cor: w1 limits}), as a consequence of the convergence of the normalized biases $(\pi_h - \pi)/h$ and $(\bar\pi_h - \pi)/h$ in the Kantorovich--Rubinstein norm (Theorem~\ref{thm: kr expansion}). For a strongly convex potential, the coefficient of SGLM is at most that of SGLD once the batch size reaches a threshold linear in $\sup_x \norm{\mathcal C(x)}$, with no explicit dependence on the dimension (Theorem~\ref{thm: batch threshold}), and Lemma~\ref{lem: w1 counterexample} shows that the order can be reversed at small batch sizes.

The coefficients also lead to approaches to bias correction in Section~\ref{sec: representors}: the explicit coefficient, Richardson--Romberg extrapolation and the shifting correction. The shifting correction replaces a sample $Z \sim \pi_h$ of SGLM by $Z + h \nabla \varphi_{\mathcal C}(Z)$ and formally removes the first order bias for all test functions at once. It is completely data-driven, since $\varphi_{\mathcal C}$ is approximated by minimizing the energy \eqref{eq: extrapolated energy} with the samples of SGLM at the step sizes $h$ and $2h$.

The numerical experiments of Section~\ref{sec: numerics} illustrate the coefficients, the bias correction and the comparison with SGLD. In a 4D Gaussian mixture, the evaluated coefficients of the numerical bias approach the exact coefficients \eqref{eq: bias coefficient} and \eqref{eq: sgld coefficient}, within $3.8\%$ at $h = 2^{-5}$. For a radial potential, the evaluated coefficients of the Wasserstein-1 errors converge to those of Corollary~\ref{cor: w1 limits} with a remainder of order $h^2$, and the coefficient of SGLM is below that of SGLD for $d = 2, 4, \dots, 256$. In a 16D product of ring distributions, the bias corrections reduce the largest absolute numerical bias of SGLM over the test functions by factors between $14.8$ and $71.1$. In Bayesian logistic regression on MNIST, SGLM has a smaller estimated MSE than SGLD for all test functions and batch sizes, and the median reduction grows with the batch size, from $34.6\%$ at $P = 1024$ to $63.3\%$ at $P = 4096$, mainly through a smaller bias, in line with the coefficient of SGLM being divided by $P$.

Four extensions are possible. The first is variance reduction: with variance-reduced gradients, the first order term of the bias of SG-UBU becomes of second order once $h$ is below a threshold \citep{lu2025mse}, and the same construction applied to \eqref{eq: sglm} would give an integrator of the overdamped dynamics with a bias of second order. The second is the remainder $C h^{\frac32} (1 + \log\frac1h)$ of \eqref{eq: w1 limits}, which is not guaranteed to be sharp in the discretization error, and the error order could be enhanced to $h^2$. The third is the underdamped Langevin dynamics, where a modified class for its stochastic gradient integrators would extend the analysis to that setting. The fourth is a rigorous study of the shifting correction, whose order $h^2$ is formal, together with the stability of the minimization of \eqref{eq: extrapolated energy}, in which the fluctuation of the empirical averages is divided by $h$.

\subsection*{Funding statement}
Lin gratefully acknowledges support from the National Science Foundation (DMS-2533878, DMS-2053746, DMS-2134209, ECCS-2328241, CBET-2347401 and OAC-2311848), and U.S.~Department of Energy (DOE) Office of Science Advanced Scientific Computing Research program under the ``Uncertainty Quantification for Multifidelity Operator Learning (MOLUcQ)'' project (Project No.~81739), DE-SC0023161, the SciDAC LEADS Institute, and DOE--Fusion Energy Science, under grant number: DE-SC0024583.

\subsection*{Acknowledgments}
The authors thank Dr. He Zhang for valuable suggestions on this paper.

\subsection*{Data availability}

The code of the numerical experiments is available at \url{https://github.com/xuda-ye-math/Poisson-SGLM}, with the commands that regenerate the numerical results.

\subsection*{Declaration of AI use}
The authors used Claude Opus 5.5 and GPT-6 Astra in deriving part of the results, in checking the correctness of the proofs, and in writing the code of the numerical experiments. Claude Opus 5.5 provided the counterexample in Lemma~\ref{lem: w1 counterexample}. The authors hold full responsibility for the content of this paper.

\bibliographystyle{plainnat}
\bibliography{references}

\appendix
\renewcommand{\theHsection}{appendix.\Alph{section}}
\renewcommand{\theHequation}{appendix.\Alph{section}.\arabic{equation}}

\section{Proof of Theorem~\ref{thm: u decay}}
\label{app: elliptic}

This appendix proves Theorem~\ref{thm: u decay}: the exponential decay of the Kolmogorov solution $u$ and of its derivatives of orders up to $4$. The decay of $u$ and of its derivatives of orders up to $3$ is proved in \citep[Theorem~1]{ye2026mse}, so essentially only the derivative of order $4$ needs to be estimated.\footnote{The paper \citep{ye2026mse} is an earlier work of the first author, Xuda Ye, and is under review at the IMA Journal of Numerical Analysis. Beyond it, the proofs of this paper use no self-cited results that are not yet peer-reviewed.} The bound \eqref{eq: poisson bound} on $\phi_h$ then follows by summation over the time steps. Appendix~\ref{sec: variations} writes out the variation processes and their pathwise bounds. Appendix~\ref{sec: semigroup} states the estimates on the semigroup $e^{t\mathcal L}$ that the proof uses and proves the short-time bound \eqref{eq: short time}. For a function $g$ with bounded derivatives of orders $1$ to $4$, Appendix~\ref{sec: fourth} bounds $\nabla^4 (e^{t\mathcal L} g)$ by differentiating the derivative formula \eqref{eq: BEL} $3$ times. Appendix~\ref{sec: elliptic conclusion} combines this bound with the decay \eqref{eq: grad decay} and the short-time bound \eqref{eq: short time}.

The proof is written for the potential $U$ and the test function $f$ of Assumptions~\ref{ass: U} and \ref{ass: f}, and it uses the derivatives of $U$ up to order $5$. Throughout this appendix, the constant $C$ depends only on $m$, $M$, $R$ and $d$, unless otherwise specified.

\subsection{Variation processes}
\label{sec: variations}

Let $X_t$ be the solution of \eqref{eq: langevin} started at $X_0 = x$, and write $\partial_w$ for the derivative with respect to the initial state $x$ along a vector $w$. For a function $g$ on $\mathbb R^d$, the derivatives of $e^{t\mathcal L} g$ in the initial state are expressed through the variation processes. The variation processes of orders $1$ to $4$ are
\begin{equation}
    J_t := \partial_x X_t , \qquad
    K_t[w] := \partial_w J_t , \qquad
    \Theta_t[w,w'] := \partial_{w'} K_t[w] , \qquad
    \Xi_t[w,w',w''] := \partial_{w''} \Theta_t[w,w'] ,
    \label{eq: variation processes}
\end{equation}
where the variations of orders $2$ to $4$ are taken along vectors $w, w', w'' \in \mathbb R^d$, so that, like $J_t$, they take values in $\mathbb R^{d \times d}$. From the definition \eqref{eq: variation processes},
\begin{equation*}
\begin{gathered}
    J_t w = \partial_w X_t , \qquad K_t[w] w' = \partial_w \partial_{w'} X_t , \\
    \Theta_t[w,w'] w'' = \partial_w \partial_{w'} \partial_{w''} X_t , \qquad \Xi_t[w,w',w''] w''' = \partial_w \partial_{w'} \partial_{w''} \partial_{w'''} X_t ,
\end{gathered}
\end{equation*}
hence $K_t[w] w'$, $\Theta_t[w,w'] w''$ and $\Xi_t[w,w',w''] w'''$ are symmetric in all their directions. Since the noise of \eqref{eq: langevin} is additive, differentiating the dynamics in its initial state gives ordinary differential equations along the path, with the initial values $J_0 = I_d$ and $K_0 = \Theta_0 = \Xi_0 = 0$:
\begin{align}
    \frac{\D J_t}{\D t} &= -\nabla^2 U(X_t)\, J_t , \label{eq: first variation}\\
    \frac{\D K_t[w]}{\D t} &= -\nabla^2 U(X_t)\, K_t[w] - \nabla^3 U(X_t)[J_t w]\, J_t , \label{eq: second variation}\\
    \frac{\D \Theta_t[w,w']}{\D t} &= -\nabla^2 U(X_t)\, \Theta_t[w,w'] - \nabla^3 U(X_t)[J_t w]\, K_t[w'] - \nabla^3 U(X_t)[J_t w']\, K_t[w] \notag\\
    &\qquad - \nabla^3 U(X_t)\bigl[ K_t[w] w' \bigr]\, J_t - \nabla^4 U(X_t)[J_t w, J_t w']\, J_t . \label{eq: third variation}
\end{align}
These equations are derived in \citep[Equations~(A.1), (A.8) and~(A.9)]{ye2026mse}. Differentiating \eqref{eq: third variation} once more, along $w''$, gives the equation of the variation $\Xi_t$ of order $4$,
\begin{align}
    & \frac{\D \Xi_t[w,w',w'']}{\D t} = -\nabla^2 U(X_t)\, \Xi_t[w,w',w'']
    - \sum_{(2,1)} \nabla^3 U(X_t)[J_t w'']\, \Theta_t[w,w']
    - \sum_{(2,1)} \nabla^3 U(X_t)\bigl[ K_t[w] w' \bigr]\, K_t[w''] \notag \\
    &\qquad - \nabla^3 U(X_t)\bigl[ \Theta_t[w,w'] w'' \bigr]\, J_t
    - \sum_{(2,1)} \nabla^4 U(X_t)[J_t w, J_t w']\, K_t[w''] \notag \\
    &\qquad - \sum_{(2,1)} \nabla^4 U(X_t)\bigl[ K_t[w] w', J_t w'' \bigr]\, J_t
    - \nabla^5 U(X_t)[J_t w, J_t w', J_t w'']\, J_t ,
    \label{eq: variations}
\end{align}
where, for an expression $F$ of the $3$ directions in which the underlined pair forms a block, the sum runs over the partitions of $\{w,w',w''\}$ into a pair and a singleton,
\begin{equation}
    \sum_{(2,1)} F(\underline{w,w'}, w'') := F(\underline{w,w'}, w'') + F(\underline{w,w''}, w') + F(\underline{w',w''}, w).
    \label{eq: sum21}
\end{equation}
In each sum $\sum_{(2,1)}$ of \eqref{eq: variations}, the summand written out has the pair $\{w,w'\}$ and the singleton $w''$.

The bounds of the next lemma are Gr\"onwall estimates, which grow exponentially in $t$. The factor $t$ in the bounds on $K_t$, $\Theta_t$ and $\Xi_t$ is what the proof of Lemma~\ref{lem: fourth from lipschitz} uses at short times.

\begin{lemma}
\label{lem: variation bounds}
    Let Assumption~\ref{ass: U} hold. Then, pathwise, for all $t \Ge 0$ and unit vectors $w, w', w''$,
    \begin{equation}
    \begin{gathered}
        \norm{J_t} \Le e^{M t} , \qquad
        \norm{K_t[w]} \Le M t\, e^{2M t} , \qquad
        \norm{\Theta_t[w,w']} \Le M(1+3M t)\, t\, e^{3M t} , \\
        \norm{\Xi_t[w,w',w'']} \Le C (1+t)^2\, t\, e^{4M t} ,
    \end{gathered}
    \label{eq: variation bounds}
    \end{equation}
    with $C$ depending only on $M$.
\end{lemma}

\begin{proof}
    The bounds on $J_t$, $K_t$ and $\Theta_t$ are \citep[Lemma~6]{ye2026mse}. For $\Xi_t[w,w',w'']$, the matrix $J_t J_s^{-1}$ solves \eqref{eq: first variation} on $[s,t]$ from $I_d$, so $\norm{J_t J_s^{-1}} \Le e^{M(t-s)}$ by Gr\"onwall's inequality, and Duhamel's principle with $\Xi_0 = 0$ gives
    \begin{equation*}
        \Xi_t[w,w',w''] = \int_0^t J_t J_s^{-1} \bigg( \frac{\D \Xi_s[w,w',w'']}{\D s} + \nabla^2 U(X_s)\, \Xi_s[w,w',w''] \bigg) \D s .
    \end{equation*}
    The factor $\frac{\D \Xi_s[w,w',w'']}{\D s} + \nabla^2 U(X_s)\, \Xi_s[w,w',w'']$ of the integrand is the sum of the $14$ terms of \eqref{eq: variations} other than $-\nabla^2 U(X_s)\, \Xi_s[w,w',w'']$. The bounds \eqref{eq: variation bounds} on $J_s$, $K_s$ and $\Theta_s$, with $\norm{\nabla^i U} \Le M$ for $i = 3, 4, 5$, bound these terms group by group: the $5$ summands of the next display bound, in this order, the sum carrying $\nabla^3 U$ and $\Theta_s$, the sum carrying $\nabla^3 U$ and $2$ factors $K_s$, the term carrying $\nabla^3 U$ and $\Theta_s[w,w'] w''$, the $2$ sums carrying $\nabla^4 U$ together, and the term carrying $\nabla^5 U$:
    \begin{equation*}
    \begin{aligned}
        & \quad \biggl\| \frac{\D \Xi_s[w,w',w'']}{\D s} + \nabla^2 U(X_s)\, \Xi_s[w,w',w''] \biggr\| \\
        & \Le 3M^2 (1+3Ms)\, s\, e^{4Ms} + 3M^3 s^2 e^{4Ms} + M^2 (1+3Ms)\, s\, e^{4Ms} + 6M^2 s\, e^{4Ms} + M e^{4Ms} \\
        & = M \bigl( 1 + 10Ms + 15M^2 s^2 \bigr) e^{4Ms} \Le 15M (1+Ms)^2 e^{4Ms} .
    \end{aligned}
    \end{equation*}
    Hence
    \begin{equation*}
        \norm{\Xi_t[w,w',w'']} \Le \int_0^t e^{M(t-s)}\, 15M (1+Ms)^2 e^{4Ms} \,\D s \Le 15M (1+Mt)^2\, t\, e^{4Mt} \Le C (1+t)^2\, t\, e^{4Mt} .
    \end{equation*}
\end{proof}

\subsection{Semigroup estimates}
\label{sec: semigroup}

The estimates of this subsection hold for every $g \in C^4(\mathbb R^d)$ with bounded derivatives of orders $1$ to $4$; the rate $\lambda > 0$ depends only on $m$, $M$, $R$ and $d$. Lemmas~\ref{lem: semigroup}, \ref{lem: BEL} and \ref{lem: variance bound} are proved in \citep[Appendix~A]{ye2026mse} under the bounds of Assumption~\ref{ass: U} of orders up to $4$, and are restated without proof, since the proof of Lemma~\ref{lem: fourth from lipschitz} uses them explicitly. Lemma~\ref{lem: short time} extends the short-time bounds of \citep[Equations~(A.11) and~(A.12)]{ye2026mse} to the derivative of order $4$ and is proved here.

\begin{lemma}
\label{lem: semigroup}
    Let Assumption~\ref{ass: U} hold. Then, for all $t \Ge 0$,
    \begin{equation}
        \bigl\| \nabla (e^{t\mathcal L} g) \bigr\|_\infty \Le C e^{-\lambda t} \norm{\nabla g}_\infty ;
        \label{eq: grad decay}
    \end{equation}
    and for $0 < t \Le 1$,
    \begin{equation}
        \bigl\| \nabla^2 (e^{t\mathcal L} g) \bigr\|_\infty \Le C\, t^{-\frac12}\, \norm{\nabla g}_\infty , \qquad
        \bigl\| \nabla^3 (e^{t\mathcal L} g) \bigr\|_\infty \Le C\, t^{-1}\, \norm{\nabla g}_\infty .
        \label{eq: smoothing 23}
    \end{equation}
\end{lemma}

The decay \eqref{eq: grad decay} is \citep[Equation~(A.7)]{ye2026mse}, obtained from the reflection coupling of the overdamped Langevin dynamics \citep{eberle2016reflection}. The smoothing bounds \eqref{eq: smoothing 23} are \citep[Lemmas~8 and~9]{ye2026mse}, obtained by differentiating the derivative formula \eqref{eq: BEL} below up to $2$ times in the initial state.

The derivative formula \citep{elworthy1994formulae}, in the form of \citep[Lemma~7]{ye2026mse}, expresses the derivative of $e^{2\tau\mathcal L} g$ along a vector $v$ through a stochastic integral of the variation process $J_s$ of order $1$.

\begin{lemma}
\label{lem: BEL}
    Let Assumption~\ref{ass: U} hold. Then, for $\tau > 0$, $\psi := e^{\tau\mathcal L} g$ and $v \in \mathbb R^d$,
    \begin{equation}
        v \cdot \nabla (e^{2\tau\mathcal L} g)(x) = \frac{1}{\sqrt2\,\tau}\, \mathbb E^x \biggl[ \psi(X_\tau) \int_0^\tau (J_s v) \cdot \D B_s \biggr] , \qquad \forall x \in \mathbb R^d .
        \label{eq: BEL}
    \end{equation}
\end{lemma}

The next lemma, \citep[Equation~(A.17)]{ye2026mse} in the proof of \citep[Lemma~8]{ye2026mse}, bounds the variance of $\psi(X_\tau)$.

\begin{lemma}
\label{lem: variance bound}
    Let Assumption~\ref{ass: U} hold. Then, for $\tau > 0$, $\psi := e^{\tau\mathcal L} g$ and $X_0 = x \in \mathbb R^d$,
    \begin{equation}
        \Var\bigl( \psi(X_\tau) \bigr) \Le 2 d\, \tau\, e^{4M\tau}\, \norm{\nabla g}_\infty^2 .
        \label{eq: variance bound}
    \end{equation}
\end{lemma}
The randomness of $\psi(X_\tau)$ comes from the Brownian motion $(B_t)_{t \Ge 0}$ on $[0,\tau]$, which drives $X_t$ from $X_0 = x$.

The last estimate of this subsection bounds the derivatives of $e^{t\mathcal L} g$ over a finite time; it covers the times $t \Le 1$ in the proof of Theorem~\ref{thm: u decay} and supplies the derivative bounds needed in the proof of Lemma~\ref{lem: fourth from lipschitz} and in the case $t \Ge 1$ of that theorem.

\begin{lemma}
\label{lem: short time}
    Let Assumption~\ref{ass: U} hold. There is $C$, depending only on $M$, such that for every $g \in C^4(\mathbb R^d)$ with bounded derivatives of orders $1$ to $4$ and all $t \Ge 0$,
    \begin{equation}
        \bigl\| \nabla^j (e^{t\mathcal L} g) \bigr\|_\infty \Le C (1+t)^3 e^{4M t} \sum_{i=1}^{j} \norm{\nabla^i g}_\infty , \qquad j = 1,2,3,4 .
        \label{eq: short time}
    \end{equation}
\end{lemma}

\begin{proof}
    For $j \Le 3$, \eqref{eq: short time} follows from the bounds \citep[Equations~(A.11) and~(A.12)]{ye2026mse},
    \begin{equation*}
    \begin{gathered}
        \bigl\| \nabla (e^{t\mathcal L} g) \bigr\|_\infty \Le e^{Mt} \norm{\nabla g}_\infty , \qquad
        \bigl\| \nabla^2 (e^{t\mathcal L} g) \bigr\|_\infty \Le e^{2Mt} \norm{\nabla^2 g}_\infty + M t\, e^{2Mt} \norm{\nabla g}_\infty , \\
        \bigl\| \nabla^3 (e^{t\mathcal L} g) \bigr\|_\infty \Le e^{3Mt} \norm{\nabla^3 g}_\infty + 3M t\, e^{3Mt} \norm{\nabla^2 g}_\infty + M (1+3Mt)\, t\, e^{3Mt} \norm{\nabla g}_\infty ,
    \end{gathered}
    \end{equation*}
    whose coefficients are each at most $C (1+t)^3 e^{4Mt}$.

    For $j = 4$, fix $x \in \mathbb R^d$ and unit vectors $w, w', w'', w'''$. By \eqref{eq: variation processes}, the derivatives of $X_t$ in its initial state along these vectors are $J_t w$, $K_t[w] w'$, $\Theta_t[w,w'] w''$ and $\Xi_t[w,w',w''] w'''$, together with the same expressions along the other subsets of $\{w,w',w'',w'''\}$. The chain rule for higher derivatives applied to $g(X_t)$, and the exchange of the derivatives with the expectation, justified by the pathwise bounds \eqref{eq: variation bounds} and the bounded derivatives of $g$, give
    \begin{align}
        & \nabla^4 (e^{t\mathcal L} g)(x)[w,w',w'',w'''] = \mathbb E^x \biggl[  \nabla^4 g(X_t)[J_t w, J_t w', J_t w'', J_t w'''] + \sum_{(2,1,1)} \nabla^3 g(X_t)\bigl[ K_t[w] w', J_t w'', J_t w''' \bigr] \notag \\
        & \quad + \sum_{(3,1)} \nabla^2 g(X_t)\bigl[ \Theta_t[w,w'] w'', J_t w''' \bigr] + \sum_{(2,2)} \nabla^2 g(X_t)\bigl[ K_t[w] w', K_t[w''] w''' \bigr] \notag \\
        & \quad + \nabla g(X_t) \cdot \Xi_t[w,w',w''] w''' \biggr] ,
        \label{eq: fourth chain rule}
    \end{align}
    where, for an expression $F$ of the $4$ directions in which each underlined group of directions forms a block, the sums run over the partitions of $\{w,w',w'',w'''\}$ into blocks of the sizes in the subscript:
    \begin{align*}
        \sum_{(2,1,1)} F(\underline{w,w'}, w'', w''') &:= F(\underline{w,w'}, w'', w''') + F(\underline{w,w''}, w', w''') + F(\underline{w,w'''}, w', w'') \\[-8pt]
        & \quad + F(\underline{w',w''}, w, w''') + F(\underline{w',w'''}, w, w'') + F(\underline{w'',w'''}, w, w') , \\[8pt]
        \sum_{(3,1)} F(\underline{w,w',w''}, w''') &:= F(\underline{w,w',w''}, w''') + F(\underline{w,w',w'''}, w'') + F(\underline{w,w'',w'''}, w') + F(\underline{w',w'',w'''}, w) , \\
        \sum_{(2,2)} F(\underline{w,w'}, \underline{w'',w'''}) &:= F(\underline{w,w'}, \underline{w'',w'''}) + F(\underline{w,w''}, \underline{w',w'''}) + F(\underline{w,w'''}, \underline{w',w''}) .
    \end{align*}

    By \eqref{eq: variation bounds}, the terms of the $5$ groups in \eqref{eq: fourth chain rule} are bounded by
    \begin{equation*}
    \begin{gathered}
        e^{4Mt} \norm{\nabla^4 g}_\infty , \qquad M t\, e^{4Mt} \norm{\nabla^3 g}_\infty , \qquad M(1+3Mt)\, t\, e^{4Mt} \norm{\nabla^2 g}_\infty , \\
        M^2 t^2 e^{4Mt} \norm{\nabla^2 g}_\infty , \qquad C (1+t)^2\, t\, e^{4Mt} \norm{\nabla g}_\infty ,
    \end{gathered}
    \end{equation*}
    each at most $C (1+t)^3 e^{4Mt}$ times the norm of the derivative of $g$ it carries. Summing the $15$ terms and taking the supremum over $x$ and over the unit vectors gives \eqref{eq: short time} for $j = 4$.
\end{proof}

\subsection{The derivative of order $4$}
\label{sec: fourth}

This subsection extends the smoothing bounds \eqref{eq: smoothing 23} to the derivative of order $4$; the proof uses Lemmas~\ref{lem: variation bounds}, \ref{lem: semigroup}, \ref{lem: BEL}, \ref{lem: variance bound} and \ref{lem: short time}.

\begin{lemma}
\label{lem: fourth from lipschitz}
    Let Assumption~\ref{ass: U} hold. For every $g \in C^4(\mathbb R^d)$ with bounded derivatives of orders $1$ to $4$ and every $0 < t \Le 1$,
    \begin{equation}
        \bigl\| \nabla^4 (e^{t\mathcal L} g)(x) \bigr\| \Le C\, t^{-\frac32}\, \norm{\nabla g}_\infty , \qquad \forall x \in \mathbb R^d .
        \label{eq: fourth from lipschitz}
    \end{equation}
\end{lemma}

\begin{proof}
    Set $\tau = t/2 \Le \frac12$ and $\psi = e^{\tau\mathcal L} g$, so that $e^{2\tau\mathcal L} g = e^{t\mathcal L} g$, and fix $x \in \mathbb R^d$ and unit vectors $v, w, w', w''$. By \eqref{eq: short time}, $\psi$ has bounded derivatives of orders $1$ to $4$. Differentiate \eqref{eq: BEL} in Lemma~\ref{lem: BEL} along $w$, $w'$ and $w''$. By \eqref{eq: variation processes}, a derivative that falls on the stochastic integral raises the order of the variation process in its integrand, from $J_s v$ to $K_s[w] v$, to $\Theta_s[w,w'] v$ and to $\Xi_s[w,w',w''] v$. Distributing the $3$ derivatives between $\psi(X_\tau)$ and the stochastic integral gives
    \begin{align}
        & \sqrt2\, \tau\, \nabla^4 (e^{t\mathcal L} g)(x)[v,w,w',w'']  = \mathbb E^x \biggl[ \partial_w \partial_{w'} \partial_{w''} \psi(X_\tau) \int_0^\tau (J_s v) \cdot \D B_s \biggr] \notag \\
        & \quad + \sum_{(2,1)} \mathbb E^x \biggl[ \partial_w \partial_{w'} \psi(X_\tau) \int_0^\tau \bigl( K_s[w''] v \bigr) \cdot \D B_s \biggr] + \sum_{(2,1)} \mathbb E^x \biggl[ \partial_{w''} \psi(X_\tau) \int_0^\tau \bigl( \Theta_s[w,w'] v \bigr) \cdot \D B_s \biggr] \notag \\
        & \quad + \mathbb E^x \biggl[ \Bigl( \psi(X_\tau) - \mathbb E^x \bigl[ \psi(X_\tau) \bigr] \Bigr) \int_0^\tau \bigl( \Xi_s[w,w',w''] v \bigr) \cdot \D B_s \biggr] ,
        \label{eq: BEL third}
    \end{align}
    where $\sum_{(2,1)}$ is the sum over the partitions of $\{w,w',w''\}$ defined in \eqref{eq: sum21}.
    
    In the term of \eqref{eq: BEL third} carrying $\Xi_s$, $\psi(X_\tau)$ is centered at its mean $\mathbb E^x \bigl[ \psi(X_\tau) \bigr] = e^{t\mathcal L} g(x)$. The centering does not change the value of this term, since the stochastic integral has mean $0$, and it makes \eqref{eq: variance bound} applicable. The exchanges of the derivatives with the expectation and with the stochastic integral are justified as in the proof of \citep[Lemma~8]{ye2026mse}. The variation processes are differentiable in the initial state along each path and bounded by \eqref{eq: variation bounds} uniformly in $x$, so the difference quotients of the integrands converge in $L^2$ over $[0,\tau]$ and the paths, and by the It\^o isometry those of the stochastic integrals converge in $L^2$. Since $\psi$ has bounded derivatives of orders $1$ to $4$, the difference quotients of the factors in front of the stochastic integrals also converge in $L^2$, and the Cauchy--Schwarz inequality passes the limits through the expectation.

    The derivatives of $\psi(X_\tau)$ in \eqref{eq: BEL third} follow from the chain rule and \eqref{eq: variation processes}:
    \begin{equation*}
    \begin{gathered}
        \partial_{w''} \psi(X_\tau) = \nabla \psi(X_\tau) \cdot J_\tau w'' , \qquad
        \partial_w \partial_{w'} \psi(X_\tau) = \nabla^2 \psi(X_\tau)[J_\tau w, J_\tau w'] + \nabla \psi(X_\tau) \cdot K_\tau[w] w' , \\
        \partial_w \partial_{w'} \partial_{w''} \psi(X_\tau) = \nabla^3 \psi(X_\tau)[J_\tau w, J_\tau w', J_\tau w''] + \sum_{(2,1)} \nabla^2 \psi(X_\tau)\bigl[ K_\tau[w] w', J_\tau w'' \bigr] + \nabla \psi(X_\tau) \cdot \Theta_\tau[w,w'] w'' .
    \end{gathered}
    \end{equation*}
    Since $\tau \Le \frac12$, \eqref{eq: short time} and \eqref{eq: smoothing 23} at time $\tau$ give
    \begin{equation*}
        \norm{\nabla\psi}_\infty \Le C\norm{\nabla g}_\infty , \qquad
        \norm{\nabla^2\psi}_\infty \Le C\tau^{-\frac12}\norm{\nabla g}_\infty , \qquad
        \norm{\nabla^3 \psi}_\infty \Le C \tau^{-1}\norm{\nabla g}_\infty ,
    \end{equation*}
    and \eqref{eq: variation bounds} gives $\norm{J_s} \Le C$ and $\norm{K_s[w]}, \norm{\Theta_s[w,w']}, \norm{\Xi_s[w,w',w'']} \Le C s$ for $s \Le \tau$. The bounds on $\nabla^j \psi$ and on the variation processes then give
    \begin{equation}
        \bigl| \partial_{w''} \psi(X_\tau) \bigr| \Le C\norm{\nabla g}_\infty , \quad
        \bigl| \partial_w \partial_{w'} \psi(X_\tau) \bigr| \Le C\tau^{-\frac12}\norm{\nabla g}_\infty , \quad
        \bigl| \partial_w \partial_{w'} \partial_{w''} \psi(X_\tau) \bigr| \Le C\tau^{-1}\norm{\nabla g}_\infty ,
        \label{eq: psi derivative bounds}
    \end{equation}
    and \eqref{eq: variance bound} bounds the centered factor of the term carrying $\Xi_s$:
    \begin{equation}
        \Var\bigl( \psi(X_\tau) \bigr)^{\frac12} = \biggl( \mathbb E^x \Bigl| \psi(X_\tau) - \mathbb E^x \bigl[ \psi(X_\tau) \bigr] \Bigr|^2 \biggr)^{\frac12} \Le C \tau^{\frac12} \norm{\nabla g}_\infty .
        \label{eq: centered factor}
    \end{equation}
    By the It\^o isometry and \eqref{eq: variation bounds},
    \begin{equation}
    \begin{gathered}
        \mathbb E^x \biggl| \int_0^\tau (J_s v) \cdot \D B_s \biggr|^2 \Le C\tau , \qquad
        \mathbb E^x \biggl| \int_0^\tau \bigl( K_s[w''] v \bigr) \cdot \D B_s \biggr|^2 \Le C\tau^3 , \\
        \mathbb E^x \biggl| \int_0^\tau \bigl( \Theta_s[w,w'] v \bigr) \cdot \D B_s \biggr|^2 \Le C\tau^3 , \qquad
        \mathbb E^x \biggl| \int_0^\tau \bigl( \Xi_s[w,w',w''] v \bigr) \cdot \D B_s \biggr|^2 \Le C\tau^3 ,
    \end{gathered}
    \label{eq: ito bounds}
    \end{equation}
    and these bounds hold for every term of the sums $\sum_{(2,1)}$ in \eqref{eq: BEL third}, since $w$, $w'$ and $w''$ are arbitrary unit vectors. The Cauchy--Schwarz inequality, with \eqref{eq: psi derivative bounds}, \eqref{eq: ito bounds} and \eqref{eq: centered factor}, bounds the terms of the $4$ groups in \eqref{eq: BEL third} by
    \begin{equation*}
        C \tau^{-1} \cdot \tau^{\frac12} , \qquad C \tau^{-\frac12} \cdot \tau^{\frac32} , \qquad C \cdot \tau^{\frac32} , \qquad C \tau^{\frac12} \cdot \tau^{\frac32}
    \end{equation*}
    times $\norm{\nabla g}_\infty$. For $\tau \Le \frac12$, each is at most $C \tau^{-\frac12} \norm{\nabla g}_\infty$, so dividing \eqref{eq: BEL third} by $\sqrt2\, \tau$ and taking the supremum over the unit vectors gives
    \begin{equation*}
        \bigl\| \nabla^4 (e^{t\mathcal L} g)(x) \bigr\| \Le C \tau^{-\frac32} \norm{\nabla g}_\infty ,
    \end{equation*}
    which is \eqref{eq: fourth from lipschitz} with $\tau = t/2$.
\end{proof}

\subsection{Conclusion of the proof}
\label{sec: elliptic conclusion}

\begin{proof}[Proof of Theorem~\ref{thm: u decay}]
    The bound on $|u(x,t)|$ in \eqref{eq: u bound} and the bounds on $\nabla^j u(x,t)$ for $1 \Le j \Le 3$ are \citep[Theorem~1]{ye2026mse}, whose assumptions on $U$ and $f$ are the bounds of Assumptions~\ref{ass: U} and \ref{ass: f} up to the orders $5$ and $3$, respectively; decreasing $\lambda$ if necessary, these bounds and \eqref{eq: grad decay} hold with the same rate $\lambda$. For the order $4$, since $\pi(f)$ is a constant, $\nabla^4 u(x,t) = \nabla^4 (e^{t\mathcal L} f)(x)$, and we split at $t = 1$. For $t \Le 1$, the short-time bound \eqref{eq: short time} with $g = f$ and Assumption~\ref{ass: f} give
    \begin{equation*}
        \norm{\nabla^4 u(x,t)} \Le C \sum_{i=1}^{4} \norm{\nabla^i f}_\infty \Le C L \Le C L\, e^{\lambda}\, e^{-\lambda t} .
    \end{equation*}
    For $t \Ge 1$, the semigroup property gives $e^{t\mathcal L} f = e^{\mathcal L} \bigl( e^{(t-1)\mathcal L} f \bigr)$, and $e^{(t-1)\mathcal L} f$ has bounded derivatives of orders $1$ to $4$, since \eqref{eq: short time} at time $t - 1$ bounds them. Hence \eqref{eq: fourth from lipschitz} at time $1$, with $e^{(t-1)\mathcal L} f$ in place of $g$, bounds $\norm{\nabla^4 u(x,t)}$ by $C \norm{\nabla (e^{(t-1)\mathcal L} f)}_\infty$, and \eqref{eq: grad decay} at time $t - 1$ bounds this gradient, so that
    \begin{equation*}
        \norm{\nabla^4 u(x,t)} \Le C\, \bigl\| \nabla (e^{(t-1)\mathcal L} f) \bigr\|_\infty \Le C\, e^{-\lambda(t-1)}\, \norm{\nabla f}_\infty \Le C L\, e^{\lambda}\, e^{-\lambda t} .
    \end{equation*}
    Together, the cases $t \Le 1$ and $t \Ge 1$ give the bound on $\nabla^4 u(x,t)$ in \eqref{eq: u bound}.

    For the discrete Poisson solution \eqref{eq: discrete poisson solution}, the geometric bounds \eqref{eq: u bound} on the summands justify differentiating the series term by term. Summing the derivative bound in \eqref{eq: u bound} over $t = kh$ gives, for $\lambda h \Le 1$ and $1 \Le j \Le 4$,
    \begin{equation*}
        \norm{\nabla^j \phi_h(x)} \Le h \sum_{k=0}^{\infty} \norm{\nabla^j u(x,kh)} \Le C L\, h \sum_{k=0}^{\infty} e^{-\lambda k h} = \frac{C L\, h}{1 - e^{-\lambda h}} \Le \frac{2 C L}{\lambda} ,
    \end{equation*}
    which is \eqref{eq: poisson bound} after $2/\lambda$ is absorbed into $C$.
\end{proof}

\section{Proofs for Section~\ref{sec: modified}}
\label{app: modified}

Throughout this appendix, the constant $C$ depends only on $m$, $M$, $R$ and $d$. The next lemma records the dissipativity of $U$ used in the moment estimates, and the Gaussian tail of $\pi$ that follows from it.

\begin{lemma}
\label{lem: dissipative}
    Let Assumption~\ref{ass: U} hold. For all $x \in \mathbb R^d$,
    \begin{equation}
        \langle x, \nabla U(x) \rangle \Ge \frac{m}{4} |x|^2 - C , \qquad U(x) - U(0) \Ge \frac{m}{8} |x|^2 - C ,
        \label{eq: dissipative bounds}
    \end{equation}
    and consequently $\pi(x) \Le C e^{-\frac{m}{8}|x|^2}$ and $\int_{\mathbb R^d} |x|^q\, \pi(\D x) \Le C_q$ for every $q \Ge 1$.
\end{lemma}

\begin{proof}
    By \citep[Lemma~4]{ye2026mse}, with $R_1 := 4R(1+M/m)$,
    \begin{equation}
        \langle \nabla U(x) - \nabla U(y), x-y \rangle \Ge -M |x-y|^2 , \quad
        \langle \nabla U(x) - \nabla U(y), x-y \rangle \Ge \frac{m}{2} |x-y|^2 \ \text{ for } |x-y| \Ge R_1 .
        \label{eq: dissipative}
    \end{equation}
    Applying \eqref{eq: dissipative} to the pair $(x,0)$ and using $|\nabla U(0)| \Le M$,
    \begin{equation*}
        \langle x, \nabla U(x) \rangle = \langle x, \nabla U(x) - \nabla U(0) \rangle + \langle x, \nabla U(0) \rangle \Ge
        \begin{cases}
            \dfrac{m}{2} |x|^2 - M|x| , & |x| \Ge R_1 , \\[1ex]
            -M|x|^2 - M|x| , & |x| < R_1 ,
        \end{cases}
    \end{equation*}
    and both cases are at least $\frac{m}{4}|x|^2 - C$, which is the drift bound in \eqref{eq: dissipative bounds}. Integrating along the segment from $0$ to $x$, with \eqref{eq: dissipative} at the pair $(rx, 0)$ and the integral split at $r|x| = R_1$ when $|x| \Ge R_1$,
    \begin{equation*}
        U(x) - U(0) = \int_0^1 \langle \nabla U(rx), x \rangle \,\D r \Ge \frac{m}{4} \bigl( |x|^2 - R_1^2 \bigr) - \frac{M}{2} R_1^2 - M |x| ,
    \end{equation*}
    which gives the lower bound on $U(x) - U(0)$ in \eqref{eq: dissipative bounds}. The upper bound $U(y) - U(0) \Le M|y| + \frac{M}{2}|y|^2$ and the lower bound on $U(x) - U(0)$ in \eqref{eq: dissipative bounds} give
    \begin{equation*}
        \pi(x) = \frac{e^{-(U(x) - U(0))}}{\int_{\mathbb R^d} e^{-(U(y) - U(0))} \D y} \Le \frac{e^{C - \frac m8 |x|^2}}{\int_{\mathbb R^d} e^{-M|y| - \frac M2 |y|^2} \D y} \Le C e^{-\frac m8 |x|^2} ,
    \end{equation*}
    and hence $\int_{\mathbb R^d} |x|^q\, \pi(\D x) \Le C_q$ for every $q \Ge 1$.
\end{proof}

\subsection{Proof of Lemma~\ref{lem: modified class}: modified potential estimates}
\label{sec: modified class proof}

\begin{proof}[Proof of Lemma~\ref{lem: modified class}]
    The proof bounds the derivatives of $V$ and of $V\theta_h$, which gives \eqref{eq: modified class}, and then chooses $\delta$ so that $\tilde U_h$ satisfies \eqref{eq: modified class bounds} and the convexity condition. Under Assumption~\ref{ass: U}, $|\nabla U(x)| \Le M(1+|x|)$ and $\norm{\nabla^j U(x)} \Le M$ for $2 \Le j \Le 7$. The derivatives of $|\nabla U|^2$ follow from the Leibniz rule,
    \begin{equation}
        \nabla^i |\nabla U|^2 = \sum_{l=0}^{i} \binom{i}{l}\, \mathrm{Sym}\, \bigl\langle \nabla^{l+1} U ,\, \nabla^{i-l+1} U \bigr\rangle , \qquad 0 \Le i \Le 5 ,
        \label{eq: Leibniz}
    \end{equation}
    where $\langle \cdot , \cdot \rangle$ contracts one index of each factor, so that both sides are tensors of order $i$. The symmetrization is needed: $\langle \nabla^{l+1} U, \nabla^{i-l+1} U \rangle$ is symmetric within its first $l$ indices and within its last $i-l$ indices, but in general not across these groups.
    
    For $i \Ge 1$, each term of \eqref{eq: Leibniz} contains the factor $\nabla U$ at most once and is therefore bounded by $C(1+|x|)$; the derivatives $\nabla^i \Delta U$, contractions of $\nabla^{i+2} U$, are bounded by $dM$ for $0 \Le i \Le 5$. Hence $V = |\nabla U|^2 - \Delta U$ of \eqref{eq: V} satisfies
    \begin{equation}
        |V(x)| \Le C \bigl( 1 + |x| \bigr)^2 , \qquad \norm{\nabla^i V(x)} \Le C \bigl( 1 + |x| \bigr) , \quad 1 \Le i \Le 5 , \qquad \forall x \in \mathbb R^d ,
        \label{eq: V growth}
    \end{equation}
    where $C$ depends only on $M$ and $d$.

    Let $\delta \in (0,1]$, fixed by \eqref{eq: delta bound} at the end of the proof, and let $h \Le \delta$ throughout. By \eqref{eq: modified potential}, $\tilde U_h - U = \frac h4\, V \theta_h$ with $\theta_h(x) = \theta(hx/\delta)$. The function $\theta_h$ is supported in $\{|x| \Le 2\delta/h\}$, where
    \begin{equation}
        \frac h\delta \bigl( 1 + |x| \bigr) \Le \frac h\delta + 2 \Le 3 ,
        \label{eq: cutoff support}
    \end{equation}
    and $\nabla^l \theta_h(x) = (h/\delta)^l\, (\nabla^l \theta)(hx/\delta)$ is bounded by $C (h/\delta)^l$ for $0 \Le l \Le 5$, with $C$ depending on the fixed function $\theta$ only. By the Leibniz rule again,
    \begin{equation}
        \nabla^i (V \theta_h) = \sum_{l=0}^{i} \binom{i}{l}\, \mathrm{Sym} \bigl( \nabla^{i-l} V \otimes \nabla^l \theta_h \bigr) , \qquad 1 \Le i \Le 5 .
        \label{eq: Leibniz cutoff}
    \end{equation}
    On the support of $\theta_h$, \eqref{eq: V growth}, $h \Le \delta$ and \eqref{eq: cutoff support} give
    \begin{equation*}
        \norm{\nabla^{i-l} V(x) \otimes \nabla^l \theta_h(x)} \Le
        \begin{cases}
            C \bigl( 1 + |x| \bigr) , & 0 \Le l \Le i - 1 , \\
            C (h/\delta)^i \bigl( 1 + |x| \bigr)^2 \Le 3C \bigl( 1 + |x| \bigr) , & l = i .
        \end{cases}
    \end{equation*}
    Hence, by \eqref{eq: Leibniz cutoff}, the derivatives of $V\theta_h$ of orders $1$ to $5$ grow at most linearly on the support of $\theta_h$,
    \begin{equation}
        \bigl\| \nabla^i (V\theta_h)(x) \bigr\| \Le C \bigl( 1 + |x| \bigr) , \qquad 1 \Le i \Le 5 , \qquad \forall x\in \mathbb R^d: |x| \Le \frac{2\delta}{h} .
        \label{eq: V cutoff growth}
    \end{equation}
    Since $\nabla^i (V\theta_h)$ vanishes for $|x| > 2\delta/h$, \eqref{eq: V cutoff growth} gives \eqref{eq: modified class}, and, by \eqref{eq: cutoff support},
    \begin{equation}
        \norm{\nabla^i \tilde U_h(x) - \nabla^i U(x)} = \frac h4\, \bigl\| \nabla^i (V\theta_h)(x) \bigr\| \Le C \delta , \qquad 1 \Le i \Le 5 , \qquad \forall x \in \mathbb R^d ,
        \label{eq: modified derivative difference}
    \end{equation}
    with $C$ independent of $\delta$ and $h$. Choose $\delta \in (0,1]$ so small that this constant satisfies
    \begin{equation}
    C\delta \Le \min\Bigl\{ \frac m2 ,\, M \Bigr\} ,
    \label{eq: delta bound}
    \end{equation}
    and let $h_0 \Le \delta$. Then \eqref{eq: modified derivative difference} gives
    \begin{equation}
        \nabla^2 \tilde U_h(x) \succeq \nabla^2 U(x) - C\delta\, I_d , \qquad \norm{\nabla^i \tilde U_h(x)} \Le \norm{\nabla^i U(x)} + C\delta , \quad 1 \Le i \Le 5 , \qquad \forall x \in \mathbb R^d ,
        \label{eq: modified derivative perturbation}
    \end{equation}
    which, with Assumption~\ref{ass: U} and \eqref{eq: delta bound}, gives \eqref{eq: modified class bounds} and the convexity condition of Lemma~\ref{lem: modified class}. Finally, $\tilde U_h \in C^5(\mathbb R^d)$, since $\theta_h \in C^\infty(\mathbb R^d)$ and $U \in C^7(\mathbb R^d)$ gives $V \in C^5(\mathbb R^d)$.
\end{proof}

\subsection{Proof of Lemma~\ref{lem: shifted bias}: modified bias estimate}
\label{sec: shifted bias proof}

\begin{proof}[Proof of Lemma~\ref{lem: shifted bias}]
    The proof bounds the distance between $\tilde\pi_h * \mathcal N(0, \frac h2 I_d)$ and $\pi$, weighted by $1 + |x|$, by $C h^2$. It writes $\tilde\pi_h$ as $\pi\, e^{-(\tilde U_h - U)}$ divided by its normalizing constant, expands this constant to the order $h$, where, since $\pi(V) = 0$, the term of order $h$ comes only from the cutoff and is of order $h^2$, and expands $\tilde\pi_h * \mathcal N(0, \frac h2 I_d)$ to the order $4$ along the Gaussian shift. Since $\tilde\pi_h \propto e^{-\tilde U_h} = e^{-U} e^{-(\tilde U_h - U)}$, where $\tilde U_h - U = \frac h4\, V\theta_h$ by \eqref{eq: modified potential}, we have
    \begin{equation}
        \tilde\pi_h(x) = \frac{\pi(x)\, e^{-(\tilde U_h - U)(x)}}{\pi\bigl( e^{-(\tilde U_h - U)} \bigr)} .
        \label{eq: modified target ratio}
    \end{equation}
    By \eqref{eq: V growth}, $|\tilde U_h - U| \Le C h (1 + |x|^2)$, while $\pi(x) \Le C e^{-\frac{m}{8}|x|^2}$ by Lemma~\ref{lem: dissipative}. Shrinking $h_0$ so that $C h_0 \Le \frac{m}{16}$, with $C$ the constant of the bound on $|\tilde U_h - U|$, we have for $h \Le h_0$,
    \begin{equation}
        e^{|\tilde U_h - U|} \Le C e^{\frac{m}{16}|x|^2} , \qquad \pi\bigl( (1+|x|)^5\, e^{|\tilde U_h - U|} \bigr) \Le C .
        \label{eq: weighted moments}
    \end{equation}
    Taylor's formula for the exponential gives
    \begin{equation}
        e^{-(\tilde U_h - U)} = 1 - \bigl( \tilde U_h - U \bigr) + \varrho_2(x) , \quad |\varrho_2(x)| \Le \frac12 \bigl( \tilde U_h - U \bigr)^2 e^{|\tilde U_h - U|} \Le C h^2 \bigl( 1 + |x| \bigr)^4 e^{|\tilde U_h - U|} ,
        \label{eq: exponential taylor}
    \end{equation}
    hence the remainder $\varrho_2(x)$ satisfies
    \begin{equation}
        \pi\bigl( (1+|x|)\, |\varrho_2(x)| \bigr) \Le C h^2\, \pi\bigl( (1+|x|)^{5} e^{|\tilde U_h - U|} \bigr) \Le C h^2 .
        \label{eq: r bound}
    \end{equation}

    Next, $\pi(V) = - \pi(\mathcal L U) = 0$ by \eqref{eq: integration by parts} with $f = 1$ and \eqref{eq: V}.
    Replacing $V$ by $V\theta_h$ changes $\pi(V)$ by an exponentially small amount: since $1 - \theta_h$ vanishes on $\{|x| \Le \delta/h\}$,
    \begin{equation}
        \pi\bigl( (1+|x|)\, |V|\, (1 - \theta_h) \bigr) \Le C\, \pi\bigl( (1+|x|)^3\, \mathbf 1_{\{|x| \Ge \delta/h\}} \bigr) \Le C e^{-\frac{m \delta^2}{16 h^2}} \Le C h^2 ,
        \label{eq: cutoff tail}
    \end{equation}
    by \eqref{eq: V growth} and $\pi(x) \Le C e^{-\frac{m}{8}|x|^2} \Le C e^{-\frac{m \delta^2}{16 h^2}}\, e^{-\frac{m}{16}|x|^2}$ on $\{|x| \Ge \delta/h\}$.
    Inserting \eqref{eq: exponential taylor} into the denominator of \eqref{eq: modified target ratio}, and using $\pi(V) = 0$, \eqref{eq: r bound} and \eqref{eq: cutoff tail}, we obtain
    \begin{equation}
        \pi\bigl( e^{-(\tilde U_h - U)} \bigr) = 1 - \frac h4\, \pi(V \theta_h) + \pi(\varrho_2) = 1 + \frac h4\, \pi\bigl( V (1 - \theta_h) \bigr) + \pi(\varrho_2) = 1 + \mathcal O(h^2) ,
        \label{eq: ratio expansion}
    \end{equation}
    where $\mathcal O(h^2)$ is bounded by $C h^2$. After shrinking $h_0$, $\pi\bigl( e^{-(\tilde U_h - U)} \bigr) \Ge 1 - C h^2 \Ge \frac12$, hence
    \begin{equation}
        \biggl| \frac{1}{\pi\bigl( e^{-(\tilde U_h - U)} \bigr)} - 1 \biggr| \Le C h^2 .
        \label{eq: ratio bound}
    \end{equation}

    Let $\chi \sim \mathcal N(0,I_d)$, so that
    \begin{equation*}
        \biggl( \tilde\pi_h * \mathcal N\Bigl( 0, \frac h2 I_d \Bigr) \biggr)(x) = \mathbb E\biggl[ \tilde\pi_h\Bigl( x - \sqrt{\frac h2}\, \chi \Bigr) \biggr] .
    \end{equation*}
    Since $\tilde U_h \in C^5$ by Lemma~\ref{lem: modified class}, $\tilde\pi_h \in C^5$. Expanding $\tilde\pi_h$ about $x$ to the order $4$ along $-\sqrt{h/2}\, \chi$,
    \begin{equation*}
        \tilde\pi_h\Bigl( x - \sqrt{\frac h2}\, \chi \Bigr) = \tilde\pi_h(x) - \sqrt{\frac h2}\, \nabla \tilde\pi_h(x) \cdot \chi + \frac h4\, \nabla^2 \tilde\pi_h(x)[\chi,\chi] - \frac16 \Bigl( \frac h2 \Bigr)^{\frac32} \nabla^3\tilde\pi_h(x)\bigl[ \chi^{\otimes3} \bigr] + \varrho_4^{\chi}(x) ,
    \end{equation*}
    with
    \begin{equation}
        |\varrho_4^{\chi}(x)| \Le \frac{h^2}{24}\, |\chi|^4 \int_0^1 (1-s)^3\, \biggl\| \nabla^4 \tilde\pi_h\Bigl( x - s\sqrt{\frac h2}\, \chi \Bigr) \biggr\| \,\D s .
        \label{eq: fourth remainder}
    \end{equation}
    Under the expectation the terms of odd order in $\chi$ vanish and the term of order $2$ has mean $\frac h4\, \Delta \tilde\pi_h(x)$, hence
    \begin{equation}
        \biggl( \tilde\pi_h * \mathcal N\Bigl( 0, \frac h2 I_d \Bigr) \biggr)(x) = \tilde\pi_h(x) + \frac h4\, \Delta \tilde\pi_h(x) + \mathbb E\bigl[ \varrho_4^{\chi}(x) \bigr] .
        \label{eq: shift of density}
    \end{equation}
    Differentiating $\tilde\pi_h \propto e^{-\tilde U_h}$ to the order $2$,
    \begin{equation}
        \Delta \tilde\pi_h = \tilde\pi_h\, \bigl( |\nabla \tilde U_h|^2 - \Delta \tilde U_h \bigr) .
        \label{eq: laplacian of density}
    \end{equation}
    By \eqref{eq: modified target ratio}, \eqref{eq: shift of density} and \eqref{eq: laplacian of density},
    \begin{equation}
        \biggl( \tilde\pi_h * \mathcal N\Bigl( 0, \frac h2 I_d \Bigr) \biggr)(x) = \frac{\pi(x)\, e^{-(\tilde U_h - U)(x)}}{\pi\bigl( e^{-(\tilde U_h - U)} \bigr)} \Bigl( 1 + \frac h4 \bigl( |\nabla \tilde U_h(x)|^2 - \Delta \tilde U_h(x) \bigr) \Bigr) + \mathbb E\bigl[ \varrho_4^{\chi}(x) \bigr] .
        \label{eq: product form}
    \end{equation}
    By \eqref{eq: modified class bounds}, \eqref{eq: modified class} and \eqref{eq: V}, $\bigl| |\nabla \tilde U_h|^2 - \Delta \tilde U_h \bigr| \Le C (1 + |x|)^2$ and
    \begin{equation}
        \bigl| |\nabla \tilde U_h|^2 - \Delta \tilde U_h - V \bigr| \Le |\nabla \tilde U_h - \nabla U|\, |\nabla \tilde U_h + \nabla U| + |\Delta \tilde U_h - \Delta U| \Le C h \bigl( 1 + |x| \bigr)^2 .
        \label{eq: modified V gap}
    \end{equation}
    Taylor's formula for the exponential, with $\tilde U_h - U = \frac h4\, V\theta_h$, the bound $|V| \Le C(1+|x|)^2$ of \eqref{eq: V growth}, the bound on $|\nabla \tilde U_h|^2 - \Delta \tilde U_h$ and the bound on $e^{|\tilde U_h - U|}$ in \eqref{eq: weighted moments} give
    \begin{equation*}
        \biggl| e^{-(\tilde U_h - U)} \Bigl( 1 + \frac h4 \bigl( |\nabla \tilde U_h|^2 - \Delta \tilde U_h \bigr) \Bigr) - 1 - \frac h4 \bigl( |\nabla \tilde U_h|^2 - \Delta \tilde U_h - V\theta_h \bigr) \biggr| \Le C h^2 \bigl( 1 + |x| \bigr)^6 e^{\frac{m}{16}|x|^2} .
    \end{equation*}
    Since $V\theta_h = V - V(1 - \theta_h)$, \eqref{eq: modified V gap} gives
    \begin{equation}
        \biggl| e^{-(\tilde U_h - U)} \Bigl( 1 + \frac h4 \bigl( |\nabla \tilde U_h|^2 - \Delta \tilde U_h \bigr) \Bigr) - 1 \biggr|
        \Le C h^2 \bigl( 1 + |x| \bigr)^6 e^{\frac{m}{16}|x|^2} + \frac h4\, |V| \bigl( 1 - \theta_h \bigr) .
        \label{eq: factor bound}
    \end{equation}
    Multiplying \eqref{eq: factor bound} by $(1+|x|)\pi(x)$ and integrating over $\mathbb R^d$, the term $C h^2 (1+|x|)^6 e^{\frac{m}{16}|x|^2}$ contributes at most $C h^2 \int_{\mathbb R^d} (1 + |x|)^{7} e^{-\frac{m}{16}|x|^2} \,\D x \Le C h^2$ by the density bound $\pi(x) \Le C e^{-\frac{m}{8}|x|^2}$ of Lemma~\ref{lem: dissipative}, and the term $\frac h4\, |V| (1 - \theta_h)$ contributes at most $C h^3$ by \eqref{eq: cutoff tail}, so
    \begin{equation}
        \int_{\mathbb R^d} \bigl( 1 + |x| \bigr) \biggl| e^{-(\tilde U_h - U)} \Bigl( 1 + \frac h4 \bigl( |\nabla \tilde U_h|^2 - \Delta \tilde U_h \bigr) \Bigr) - 1 \biggr|\, \pi(x) \,\D x \Le C h^2 .
        \label{eq: factor integral}
    \end{equation}

    Differentiating $\tilde\pi_h \propto e^{-\tilde U_h}$ to the order $4$ and using \eqref{eq: modified class bounds} with $|\nabla \tilde U_h(y)| \Le 2M(1+|y|)$,
    \begin{equation*}
        \norm{\nabla^4 \tilde\pi_h(y)} \Le C \bigl( 1 + |y| \bigr)^4 \tilde\pi_h(y) , \qquad y \in \mathbb R^d .
    \end{equation*}
    By \eqref{eq: modified target ratio}, \eqref{eq: weighted moments}, the lower bound $\pi( e^{-(\tilde U_h - U)} ) \Ge \frac12$ and the density bound of Lemma~\ref{lem: dissipative}, $\tilde\pi_h\bigl( (1 + |x|)^q \bigr) \Le 2C\, \pi\bigl( (1 + |x|)^q e^{\frac{m}{16}|x|^2} \bigr) \Le C_q$ for $q \Ge 1$. For each fixed $s \in [0,1]$, the substitution $y = x - s\sqrt{h/2}\,\chi$ and the bound on $\tilde\pi_h\bigl( (1 + |x|)^q \bigr)$ with $q = 5$ give
    \begin{equation*}
        \int_{\mathbb R^d} \bigl( 1 + |x| \bigr) \biggl( 1 + \Bigl| x - s\sqrt{\frac h2}\, \chi \Bigr| \biggr)^4 \tilde\pi_h\Bigl( x - s\sqrt{\frac h2}\, \chi \Bigr) \,\D x
        \Le \bigl( 1 + |\chi| \bigr) \int_{\mathbb R^d} \bigl( 1 + |y| \bigr)^5 \tilde\pi_h(y) \,\D y \Le C \bigl( 1 + |\chi| \bigr) .
    \end{equation*}
    Hence, by \eqref{eq: fourth remainder}, Tonelli's theorem and $\mathbb E\bigl[ |\chi|^4 (1+|\chi|) \bigr] \Le C$,
    \begin{equation}
    \begin{aligned}
        & \quad \int_{\mathbb R^d} \bigl( 1 + |x| \bigr) \mathbb E|\varrho_4^{\chi}(x)| \,\D x \\
        & \Le C h^2\, \mathbb E\Biggl[ |\chi|^4 \int_0^1 (1-s)^3 \int_{\mathbb R^d} \bigl( 1 + |x| \bigr) \biggl( 1 + \Bigl| x - s\sqrt{\frac h2}\, \chi \Bigr| \biggr)^4 \tilde\pi_h\Bigl( x - s\sqrt{\frac h2}\, \chi \Bigr) \,\D x \,\D s \Biggr] \Le C h^2 .
    \end{aligned}
    \label{eq: remainder integral}
    \end{equation}

    By \eqref{eq: product form}, the difference between $\tilde\pi_h * \mathcal N(0, \frac h2 I_d)$ and $\pi$ splits into the error $\mathcal T_1$ of the normalizing constant, the error $\mathcal T_2$ of the factor and the remainder $\mathcal T_3$ of the expansion along the Gaussian shift,
    \begin{align}
        & \biggl( \tilde\pi_h * \mathcal N\Bigl( 0, \frac h2 I_d \Bigr) \biggr)(x) - \pi(x) \notag \\
        & = \underbrace{\biggl( \frac{1}{\pi\bigl( e^{-(\tilde U_h - U)} \bigr)} - 1 \biggr) \pi(x)}_{=: \mathcal T_1} + \underbrace{\frac{\pi(x)}{\pi\bigl( e^{-(\tilde U_h - U)} \bigr)} \biggl( e^{-(\tilde U_h - U)(x)} \Bigl( 1 + \frac h4 \bigl( |\nabla \tilde U_h(x)|^2 - \Delta \tilde U_h(x) \bigr) \Bigr) - 1 \biggr)}_{=: \mathcal T_2} + \underbrace{\mathbb E\bigl[ \varrho_4^{\chi}(x) \bigr]}_{=: \mathcal T_3} .
        \label{eq: density split}
    \end{align}
    Integrated against $1 + |x|$, each term of \eqref{eq: density split} is at most $C h^2$: $\mathcal T_1$ by \eqref{eq: ratio bound} and $\pi(1+|x|) \Le C$, $\mathcal T_2$ by the lower bound $\pi\bigl( e^{-(\tilde U_h - U)} \bigr) \Ge \frac12$ and \eqref{eq: factor integral}, and $\mathcal T_3$ by \eqref{eq: remainder integral}, which gives
    \begin{equation}
        \int_{\mathbb R^d} \bigl( 1 + |x| \bigr)\, \Biggl| \biggl( \tilde\pi_h * \mathcal N\Bigl( 0, \frac h2 I_d \Bigr) \biggr)(x) - \pi(x) \Biggr| \,\D x \Le C h^2 .
        \label{eq: density comparison}
    \end{equation}
    Finally, the symmetry of $\mathcal N(0, \frac h2 I_d)$ gives
    \begin{equation*}
        \tilde\pi_h(\tilde f_h) = \int_{\mathbb R^d} f(x)\, \biggl( \tilde\pi_h * \mathcal N\Bigl( 0, \frac h2 I_d \Bigr) \biggr)(x) \,\D x ,
    \end{equation*}
    and since both $\tilde\pi_h * \mathcal N(0, \frac h2 I_d)$ and $\pi$ are probability densities, $f$ may be replaced by $f - f(0)$, with $|f(x) - f(0)| \Le \norm{\nabla f}_\infty |x|$; \eqref{eq: density comparison} then gives \eqref{eq: shifted bias}.
\end{proof}

\section{Proofs for Section~\ref{sec: local}}
\label{app: local}

Statements about a single step are made at the index $0$, with the step $D_0 = Y_1 - Y_0$ of \eqref{eq: chain step} and $\chi_0 \sim \mathcal N(0, I_d)$ independent of $\mathcal F_0$, so that under $\mathbb E_0$ the state $Y_0$ is fixed and
\begin{equation*}
    \mathbb E_0[\chi_0] = 0 , \qquad \mathbb E_0\bigl[ \chi_0^{\otimes 2} \bigr] = I_d , \qquad \mathbb E_0\bigl[ \chi_0^{\otimes 3} \bigr] = 0 , \qquad \mathbb E_0 |\chi_0|^p \Le C , \quad p \Ge 1 .
\end{equation*}
The displacement of the modified dynamics started at $Y_0$ is written $\tilde W_s := \tilde Y_0(s) - Y_0$, $0 \Le s \Le h$. Throughout this appendix, the constant $C$ depends only on $m$, $M$, $R$ and $d$, and also on $p$ where a moment order $p$ appears. The largest moment order used is $p = 12$, for the bounds \eqref{eq: moment bounds} in the proof of Lemma~\ref{lem: moments}; all other proofs in this appendix use $p \Le 6$.

\subsection{Proof of Lemma~\ref{lem: moments}}
\label{sec: moments proof}

The proof uses the bounds \eqref{eq: exact sizes} on one step of the modified dynamics. Lemma~\ref{lem: exact moments} states them first, together with the moments of that step, which are used in Appendix~\ref{sec: matching proof}.

\begin{lemma}
\label{lem: exact moments}
    Let Assumption~\ref{ass: U} hold and let $\tilde U_h$ be the modified potential \eqref{eq: modified potential}. Then, for $0 \Le s \Le h$ and $p \Ge 1$, with $C$ depending only on $M$, $d$ and $p$ and $\pnorm{\cdot}_p$ taken under $\mathbb E_0$,
    \begin{equation}
        \pnorm{\tilde Y_0(s) - Y_0}_p \Le C \sqrt h \bigl( 1 + |Y_0| \bigr) , \qquad
        \pnorm{\tilde Y_0(h) - Y_0 - \sqrt{2h}\, \chi_0}_p \Le C h \bigl( 1 + |Y_0| \bigr) ,
        \label{eq: exact sizes}
    \end{equation}
    and, with $C$ depending only on $M$ and $d$, all derivatives evaluated at $Y_0$ and $1 \Le i,j,k \Le d$,
    \begin{align}
        \Bigl| \mathbb E_0\bigl[ \tilde Y_0(h) - Y_0 \bigr] + h\, \nabla \tilde U_h - \frac{h^2}{2} \bigl( \nabla^2 \tilde U_h\, \nabla \tilde U_h - \nabla \Delta \tilde U_h \bigr) \Bigr| &\Le C h^3 \bigl( 1+|Y_0| \bigr)^3 , \label{eq: exact m1}\\
        \Bigl\| \mathbb E_0\bigl[ (\tilde Y_0(h) - Y_0)^{\otimes 2} \bigr] - 2h\, I_d - h^2 \bigl( \nabla \tilde U_h^{\otimes 2} - 2\, \nabla^2 \tilde U_h \bigr) \Bigr\| &\Le C h^3 \bigl( 1+|Y_0| \bigr)^4 , \label{eq: exact m2}\\
        \Bigl| \bigl( \mathbb E_0\bigl[ (\tilde Y_0(h) - Y_0)^{\otimes 3} \bigr] \bigr)_{ijk} + 2h^2 \bigl( \delta_{ij}\, \partial_k \tilde U_h + \delta_{ik}\, \partial_j \tilde U_h + \delta_{jk}\, \partial_i \tilde U_h \bigr) \Bigr| &\Le C h^3 \bigl( 1+|Y_0| \bigr)^{6} . \label{eq: exact m3}
    \end{align}
    The same statements hold with $U$ in place of $\tilde U_h$ and $X_s$ in place of $\tilde Y_0(s)$, for the exact dynamics \eqref{eq: langevin} started at $X_0 = Y_0$, with the same constants.
\end{lemma}

\begin{proof}
    By Lemma~\ref{lem: modified class}, $|\nabla \tilde U_h(0)| \Le 2M$ and $\norm{\nabla^2 \tilde U_h} \Le 2M$, so that $|\nabla \tilde U_h(x)| \Le 2M(1+|x|)$ for every $x$. The displacement of the modified dynamics \eqref{eq: modified langevin} started at $Y_0$ is
    \begin{equation*}
        \tilde W_s = \tilde Y_0(s) - Y_0 = -\int_0^s \nabla \tilde U_h\bigl( Y_0 + \tilde W_r \bigr) \D r + \sqrt2\, B_s , \qquad 0 \Le s \Le h ,
    \end{equation*}
    where $(B_t)_{t \Ge 0}$ is the Brownian motion of the step, so that $\sqrt2\, B_h = \sqrt{2h}\, \chi_0$. Hence
    \begin{equation*}
        |\tilde W_s| \Le 2M \int_0^s \bigl( 1 + |Y_0| + |\tilde W_r| \bigr) \D r + \sqrt2\, \sup_{r \Le h} |B_r| ,
    \end{equation*}
    and Gr\"onwall's inequality gives
    \begin{equation*}
        \sup_{s \Le h} |\tilde W_s| \Le e^{2Mh} \Bigl( 2Mh \bigl( 1 + |Y_0| \bigr) + \sqrt2\, \sup_{r \Le h} |B_r| \Bigr) .
    \end{equation*}
    Since $\mathbb E \big[\sup_{r \Le h} |B_r|^p \big] \Le C h^{\frac p2}$ by the Burkholder--Davis--Gundy inequality and $h \Le 1$, the bound on $\tilde Y_0(s) - Y_0$ in \eqref{eq: exact sizes} follows. For the bound on $\tilde Y_0(h) - Y_0 - \sqrt{2h}\, \chi_0$,
    \begin{equation*}
        \bigl| \tilde W_h - \sqrt{2h}\, \chi_0 \bigr| = \biggl| \int_0^h \nabla \tilde U_h\bigl( Y_0 + \tilde W_r \bigr) \D r \biggr| \Le 2Mh \Bigl( 1 + |Y_0| + \sup_{r \Le h} |\tilde W_r| \Bigr) ,
    \end{equation*}
    and the Gr\"onwall bound on $\sup_{s \Le h} |\tilde W_s|$ above applies. This proves the bound on $\tilde Y_0(h) - Y_0 - \sqrt{2h}\, \chi_0$ in \eqref{eq: exact sizes}. In particular, $\mathbb E_0 |\tilde Y_0(s)|^q \Le C ( 1 + |Y_0| )^q$ for $s \Le h$ and every $q \Ge 1$.

    Next we consider the moments of $\tilde W_h = \tilde Y_0(h) - Y_0$ shown in \eqref{eq: exact m1}--\eqref{eq: exact m3}. Since $\tilde Y_0(h)$ is driven by the modified dynamics \eqref{eq: modified langevin} with the generator $\tilde{\mathcal L}_h \psi = \Delta \psi - \nabla \tilde U_h \cdot \nabla \psi$ of \eqref{eq: modified generator}, the estimate of these moments relies on a Taylor expansion with respect to $\tilde{\mathcal L}_h$.
    
    Let $\ell(z) := z - Y_0$, a linear function of $z \in \mathbb R^d$ with components $\ell_i$, so that $\tilde W_h = \ell(\tilde Y_0(h))$ and the entries of $\mathbb E_0\bigl[ \tilde W_h^{\otimes n} \bigr]$, $n = 1, 2, 3$, are $\mathbb E_0\bigl[ \psi(\tilde Y_0(h)) \bigr]$ for the polynomials $\psi = \ell_i$, $\ell_i \ell_j$, $\ell_i \ell_j \ell_k$. Here and in the rest of the proof, $z$ denotes the variable of $\ell$ and of these polynomials.
    
    For each such $\psi$, since $\tilde U_h \in C^5(\mathbb R^d)$ by Lemma~\ref{lem: modified class}, the functions $\tilde{\mathcal L}_h \psi$ and $\tilde{\mathcal L}_h^2 \psi$ are of class $C^2$ and $\tilde{\mathcal L}_h^3 \psi$ is continuous. By \eqref{eq: modified class bounds}, $\tilde{\mathcal L}_h \psi$, $\tilde{\mathcal L}_h^2 \psi$, $\tilde{\mathcal L}_h^3 \psi$ and the gradients of $\psi$, $\tilde{\mathcal L}_h \psi$, $\tilde{\mathcal L}_h^2 \psi$ have polynomial growth, so that, by $\mathbb E_0 |\tilde Y_0(s)|^q \Le C ( 1 + |Y_0| )^q$, the stochastic integral in It\^o's formula is a martingale. Dynkin's formula, applied to $\psi$, then to $\tilde{\mathcal L}_h \psi$ and $\tilde{\mathcal L}_h^2 \psi$ inside the integrals, gives with Fubini's theorem the Taylor expansion to the order $2$ with integral remainder, where $\psi(Y_0) = 0$ is omitted since $\ell(Y_0) = 0$,
    \begin{equation}
        \mathbb E_0\bigl[ \psi(\tilde Y_0(h)) \bigr] = h\, (\tilde{\mathcal L}_h\psi)(Y_0) + \frac{h^2}{2}\, (\tilde{\mathcal L}_h^2\psi)(Y_0)
        + \frac12 \int_0^h (h-s)^2\, \mathbb E_0\bigl[ (\tilde{\mathcal L}_h^3 \psi)(\tilde Y_0(s)) \bigr] \D s .
        \label{eq: ito taylor}
    \end{equation}
    Here $\tilde{\mathcal L}_h \psi(Y_0)$ and $\tilde{\mathcal L}_h^2 \psi(Y_0)$ are computed with the product rule
    \begin{equation}
        \tilde{\mathcal L}_h(g_1 g_2) = g_1 \tilde{\mathcal L}_h g_2 + g_2 \tilde{\mathcal L}_h g_1 + 2\, \nabla g_1 \cdot \nabla g_2 , \qquad g_1, g_2 \in C^2(\mathbb R^d) .
        \label{eq: generator product}
    \end{equation}

    Next we evaluate the explicit $\tilde{\mathcal L}_h$ and $\tilde{\mathcal L}_h^2$ operations on $\psi = \ell_i$, $\ell_i \ell_j$, $\ell_i \ell_j \ell_k$, to give the explicit $h$ and $h^2$ coefficients in the moment estimates \eqref{eq: exact m1}--\eqref{eq: exact m3}.
    \begin{enumerate}
    \item For $\psi = \ell_i$, the linearity of $\ell_i$ gives $\tilde{\mathcal L}_h \ell_i = -\partial_i \tilde U_h$, and
    \begin{equation*}
        \tilde{\mathcal L}_h^2 \ell_i = -\tilde{\mathcal L}_h \partial_i \tilde U_h = -\partial_i \Delta \tilde U_h + \sum_{j} \partial_j \tilde U_h\, \partial_i \partial_j \tilde U_h = \bigl( \nabla^2 \tilde U_h\, \nabla \tilde U_h - \nabla \Delta \tilde U_h \bigr)_i .
    \end{equation*}
    At $z = Y_0$, $h\, \tilde{\mathcal L}_h \ell_i + \frac{h^2}{2}\, \tilde{\mathcal L}_h^2 \ell_i$ is the component $i$ of $-h\, \nabla \tilde U_h + \frac{h^2}{2} \bigl( \nabla^2 \tilde U_h\, \nabla \tilde U_h - \nabla \Delta \tilde U_h \bigr)$, the approximation of $\mathbb E_0[\tilde W_h]$ in \eqref{eq: exact m1}.

    \item For $\psi = \ell_i \ell_j$, \eqref{eq: modified generator} gives
    \begin{equation*}
        \tilde{\mathcal L}_h(\ell_i \ell_j) = 2\delta_{ij} - \partial_i \tilde U_h\, \ell_j - \partial_j \tilde U_h\, \ell_i ,
    \end{equation*}
    which equals $2\delta_{ij}$ at $z = Y_0$. Since $\tilde{\mathcal L}_h$ annihilates the constant $2\delta_{ij}$, $\tilde{\mathcal L}_h^2(\ell_i \ell_j) = -\tilde{\mathcal L}_h(\partial_i \tilde U_h\, \ell_j) - \tilde{\mathcal L}_h(\partial_j \tilde U_h\, \ell_i)$, and \eqref{eq: generator product} with $g_1 = \partial_i \tilde U_h$ and $g_2 = \ell_j$ gives, at $z = Y_0$, where $\ell_j = 0$,
    \begin{equation*}
        \tilde{\mathcal L}_h\bigl( \partial_i \tilde U_h\, \ell_j \bigr)(Y_0) = \partial_i \tilde U_h\, \tilde{\mathcal L}_h \ell_j + 2\, \nabla \partial_i \tilde U_h \cdot \nabla \ell_j = -\partial_i \tilde U_h\, \partial_j \tilde U_h + 2\, \partial_i \partial_j \tilde U_h .
    \end{equation*}
    This identity and its copy with $i$ and $j$ exchanged give
    \begin{equation*}
        \tilde{\mathcal L}_h^2(\ell_i \ell_j)(Y_0) = 2\, \partial_i \tilde U_h\, \partial_j \tilde U_h - 4\, \partial_i \partial_j \tilde U_h ,
    \end{equation*}
    so that, at $z = Y_0$, $h\, \tilde{\mathcal L}_h(\ell_i \ell_j) + \frac{h^2}{2}\, \tilde{\mathcal L}_h^2(\ell_i \ell_j)$ is the entry $(i,j)$ of $2h\, I_d + h^2 \bigl( \nabla \tilde U_h^{\otimes 2} - 2\, \nabla^2 \tilde U_h \bigr)$, the approximation of $\mathbb E_0[\tilde W_h^{\otimes 2}]$ in \eqref{eq: exact m2}.

    \item For $\psi = \ell_i \ell_j \ell_k$,
    \begin{equation*}
        \tilde{\mathcal L}_h\bigl( \ell_i \ell_j \ell_k \bigr) = 2 \bigl( \delta_{ij}\, \ell_k + \delta_{ik}\, \ell_j + \delta_{jk}\, \ell_i \bigr) - \bigl( \partial_i \tilde U_h\, \ell_j \ell_k + \partial_j \tilde U_h\, \ell_i \ell_k + \partial_k \tilde U_h\, \ell_i \ell_j \bigr) ,
    \end{equation*}
    which vanishes at $z = Y_0$. In $\tilde{\mathcal L}_h^2(\ell_i \ell_j \ell_k)(Y_0)$, the terms $2\delta_{ij}\, \ell_k$, $2\delta_{ik}\, \ell_j$, $2\delta_{jk}\, \ell_i$ contribute $-2\delta_{ij}\, \partial_k \tilde U_h$, $-2\delta_{ik}\, \partial_j \tilde U_h$, $-2\delta_{jk}\, \partial_i \tilde U_h$, since $\tilde{\mathcal L}_h \ell_k = -\partial_k \tilde U_h$. By \eqref{eq: generator product} with $g_1 = \partial_i \tilde U_h$ and $g_2 = \ell_j \ell_k$, whose value and gradient vanish at $Y_0$, the term $-\partial_i \tilde U_h\, \ell_j \ell_k$ contributes $-\partial_i \tilde U_h\, \tilde{\mathcal L}_h(\ell_j \ell_k)(Y_0) = -2\delta_{jk}\, \partial_i \tilde U_h$, and in the same way the terms $-\partial_j \tilde U_h\, \ell_i \ell_k$ and $-\partial_k \tilde U_h\, \ell_i \ell_j$ contribute $-2\delta_{ik}\, \partial_j \tilde U_h$ and $-2\delta_{ij}\, \partial_k \tilde U_h$. Hence
    \begin{equation*}
        \tilde{\mathcal L}_h^2 \bigl( \ell_i \ell_j \ell_k \bigr)(Y_0) = -4 \bigl( \delta_{ij}\, \partial_k \tilde U_h + \delta_{ik}\, \partial_j \tilde U_h + \delta_{jk}\, \partial_i \tilde U_h \bigr) ,
    \end{equation*}
    so that $\frac{h^2}{2}\, \tilde{\mathcal L}_h^2(\ell_i \ell_j \ell_k)(Y_0)$ is the entry $(i,j,k)$ of the approximation of $\mathbb E_0[\tilde W_h^{\otimes 3}]$ in \eqref{eq: exact m3}.
    \end{enumerate}

    It remains to bound the integral term of \eqref{eq: ito taylor}. Each of $\psi$, $\tilde{\mathcal L}_h \psi$, $\tilde{\mathcal L}_h^2 \psi$ and $\tilde{\mathcal L}_h^3 \psi$ is a linear combination of products of components of $\ell$, components of $\nabla \tilde U_h$ and derivatives of $\tilde U_h$ of orders $2$ to $5$; call the degree of such a product the number of its factors that are components of $\ell$ or of $\nabla \tilde U_h$. By the linear growth of $\nabla \tilde U_h$ and \eqref{eq: modified class bounds}, a product of degree $q$ is bounded by $C (1 + |z| + |Y_0|)^q$. In $\tilde{\mathcal L}_h$, the Laplacian does not raise the degree, and $-\nabla \tilde U_h \cdot \nabla$ raises it by at most $1$, and only where it differentiates a derivative of $\tilde U_h$ of order at least $2$. No such derivative occurs in $\psi$ or $\tilde{\mathcal L}_h \psi$, so for $\psi$ a product of $n$ components of $\ell$ the products in $\tilde{\mathcal L}_h \psi$ and $\tilde{\mathcal L}_h^2 \psi$ have degree at most $n$, those in $\tilde{\mathcal L}_h^3 \psi$ have degree at most $n + 1$, and
    \begin{equation*}
        \bigl| (\tilde{\mathcal L}_h^3 \psi)(z) \bigr| \Le C \bigl( 1 + |z| + |Y_0| \bigr)^{n+1} .
    \end{equation*}
    With $\mathbb E_0 (1 + |\tilde Y_0(s)| + |Y_0|)^{n+1} \Le C(1+|Y_0|)^{n+1}$ for $s \Le h$, the integral term is bounded by $C h^3 (1+|Y_0|)^{n+1}$, which is at most the remainders stated in \eqref{eq: exact m1}--\eqref{eq: exact m3}.

    For the exact Langevin dynamics \eqref{eq: langevin}, the same argument applies with $U$ and $X_s$ in place of $\tilde U_h$ and $\tilde Y_0(s)$, since by Assumption~\ref{ass: U} the derivatives of $U$ satisfy \eqref{eq: modified class bounds} with $M$ in place of $2M$.
\end{proof}

\begin{proof}[Proof of Lemma~\ref{lem: moments}]
    We prove the bounds with the exponents $2q$ and $q$ in place of $12$ and $6$ for an integer $q \Ge 1$, the case $q = 6$ being \eqref{eq: moment bounds}; in this proof $C$ may depend on $q$ as well. The proof is a Lyapunov argument: a bound on the displacement over one step and a drift condition give a contraction of $\mathbb E_k |Y_{k+1}|^{2q}$ against $|Y_k|^{2q}$, which is iterated over $k$. By \eqref{eq: lm markov}, the step of the chain at the index $k$ is
    \begin{equation*}
        D_k := Y_{k+1} - Y_k = \sqrt{2h}\, \chi_k - h\, \nabla U(Z_k) , \qquad Z_k = Y_k + \sqrt{\frac h2}\, \chi_k ,
    \end{equation*}
    and Assumption~\ref{ass: U} gives $|\nabla U(Z_k)| \Le M ( 1 + |Y_k| + \sqrt{h/2}\, |\chi_k| )$.

    We first bound the displacement over one step. Since $\mathbb E_k |\sqrt{2h}\, \chi_k|^{p} \Le C h^{\frac p2}$, the bound on $|\nabla U(Z_k)|$ gives, for every $p \Ge 1$,
    \begin{equation}
        \mathbb E_k |D_k|^{p} \Le C \bigl( \mathbb E_k |\sqrt{2h}\, \chi_k|^{p} + h^{p}\, \mathbb E_k |\nabla U(Z_k)|^{p} \bigr)
        \Le C h^{\frac p2} + C h^{p} \bigl( 1 + |Y_k| \bigr)^{p} ,
        \label{eq: step displacement}
    \end{equation}
    and, since \eqref{eq: exact sizes} at the index $k$ gives $\pnorm{\tilde Y_k(h) - Y_k - \sqrt{2h}\,\chi_k}_p \Le C h (1 + |Y_k|)$ under $\mathbb E_k$,
    \begin{equation}
        \mathbb E_k |\tilde Y_k(h) - Y_k|^{p} \Le C h^{\frac p2} + C h^{p} \bigl( 1 + |Y_k| \bigr)^{p} .
        \label{eq: exact displacement}
    \end{equation}

    Next we derive the drift condition. Splitting $\langle Y_k, \nabla U(Z_k)\rangle = \langle Z_k, \nabla U(Z_k)\rangle - \sqrt{h/2}\, \langle \chi_k, \nabla U(Z_k)\rangle$ and applying the drift bound in \eqref{eq: dissipative bounds} at $Z_k$ and the bound on $|\nabla U(Z_k)|$,
    \begin{equation*}
        \langle Y_k, \nabla U(Z_k) \rangle \Ge \frac m4 |Z_k|^2 - C - M \sqrt{\frac h2}\, |\chi_k| \Bigl( 1 + |Y_k| + \sqrt{\frac h2}\, |\chi_k| \Bigr) .
    \end{equation*}
    Taking $\mathbb E_k$, with $\mathbb E_k |Z_k|^2 = |Y_k|^2 + \frac{dh}{2}$ and $\mathbb E_k |\chi_k|^2 = d$,
    \begin{equation}
        \langle Y_k, \mathbb E_k [\nabla U(Z_k)] \rangle \Ge \frac m4 |Y_k|^2 - C - C \sqrt h \bigl( 1 + |Y_k| \bigr) \Ge \frac{m}{8} |Y_k|^2 - C ,
        \label{eq: drift condition}
    \end{equation}
    since $C \sqrt h\, |Y_k| \Le \frac m8 |Y_k|^2 + C$ by Young's inequality and $h \Le 1$.

    Next we combine \eqref{eq: step displacement} and \eqref{eq: drift condition} into the Lyapunov estimate. Let $S_k := |Y_{k+1}|^2 - |Y_k|^2 = 2\langle Y_k, D_k\rangle + |D_k|^2$. Since $\chi_k$ is centered and independent of $\mathcal F_k$, \eqref{eq: drift condition} and \eqref{eq: step displacement} with $p = 2$ give
    \begin{equation*}
        \mathbb E_k S_k = -2h \langle Y_k, \mathbb E_k [\nabla U(Z_k)] \rangle + \mathbb E_k |D_k|^2 \Le -\frac m4 h |Y_k|^2 + C h + C h^2 \bigl( 1 + |Y_k| \bigr)^2 \Le -\frac m8 h |Y_k|^2 + C h
    \end{equation*}
    for $h \Le h_0$ with $C h_0 \Le \frac m8$, while $|S_k| \Le 2|Y_k||D_k| + |D_k|^2$ and \eqref{eq: step displacement} give, for $1 \Le j \Le q$,
    \begin{equation*}
        \mathbb E_k |S_k|^{j} \Le C \bigl( |Y_k|^j\, \mathbb E_k |D_k|^j + \mathbb E_k |D_k|^{2j} \bigr) \Le C h^{\frac j2} |Y_k|^{j} + C h^{j} \bigl( 1 + |Y_k| \bigr)^{2j} .
    \end{equation*}
    Expanding
    \begin{equation*}
        |Y_{k+1}|^{2q} = \bigl( |Y_k|^2 + S_k \bigr)^q = \sum_{j=0}^{q} \binom{q}{j} |Y_k|^{2(q-j)} S_k^{j}
    \end{equation*}
    and inserting the bounds on $\mathbb E_k S_k$ and $\mathbb E_k |S_k|^{j}$, with $h \Le 1$,
    \begin{equation*}
        q\, |Y_k|^{2q-2}\, \mathbb E_k S_k \Le -\frac{qm}{8} h |Y_k|^{2q} + C h |Y_k|^{2q-2} , \qquad
        |Y_k|^{2(q-j)}\, \mathbb E_k |S_k|^{j} \Le C h |Y_k|^{2q - j} + C h^2 (1+|Y_k|)^{2q} , \quad j \Ge 2 .
    \end{equation*}
    By Young's inequality, the terms of degree below $2q$ in $|Y_k|$ satisfy
    \begin{equation*}
        C h |Y_k|^{2q-2} + C h \sum_{j=2}^{q} |Y_k|^{2q-j} \Le \frac{qm}{32}\, h |Y_k|^{2q} + C h ,
    \end{equation*}
    and, for $h \Le h_0$ with $C h_0$ small enough, the terms of order $h^2$ satisfy
    \begin{equation*}
        C h^2 \bigl( 1+|Y_k| \bigr)^{2q} \Le \frac{qm}{32}\, h |Y_k|^{2q} + C h .
    \end{equation*}
    With these bounds on the terms of degree below $2q$ and of order $h^2$, the expansion of $\mathbb E_k |Y_{k+1}|^{2q}$ gives
    \begin{equation*}
        \mathbb E_k |Y_{k+1}|^{2q} \Le |Y_k|^{2q} - \frac{qm}{8}\, h |Y_k|^{2q} + 2 \cdot \frac{qm}{32}\, h |Y_k|^{2q} + C h ,
    \end{equation*}
    that is,
    \begin{equation}
        \mathbb E_k |Y_{k+1}|^{2q} \Le \Bigl( 1 - \frac{qm}{16} h \Bigr) |Y_k|^{2q} + C h ,
        \label{eq: lyapunov}
    \end{equation}
    with $h_0 \Le 16/(qm)$ making the factor lie in $[0,1)$. Taking expectations in \eqref{eq: lyapunov} and iterating over $k$ gives the bound on $\mathbb E|Y_k|^{2q}$ in \eqref{eq: moment bounds}, and averaging \eqref{eq: step displacement} and \eqref{eq: exact displacement} with $p = 2q$ over the law of $Y_k$ gives the displacement bounds in \eqref{eq: moment bounds}:
    \begin{equation*}
        \mathbb E |Y_k|^{2q} \Le |Y_0|^{2q} + \frac{16 C}{qm} , \qquad
        \max\bigl\{ \mathbb E |D_k|^{2q} ,\; \mathbb E |\tilde Y_k(h) - Y_k|^{2q} \bigr\} \Le C h^{q} + C h^{2q}\, \mathbb E \bigl( 1 + |Y_k| \bigr)^{2q} \Le C h^{q} \bigl( 1 + |Y_0| \bigr)^{2q} .
    \end{equation*}

    The argument carries over to the stochastic gradient chain \eqref{eq: sglm markov}, with $\hat{\mathbb E}_k$ of \eqref{eq: sg filtration} in place of $\mathbb E_k$, once the displacement bounds and the drift condition are verified for it. The step \eqref{eq: sg chain step} of this chain, at the index $k$, is
    \begin{equation*}
        D_k = Y_{k+1} - Y_k = \sqrt{2h}\, \chi_k - h\, \nabla \hat U(Z_k; \xi_k) .
    \end{equation*}
    By the moment bound of Assumption~\ref{ass: sg} and Jensen's inequality, for $1 \Le p \Le 12$,
    \begin{equation*}
        \hat{\mathbb E}_k |\nabla \hat U(Z_k; \xi_k)|^p \Le M^p\, \hat{\mathbb E}_k \bigl( 1 + |Z_k| \bigr)^p \Le C \bigl( 1 + |Y_k| \bigr)^p ,
    \end{equation*}
    which is the bound used for $\nabla U(Z_k)$ in \eqref{eq: step displacement}, with the same constant, so that
    \begin{equation*}
        \hat{\mathbb E}_k |D_k|^{p} \Le C h^{\frac p2} + C h^{p} \bigl( 1 + |Y_k| \bigr)^{p} , \qquad 1 \Le p \Le 12 .
    \end{equation*}
    The bound \eqref{eq: exact displacement} holds with $\hat{\mathbb E}_k$ in place of $\mathbb E_k$, since $\tilde Y_k(h) - Y_k$ is driven by $(B_{kh+t} - B_{kh})_{0 \Le t \Le h}$, which is independent of $\hat{\mathcal F}_k$. Since the estimator is unbiased and $\xi_k$ is independent of $\hat{\mathcal F}_{k}$ and $\chi_k$, the derivation of \eqref{eq: drift condition} with $\hat{\mathbb E}_k$ in place of $\mathbb E_k$ gives
    \begin{equation*}
        \bigl\langle Y_k, \hat{\mathbb E}_k [\nabla \hat U(Z_k; \xi_k)] \bigr\rangle = \bigl\langle Y_k, \hat{\mathbb E}_k [\nabla U(Z_k)] \bigr\rangle \Ge \frac{m}{8} |Y_k|^2 - C .
    \end{equation*}
    The Lyapunov estimate uses only $\hat{\mathbb E}_k[\chi_k] = 0$, the drift condition and the bound on $\hat{\mathbb E}_k |D_k|^{p}$ with $p \Le 2q$, so for $q \Le 6$ it gives \eqref{eq: lyapunov} with $\hat{\mathbb E}_k$ in place of $\mathbb E_k$ and the same constants. Averaging the bound on $\hat{\mathbb E}_k |D_k|^{p}$ and \eqref{eq: exact displacement} under $\hat{\mathbb E}_k$, with $p = 2q \Le 12$, over the law of $Y_k$ then gives \eqref{eq: moment bounds} for \eqref{eq: sglm markov}.
\end{proof}

\subsection{Proof of Lemma~\ref{lem: sizes}}
\label{sec: sizes proof}

\begin{proof}[Proof of Lemma~\ref{lem: sizes}]
    We prove the bounds in $\pnorm{\cdot}_p$ for every $p \Ge 1$, with $C$ depending also on $p$; the case $p = 4$ is \eqref{eq: sizes}. The bounds on $Y_1 - Y_0$ and on $\tilde Y_0(s) - Y_0$ in \eqref{eq: sizes} are direct. For the step $D_0 = Y_1 - Y_0 = \sqrt{2h}\, \chi_0 - h\, \nabla U(Z_0)$ of \eqref{eq: chain step}, the moments of $\chi_0$ and Assumption~\ref{ass: U} give
    \begin{equation*}
        \pnorm{\sqrt{2h}\, \chi_0}_p \Le C\sqrt h , \qquad \pnorm{h\, \nabla U(Z_0)}_p \Le M h \Bigl( 1 + |Y_0| + \sqrt{\frac h2}\, \pnorm{\chi_0}_p \Bigr) \Le C h (1+|Y_0|) ,
    \end{equation*}
    and the bound on $\tilde Y_0(s) - Y_0$ is that in \eqref{eq: exact sizes}.

    Next we bound the local error $Y_1 - \tilde Y_0(h) = D_0 - \tilde W_h$. Subtracting the displacement $\tilde W_h = \sqrt{2h}\, \chi_0 - \int_0^h \nabla \tilde U_h(Y_0 + \tilde W_r)\,\D r$ of \eqref{eq: modified langevin} from $D_0$ gives
    \begin{equation}
        Y_1 - \tilde Y_0(h) = \underbrace{\int_0^h \Bigl( \nabla U\bigl( Y_0 + \tilde W_r \bigr) - \nabla U( Z_0 ) \Bigr) \D r}_{=: \mathcal T_1} + \underbrace{\int_0^h \Bigl( \nabla \tilde U_h\bigl( Y_0 + \tilde W_r \bigr) - \nabla U\bigl( Y_0 + \tilde W_r \bigr) \Bigr) \D r}_{=: \mathcal T_2} .
        \label{eq: local error explicit}
    \end{equation}
    The integrands of $\mathcal T_1$ and $\mathcal T_2$ satisfy
    \begin{equation*}
        \bigl| \nabla U( Y_0 + \tilde W_r ) - \nabla U( Z_0 ) \bigr| \Le M \Bigl| \tilde W_r - \sqrt{\frac h2}\, \chi_0 \Bigr| , \qquad
        \bigl| ( \nabla \tilde U_h - \nabla U )( Y_0 + \tilde W_r ) \bigr| \Le C h \bigl( 1 + |Y_0| + |\tilde W_r| \bigr) ,
    \end{equation*}
    by $\norm{\nabla^2 U} \Le M$ of Assumption~\ref{ass: U} with $Z_0 = Y_0 + \sqrt{h/2}\, \chi_0$, and by \eqref{eq: modified class}, respectively. Minkowski's inequality for the integrals over $r \in [0,h]$ and the bound $\pnorm{\tilde W_r}_p \Le C \sqrt h\, (1 + |Y_0|)$ of \eqref{eq: exact sizes} then give
    \begin{equation*}
    \begin{aligned}
        \pnorm{\mathcal T_1}_p & \Le M \int_0^h \Bigl( \pnorm{\tilde W_r}_p + \sqrt{\frac h2}\, \pnorm{\chi_0}_p \Bigr) \D r \Le C h^{\frac32} (1+|Y_0|) , \\
        \pnorm{\mathcal T_2}_p & \Le C h \int_0^h \bigl( 1 + |Y_0| + \pnorm{\tilde W_r}_p \bigr) \D r \Le C h^2 (1+|Y_0|) .
    \end{aligned}
    \end{equation*}
    Adding the bounds on $\pnorm{\mathcal T_1}_p$ and $\pnorm{\mathcal T_2}_p$ gives $\pnorm{Y_1 - \tilde Y_0(h)}_p \Le C h^{\frac32} (1+|Y_0|)$, as stated in \eqref{eq: sizes} for $p = 4$.
\end{proof}

\subsection{Proof of Lemma~\ref{lem: moment matching}}
\label{sec: matching proof}

The proof compares the moments of $\tilde W_h$ in Lemma~\ref{lem: exact moments} with the following expansion of the moments of the step $Y_1 - Y_0$ of the chain.

\begin{lemma}
\label{lem: step moments}
    Let Assumption~\ref{ass: U} hold and let $h \in (0,1]$. Then, with all derivatives of $U$ evaluated at $Y_0$ and $C$ depending only on $M$ and $d$, the moments of the step \eqref{eq: chain step} of the chain satisfy, for $1 \Le i,j,k \Le d$,
    \begin{align}
        \Bigl| \mathbb E_0[Y_1 - Y_0] + h\, \nabla U + \frac{h^2}{4}\, \nabla \Delta U \Bigr| &\Le C h^3 , \label{eq: step m1}\\
        \Bigl\| \mathbb E_0\bigl[ (Y_1 - Y_0)^{\otimes 2} \bigr] - 2h\, I_d - h^2 \bigl( \nabla U^{\otimes 2} - 2\, \nabla^2 U \bigr) \Bigr\| &\Le C h^3 \bigl( 1+|Y_0| \bigr) , \label{eq: step m2}\\
        \Bigl| \bigl( \mathbb E_0\bigl[ (Y_1 - Y_0)^{\otimes 3} \bigr] \bigr)_{ijk} + 2h^2 \bigl( \delta_{ij}\, \partial_k U + \delta_{ik}\, \partial_j U + \delta_{jk}\, \partial_i U \bigr) \Bigr| &\Le C h^3 \bigl( 1+|Y_0| \bigr)^{3} . \label{eq: step m3}
    \end{align}
\end{lemma}

\begin{proof}
    The step \eqref{eq: chain step} reads $Y_1 - Y_0 = \sqrt{2h}\, \chi_0 - h\, \nabla U(Z_0)$ with $Z_0 - Y_0 = \sqrt{h/2}\, \chi_0$, so its moments follow from expansions of $\nabla U(Z_0)$ about $Y_0$ in powers of $\chi_0$. In this proof, $\mathcal O(h^q)$ denotes a tensor of norm at most $C h^q$. Taylor's formula to the order $3$, with the remainder controlled by $\norm{\nabla^5 U} \Le M$, gives
    \begin{equation}
        \nabla U(Z_0) = \nabla U + \sqrt{\frac h2}\, \nabla^2 U\, \chi_0 + \frac h4\, \nabla^3 U\bigl[ \chi_0^{\otimes 2} \bigr] + \frac16 \Bigl( \frac h2 \Bigr)^{\frac32} \nabla^4 U\bigl[ \chi_0^{\otimes 3} \bigr] + \varrho_4 , \qquad |\varrho_4| \Le \frac{M}{24} \Bigl( \frac h2 \Bigr)^2 |\chi_0|^4 ,
        \label{eq: gradient taylor}
    \end{equation}
    and Taylor's formula to the order $1$, with the remainder controlled by $\norm{\nabla^3 U} \Le M$, gives
    \begin{equation}
        \nabla U(Z_0) = \nabla U + \sqrt{\frac h2}\, \nabla^2 U\, \chi_0 + \varrho_2 , \qquad |\varrho_2| \Le \frac M2\, \frac h2\, |\chi_0|^2 .
        \label{eq: gradient taylor crude}
    \end{equation}
    The expansion \eqref{eq: gradient taylor} is used where $\nabla U(Z_0)$ carries a factor $h$ or $h^{\frac32}$, namely in the mean and in the cross term of the moment of order $2$, and \eqref{eq: gradient taylor crude} where it carries a factor $h^2$ or $h^{\frac52}$. In addition, $|\nabla U(x)| \Le M(1+|x|)$ gives $|\nabla U(Z_0)| \Le M ( 1 + |Y_0| + \sqrt{h/2}\, |\chi_0| )$, and under $\mathbb E_0$ every term of odd degree in $\chi_0$ in \eqref{eq: gradient taylor} and \eqref{eq: gradient taylor crude}, and in their products with powers of $\chi_0$, has mean zero.

    Next we compute the moments of orders $1$, $2$, $3$ of $Y_1 - Y_0$ from \eqref{eq: gradient taylor} and \eqref{eq: gradient taylor crude}, to give the explicit $h$ and $h^2$ coefficients in \eqref{eq: step m1}--\eqref{eq: step m3}.
    \begin{enumerate}
    \item For the mean, the terms of degrees $1$ and $3$ in $\chi_0$ of \eqref{eq: gradient taylor} have mean zero, and $\nabla^3 U[\mathbb E_0 \chi_0^{\otimes 2}] = \nabla^3 U[I_d] = \nabla \Delta U$ by the contraction convention of Section~\ref{subsec: notation}, so that
    \begin{equation*}
        \mathbb E_0\bigl[ \nabla U(Z_0) \bigr] = \nabla U + \frac h4\, \nabla \Delta U + \mathbb E_0[\varrho_4] , \qquad |\mathbb E_0[\varrho_4]| \Le C h^2 .
    \end{equation*}
    Multiplied by $-h$, this is \eqref{eq: step m1}, since $\mathbb E_0[Y_1 - Y_0] = -h\, \mathbb E_0[\nabla U(Z_0)]$.

    \item For the moment of order $2$,
    \begin{equation*}
        (Y_1 - Y_0)^{\otimes 2} = 2h\, \chi_0^{\otimes 2} - 2h\sqrt{2h}\, \mathrm{Sym}\bigl( \chi_0 \otimes \nabla U(Z_0) \bigr) + h^2\, \nabla U(Z_0)^{\otimes 2} .
    \end{equation*}
    The term $2h\, \chi_0^{\otimes 2}$ has conditional mean $2h I_d$. For the cross term $-2h\sqrt{2h}\, \mathrm{Sym}( \chi_0 \otimes \nabla U(Z_0) )$, insert \eqref{eq: gradient taylor} into $\chi_0 \otimes \nabla U(Z_0)$. The terms $\chi_0 \otimes \nabla U$ and $\chi_0 \otimes \nabla^3 U[\chi_0^{\otimes 2}]$ have odd degree and mean zero, $\mathbb E_0[\chi_0 \otimes \nabla^2 U\, \chi_0] = \nabla^2 U$ since $\mathbb E_0[\chi_0^{\otimes 2}] = I_d$ and $\nabla^2 U$ is symmetric, and the terms $\frac16 (\frac h2)^{\frac32}\, \chi_0 \otimes \nabla^4 U[\chi_0^{\otimes 3}]$ and $\chi_0 \otimes \varrho_4$ have conditional means of norm at most $C h^{\frac32}$. Hence
    \begin{equation*}
        \mathbb E_0\bigl[ \chi_0 \otimes \nabla U(Z_0) \bigr] = \sqrt{\frac h2}\, \nabla^2 U + \mathcal O\bigl( h^{\frac32} \bigr) , \qquad
        -2h\sqrt{2h}\, \mathrm{Sym}\, \mathbb E_0\bigl[ \chi_0 \otimes \nabla U(Z_0) \bigr] = -2h^2\, \nabla^2 U + \mathcal O( h^3 ) .
    \end{equation*}
    For the term $h^2\, \nabla U(Z_0)^{\otimes 2}$, \eqref{eq: gradient taylor crude} gives
    \begin{equation}
    \begin{aligned}
        \nabla U(Z_0)^{\otimes 2} & = \nabla U^{\otimes 2} + 2\, \mathrm{Sym}\Bigl( \nabla U \otimes \sqrt{\frac h2}\, \nabla^2 U\, \chi_0 \Bigr) + 2\, \mathrm{Sym}\bigl( \nabla U \otimes \varrho_2 \bigr) \\
        & \quad + \frac h2 \bigl( \nabla^2 U\, \chi_0 \bigr)^{\otimes 2} + 2\, \mathrm{Sym}\Bigl( \sqrt{\frac h2}\, \nabla^2 U\, \chi_0 \otimes \varrho_2 \Bigr) + \varrho_2^{\otimes 2} .
    \end{aligned}
    \label{eq: gradient square}
    \end{equation}
    The term $2\, \mathrm{Sym}( \nabla U \otimes \sqrt{h/2}\, \nabla^2 U\, \chi_0 )$ has mean zero. By $|\nabla U| \Le M(1+|Y_0|)$ and the bound on $\varrho_2$ in \eqref{eq: gradient taylor crude}, the conditional means of the remaining $4$ terms of \eqref{eq: gradient square} satisfy
    \begin{equation*}
    \begin{aligned}
        & \Bigl\| \mathbb E_0\bigl[ 2\, \mathrm{Sym}(\nabla U \otimes \varrho_2) \bigr] \Bigr\| \Le C h (1+|Y_0|) , \qquad
        \biggl\| \mathbb E_0\Bigl[ \frac h2 (\nabla^2 U\, \chi_0)^{\otimes 2} \Bigr] \biggr\| \Le C h , \\
        & \Bigl\| \mathbb E_0\bigl[ 2\, \mathrm{Sym}\bigl(\sqrt{h/2}\, \nabla^2 U\, \chi_0 \otimes \varrho_2\bigr) \bigr] \Bigr\| \Le C h^{\frac32} , \qquad
        \Bigl\| \mathbb E_0\bigl[ \varrho_2^{\otimes 2} \bigr] \Bigr\| \Le C h^2 ,
    \end{aligned}
    \end{equation*}
    so that
    \begin{equation*}
        \Bigl\| h^2\, \mathbb E_0\bigl[ \nabla U(Z_0)^{\otimes 2} \bigr] - h^2\, \nabla U^{\otimes 2} \Bigr\| \Le C h^3 (1+|Y_0|) .
    \end{equation*}
    Adding the conditional means of the $3$ terms of $(Y_1 - Y_0)^{\otimes 2}$, namely $2h I_d$, $-2h^2\, \nabla^2 U + \mathcal O(h^3)$ and $h^2\, \nabla U^{\otimes 2}$ up to $C h^3 (1+|Y_0|)$, gives \eqref{eq: step m2}.

    \item For the moment of order $3$,
    \begin{equation*}
        (Y_1 - Y_0)^{\otimes 3} = (2h)^{\frac32}\, \chi_0^{\otimes 3} - 6h^2\, \mathrm{Sym}\bigl( \chi_0^{\otimes 2} \otimes \nabla U(Z_0) \bigr) + 3h^2 \sqrt{2h}\, \mathrm{Sym}\bigl( \chi_0 \otimes \nabla U(Z_0)^{\otimes 2} \bigr) - h^3\, \nabla U(Z_0)^{\otimes 3} .
    \end{equation*}
    The term $(2h)^{\frac32}\, \chi_0^{\otimes 3}$ has mean zero. For the term $-6h^2\, \mathrm{Sym}( \chi_0^{\otimes 2} \otimes \nabla U(Z_0) )$, insert \eqref{eq: gradient taylor crude} into the entry $\mathbb E_0[\chi_{0,i}\, \chi_{0,j}\, \nabla U(Z_0)_k]$: the term $\nabla U$ gives $\delta_{ij}\, \partial_k U$, the term linear in $\chi_0$ has mean zero, and the remainder $\varrho_2$ contributes at most $C h$, so that
    \begin{equation*}
        \mathbb E_0\bigl[ \chi_{0,i}\, \chi_{0,j}\, \nabla U(Z_0)_k \bigr] = \delta_{ij}\, \partial_k U + \mathcal O(h) ,
    \end{equation*}
    and, since $\delta_{ij}\, \partial_k U$ is symmetric in $i$ and $j$, its symmetrization is $\frac13 ( \delta_{ij}\, \partial_k U + \delta_{ik}\, \partial_j U + \delta_{jk}\, \partial_i U )$, so that
    \begin{equation*}
        -6h^2\, \mathrm{Sym}\, \mathbb E_0\bigl[ \chi_0^{\otimes 2} \otimes \nabla U(Z_0) \bigr]_{ijk} = -2h^2 \bigl( \delta_{ij}\, \partial_k U + \delta_{ik}\, \partial_j U + \delta_{jk}\, \partial_i U \bigr) + \mathcal O(h^3) .
    \end{equation*}
    For the term $3h^2 \sqrt{2h}\, \mathrm{Sym}( \chi_0 \otimes \nabla U(Z_0)^{\otimes 2} )$, insert \eqref{eq: gradient square} into $\chi_0 \otimes \nabla U(Z_0)^{\otimes 2}$: the term $\chi_0 \otimes \nabla U^{\otimes 2}$ has mean zero, and every other term carries a factor $\sqrt{h/2}\, \nabla^2 U\, \chi_0$ or $\varrho_2$, so its conditional mean has norm at most $C \sqrt h (1+|Y_0|)$. For the term $h^3\, \nabla U(Z_0)^{\otimes 3}$, the bound on $|\nabla U(Z_0)|$ suffices. Hence
    \begin{equation*}
    \begin{aligned}
        & \Bigl\| 3h^2 \sqrt{2h}\, \mathrm{Sym}\, \mathbb E_0\bigl[ \chi_0 \otimes \nabla U(Z_0)^{\otimes 2} \bigr] \Bigr\| \Le C h^3 (1+|Y_0|) , \\
        & \Bigl\| h^3\, \mathbb E_0\bigl[ \nabla U(Z_0)^{\otimes 3} \bigr] \Bigr\| \Le h^3\, \mathbb E_0 |\nabla U(Z_0)|^3 \Le C h^3 (1+|Y_0|)^3 .
    \end{aligned}
    \end{equation*}
    Adding the conditional means of the $4$ terms of $(Y_1 - Y_0)^{\otimes 3}$ gives \eqref{eq: step m3}.
    \end{enumerate}
\end{proof}

\begin{proof}[Proof of Lemma~\ref{lem: moment matching}]
    The moments of $\tilde W_h = \tilde Y_0(h) - Y_0$ are expanded in \eqref{eq: exact m1}--\eqref{eq: exact m3} with coefficients in $\tilde U_h$, and those of $D_0 = Y_1 - Y_0$ in \eqref{eq: step m1}--\eqref{eq: step m3} with coefficients in $U$; Lemma~\ref{lem: step moments} applies since $h \Le h_0 \Le 1$, and all derivatives are evaluated at $Y_0$. Next we compare the two expansions. Their coefficients differ through $\tilde U_h - U$, and by \eqref{eq: modified class} and the linear growth of $\nabla U$, and of $\nabla \tilde U_h$ from Lemma~\ref{lem: modified class},
    \begin{equation}
        |\nabla \tilde U_h - \nabla U| \Le C h \bigl( 1+|Y_0| \bigr) , \quad \norm{\nabla \tilde U_h^{\otimes 2} - \nabla U^{\otimes 2}} \Le C h \bigl( 1+|Y_0| \bigr)^2 , \quad \max_{2 \Le i \Le 3} \norm{\nabla^i \tilde U_h - \nabla^i U} \Le C h \bigl( 1+|Y_0| \bigr) .
        \label{eq: coefficient replacement}
    \end{equation}
    
    Next we prove \eqref{eq: moment matching} for $j = 1, 2, 3$ with \eqref{eq: coefficient replacement}, first for $j = 2, 3$ and then for $j = 1$.
    \begin{enumerate}
    \item For $j = 2$ and $j = 3$, replacing $\tilde U_h$ by $U$ in the terms of order $h^2$ of \eqref{eq: exact m2} and \eqref{eq: exact m3} changes them in norm by at most $C h^3 (1+|Y_0|)^2$, by \eqref{eq: coefficient replacement}. After this replacement, the approximations of $\mathbb E_0[\tilde W_h^{\otimes j}]$ in \eqref{eq: exact m2} and \eqref{eq: exact m3} are those of $\mathbb E_0[D_0^{\otimes j}]$ in \eqref{eq: step m2} and \eqref{eq: step m3}, which proves \eqref{eq: moment matching} for $j = 2$ and $j = 3$.
    \item For $j = 1$, the term $-h\, \nabla \tilde U_h$ of the approximation of $\mathbb E_0[\tilde W_h]$ in \eqref{eq: exact m1} is kept, since \eqref{eq: coefficient replacement} bounds $h(\nabla \tilde U_h - \nabla U)$ only by $C h^2 (1+|Y_0|)$. In the term of order $h^2$, the replacement of $\tilde U_h$ by $U$ introduces an additional error of at most $C h^3 (1+|Y_0|)^2$ in norm, by \eqref{eq: coefficient replacement}, the linear growth of $\nabla \tilde U_h$ and $\norm{\nabla^2 U} \Le M$, the bound on $\nabla^3 \tilde U_h - \nabla^3 U$ in \eqref{eq: coefficient replacement} controlling $\nabla \Delta \tilde U_h - \nabla \Delta U$. Subtracting \eqref{eq: step m1} then gives
    \begin{equation}
        \Bigl| \mathbb E_0[\tilde W_h] - \mathbb E_0[D_0] + h \bigl( \nabla \tilde U_h - \nabla U \bigr) - \frac{h^2}{2} \bigl( \nabla^2 U\, \nabla U - \nabla \Delta U \bigr) - \frac{h^2}{4}\, \nabla \Delta U \Bigr| \Le C h^3 \bigl( 1+|Y_0| \bigr)^3 .
        \label{eq: m1 difference}
    \end{equation}
    It remains to bound the terms of \eqref{eq: m1 difference} other than $\mathbb E_0[\tilde W_h] - \mathbb E_0[D_0]$ by $C h^3 (1+|Y_0|)^2$, and here the cutoff $\theta_h$ of \eqref{eq: modified potential} enters. On the ball $\{|Y_0| \Le \delta/h\}$, where $\theta_h = 1$ and $\nabla \theta_h = 0$, $\nabla \tilde U_h - \nabla U = \frac h4 \nabla V = \frac h4 ( 2\, \nabla^2 U\, \nabla U - \nabla \Delta U )$ by \eqref{eq: V}, so that
    \begin{equation}
    \begin{aligned}
        & \quad h \bigl( \nabla \tilde U_h - \nabla U \bigr) - \frac{h^2}{2} \bigl( \nabla^2 U\, \nabla U - \nabla \Delta U \bigr) - \frac{h^2}{4}\, \nabla \Delta U \\
        & = \frac{h^2}{4} \bigl( 2\, \nabla^2 U\, \nabla U - \nabla \Delta U \bigr) - \frac{h^2}{2}\, \nabla^2 U\, \nabla U + \frac{h^2}{4}\, \nabla \Delta U = 0 .
    \end{aligned}
    \label{eq: m1 cancellation}
    \end{equation}
    Outside that ball, \eqref{eq: coefficient replacement} and Assumption~\ref{ass: U} give
    \begin{equation*}
        h\, \bigl| \nabla \tilde U_h(Y_0) - \nabla U(Y_0) \bigr| \Le C h^2 (1+|Y_0|) , \qquad
        \Bigl| \frac{h^2}{2} \bigl( \nabla^2 U\, \nabla U - \nabla \Delta U \bigr) + \frac{h^2}{4}\, \nabla \Delta U \Bigr| \Le C h^2 (1+|Y_0|) ,
    \end{equation*}
    and $1 \Le h|Y_0|/\delta$ there, so the left side of \eqref{eq: m1 cancellation} is at most
    \begin{equation*}
        C h^2 (1+|Y_0|) \cdot \frac{h |Y_0|}{\delta} \Le C h^3 (1+|Y_0|)^2
    \end{equation*}
    in norm. On the ball and outside it, \eqref{eq: m1 difference} therefore gives
    \begin{equation*}
        \bigl| \mathbb E_0[\tilde W_h] - \mathbb E_0[D_0] \bigr| \Le C h^3 (1+|Y_0|)^3 \Le C h^3 (1+|Y_0|)^6 ,
    \end{equation*}
    which proves \eqref{eq: moment matching} for $j = 1$.
    \end{enumerate}
\end{proof}

\subsection{Proof of Lemma~\ref{lem: weak}}
\label{sec: weak proof}

\begin{proof}[Proof of Lemma~\ref{lem: weak}]
    The proof expands $G(Y_1) - G(\tilde Y_0(h))$ about $Y_0$ to the order $3$, so that the terms of orders $1$ to $3$ are differences of the moments of $D_0$ and $\tilde W_h$, bounded by Lemma~\ref{lem: moment matching}, and the remainder is bounded by Lemma~\ref{lem: sizes}.
    
    Work at the index $0$, with $D_0$ and $\tilde W_h$ as in the opening of Appendix~\ref{app: local}. Write $Q_\varsigma := \varsigma Y_1 + (1 - \varsigma) \tilde Y_0(h)$ for $\varsigma \in [0,1]$, the interpolated position of \citep{ye2026mse}, so that $\frac{\D}{\D\varsigma} Q_\varsigma = Y_1 - \tilde Y_0(h)$, the local error, and
    \begin{equation}
        G(Y_1) - G(\tilde Y_0(h)) = \int_0^1 \bigl\langle \nabla G(Q_\varsigma),\, Y_1 - \tilde Y_0(h) \bigr\rangle \,\D\varsigma .
        \label{eq: weak start}
    \end{equation}
    Expand $\nabla G$ about $Y_0$ to the order $2$ with the integral remainder,
    \begin{equation*}
        \nabla G(Y_0 + v) = \nabla G(Y_0) + \nabla^2 G(Y_0)\, v + \frac12\, \nabla^3 G(Y_0)\bigl[ v^{\otimes2} \bigr] + \varrho_3(v) ,
        \qquad |\varrho_3(v)| \Le \frac16\, \norm{\nabla^4 G}_\infty\, |v|^3 ,
    \end{equation*}
    and insert it into \eqref{eq: weak start} with $v = Q_\varsigma - Y_0$. Since $\nabla^2 G(Y_0)$ and $\nabla^3 G(Y_0)$ are symmetric, $\nabla^j G(Y_0) : T = \nabla^j G(Y_0) : \mathrm{Sym}\, T$ for every tensor $T$ of order $j = 2, 3$. Moreover,
    \begin{equation*}
        \frac{\D}{\D\varsigma} \bigl( Q_\varsigma - Y_0 \bigr)^{\otimes j} = j\, \mathrm{Sym}\bigl( (Q_\varsigma - Y_0)^{\otimes (j-1)} \otimes ( Y_1 - \tilde Y_0(h) ) \bigr) , \qquad j = 2, 3 ,
    \end{equation*}
    and integrating over $\varsigma \in [0,1]$, with $Q_0 - Y_0 = \tilde W_h$ and $Q_1 - Y_0 = D_0$,
    \begin{equation*}
    \begin{aligned}
        \mathrm{Sym} \biggl( \int_0^1 ( Q_\varsigma - Y_0 ) \,\D\varsigma \otimes \bigl( Y_1 - \tilde Y_0(h) \bigr) \biggr) & = \frac12 \Bigl( D_0^{\otimes2} - \tilde W_h^{\otimes2} \Bigr) , \\
        \mathrm{Sym} \biggl( \int_0^1 ( Q_\varsigma - Y_0 )^{\otimes2} \,\D\varsigma \otimes \bigl( Y_1 - \tilde Y_0(h) \bigr) \biggr) & = \frac13 \Bigl( D_0^{\otimes3} - \tilde W_h^{\otimes3} \Bigr) .
    \end{aligned}
    \end{equation*}
    These identities turn \eqref{eq: weak start} into
    \begin{equation}
        G(Y_1) - G(\tilde Y_0(h)) = \sum_{j=1}^{3} \frac{1}{j!}\, \nabla^j G(Y_0) : \Bigl( D_0^{\otimes j} - \tilde W_h^{\otimes j} \Bigr) + \varrho_4 ,
        \label{eq: G expansion full}
    \end{equation}
    where the remainder $\varrho_4 := \int_0^1 \bigl\langle \varrho_3(Q_\varsigma - Y_0),\, Y_1 - \tilde Y_0(h) \bigr\rangle \,\D\varsigma$ satisfies
    \begin{equation*}
        |\varrho_4| \Le \frac16\, \norm{\nabla^4 G}_\infty\, \bigl| Y_1 - \tilde Y_0(h) \bigr|\, \max\bigl\{ |D_0| ,\, |\tilde W_h| \bigr\}^3 .
    \end{equation*}
    In \eqref{eq: G expansion full}, $G(Y_1)$ and $G(\tilde Y_0(h))$ are both expanded about the initial value $Y_0$, which is why Lemma~\ref{lem: moment matching} estimates the moments of $D_0$ and $\tilde W_h$ up to the order $3$. Taking $\mathbb E_0$, \eqref{eq: moment matching} bounds the terms of orders $1$ to $3$ in \eqref{eq: G expansion full}:
    \begin{equation*}
        \biggl| \sum_{j=1}^{3} \frac{1}{j!}\, \nabla^j G(Y_0) : \mathbb E_0\bigl[ D_0^{\otimes j} - \tilde W_h^{\otimes j} \bigr] \biggr|
        \Le C h^3 \sum_{j=1}^{3} \norm{\nabla^j G}_\infty \bigl( 1 + |Y_0| \bigr)^{6} .
    \end{equation*}
    For the remainder $\varrho_4$, \eqref{eq: sizes} gives
    \begin{equation*}
        \pnorm{Y_1 - \tilde Y_0(h)}_4 \Le C h^{\frac32}(1+|Y_0|) , \qquad \pnorm{D_0}_4 + \pnorm{\tilde W_h}_4 \Le C \sqrt h (1+|Y_0|) ,
    \end{equation*}
    and H\"older's inequality gives
    \begin{equation*}
        \mathbb E_0 |\varrho_4|
        \Le \frac16 \norm{\nabla^4 G}_\infty\, \pnorm{Y_1 - \tilde Y_0(h)}_4 \bigl( \pnorm{D_0}_{4} + \pnorm{\tilde W_h}_{4} \bigr)^3
        \Le C \norm{\nabla^4 G}_\infty\, h^3 \bigl( 1+|Y_0| \bigr)^{4} .
    \end{equation*}
    Adding the bounds on the terms of orders $1$ to $3$ and on $\varrho_4$ proves \eqref{eq: weak} at the index $0$, and by time homogeneity at every $k$.
\end{proof}

\subsection{Proof of Lemma~\ref{lem: sg weak}}
\label{sec: sg weak proof}

The proof uses two bounds on the covariance of the gradient noise, which also serve Appendices~\ref{app: bias} to \ref{app: w1}. We write the gradient noise and its covariance \eqref{eq: covariance} as
\begin{equation*}
    \eta(x;\xi) := \nabla \hat U(x;\xi) - \nabla U(x) , \qquad \mathcal C(x) = \mathbb E_\xi\bigl[ \eta(x;\xi) \otimes \eta(x;\xi) \bigr] .
\end{equation*}
Let
\begin{equation}
    S(x;\xi) := \sup_{|x'-x| \Le 1} \norm{\nabla^2 \hat U(x';\xi)} , \qquad \mathbb E_\xi \big[S(x;\xi)^2 \big] \Le M^2 ,
    \label{eq: hessian sup}
\end{equation}
the bound being the Hessian condition of Assumption~\ref{ass: sg}. Let $|x - x'| \Le 1$. Then
\begin{equation*}
    \eta(x;\xi) - \eta(x';\xi) = \int_0^1 \bigl( \nabla^2 \hat U(x_s;\xi) - \nabla^2 U(x_s) \bigr) (x - x') \D s , \qquad x_s := x' + s(x - x') ,
\end{equation*}
so the definition of $S$ and $\norm{\nabla^2 U} \Le M$ give
\begin{equation}
    |\eta(x;\xi) - \eta(x';\xi)| \Le \bigl( S(x;\xi) + M \bigr) |x - x'| .
    \label{eq: noise lipschitz}
\end{equation}
The identity
\begin{equation*}
    \eta(x;\xi)^{\otimes 2} - \eta(x';\xi)^{\otimes 2} = \eta(x;\xi) \otimes \bigl( \eta(x;\xi) - \eta(x';\xi) \bigr) + \bigl( \eta(x;\xi) - \eta(x';\xi) \bigr) \otimes \eta(x';\xi)
\end{equation*}
gives
\begin{equation*}
    \bigl\| \eta(x;\xi)^{\otimes 2} - \eta(x';\xi)^{\otimes 2} \bigr\| \Le \bigl( |\eta(x;\xi)| + |\eta(x';\xi)| \bigr) |\eta(x;\xi) - \eta(x';\xi)| .
\end{equation*}
Minkowski's inequality, with Assumption~\ref{ass: sg} for $\eta(x;\xi)$ and $\eta(x';\xi)$ and with \eqref{eq: hessian sup} for $S$, gives
\begin{equation}
    \mathbb E_\xi \Big[\bigl( |\eta(x;\xi)| + |\eta(x';\xi)| \bigr)^2 \Big]  \Le 4\Sigma , \qquad \mathbb E_\xi \big[( S(x;\xi) + M )^2\big] \Le 4M^2 ,
    \label{eq: noise moments}
\end{equation}
so \eqref{eq: covariance}, the bound on $\eta(x;\xi)^{\otimes 2} - \eta(x';\xi)^{\otimes 2}$, \eqref{eq: noise lipschitz} and the Cauchy--Schwarz inequality with \eqref{eq: noise moments} give
\begin{equation*}
\begin{aligned}
    \norm{\mathcal C(x) - \mathcal C(x')} & = \Bigl\| \mathbb E_\xi\bigl[ \eta(x;\xi)^{\otimes 2} - \eta(x';\xi)^{\otimes 2} \bigr] \Bigr\| \Le \mathbb E_\xi\Bigl[ \bigl( |\eta(x;\xi)| + |\eta(x';\xi)| \bigr) |\eta(x;\xi) - \eta(x';\xi)| \Bigr] \\
    & \Le \mathbb E_\xi\Bigl[ \bigl( |\eta(x;\xi)| + |\eta(x';\xi)| \bigr) \bigl( S(x;\xi) + M \bigr) \Bigr]\, |x - x'| \Le 2\Sigma^{\frac12} \cdot 2M\, |x - x'| .
\end{aligned}
\end{equation*}
For $|x - x'| > 1$, since $\mathcal C(x)$ and $\mathcal C(x')$ are positive semidefinite, \eqref{eq: covariance} gives
\begin{equation*}
    \norm{\mathcal C(x) - \mathcal C(x')} \Le \mathrm{Tr}\bigl[ \mathcal C(x) \bigr] + \mathrm{Tr}\bigl[ \mathcal C(x') \bigr] \Le 2\Sigma \Le 2\Sigma\, |x - x'| .
\end{equation*}
Since $\Sigma \Le M^2$, both cases give
\begin{equation}
    \norm{\mathcal C(x) - \mathcal C(x')} \Le 4 M^2 |x - x'| , \qquad x, x' \in \mathbb R^d .
    \label{eq: covariance lipschitz}
\end{equation}
Writing $\mathcal C(x) = \sum_{i=1}^d \lambda_i\, u_i \otimes u_i$ with $\lambda_i \Ge 0$ and $(u_i)$ orthonormal, \eqref{eq: covariance} gives
\begin{equation}
    |A : \mathcal C(x)| = \Bigl| \sum_{i=1}^d \lambda_i \langle u_i, A u_i \rangle \Bigr| \Le \norm{A}\, \mathrm{Tr}\bigl[ \mathcal C(x) \bigr] \Le \Sigma \norm{A} , \qquad A \in \mathbb R^{d \times d} , \ x \in \mathbb R^d .
    \label{eq: trace bound}
\end{equation}

\begin{proof}[Proof of Lemma~\ref{lem: sg weak}]
    Work at the index $0$, where $\hat{\mathbb E}_0$ and $\mathbb E_0$ agree on functions of $Y_0$ and $(B_s)_{s \Le h}$, both fixing $Y_0$ and averaging over $(B_s)_{s \Le h}$. With $D_0 = \sqrt{2h}\, \chi_0 - h\, \nabla U(Z_0)$, the right side of \eqref{eq: chain step}, the step \eqref{eq: sg chain step} adds the gradient noise $\eta_0$ at $Z_0$, multiplied by $-h$,
    \begin{equation*}
        Y_1 = Y_0 + D_0 - h \eta_0 , \qquad \eta_0 := \eta(Z_0;\xi_0) = \nabla \hat U(Z_0;\xi_0) - \nabla U(Z_0) .
    \end{equation*}
    The proof expands $G(Y_1)$ about $Y_0 + D_0$ to the order $2$ in $h\eta_0$, compares $G(Y_0 + D_0)$ with $G(\tilde Y_0(h))$ by Lemma~\ref{lem: weak}, and finally replaces $\nabla^2 G(Y_0 + D_0)$ by $\nabla^2 G(Y_0)$ in the term of order $h^2$.
    
    Let $\hat{\mathcal H}_0 := \sigma\bigl( \hat{\mathcal F}_0, (B_s)_{s \Le h} \bigr)$. Then $Y_0$, $D_0$, $Z_0$ and $\tilde Y_0(h)$ are $\hat{\mathcal H}_0$-measurable, and since $\xi_0$ is independent of $\hat{\mathcal H}_0$, Assumption~\ref{ass: sg} and \eqref{eq: covariance} give
    \begin{equation*}
        \mathbb E\bigl[ \eta_0 \,\big|\, \hat{\mathcal H}_0 \bigr] = 0 , \qquad \mathbb E\bigl[ \eta_0 \otimes \eta_0 \,\big|\, \hat{\mathcal H}_0 \bigr] = \mathcal C(Z_0) .
    \end{equation*}
    Taylor's formula for $G$ at $Y_0 + D_0$, to the order $2$ in $h\eta_0$, gives
    \begin{equation*}
        G(Y_1) = G(Y_0 + D_0) - h\, \nabla G(Y_0+D_0) \cdot \eta_0 + \frac{h^2}{2}\, \nabla^2 G(Y_0+D_0)[\eta_0, \eta_0] + \varrho_3^{\eta_0} , \qquad |\varrho_3^{\eta_0}| \Le \frac{h^3}{6}\, \norm{\nabla^3 G}_\infty\, |\eta_0|^3 .
    \end{equation*}
    Under $\mathbb E[\,\cdot \,|\, \hat{\mathcal H}_0]$, the term linear in $\eta_0$ vanishes and the term quadratic in $\eta_0$ becomes $\frac{h^2}{2}\, \nabla^2 G(Y_0+D_0) : \mathcal C(Z_0)$. For the remainder, $|\eta_0| \Le |\nabla \hat U(Z_0;\xi_0)| + |\nabla U(Z_0)|$, and the moment bound of Assumption~\ref{ass: sg}, applied given $\hat{\mathcal H}_0$, with $|\nabla U(x)| \Le M(1+|x|)$ and $h \Le 1$, gives
    \begin{equation*}
        \hat{\mathbb E}_0 |\eta_0|^3 \Le C\, \hat{\mathbb E}_0 \bigl( 1 + |Z_0| \bigr)^3 \Le C \bigl( 1 + |Y_0| \bigr)^3 .
    \end{equation*}
    Hence
    \begin{equation}
        \Bigl| \hat{\mathbb E}_0\bigl[ G(Y_1) \bigr] - \mathbb E_0\bigl[ G(Y_0 + D_0) \bigr] - \frac{h^2}{2}\, \hat{\mathbb E}_0\bigl[ \nabla^2 G(Y_0 + D_0) : \mathcal C(Z_0) \bigr] \Bigr| \Le C h^3 \norm{\nabla^3 G}_\infty \bigl( 1+|Y_0| \bigr)^3 .
        \label{eq: sg noise expansion}
    \end{equation}
    Next, $Y_0 + D_0$ is the state after one step of \eqref{eq: lm markov} from $Y_0$, so $\mathbb E_0[G(Y_0+D_0)] - \mathbb E_0[G(\tilde Y_0(h))]$ is the local weak error of the chain with the full gradient, and \eqref{eq: weak} at $k = 0$ gives
    \begin{equation}
        \bigl| \mathbb E_0\bigl[ G(Y_0 + D_0) - G(\tilde Y_0(h)) \bigr] \bigr| \Le C h^3 \sum_{j=1}^{4} \norm{\nabla^j G}_\infty \bigl( 1+|Y_0| \bigr)^6 .
        \label{eq: sg exact part}
    \end{equation}

    It remains to replace $\nabla^2 G(Y_0 + D_0)$ by $\nabla^2 G(Y_0)$ in the term carrying $\mathcal C(Z_0)$, at the cost $C h^3$. By Taylor's formula,
    \begin{equation*}
        \nabla^2 G(Y_0 + D_0) = \nabla^2 G(Y_0) + \nabla^3 G(Y_0)[D_0] + \varrho_2^{D_0} , \qquad \norm{\varrho_2^{D_0}} \Le \frac12 \norm{\nabla^4 G}_\infty |D_0|^2 ,
    \end{equation*}
    and, with $D_0 = \sqrt{2h}\, \chi_0 - h\, \nabla U(Z_0)$,
    \begin{equation}
    \begin{aligned}
        & \bigl( \nabla^2 G(Y_0 + D_0) - \nabla^2 G(Y_0) \bigr) : \mathcal C(Z_0) \\
        & \quad = \varrho_2^{D_0} : \mathcal C(Z_0) - h\, \nabla^3 G(Y_0)[\nabla U(Z_0)] : \mathcal C(Z_0) + \sqrt{2h}\, \nabla^3 G(Y_0)[\chi_0] : \mathcal C(Z_0) .
    \end{aligned}
    \label{eq: sg hessian split}
    \end{equation}
    Assumption~\ref{ass: U} gives
    \begin{equation*}
        |\nabla U(Z_0)| \Le M \bigl( 1 + |Z_0| \bigr) \Le M \Bigl( 1 + |Y_0| + \sqrt{\frac h2}\, |\chi_0| \Bigr) .
    \end{equation*}
    By \eqref{eq: trace bound}, $\Sigma \Le M^2$, the bound on $|\nabla U(Z_0)|$ and \eqref{eq: sizes} for the step $D_0$ of \eqref{eq: chain step}, the terms of \eqref{eq: sg hessian split} carrying $\varrho_2^{D_0}$ and $\nabla U(Z_0)$ contribute
    \begin{equation}
    \begin{aligned}
        & \bigl| \hat{\mathbb E}_0[\varrho_2^{D_0} : \mathcal C(Z_0)] \bigr| \Le \frac12 \norm{\nabla^4 G}_\infty\, \Sigma\, \hat{\mathbb E}_0 |D_0|^2 \Le C h \norm{\nabla^4 G}_\infty (1+|Y_0|)^2 , \\
        & h\, \bigl| \hat{\mathbb E}_0\bigl[ \nabla^3 G(Y_0)[\nabla U(Z_0)] : \mathcal C(Z_0) \bigr] \bigr| \Le \Sigma h \norm{\nabla^3 G}_\infty\, \hat{\mathbb E}_0 |\nabla U(Z_0)| \Le C h \norm{\nabla^3 G}_\infty (1+|Y_0|) .
    \end{aligned}
    \label{eq: sg remainder terms}
    \end{equation}
    For the term of \eqref{eq: sg hessian split} carrying $\chi_0$, split
    \begin{equation}
        \nabla^3 G(Y_0)[\chi_0] : \mathcal C(Z_0) = \nabla^3 G(Y_0)[\chi_0] : \mathcal C(Y_0) + \nabla^3 G(Y_0)[\chi_0] : \bigl( \mathcal C(Z_0) - \mathcal C(Y_0) \bigr) .
        \label{eq: sg chi split}
    \end{equation}
    Since $\chi_0$ has mean $0$ and is independent of $\hat{\mathcal F}_0$, while $\nabla^3 G(Y_0)$ and $\mathcal C(Y_0)$ are $\hat{\mathcal F}_0$-measurable,
    \begin{equation}
        \hat{\mathbb E}_0 \bigl[ \nabla^3 G(Y_0)[\chi_0] : \mathcal C(Y_0) \bigr] = 0 ,
        \label{eq: sg centered term}
    \end{equation}
    and, since $Z_0 - Y_0 = \sqrt{h/2}\,\chi_0$, \eqref{eq: covariance lipschitz} gives
    \begin{equation*}
        \norm{\mathcal C(Z_0) - \mathcal C(Y_0)} \Le 4 M^2 \sqrt{\frac h2}\, |\chi_0| ,
    \end{equation*}
    and the facts
    \begin{equation*}
        |A : B| \Le d \norm{A} \norm{B} , \qquad \norm{\nabla^3 G(Y_0)[\chi_0]} \Le \norm{\nabla^3 G}_\infty |\chi_0| , \qquad \hat{\mathbb E}_0 |\chi_0|^2 = d
    \end{equation*}
    give
    \begin{equation}
    \begin{aligned}
        & \sqrt{2h}\, \Bigl| \hat{\mathbb E}_0 \bigl[ \nabla^3 G(Y_0)[\chi_0] : \bigl( \mathcal C(Z_0) - \mathcal C(Y_0) \bigr) \bigr] \Bigr| \Le \sqrt{2h}\, d\, \norm{\nabla^3 G}_\infty\, \hat{\mathbb E}_0\bigl[ |\chi_0|\, \norm{\mathcal C(Z_0) - \mathcal C(Y_0)} \bigr] \\
        & \quad \Le \sqrt{2h}\, d\, \norm{\nabla^3 G}_\infty \cdot 4 M^2 \sqrt{\frac h2}\, \hat{\mathbb E}_0 |\chi_0|^2 = 4 d^2 M^2 h \norm{\nabla^3 G}_\infty .
    \end{aligned}
    \label{eq: sg covariance shift}
    \end{equation}
    Inserting \eqref{eq: sg remainder terms}, \eqref{eq: sg centered term} and \eqref{eq: sg covariance shift} into \eqref{eq: sg hessian split} and \eqref{eq: sg chi split}, and multiplying by $\frac{h^2}{2}$, gives
    \begin{equation}
        \frac{h^2}{2} \Bigl| \hat{\mathbb E}_0\bigl[ \bigl( \nabla^2 G(Y_0 + D_0) - \nabla^2 G(Y_0) \bigr) : \mathcal C(Z_0) \bigr] \Bigr| \Le C h^3 \bigl( \norm{\nabla^3 G}_\infty + \norm{\nabla^4 G}_\infty \bigr) (1+|Y_0|)^2 .
        \label{eq: sg hessian shift}
    \end{equation}
    Finally, the left side of \eqref{eq: sg weak} at the index $0$ splits as
    \begin{align}
        & \quad \hat{\mathbb E}_0\bigl[ G(Y_1) - G(\tilde Y_0(h)) \bigr] - \frac{h^2}{2}\, \nabla^2 G(Y_0) : \hat{\mathbb E}_0\bigl[ \mathcal C(Z_0) \bigr] \notag \\
        & = \Bigl( \hat{\mathbb E}_0\bigl[ G(Y_1) \bigr] - \mathbb E_0\bigl[ G(Y_0 + D_0) \bigr] - \frac{h^2}{2}\, \hat{\mathbb E}_0\bigl[ \nabla^2 G(Y_0 + D_0) : \mathcal C(Z_0) \bigr] \Bigr) + \mathbb E_0\bigl[ G(Y_0 + D_0) - G(\tilde Y_0(h)) \bigr] \notag \\
        & \quad + \frac{h^2}{2}\, \hat{\mathbb E}_0\bigl[ \bigl( \nabla^2 G(Y_0 + D_0) - \nabla^2 G(Y_0) \bigr) : \mathcal C(Z_0) \bigr] .
        \label{eq: sg weak split}
    \end{align}
	On the right side of \eqref{eq: sg weak split}, \eqref{eq: sg noise expansion} bounds the first bracket, \eqref{eq: sg exact part} the second term, and \eqref{eq: sg hessian shift} the third term. Adding these bounds gives
    \begin{equation*}
    \begin{aligned}
        & \quad \Bigl| \hat{\mathbb E}_0\bigl[ G(Y_1) - G(\tilde Y_0(h)) \bigr] - \frac{h^2}{2}\, \nabla^2 G(Y_0) : \hat{\mathbb E}_0\bigl[ \mathcal C(Z_0) \bigr] \Bigr| \\
        & \Le C h^3 \norm{\nabla^3 G}_\infty (1+|Y_0|)^3 + C h^3 \sum_{j=1}^{4} \norm{\nabla^j G}_\infty (1+|Y_0|)^6 + C h^3 \bigl( \norm{\nabla^3 G}_\infty + \norm{\nabla^4 G}_\infty \bigr) (1+|Y_0|)^2 \\
        & \Le C h^3 \sum_{j=1}^{4} \norm{\nabla^j G}_\infty (1+|Y_0|)^6 ,
    \end{aligned}
    \end{equation*}
    which is \eqref{eq: sg weak} at the index $0$; by time homogeneity, \eqref{eq: sg weak} holds at every $k$.
\end{proof}

\section{Proofs of the MSE analysis}
\label{app: mse}

Throughout this appendix, $C$ depends only on $m$, $M$, $R$ and $d$, and $h \Le h_0$ with $h_0 \Le \min\{\lambda^{-1}, \frac12\}$, $\lambda$ the rate of Theorem~\ref{thm: u decay} for $\tilde U_h$, so that Corollary~\ref{cor: modified poisson} applies. By \eqref{eq: modified poisson bound}, $\tilde\phi_h$ is $CL$-Lipschitz.

\subsection{Proof of Theorem~\ref{thm: lm mse}: MSE of LM}
\label{sec: lm mse proof}

\begin{proof}[Proof of Theorem~\ref{thm: lm mse}]
    Let first $Nh \Ge 1$, so that $N \Ge 2$ since $h \Le \frac12$. The first step of \eqref{eq: lm} gives
    \begin{equation*}
        Z_k = Y_k + \sqrt{\frac h2}\,\chi_k , \quad k \Ge 1 , \qquad Y_1 := Z_0 - h \nabla U(Z_0) + \sqrt{\frac h2}\,\chi_0 ,
    \end{equation*}
    with $(Y_k)_{k \Ge 1}$ the chain \eqref{eq: lm markov} started from $Y_1$. Here $Y_1$ is independent of $(B_t - B_h)_{t \Ge h}$, and $\mathbb E |Y_1|^{12} \Le C (1 + |Z_0|)^{12}$. Splitting off the term $k = 0$ in \eqref{eq: mse},
    \begin{equation*}
        \MSE(N,h) \Le \frac{2}{N^2} \bigl| f(Z_0) - \pi(f) \bigr|^2 + 2 \Bigl( \frac{N-1}{N} \Bigr)^2\, \mathbb E \Biggl[ \biggl| \frac{1}{N-1} \sum_{k=1}^{N-1} f(Z_k) - \pi(f) \biggr|^2 \Biggr] ,
    \end{equation*}
    where
    \begin{equation*}
        \frac{2}{N^2} \bigl| f(Z_0) - \pi(f) \bigr|^2 \Le \frac{C L^2}{N^2} \bigl( 1+|Z_0| \bigr)^2 \Le \frac{C L^2}{Nh} \bigl( 1+|Z_0| \bigr)^2 ,
    \end{equation*}
    and the expectation of the time average over $Z_1, \dots, Z_{N-1}$ is the MSE of $N-1$ samples of a chain of the same form started from $Y_1$, with $\frac{1}{(N-1)h} \Le \frac{2}{Nh}$ for $N \Ge 2$. Conditioning on $Y_1$, it therefore suffices to prove \eqref{eq: lm mse} for the samples $Z_k = Y_k + \sqrt{h/2}\,\chi_k$, $k \Ge 0$, of the chain \eqref{eq: lm markov} started from a given $Y_0$, with $Y_0$ in place of $Z_0$, and then to apply it with $Y_1$ in place of $Y_0$ and take the expectation, with $\mathbb E(1 + |Y_1|)^2 \Le C (1 + |Z_0|)^2$ and $\mathbb E(1 + |Y_1|)^{12} \Le C (1 + |Z_0|)^{12}$. By \eqref{eq: mse split modified} and $\frac{h}{N} \Le \frac{1}{Nh}$ for $h \Le 1$, it suffices to prove
    \begin{equation}
        \mathbb E \Biggl[ \biggl| \frac1N \sum_{k=0}^{N-1} \bigl( \tilde f_h(Y_k) - \tilde\pi_h(\tilde f_h) \bigr) \biggr|^2 \Biggr] \Le C L^2 \Bigl( \bigl( 1 + |Y_0| \bigr)^2 \frac{1}{Nh} + \bigl( 1 + |Y_0| \bigr)^{12} h^4 \Bigr) .
        \label{eq: markov mse}
    \end{equation}
    By \eqref{eq: time average formula} and Cauchy's inequality, the left side of \eqref{eq: markov mse} is at most
    \begin{equation}
        \frac{3}{N^2h^2}\, \mathbb E\bigl| \tilde\phi_h(Y_0) - \tilde\phi_h(Y_N) \bigr|^2
        + \frac{3}{N^2}\, \mathbb E\biggl| \sum_{k=0}^{N-1} \Upsilon_k \biggr|^2
        + \frac{3}{N^2}\, \mathbb E\biggl| \sum_{k=0}^{N-1} \Psi_k \biggr|^2 .
        \label{eq: three terms}
    \end{equation}
    The terms $\Upsilon_k$ and $\Psi_k$ are $\mathcal F_{k+1}$-measurable, since $\tilde Y_k(h)$ is driven by the Brownian motion on $[kh, (k+1)h]$, and $\mathbb E_k \Psi_k = 0$ by \eqref{eq: modified poisson}. Hence the $\Psi_k$, and the centered terms $\Upsilon_k - \mathbb E_k \Upsilon_k$, are martingale differences with respect to $(\mathcal F_k)_{k \Ge 0}$, and their sums have additive second moments. Lemmas~\ref{lem: moments} and \ref{lem: sizes} also give
    \begin{equation}
        \mathbb E|Y_k|^2 \Le C \bigl( 1+|Y_0| \bigr)^2 , \qquad \mathbb E|Y_{k+1} - Y_k|^2 + \mathbb E|\tilde Y_k(h) - Y_k|^2 \Le C h \bigl( 1+|Y_0| \bigr)^2 .
        \label{eq: markov second moments}
    \end{equation}

    For the boundary term of \eqref{eq: three terms}, the Lipschitz bound on $\tilde\phi_h$ and \eqref{eq: markov second moments} give
    \begin{equation}
        \frac{3}{N^2h^2}\, \mathbb E\bigl| \tilde\phi_h(Y_0) - \tilde\phi_h(Y_N) \bigr|^2 \Le \frac{C L^2}{N^2h^2}\, \mathbb E|Y_N - Y_0|^2 \Le \frac{C L^2}{N^2 h^2} \bigl( 1+|Y_0| \bigr)^2 \Le \frac{C L^2}{Nh} \bigl( 1+|Y_0| \bigr)^2 ,
    \label{eq: boundary term bound}
    \end{equation}
    using $Nh \Ge \frac12$. This holds for the $N-1$ samples after the first step, since $(N-1)h \Ge 1 - h \Ge \frac12$.

    For the discretization term of \eqref{eq: three terms}, carrying $\Upsilon_k$ with $h \Upsilon_k = \tilde\phi_h(Y_{k+1}) - \tilde\phi_h(\tilde Y_k(h))$ by \eqref{eq: Upsilon Psi}, Lemma~\ref{lem: weak} with $G = \tilde\phi_h$ and \eqref{eq: modified poisson bound} give
    \begin{equation*}
        |\mathbb E_k \Upsilon_k| \Le C L h^2 \bigl( 1+|Y_k| \bigr)^6 ,
    \end{equation*}
    so by Cauchy's inequality and Lemma~\ref{lem: moments},
    \begin{equation*}
        \mathbb E\biggl| \sum_{k=0}^{N-1} \mathbb E_k \Upsilon_k \biggr|^2 \Le N \sum_{k=0}^{N-1} \mathbb E\bigl| \mathbb E_k \Upsilon_k \bigr|^2 \Le C L^2 N^2 h^4 \bigl( 1 + |Y_0| \bigr)^{12} .
    \end{equation*}
    For its fluctuation, the Lipschitz bound on $\tilde\phi_h$ and Lemma~\ref{lem: sizes} give
    \begin{equation*}
        h |\Upsilon_k| \Le C L\, \bigl| Y_{k+1} - \tilde Y_k(h) \bigr| , \qquad \mathbb E_k |\Upsilon_k|^2 \Le C L^2 h \bigl( 1+|Y_k| \bigr)^2 ,
    \end{equation*}
    so by the orthogonality of martingale differences,
    \begin{equation*}
        \mathbb E\biggl| \sum_{k=0}^{N-1} \bigl( \Upsilon_k - \mathbb E_k \Upsilon_k \bigr) \biggr|^2 = \sum_{k=0}^{N-1} \mathbb E\bigl| \Upsilon_k - \mathbb E_k \Upsilon_k \bigr|^2 \Le C L^2 N h \bigl( 1 + |Y_0| \bigr)^2 .
    \end{equation*}
    The discretization term of \eqref{eq: three terms} is therefore at most
    \begin{equation}
        \frac{3}{N^2}\, \mathbb E\biggl| \sum_{k=0}^{N-1} \Upsilon_k \biggr|^2 \Le C L^2 \Bigl( \bigl( 1+|Y_0| \bigr)^{12} h^4 + \bigl( 1+|Y_0| \bigr)^2 \frac hN \Bigr) .
        \label{eq: discretization term bound}
    \end{equation}

    For the martingale term of \eqref{eq: three terms}, carrying $\Psi_k$ with $h \Psi_k = \tilde\phi_h(\tilde Y_k(h)) - \mathbb E_k[\tilde\phi_h(\tilde Y_k(h))]$ by \eqref{eq: Upsilon Psi} and \eqref{eq: modified poisson}, the Lipschitz bound and \eqref{eq: markov second moments} give
    \begin{equation*}
        h^2\, \mathbb E|\Psi_k|^2 \Le C L^2\, \mathbb E\bigl| \tilde Y_k(h) - Y_k \bigr|^2 \Le C L^2 h \bigl( 1+|Y_0| \bigr)^2 ,
    \end{equation*}
    and by orthogonality,
    \begin{equation}
        \frac{3}{N^2}\, \mathbb E\biggl| \sum_{k=0}^{N-1} \Psi_k \biggr|^2 = \frac{3}{N^2} \sum_{k=0}^{N-1} \mathbb E|\Psi_k|^2 \Le \frac{C L^2}{Nh} \bigl( 1+|Y_0| \bigr)^2 .
        \label{eq: martingale term bound}
    \end{equation}

    Adding \eqref{eq: boundary term bound}, \eqref{eq: discretization term bound} and \eqref{eq: martingale term bound}, and using $\frac{h}{N} \Le \frac{1}{Nh}$ for $h \Le 1$, proves \eqref{eq: markov mse}.

    Let finally $Nh < 1$. By Assumption~\ref{ass: f} and the moments of $\pi$ in Lemma~\ref{lem: dissipative}, $|f(x) - \pi(f)| \Le C L (1 + |x|)$, and $\mathbb E|Z_k|^2 \Le C (1 + |Z_0|)^2$ for every $k \Ge 0$, by the moment bound on $Y_1$ and Lemma~\ref{lem: moments}. Hence \eqref{eq: mse} and Cauchy's inequality give
    \begin{equation*}
        \MSE(N,h) \Le \frac1N \sum_{k=0}^{N-1} \mathbb E\bigl| f(Z_k) - \pi(f) \bigr|^2 \Le C L^2 \bigl( 1 + |Z_0| \bigr)^2 \Le \frac{C L^2}{Nh} \bigl( 1 + |Z_0| \bigr)^2 ,
    \end{equation*}
    which is at most the right side of \eqref{eq: lm mse}.
\end{proof}

\subsection{Proof of Theorem~\ref{thm: sglm mse}: MSE of SGLM}
\label{sec: sglm mse proof}

\begin{proof}[Proof of Theorem~\ref{thm: sglm mse}]
    As in the proof of Theorem~\ref{thm: lm mse}, with $Y_1 := Z_0 - h \nabla \hat U(Z_0;\xi_0) + \sqrt{h/2}\,\chi_0$, it suffices to bound the MSE of the samples $Z_k = Y_k + \sqrt{h/2}\,\chi_k$, $k \Ge 0$, of the chain \eqref{eq: sglm markov} started from a given $Y_0$ by the right side of \eqref{eq: sglm mse} with $Y_0$ in place of $Z_0$, and then to apply it with $Y_1$ in place of $Y_0$ and take the expectation, with $\mathbb E(1 + |Y_1|)^2 \Le C (1 + |Z_0|)^2$ and $\mathbb E(1 + |Y_1|)^{12} \Le C (1 + |Z_0|)^{12}$ by Assumption~\ref{ass: sg}. The bound \eqref{eq: mse split modified} and the formula \eqref{eq: time average formula} involve the modified test function, the modified dynamics and $\tilde\phi_h$ only, and hold for this chain with $\hat{\mathcal F}_k$ and $\hat{\mathbb E}_k$ in place of $\mathcal F_k$ and $\mathbb E_k$. The proof of Theorem~\ref{thm: lm mse} is repeated with two changes.

    First, the moment bounds of Lemma~\ref{lem: moments} hold for \eqref{eq: sglm markov} with the same constants, so the boundary term and the martingale term are bounded by $C L^2 (1+|Y_0|)^2 / (Nh)$ as before.

    Second, in the discretization term, Lemma~\ref{lem: sg weak} with $G = \tilde\phi_h$, \eqref{eq: modified poisson bound} and \eqref{eq: trace bound} give
    \begin{equation*}
        \bigl| \hat{\mathbb E}_k \Upsilon_k \bigr| \Le C L \Bigl( \Sigma h + h^2 \bigl( 1 + |Y_k| \bigr)^6 \Bigr) ,
    \end{equation*}
    since $\frac{h^2}{2} |\nabla^2 \tilde\phi_h(Y_k) : \hat{\mathbb E}_k[\mathcal C(Z_k)]| \Le C L \Sigma h^2$. Cauchy's inequality and Lemma~\ref{lem: moments} for the term $h^2 (1 + |Y_k|)^6$ give
    \begin{equation*}
        \mathbb E\biggl| \sum_{k=0}^{N-1} \hat{\mathbb E}_k \Upsilon_k \biggr|^2 \Le C L^2 N^2 \Bigl( \Sigma^2 h^2 + h^4 \bigl( 1 + |Y_0| \bigr)^{12} \Bigr) .
    \end{equation*}
    For its fluctuation, the step of \eqref{eq: sglm markov} differs from that of \eqref{eq: lm markov} by $-h \bigl( \nabla \hat U(Z_k;\xi_k) - \nabla U(Z_k) \bigr)$, and
    \begin{equation*}
        \hat{\mathbb E}_k \bigl| h \bigl( \nabla \hat U(Z_k;\xi_k) - \nabla U(Z_k) \bigr) \bigr|^2 \Le \Sigma h^2 \Le M^2 h^2 .
    \end{equation*}
    With the Lipschitz bound on $\tilde\phi_h$ and Lemma~\ref{lem: sizes}, this gives
    \begin{equation*}
        \hat{\mathbb E}_k |\Upsilon_k|^2 \Le C L^2 \Bigl( h \bigl( 1+|Y_k| \bigr)^2 + \Sigma \Bigr) \Le C L^2 \bigl( 1+|Y_k| \bigr)^2 ,
    \end{equation*}
    and, by orthogonality and $h \Le 1$,
    \begin{equation*}
        \frac{3}{N^2}\, \mathbb E\biggl| \sum_{k=0}^{N-1} \bigl( \Upsilon_k - \hat{\mathbb E}_k \Upsilon_k \bigr) \biggr|^2 \Le \frac{C L^2}{N} \bigl( 1+|Y_0| \bigr)^2 \Le \frac{C L^2}{Nh} \bigl( 1+|Y_0| \bigr)^2 .
    \end{equation*}
    Adding the bounds proves \eqref{eq: sglm mse} for $Nh \Ge 1$. For $Nh < 1$, it follows as at the end of the proof of Theorem~\ref{thm: lm mse}, since the moment bounds of Lemma~\ref{lem: moments} hold for \eqref{eq: sglm markov} and $\mathbb E(1 + |Y_1|)^2 \Le C (1 + |Z_0|)^2$.

    For the average of $P$ independent estimators, the average is unbiased, its moment of order $12$ satisfies the bound of Assumption~\ref{ass: sg} with the same $M$ by Minkowski's inequality, and its variance is at most $\Sigma / P \Le M^2$. Minkowski's inequality in $L^2$ also gives the Hessian bound of Assumption~\ref{ass: sg} with the same $M$. Assumption~\ref{ass: sg} therefore holds with $\Sigma / P$ in place of $\Sigma$, and \eqref{eq: sglm mse} gives the term $\Sigma^2 h^2 / P^2$.
\end{proof}

\section{Proofs of the numerical bias estimate}
\label{app: bias}

Throughout this appendix, $(Y_k)_{k \Ge 0}$ is the chain \eqref{eq: sglm markov}, $\mathcal P_h$ its transition kernel, and $h \Le h_0$ with $M h_0 \Le 1$. The local Hessian bound \eqref{eq: hessian sup} and the bounds \eqref{eq: covariance lipschitz} and \eqref{eq: trace bound} on the covariance $\mathcal C(x)$ of the gradient noise are those of Appendix~\ref{sec: sg weak proof}.

\subsection{Proof of Lemma~\ref{lem: invariant}: invariant distribution}
\label{sec: invariant proof}

\begin{proof}[Proof of Lemma~\ref{lem: invariant}]
    The proof obtains the existence and uniqueness of $\mu_h$ from Harris's ergodic theorem, deduces the moment bound of $\mu_h$ from Lemma~\ref{lem: moments}, and proves the limit of time averages first for the chain \eqref{eq: sglm markov} and then, in \eqref{eq: ergodic limit}, for the samples of \eqref{eq: sglm}. Harris's ergodic theorem in the form of \citep[Theorem~1.2]{hairer2008yet} assumes that a measurable $W \colon \mathbb R^d \to [0,\infty)$, constants $\gamma, \alpha_h \in (0,1)$, $K \Ge 0$, a level $R_W > 0$ and a probability measure $\nu$ satisfy the drift condition, the minorization on the level set $\{W \Le R_W\}$, and the condition on the level,
    \begin{equation}
        (\mathcal P_h W)(y) \Le \gamma W(y) + K , \quad y \in \mathbb R^d , \qquad
        \mathcal P_h(y, \cdot) \Ge \alpha_h\, \nu , \quad W(y) \Le R_W , \qquad
        R_W > \frac{2K}{1-\gamma} .
        \label{eq: harris conditions}
    \end{equation}
    That theorem then gives a unique invariant distribution $\mu_h$ and constants $C_h > 0$ and $\bar\gamma \in (0,1)$, which may depend on $h$, with
    \begin{equation}
        \sup_{y \in \mathbb R^d} \frac{\bigl|(\mathcal P_h^n \varphi)(y) - \mu_h(\varphi)\bigr|}{1 + W(y)} \Le C_h\, \bar\gamma^{\,n} \sup_{y \in \mathbb R^d} \frac{|\varphi(y) - \mu_h(\varphi)|}{1 + W(y)} , \qquad n \Ge 0 ,
        \label{eq: harris}
    \end{equation}
    for every measurable $\varphi$ with $\sup_{y \in \mathbb R^d} |\varphi(y)|/(1 + W(y)) < \infty$.

    For the drift condition in \eqref{eq: harris conditions}, take $W(y) := |y|^2$. The Lyapunov estimate \eqref{eq: lyapunov} with $q = 1$, which holds for \eqref{eq: sglm markov} as shown in the proof of Lemma~\ref{lem: moments}, gives
    \begin{equation*}
        (\mathcal P_h W)(y) \Le \Bigl( 1 - \frac{m}{16} h \Bigr) W(y) + C h , \qquad y \in \mathbb R^d .
    \end{equation*}
    This is the drift condition in \eqref{eq: harris conditions} with
    \begin{equation*}
        \gamma := 1 - \frac{m}{16} h , \qquad K := C h ,
    \end{equation*}
    where $\gamma \in (0,1)$ since $m \Le M$ by Assumption~\ref{ass: U} and $Mh \Le 1$. With the constant $C$ in $K = C h$, set
    \begin{equation*}
        R_W := \frac{64 C}{m} > \frac{32 C}{m} = \frac{2K}{1 - \gamma} ,
    \end{equation*}
    which depends only on $m$, $M$, $R$ and $d$.

    Next we prove the minorization on the level set $\{W \Le R_W\}$, with a constant $\alpha_h$ that may depend on $h$. Let $h_0$ also satisfy $4 M h_0 (2 + R_W^{\frac12}) \Le 1$, let $\mathbb B$ be the closed unit ball and $\nu$ the uniform distribution on $\mathbb B$, and fix $y$ with $|y|^2 \Le R_W$. From $Y_0 = y$, the step \eqref{eq: sglm markov} reads
    \begin{equation}
        Y_1 = T_{\xi_0}\Bigl( \frac y2 + \sqrt{\frac h2}\,\chi_0 \Bigr) , \qquad T_\xi(v) := 2v - h \nabla \hat U\Bigl( \frac y2 + v; \xi \Bigr) , \qquad v \in \mathbb R^d .
        \label{eq: integrator map}
    \end{equation}
    By Markov's inequality, \eqref{eq: hessian sup} and the moment bound of Assumption~\ref{ass: sg}, the event $E_0 := \{ \xi : S(\frac y2;\xi) \Le 4M ,\ |\nabla \hat U(\frac y2;\xi)| \Le 4M (1 + \frac{|y|}{2}) \}$ satisfies $\mathbb P_\xi(E_0^c) \Le \frac1{16} + 4^{-12} \Le \frac12$. For $\xi \in E_0$ and $v, v' \in \mathbb B$,
    \begin{equation*}
        h \Bigl| \nabla \hat U\Bigl( \frac y2 + v;\xi \Bigr) - \nabla \hat U\Bigl( \frac y2 + v';\xi \Bigr) \Bigr| \Le 4Mh\, |v - v'| , \qquad
        h \Bigl| \nabla \hat U\Bigl( \frac y2 + v;\xi \Bigr) \Bigr| \Le 4Mh \Bigl( 2 + \frac{|y|}{2} \Bigr) \Le 1 .
    \end{equation*}
    Hence $T_\xi$ is Lipschitz with constant $3$ on $\mathbb B$, and for every $w \in \mathbb B$ the map
    \begin{equation*}
        v \mapsto \frac12 \Bigl( w + h \nabla \hat U\Bigl( \frac y2 + v;\xi \Bigr) \Bigr)
    \end{equation*}
    maps $\mathbb B$ into itself with Lipschitz constant $\frac12$, so, by the contraction mapping theorem, it has a fixed point, which solves $T_\xi(v) = w$; hence $\mathbb B \subset T_\xi(\mathbb B)$. For a Borel set $A$, the set $A \cap \mathbb B$ is therefore contained in $T_\xi(T_\xi^{-1}(A) \cap \mathbb B)$, whose outer Lebesgue measure is at most $3^d\, \mathrm{Leb}(T_\xi^{-1}(A) \cap \mathbb B)$ by the Lipschitz bound for the $d$-dimensional Hausdorff measure, which equals the Lebesgue measure $\mathrm{Leb}$ on $\mathbb R^d$; on $\mathbb B$, whose points lie within $1 + \frac12 R_W^{\frac12}$ of $\frac y2$, the density of $\frac y2 + \sqrt{h/2}\,\chi_0 \sim \mathcal N(\frac y2, \frac h2 I_d)$ is at least $\beta_h := (\pi h)^{-\frac d2} \exp\bigl( -(1 + \frac12 R_W^{\frac12})^2 / h \bigr)$. Conditioning \eqref{eq: integrator map} on $\xi_0$, which is independent of $\chi_0$, gives
    \begin{equation*}
        \mathcal P_h(y, A) \Ge \mathbb E_\xi\biggl[ \mathbf 1_{\{\xi \in E_0\}}\, \mathbb P\Bigl( \frac y2 + \sqrt{\frac h2}\,\chi_0 \in T_\xi^{-1}(A) \cap \mathbb B \Bigr) \biggr] \Ge \frac{\beta_h}{2 \cdot 3^d}\, \mathrm{Leb}(A \cap \mathbb B) = \frac{\beta_h\, \mathrm{Leb}(\mathbb B)}{2 \cdot 3^d}\, \nu(A) ,
    \end{equation*}
    which is the minorization with $\alpha_h := \min\{ \frac12 , \beta_h\, \mathrm{Leb}(\mathbb B) / (2 \cdot 3^d) \}$, and \eqref{eq: harris} holds.

    For the moment bound, let $\varphi_\ell(y) := \min\{|y|^{12}, \ell\}$, $\ell \Ge 1$. Since $0 \Le \varphi_\ell \Le \ell$, \eqref{eq: harris} with $\varphi = \varphi_\ell$ and $n = k$ gives, for the chain from $Y_0 = 0$,
    \begin{equation*}
        \bigl| \mathbb E\big[ \varphi_\ell(Y_k) \big] - \mu_h(\varphi_\ell) \bigr| = \bigl| (\mathcal P_h^k \varphi_\ell)(0) - \mu_h(\varphi_\ell) \bigr| \Le \ell\, C_h\, \bar\gamma^{\,k} \to 0 , \qquad k \to \infty .
    \end{equation*}
    With Lemma~\ref{lem: moments} at $Y_0 = 0$, this gives
    \begin{equation*}
        \mu_h(\varphi_\ell) = \lim_{k \to \infty} \mathbb E\big[ \varphi_\ell(Y_k)\big] \Le \sup_{k \Ge 0} \mathbb E |Y_k|^{12} \Le C .
    \end{equation*}
    Monotone convergence as $\ell \to \infty$ gives $\mu_h(|y|^{12}) \Le C$.

    Next we prove the limit of time averages for the chain \eqref{eq: sglm markov}. By the minorization, the level set $\{W \Le R_W\}$ is petite \citep[Section~3.1]{meyn1992stability}. Outside this set, $W > R_W > \frac{2K}{1-\gamma}$, so the drift condition in \eqref{eq: harris conditions} gives $\mathcal P_h W \Le \frac{1+\gamma}{2}\, W$. For $Y_0 = y$ with $W(y) > R_W$ and the entrance time $\tau$ of the chain into $\{W \Le R_W\}$, the process $(\frac{1+\gamma}{2})^{-(k \wedge \tau)}\, W(Y_{k \wedge \tau})$ is therefore a nonnegative supermartingale; it exceeds $(\frac{1+\gamma}{2})^{-k} R_W$ on $\{\tau > k\}$ and has mean at most $W(y)$, so
    \begin{equation*}
        \mathbb P(\tau > k) \Le \Bigl( \frac{1+\gamma}{2} \Bigr)^{k} \frac{W(y)}{R_W} \to 0 , \qquad k \to \infty .
    \end{equation*}
    By the Markov property at the time $1$, the chain started from any initial value returns to $\{W \Le R_W\}$ almost surely, so it is Harris recurrent by \citep[Theorem~4.1(i)]{meyn1992stability}, and, since it has the invariant probability $\mu_h$, positive Harris recurrent \citep[Section~1.3]{meyn1992stability}. Let $g$ be measurable and of linear growth, so that $\mu_h(|g|) < \infty$ by the moment bound. The strong law of large numbers for positive Harris recurrent chains \citep[Theorem~17.0.1(i)]{meyn1993markov} gives, for every initial value $Y_0$,
    \begin{equation}
        \frac1N \sum_{k=0}^{N-1} g(Y_k) \to \mu_h(g) \qquad \text{almost surely} .
        \label{eq: chain limit}
    \end{equation}

    For the samples of \eqref{eq: sglm} from the given $Z_0$ of \eqref{eq: ergodic limit}, start \eqref{eq: sglm markov}, as in Section~\ref{subsec: mse analysis}, at
    \begin{equation*}
        Y_1 = Z_0 - h\nabla \hat U(Z_0;\xi_0) + \sqrt{\frac h2}\,\chi_0 , \qquad \text{so that} \qquad Z_k = Y_k + \sqrt{\frac h2}\,\chi_k , \quad k \Ge 1 .
    \end{equation*}
    Since $Y_1$ is independent of $(\chi_k, \xi_k)_{k \Ge 1}$, the sequence $(Y_{j+1})_{j \Ge 0}$ is the chain \eqref{eq: sglm markov} started at $Y_1$, and \eqref{eq: chain limit}, which holds for every initial value, gives by Fubini's theorem, almost surely,
    \begin{equation}
        \frac1N \sum_{j=0}^{N-1} g(Y_{j+1}) \to \mu_h(g) , \qquad N \to \infty .
        \label{eq: shifted chain limit}
    \end{equation}
    For $k \Ge 1$, as in Section~\ref{subsec: modified test function}, $\mathfrak M_k = f(Z_k) - \tilde f_h(Y_k)$, where $\hat{\mathbb E}_k[\mathfrak M_k] = 0$ by \eqref{eq: heat shift} and $\mathbb E[\mathfrak M_k^2] \Le \frac{dh}{2} L^2$ by \eqref{eq: martingale variance}, both with $\hat{\mathbb E}_k$ of \eqref{eq: sg filtration} in place of $\mathbb E_k$, since $Y_k$ is $\hat{\mathcal F}_k$-measurable and $\chi_k$ is independent of $\hat{\mathcal F}_k$; the $\mathfrak M_k$ are therefore $\hat{\mathcal F}_{k+1}$-measurable martingale differences. By \eqref{eq: shifted chain limit} with $g = \tilde f_h$, which is Lipschitz with constant $L$ by Section~\ref{subsec: modified test function}, hence of linear growth,
    \begin{equation*}
        \frac1N \sum_{k=1}^{N-1} \tilde f_h(Y_k) \to \mu_h(\tilde f_h) \qquad \text{almost surely} .
    \end{equation*}
    By the bound on $\mathbb E [ \mathfrak M_k^2 ]$, the martingale $\sum_{k=1}^{N} \mathfrak M_k / k$ is bounded in $L^2$ and converges almost surely, and Kronecker's lemma gives
    \begin{equation*}
        \frac1N \sum_{k=1}^{N-1} \mathfrak M_k \to 0 \qquad \text{almost surely} .
    \end{equation*}
    Let $Y \sim \mu_h$ and $\chi \sim \mathcal N(0, I_d)$ be independent. Fubini's theorem, \eqref{eq: heat shift} and the definition of $\pi_h$ in Lemma~\ref{lem: invariant} give
    \begin{equation}
        \mu_h(\tilde f_h) = \mathbb E \biggl[ f\Big( Y + \sqrt{\frac h2}\,\chi \Big) \biggr] = \pi_h(f) .
        \label{eq: heat shift invariant}
    \end{equation}
    Combining the limits of the averages of $\tilde f_h(Y_k)$ and of $\mathfrak M_k$ with \eqref{eq: heat shift invariant} and with $f(Z_0)/N \to 0$, since $Z_0$ is given,
    \begin{equation*}
        \frac1N \sum_{k=0}^{N-1} f(Z_k) = \frac{f(Z_0)}{N} + \frac1N \sum_{k=1}^{N-1} \tilde f_h(Y_k) + \frac1N \sum_{k=1}^{N-1} \mathfrak M_k \to \pi_h(f) \qquad \text{almost surely} ,
    \end{equation*}
    which is \eqref{eq: ergodic limit}.
\end{proof}

\subsection{Proof of Theorem~\ref{thm: numerical bias}: numerical bias}
\label{sec: numerical bias proof}

\begin{proof}[Proof of Theorem~\ref{thm: numerical bias}]
    The proof writes $\mu_h(\tilde f_h) - \tilde\pi_h(\tilde f_h)$ through the modified discrete Poisson equation \eqref{eq: modified poisson} as the local weak error of Lemma~\ref{lem: sg weak} averaged over $\mu_h$, and then adds the bound of Lemma~\ref{lem: shifted bias} on $\tilde\pi_h(\tilde f_h) - \pi(f)$.

    Let $Y_0 \sim \mu_h$, the invariant distribution of \eqref{eq: sglm markov} given by Lemma~\ref{lem: invariant}, independent of $(B_t)_{t \Ge 0}$ and $(\xi_k)_{k \Ge 0}$. Then $Y_1 \sim \mu_h$. Add $Y_0$ to the generators of $\hat{\mathcal F}_k$ in \eqref{eq: sg filtration}, so that $\hat{\mathcal F}_0 = \sigma(Y_0)$, while $\chi_k$ and $\xi_k$ remain independent of $\hat{\mathcal F}_k$. The process $\tilde Y_0(s)$ evolves by the modified dynamics \eqref{eq: modified langevin}, with the generator $\tilde{\mathcal L}_h$ of \eqref{eq: modified generator}, from $\tilde Y_0(0) = Y_0$, so that
    \begin{equation*}
        (e^{h\tilde{\mathcal L}_h} \tilde\phi_h)(Y_0) = \hat{\mathbb E}_0\bigl[ \tilde\phi_h(\tilde Y_0(h)) \bigr] .
    \end{equation*}
    The modified discrete Poisson equation \eqref{eq: modified poisson} at $Y_0$ therefore reads
    \begin{equation*}
        \tilde f_h(Y_0) - \tilde\pi_h(\tilde f_h) = \frac1h\, \hat{\mathbb E}_0\bigl[ \tilde\phi_h(Y_0) - \tilde\phi_h(\tilde Y_0(h)) \bigr] .
    \end{equation*}
    Taking expectations, with $\mathbb E\, \tilde f_h(Y_0) = \mu_h(\tilde f_h)$, gives
    \begin{equation}
        \mu_h(\tilde f_h) - \tilde\pi_h(\tilde f_h) = \frac1h\, \mathbb E\bigl[ \tilde\phi_h(Y_0) - \tilde\phi_h(\tilde Y_0(h)) \bigr] = \frac1h\, \mathbb E\bigl[ \tilde\phi_h(Y_1) - \tilde\phi_h(\tilde Y_0(h)) \bigr] ,
        \label{eq: bias identity}
    \end{equation}
    where $\mathbb E\, \tilde\phi_h(Y_0) = \mathbb E\, \tilde\phi_h(Y_1)$ since $Y_1 \sim \mu_h$. All expectations in \eqref{eq: bias identity} are finite, since $\tilde\phi_h$ is Lipschitz with constant $C L$ by \eqref{eq: modified poisson bound}, $\mu_h(|y|^{12}) \Le C$ by Lemma~\ref{lem: invariant}, and $\mathbb E|\tilde Y_0(h) - Y_0| \Le C \sqrt h\, ( 1 + \mathbb E|Y_0| )$ by the bound on $\tilde Y_0(s) - Y_0$ in \eqref{eq: exact sizes}.

    Next we expand $\mathbb E\bigl[ \tilde\phi_h(Y_1) - \tilde\phi_h(\tilde Y_0(h)) \bigr]$ in \eqref{eq: bias identity} by Lemma~\ref{lem: sg weak}. By \eqref{eq: modified poisson bound}, the derivatives of $\tilde\phi_h$ of orders $1$ to $4$ are bounded by $C L$, so \eqref{eq: sg weak} with $G = \tilde\phi_h$ at $k = 0$ gives
    \begin{equation*}
        \Bigl| \hat{\mathbb E}_0\bigl[ \tilde\phi_h(Y_1) - \tilde\phi_h(\tilde Y_0(h)) \bigr] - \frac{h^2}{2}\, \nabla^2 \tilde\phi_h(Y_0) : \hat{\mathbb E}_0\bigl[ \mathcal C(Z_0) \bigr] \Bigr| \Le C L h^3 \bigl( 1 + |Y_0| \bigr)^6 .
    \end{equation*}
    Take expectations in the bound on $\hat{\mathbb E}_0\bigl[ \tilde\phi_h(Y_1) - \tilde\phi_h(\tilde Y_0(h)) \bigr]$. Its term of order $h^2$ carries the $\hat{\mathcal F}_0$-measurable factor $\nabla^2 \tilde\phi_h(Y_0)$, so the tower property gives
    \begin{equation*}
        \mathbb E\bigl[ \nabla^2 \tilde\phi_h(Y_0) : \hat{\mathbb E}_0[\mathcal C(Z_0)] \bigr] = \mathbb E\bigl[ \nabla^2 \tilde\phi_h(Y_0) : \mathcal C(Z_0) \bigr] ,
    \end{equation*}
    and hence
    \begin{equation*}
        \Bigl| \mathbb E\bigl[ \tilde\phi_h(Y_1) - \tilde\phi_h(\tilde Y_0(h)) \bigr] - \frac{h^2}{2}\, \mathbb E\bigl[ \nabla^2 \tilde\phi_h(Y_0) : \mathcal C(Z_0) \bigr] \Bigr| \Le C L h^3\, \mathbb E \bigl( 1 + |Y_0| \bigr)^6 .
    \end{equation*}
    Dividing by $h$ and inserting \eqref{eq: bias identity},
    \begin{equation}
        \Bigl| \mu_h(\tilde f_h) - \tilde\pi_h(\tilde f_h) - \frac h2\, \mathbb E\bigl[ \nabla^2 \tilde\phi_h(Y_0) : \mathcal C(Z_0) \bigr] \Bigr| \Le C L h^2\, \mathbb E \bigl( 1 + |Y_0| \bigr)^6 \Le C L h^2 ,
        \label{eq: shifted numerical bias}
    \end{equation}
    by the moment bound $\mu_h(|y|^{12}) \Le C$ of Lemma~\ref{lem: invariant}. By \eqref{eq: heat shift invariant}, $\pi_h(f) = \mu_h(\tilde f_h)$, so
    \begin{equation}
    \begin{aligned}
        & \quad \pi_h(f) - \pi(f) - \frac h2\, \mathbb E\bigl[ \nabla^2 \tilde\phi_h(Y_0) : \mathcal C(Z_0) \bigr] \\
        & = \Bigl( \mu_h(\tilde f_h) - \tilde\pi_h(\tilde f_h) - \frac h2\, \mathbb E\bigl[ \nabla^2 \tilde\phi_h(Y_0) : \mathcal C(Z_0) \bigr] \Bigr) + \bigl( \tilde\pi_h(\tilde f_h) - \pi(f) \bigr) .
    \end{aligned}
        \label{eq: bias split}
    \end{equation}
    The first bracket of \eqref{eq: bias split} is bounded by \eqref{eq: shifted numerical bias}, and the second by \eqref{eq: shifted bias} of Lemma~\ref{lem: shifted bias}, which proves \eqref{eq: numerical bias}. Finally, by \eqref{eq: trace bound} and \eqref{eq: modified poisson bound},
    \begin{equation*}
        \frac h2 \Bigl| \mathbb E\bigl[ \nabla^2 \tilde\phi_h(Y_0) : \mathcal C(Z_0) \bigr] \Bigr| \Le \frac{\Sigma h}{2}\, \norm{\nabla^2 \tilde\phi_h}_\infty \Le C L \Sigma h ,
    \end{equation*}
    which with \eqref{eq: numerical bias} gives
    \begin{equation*}
        |\pi_h(f) - \pi(f)| \Le C L (\Sigma h + h^2) .
    \end{equation*}
\end{proof}

\subsection{Proof of Corollary~\ref{cor: bias coefficient}: coefficient of the numerical bias}
\label{sec: bias coefficient proof}

\begin{proof}[Proof of Corollary~\ref{cor: bias coefficient}]
    With $Y_0 \sim \mu_h$ and $Z_0$ as in Theorem~\ref{thm: numerical bias}, let
    \begin{equation}
        a_h := \mathbb E\bigl[ \nabla^2 \tilde\phi_h(Y_0) : \mathcal C(Z_0) \bigr] .
        \label{eq: coefficient mean}
    \end{equation}
    Dividing \eqref{eq: numerical bias} by $h$ gives
    \begin{equation*}
        \biggl| \frac{\pi_h(f) - \pi(f)}{h} - \frac12\, a_h \biggr| \Le C L h ,
    \end{equation*}
    so it suffices to prove
    \begin{equation*}
        a_h \to \pi\bigl( \nabla^2 \phi_0(x) : \mathcal C(x) \bigr) , \qquad h \to 0 .
    \end{equation*}
    Since $Z_0 = Y_0 + \sqrt{h/2}\,\chi_0$ with $\chi_0 \sim \mathcal N(0, I_d)$ independent of $Y_0 \sim \mu_h$, conditioning on $Y_0$ in \eqref{eq: coefficient mean} gives
    \begin{equation}
        a_h = \mathbb E\Bigl[ \nabla^2 \tilde\phi_h(Y_0) : \mathbb E\bigl[ \mathcal C(Z_0) \,\big|\, Y_0 \bigr] \Bigr] = \mu_h(A_h) ,
        \label{eq: coefficient conditioning}
    \end{equation}
    with
    \begin{equation}
        A_h(y) := \nabla^2 \tilde\phi_h(y) : \tilde{\mathcal C}_h(y) , \qquad \tilde{\mathcal C}_h(y) := \mathbb E\bigg[ \mathcal C\Big( y + \sqrt{\frac h2}\, \chi \Big) \bigg] , \quad \chi \sim \mathcal N(0, I_d) .
        \label{eq: coefficient integrand}
    \end{equation}
    
    First, $\mu_h \to \pi$ weakly as $h\to0$. The family $(\mu_h)$ is tight: since $\mu_h(|y|^{12}) \Le C$ uniformly in $h$ by Lemma~\ref{lem: invariant}, Markov's inequality gives
    \begin{equation}
        \mu_h\bigl( |y| > \rho \bigr) \Le \rho^{-12}\, \mu_h\bigl( |y|^{12} \bigr) \Le C \rho^{-12} , \qquad \forall \rho > 0 .
        \label{eq: mu tail}
    \end{equation}
    Let $\varphi \in C_c^\infty(\mathbb R^d)$; it satisfies Assumption~\ref{ass: f} with some constant $L_\varphi$. By \eqref{eq: heat shift invariant} and \eqref{eq: shift expansion},
    \begin{equation*}
        \pi_h(\varphi) = \mu_h\bigl( e^{\frac h4\Delta} \varphi \bigr) , \qquad \bigl| (e^{\frac h4\Delta} \varphi)(y) - \varphi(y) \bigr| \Le C L_\varphi h , \quad y \in \mathbb R^d ,
    \end{equation*}
    and with Theorem~\ref{thm: numerical bias} applied to $\varphi$ these give
    \begin{equation*}
        |\mu_h(\varphi) - \pi(\varphi)| \Le \Bigl| \mu_h\bigl( \varphi - e^{\frac h4\Delta} \varphi \bigr) \Bigr| + |\pi_h(\varphi) - \pi(\varphi)| \Le C L_\varphi h + C L_\varphi (\Sigma h + h^2) .
    \end{equation*}
    Hence $\mu_h(\varphi) \to \pi(\varphi)$ for every $\varphi \in C_c^\infty(\mathbb R^d)$. Let $b$ be a bounded continuous function on $\mathbb R^d$, and let $\rho > 0$ and $\varepsilon > 0$. With the cutoff $\theta$ of \eqref{eq: modified potential}, which satisfies $\mathbf 1_{\{|y| \Le \rho\}} \Le \theta(y/\rho) \Le 1$, choose a mollification $\varphi \in C_c^\infty(\mathbb R^d)$ of $b\,\theta(\cdot/\rho)$ with $\norm{b\,\theta(\cdot/\rho) - \varphi}_\infty \Le \varepsilon$. The split
    \begin{equation*}
        b = b\,\bigl( 1 - \theta(\cdot/\rho) \bigr) + \bigl( b\,\theta(\cdot/\rho) - \varphi \bigr) + \varphi , \qquad \Bigl| b\,\bigl( 1 - \theta(\cdot/\rho) \bigr) \Bigr| \Le \norm{b}_\infty\, \mathbf 1_{\{|y| > \rho\}} ,
    \end{equation*}
    taken under $\mu_h$ and under $\pi$, gives
    \begin{equation*}
        |\mu_h(b) - \pi(b)| \Le 2\varepsilon + |\mu_h(\varphi) - \pi(\varphi)| + \norm{b}_\infty \bigl( \mu_h(|y| > \rho) + \pi(|y| > \rho) \bigr) ,
    \end{equation*}
    and \eqref{eq: mu tail}, with $\mu_h(\varphi) \to \pi(\varphi)$, gives
    \begin{equation*}
        \limsup_{h \to 0}\, |\mu_h(b) - \pi(b)| \Le 2\varepsilon + \norm{b}_\infty \bigl( C \rho^{-12} + \pi(|y| > \rho) \bigr) ,
    \end{equation*}
    which tends to $0$ as $\varepsilon \to 0$ and $\rho \to \infty$. So $\mu_h \to \pi$ weakly.

    Second, $\nabla^2 \tilde\phi_h \to \nabla^2 \phi_0$ uniformly on compact sets as $h \to 0$. Let
    \begin{equation}
        \psi_h(x) := \tilde\phi_h(x) - \tilde\phi_h(0) , \qquad x \in \mathbb R^d ,
        \label{eq: psi h}
    \end{equation}
    so that $\psi_h(0) = 0$. By \eqref{eq: modified poisson bound}, the derivatives of $\psi_h$ of orders $1$ to $4$ are bounded by $C L$ uniformly in $h$; hence, for $x, x' \in \mathbb R^d$ and $j = 0, 1, 2, 3$,
    \begin{equation*}
        \norm{\nabla^j \psi_h(x) - \nabla^j \psi_h(x')} \Le C L |x - x'| , \qquad |\psi_h(x)| \Le C L |x| ,
    \end{equation*}
    so, on every ball $\{|x| \Le \rho\}$ and uniformly in $h$,
    \begin{equation*}
        \norm{\nabla^j \psi_h(x)} \Le C L (1 + \rho) , \qquad j = 0, 1, 2, 3 .
    \end{equation*}
    The Arzel\`a--Ascoli theorem gives, for every sequence $h_n \to 0$, a subsequence, still written $h_n$, along which each $\nabla^j \psi_{h_n}$ converges uniformly on compact sets. Since the uniform limit of derivatives is the derivative of the limit, the limits of $\nabla^j \psi_{h_n}$ are $\nabla^j \psi$ for a function $\psi \in C^3(\mathbb R^d)$ with $\nabla \psi$, $\nabla^2 \psi$, $\nabla^3 \psi$ bounded by $C L$, and the uniform convergence on compact sets reads
    \begin{equation}
        \lim_{n \to \infty} \sup_{|x| \Le \rho}\Big\{\max_{0 \Le j \Le 3} \norm{\nabla^j \psi_{h_n}(x) - \nabla^j \psi(x)} \Big\} = 0 , \qquad \forall \rho > 0 .
        \label{eq: psi convergence}
    \end{equation}
    
    To identify the limit $\psi$, fix a deterministic initial point $Y_0 \in \mathbb R^d$. By the modified discrete Poisson equation \eqref{eq: modified poisson} at $Y_0$, with $(e^{h\tilde{\mathcal L}_h} \tilde\phi_h)(Y_0) = \mathbb E[\tilde\phi_h(\tilde Y_0(h))]$,
    \begin{equation}
        \mathbb E\bigl[ \tilde\phi_h(\tilde Y_0(h)) \bigr] - \tilde\phi_h(Y_0) = -h \bigl( \tilde f_h(Y_0) - \tilde\pi_h(\tilde f_h) \bigr) .
        \label{eq: modified poisson expectation}
    \end{equation}
    Since $\tilde Y_0(h)$ is driven by the modified dynamics \eqref{eq: modified langevin} from $Y_0$, with the generator $\tilde{\mathcal L}_h$ of \eqref{eq: modified generator}, It\^o's formula for $\tilde\phi_h$, whose stochastic integral is a martingale since $\nabla\tilde\phi_h$ is bounded, expands the difference on the left side of \eqref{eq: modified poisson expectation} as
    \begin{equation}
        \mathbb E\bigl[ \tilde\phi_h(\tilde Y_0(h)) \bigr] - \tilde\phi_h(Y_0) = \int_0^h \mathbb E\bigl[ \tilde g_h(\tilde Y_0(s)) \bigr] \D s ,
        \label{eq: modified dynkin}
    \end{equation}
    where, with $\psi_h$ of \eqref{eq: psi h},
    \begin{equation}
        \tilde g_h := \tilde{\mathcal L}_h \tilde\phi_h = \Delta \psi_h - \nabla \tilde U_h \cdot \nabla \psi_h .
        \label{eq: modified g}
    \end{equation}
    Compared with the modified discrete Poisson equation \eqref{eq: modified poisson}, $-\tilde g_h$ is the counterpart of $\tilde f_h - \tilde\pi_h(\tilde f_h)$, with the difference quotient $(\tilde\phi_h - e^{h\tilde{\mathcal L}_h} \tilde\phi_h)/h$ replaced by $-\tilde{\mathcal L}_h \tilde\phi_h$. Next we bound $\tilde g_h(Y_0) + \tilde f_h(Y_0) - \tilde\pi_h(\tilde f_h)$. Differentiating $\tilde g_h$ in \eqref{eq: modified g},
    \begin{equation*}
        \nabla \tilde g_h = \nabla \Delta \psi_h - \nabla^2 \tilde U_h\, \nabla \psi_h - \nabla^2 \psi_h\, \nabla \tilde U_h ,
    \end{equation*}
    and \eqref{eq: modified class bounds} gives
    \begin{equation*}
        \norm{\nabla^2 \tilde U_h(z)} \Le C , \qquad |\nabla \tilde U_h(z)| \Le C (1 + |z|) , \qquad z \in \mathbb R^d .
    \end{equation*}
    Hence the bounds of \eqref{eq: modified poisson bound} on the derivatives of $\psi_h$ of orders $1$ to $3$ give, for $z \in \mathbb R^d$ and $0 \Le s \Le h$,
    \begin{equation*}
        |\nabla \tilde g_h(z)| \Le C L (1 + |z|) , \qquad |\tilde g_h(\tilde Y_0(s)) - \tilde g_h(Y_0)| \Le C L \bigl( 1 + |Y_0| + |\tilde Y_0(s) - Y_0| \bigr) |\tilde Y_0(s) - Y_0| .
    \end{equation*}
    The bound on $\tilde Y_0(s) - Y_0$ in \eqref{eq: exact sizes} gives, for $0 \Le s \Le h \Le 1$,
    \begin{equation*}
        \mathbb E |\tilde Y_0(s) - Y_0| \Le C \sqrt h (1+|Y_0|) , \qquad \mathbb E |\tilde Y_0(s) - Y_0|^2 \Le C h (1+|Y_0|)^2 .
    \end{equation*}
    Taking the expectation of the bound on $|\tilde g_h(\tilde Y_0(s)) - \tilde g_h(Y_0)|$ and inserting the bounds on $\mathbb E |\tilde Y_0(s) - Y_0|$ and $\mathbb E |\tilde Y_0(s) - Y_0|^2$ gives
    \begin{equation*}
        \Bigl| \mathbb E\bigl[ \tilde g_h(\tilde Y_0(s)) \bigr] - \tilde g_h(Y_0) \Bigr| \Le C L (1+|Y_0|)\, \mathbb E |\tilde Y_0(s) - Y_0| + C L\, \mathbb E |\tilde Y_0(s) - Y_0|^2 \Le C L \sqrt h\, (1+|Y_0|)^2 .
    \end{equation*}
    By \eqref{eq: modified poisson expectation} and \eqref{eq: modified dynkin}, averaging the bound on $\mathbb E\bigl[ \tilde g_h(\tilde Y_0(s)) \bigr] - \tilde g_h(Y_0)$ over $s \in [0,h]$ gives
    \begin{equation}
        \bigl| \tilde g_h(Y_0) + \tilde f_h(Y_0) - \tilde\pi_h(\tilde f_h) \bigr| = \biggl| \frac1h \int_0^h \Bigl( \mathbb E\bigl[ \tilde g_h(\tilde Y_0(s)) \bigr] - \tilde g_h(Y_0) \Bigr) \D s \biggr| \Le C L \sqrt h\, (1+|Y_0|)^2 , \quad Y_0 \in \mathbb R^d .
        \label{eq: g consistency}
    \end{equation}
    To pass to the limit in \eqref{eq: g consistency} along $h_n$, we write, with $\mathcal L$ as in \eqref{eq: generator} and $\tilde g_h$ as in \eqref{eq: modified g},
    \begin{equation}
        \tilde g_{h_n} - \mathcal L \psi = \bigl( \Delta \psi_{h_n} - \Delta \psi \bigr) - \bigl( \nabla \tilde U_{h_n} - \nabla U \bigr) \cdot \nabla \psi_{h_n} - \nabla U \cdot \bigl( \nabla \psi_{h_n} - \nabla \psi \bigr) ,
        \label{eq: g limit split}
    \end{equation}
    and, for every $x \in \mathbb R^d$, with $|\nabla \psi_{h_n}| \Le C L$ by \eqref{eq: modified poisson bound}, \eqref{eq: g limit split} gives
    \begin{equation}
        \bigl| \tilde g_{h_n}(x) - \mathcal L \psi(x) \bigr| \Le \bigl| \Delta \psi_{h_n}(x) - \Delta \psi(x) \bigr| + C L\, \bigl| \nabla \tilde U_{h_n}(x) - \nabla U(x) \bigr| + |\nabla U(x)|\, \bigl| \nabla \psi_{h_n}(x) - \nabla \psi(x) \bigr| .
        \label{eq: g limit bound}
    \end{equation}
    On the right side of \eqref{eq: g limit bound}, the terms in $\psi_{h_n} - \psi$ vanish by \eqref{eq: psi convergence} at $j = 1$ and $j = 2$, and the term in $\nabla \tilde U_{h_n} - \nabla U$ vanishes by \eqref{eq: modified class}:
    \begin{equation*}
    \begin{gathered}
        \lim_{n \to \infty} \bigl| \Delta \psi_{h_n}(x) - \Delta \psi(x) \bigr| \Le d \lim_{n \to \infty} \norm{\nabla^2 \psi_{h_n}(x) - \nabla^2 \psi(x)} = 0 , \\
        \lim_{n \to \infty} \bigl| \nabla \psi_{h_n}(x) - \nabla \psi(x) \bigr| = 0 , \qquad \lim_{n \to \infty} \bigl| \nabla \tilde U_{h_n}(x) - \nabla U(x) \bigr| \Le \lim_{n \to \infty} C h_n (1+|x|) = 0 .
    \end{gathered}
    \end{equation*}
    Hence $\lim_{n \to \infty} \tilde g_{h_n}(x) = \mathcal L \psi(x)$, and \eqref{eq: shift expansion} and \eqref{eq: shifted bias} give
    \begin{equation*}
        \lim_{n \to \infty} \bigl( \tilde f_{h_n}(x) - \tilde\pi_{h_n}(\tilde f_{h_n}) \bigr) = f(x) - \pi(f) .
    \end{equation*}
    Letting $n \to \infty$ in \eqref{eq: g consistency} at $h = h_n$ shows that $\psi$ solves the continuous Poisson equation \eqref{eq: continuous poisson}.
    
    It remains to show that $\nabla^2 \psi = \nabla^2 \phi_0$. The continuous Poisson equation \eqref{eq: continuous poisson} has the solution
    \begin{equation*}
        \phi_0 = \int_0^\infty u(\cdot, t)\, \D t .
    \end{equation*}
    The integral is well defined: by \eqref{eq: u bound}, $u(\cdot,t)$ decays like $e^{-\lambda t}$ uniformly on compact sets, and its derivatives of orders $1$ to $4$ decay like $e^{-\lambda t}$ uniformly on $\mathbb R^d$. Hence
    \begin{equation}
        \norm{\nabla^j \phi_0}_\infty \Le C L , \qquad j = 1, \dots, 4 ,
        \label{eq: continuous poisson bound}
    \end{equation}
    and, since $\partial_t u = \mathcal L u$ and $u(\cdot,0) = f - \pi(f)$, the integral solves \eqref{eq: continuous poisson}. The difference
    \begin{equation*}
        w := \psi - \phi_0 \in C^3(\mathbb R^d) \qquad \text{satisfies} \qquad \mathcal L w = 0 , \qquad \norm{\nabla w}_\infty < \infty .
    \end{equation*}
    Since $\mathcal L w = 0$ and $\nabla w$ is bounded, It\^o's formula gives $(e^{t\mathcal L} w)(x) = \mathbb E^x[w(X_t)] = w(x)$ for $t \Ge 0$.
    The contraction of \citep[Lemma~5]{ye2026mse}, which holds under Assumption~\ref{ass: U} and requires of its initial function only a bounded gradient, gives
    \begin{equation*}
        \norm{\nabla w}_\infty = \norm{\nabla (e^{t\mathcal L} w)}_\infty \Le C e^{-\lambda t} \norm{\nabla w}_\infty \to 0 , \qquad t \to \infty .
    \end{equation*}
    So $w$ is constant and $\nabla^2 \psi = \nabla^2 \phi_0$. The argument for $w$ uses only $\mathcal L w = 0$ and the bounded gradient of $w$, so \eqref{eq: continuous poisson} determines $\phi_0$ up to an additive constant among the functions with a bounded gradient, as stated in Section~\ref{subsec: poisson}. Since $\nabla^2 \psi_h = \nabla^2 \tilde\phi_h$ and $\nabla^2 \psi = \nabla^2 \phi_0$, \eqref{eq: psi convergence} at $j = 2$ shows that every sequence $h_n \to 0$ has a subsequence along which
    \begin{equation*}
        \lim_{n \to \infty} \sup_{|x| \Le \rho} \norm{\nabla^2 \tilde\phi_{h_n}(x) - \nabla^2 \phi_0(x)} = 0 , \qquad \forall \rho > 0 ,
    \end{equation*}
    and hence $\nabla^2 \tilde\phi_h \to \nabla^2 \phi_0$ uniformly on compact sets as $h \to 0$.

    Third, $\mu_h(A_h) \to \pi(A_0)$, with $A_h$ of \eqref{eq: coefficient integrand} and
    \begin{equation}
        A_0(y) := \nabla^2 \phi_0(y) : \mathcal C(y) , \qquad y \in \mathbb R^d .
        \label{eq: coefficient integrand limit}
    \end{equation}
    By \eqref{eq: covariance lipschitz} and the definition of $\tilde{\mathcal C}_h$ in \eqref{eq: coefficient integrand}, for every $y \in \mathbb R^d$,
    \begin{equation}
        \norm{\tilde{\mathcal C}_h(y) - \mathcal C(y)} \Le \mathbb E \biggl\| \mathcal C\Bigl( y + \sqrt{\frac h2}\, \chi \Bigr) - \mathcal C(y) \biggr\| \Le 4 M^2 \sqrt{\frac h2}\, \mathbb E|\chi| .
        \label{eq: covariance smoothing}
    \end{equation}
    By \eqref{eq: coefficient integrand} and \eqref{eq: coefficient integrand limit},
    \begin{equation}
        A_h(y) - A_0(y) = \nabla^2 \tilde\phi_h(y) : \bigl( \tilde{\mathcal C}_h(y) - \mathcal C(y) \bigr) + \bigl( \nabla^2 \tilde\phi_h(y) - \nabla^2 \phi_0(y) \bigr) : \mathcal C(y) .
        \label{eq: coefficient integrand split}
    \end{equation}
    For the term carrying $\tilde{\mathcal C}_h(y) - \mathcal C(y)$ in \eqref{eq: coefficient integrand split}, $|A : B| \Le d \norm{A} \norm{B}$, \eqref{eq: modified poisson bound}, \eqref{eq: covariance smoothing} and $\mathbb E|\chi| \Le \sqrt d$ give
    \begin{equation*}
        \Bigl| \nabla^2 \tilde\phi_h(y) : \bigl( \tilde{\mathcal C}_h(y) - \mathcal C(y) \bigr) \Bigr| \Le d\, \norm{\nabla^2 \tilde\phi_h}_\infty \norm{\tilde{\mathcal C}_h(y) - \mathcal C(y)} \Le d \cdot C L \cdot 4 M^2 \sqrt{\frac h2}\, \sqrt d \Le C L \sqrt h ,
    \end{equation*}
    and for the term carrying $\nabla^2 \tilde\phi_h(y) - \nabla^2 \phi_0(y)$, \eqref{eq: trace bound} gives
    \begin{equation*}
        \Bigl| \bigl( \nabla^2 \tilde\phi_h(y) - \nabla^2 \phi_0(y) \bigr) : \mathcal C(y) \Bigr| \Le \Sigma\, \norm{\nabla^2 \tilde\phi_h(y) - \nabla^2 \phi_0(y)} .
    \end{equation*}
    Hence
    \begin{equation}
        |A_h(y) - A_0(y)| \Le C L \sqrt h + \Sigma\, \norm{\nabla^2 \tilde\phi_h(y) - \nabla^2 \phi_0(y)} , \qquad y \in \mathbb R^d .
        \label{eq: coefficient integrand difference}
    \end{equation}
    With the uniform convergence of $\nabla^2 \tilde\phi_h$ to $\nabla^2 \phi_0$ on compact sets, \eqref{eq: coefficient integrand difference} gives $A_h \to A_0$ uniformly on compact sets. Next we bound $A_h$ and $A_0$ uniformly in $y$ and $h$. Since $A : \tilde{\mathcal C}_h(y) = \mathbb E\bigl[ A : \mathcal C(y + \sqrt{h/2}\,\chi) \bigr]$ by \eqref{eq: coefficient integrand}, \eqref{eq: trace bound} gives
    \begin{equation*}
        |A : \tilde{\mathcal C}_h(y)| \Le \Sigma \norm{A} , \qquad A \in \mathbb R^{d \times d} , \ y \in \mathbb R^d .
    \end{equation*}
    With \eqref{eq: modified poisson bound} and \eqref{eq: continuous poisson bound}, the bound on $A : \tilde{\mathcal C}_h(y)$ and \eqref{eq: trace bound} give
    \begin{equation*}
        |A_h(y)| + |A_0(y)| \Le \Sigma \bigl( \norm{\nabla^2 \tilde\phi_h}_\infty + \norm{\nabla^2 \phi_0}_\infty \bigr) \Le C L \Sigma , \qquad y \in \mathbb R^d ,
    \end{equation*}
    and $A_0$ is continuous, since $\nabla^2 \phi_0$ is continuous and $x \mapsto \mathcal C(x)$ is Lipschitz by \eqref{eq: covariance lipschitz}. For every $\rho > 0$,
    \begin{equation*}
        |\mu_h(A_h) - \pi(A_0)| \Le \Bigl| \mu_h\bigl( (A_h - A_0)\, \mathbf 1_{\{|y| \Le \rho\}} \bigr) \Bigr| + \Bigl| \mu_h\bigl( (A_h - A_0)\, \mathbf 1_{\{|y| > \rho\}} \bigr) \Bigr| + |\mu_h(A_0) - \pi(A_0)| ,
    \end{equation*}
    and the bound $|A_h| + |A_0| \Le C L \Sigma$ on the term carrying $\mathbf 1_{\{|y| > \rho\}}$, followed by \eqref{eq: mu tail}, gives
    \begin{equation*}
    \begin{aligned}
        |\mu_h(A_h) - \pi(A_0)| & \Le \sup_{|y| \Le \rho} |A_h(y) - A_0(y)| + C L \Sigma\, \mu_h\bigl( |y| > \rho \bigr) + |\mu_h(A_0) - \pi(A_0)| \\
        & \Le \sup_{|y| \Le \rho} |A_h(y) - A_0(y)| + C L \Sigma\, \rho^{-12} + |\mu_h(A_0) - \pi(A_0)| .
    \end{aligned}
    \end{equation*}
    As $h \to 0$, the uniform convergence of $A_h$ to $A_0$ on compact sets and the weak convergence $\mu_h \to \pi$, applied to the bounded continuous $A_0$, give
    \begin{equation*}
        \sup_{|y| \Le \rho} |A_h(y) - A_0(y)| \to 0 , \qquad |\mu_h(A_0) - \pi(A_0)| \to 0 ,
    \end{equation*}
    so, with $a_h = \mu_h(A_h)$ by \eqref{eq: coefficient conditioning} and $\pi(A_0) = \pi\bigl( \nabla^2 \phi_0(x) : \mathcal C(x) \bigr)$ by \eqref{eq: coefficient integrand limit},
    \begin{equation*}
        \limsup_{h \to 0}\, \Bigl| a_h - \pi\bigl( \nabla^2 \phi_0(x) : \mathcal C(x) \bigr) \Bigr| \Le C L \Sigma\, \rho^{-12} , \qquad \rho > 0 .
    \end{equation*}
    Letting $\rho \to \infty$ gives $a_h \to \pi\bigl( \nabla^2 \phi_0(x) : \mathcal C(x) \bigr)$, which proves \eqref{eq: bias coefficient}.
\end{proof}

\section{Proofs for the numerical bias comparison with SGLD}
\label{app: sgld}

Throughout this appendix, $(Z_k)_{k \Ge 0}$ is the chain \eqref{eq: sgld} and $\mathcal P_h$ its transition kernel. Since $(\chi_k)_{k \Ge 1}$ is independent of $(\xi_k)_{k \Ge 0}$, the law of the chain is unchanged if its Gaussian increments are taken as $\sqrt{2h}\, \chi_{k+1} = \sqrt2\, (B_{(k+1)h} - B_{kh})$, with $(B_t)_{t \Ge 0}$ a Brownian motion independent of $(\xi_k)_{k \Ge 0}$; then, for a given $Z_0$, $Z_k$ is $\hat{\mathcal F}_k$-measurable, with $\hat{\mathcal F}_k$ and $\hat{\mathbb E}_k$ of \eqref{eq: sg filtration}. Statements about a single step are made at the index $0$, with the solution $X_s$ of \eqref{eq: langevin} from $X_0 = Z_0$, driven by the same Brownian motion. Then
\begin{equation}
    Z_1 = Z_0 + D_0 - h\, \eta_0 , \qquad D_0 := \sqrt2\, B_h - h\, \nabla U(Z_0) , \qquad \eta_0 := \eta(Z_0;\xi_0) = \nabla \hat U(Z_0;\xi_0) - \nabla U(Z_0) ,
    \label{eq: sgld step}
\end{equation}
where $Z_0 + D_0$ is the state after one step of EM from $Z_0$, and $\eta$ is the gradient noise of Appendix~\ref{sec: sg weak proof}.

\subsection{Invariant distribution of SGLD}
\label{sec: sgld invariant proof}

\begin{lemma}
\label{lem: sgld invariant}
    Let Assumptions~\ref{ass: U} and \ref{ass: sg} hold. For every $h \in (0,h_0]$, \eqref{eq: sgld} has a unique invariant distribution $\bar\pi_h$, with $\bar\pi_h(|x|^{12}) \Le C$. For every $f$ satisfying Assumption~\ref{ass: f} and every given $Z_0 \in \mathbb R^d$, 
    \begin{equation*}
        \lim_{N \to \infty} \frac1N \sum_{k=0}^{N-1} f(Z_k) = \bar\pi_h(f) \qquad \text{almost surely} .
    \end{equation*}
\end{lemma}

\begin{proof}
    The proof follows that of Lemma~\ref{lem: invariant} in Appendix~\ref{sec: invariant proof}, and we indicate the changes only. The step of \eqref{eq: sgld} from $Z_k$ is $\sqrt{2h}\, \chi_{k+1} - h\, \nabla \hat U(Z_k;\xi_k)$, with the stochastic gradient evaluated at $Z_k$ itself. The proof of Lemma~\ref{lem: moments} for \eqref{eq: sglm markov} therefore applies with $Z_k$ and $\chi_{k+1}$ in place of $Y_k$ and $\chi_k$, and with the Gaussian shift $\sqrt{h/2}\, \chi_k$ replaced by $0$; in \eqref{eq: drift condition} the drift bound in \eqref{eq: dissipative bounds} applies at $Z_k$ directly. For $q \Le 6$ that proof gives \eqref{eq: lyapunov} and the bound on $\mathbb E |Y_k|^{12}$ in \eqref{eq: moment bounds} for $Z_k$,
    \begin{equation*}
        \hat{\mathbb E}_k |Z_{k+1}|^{2q} \Le \Bigl( 1 - \frac{qm}{16} h \Bigr) |Z_k|^{2q} + C h , \qquad \mathbb E |Z_k|^{12} \Le C \bigl( 1 + |Z_0| \bigr)^{12} ,
    \end{equation*}
    and, as in Appendix~\ref{sec: invariant proof}, the drift condition of Harris's ergodic theorem with $W(x) = |x|^2$ and a level $R_W$ depending only on $m$, $M$, $R$ and $d$.

    The minorization is simpler than in Appendix~\ref{sec: invariant proof} and does not use the Hessian bound of Assumption~\ref{ass: sg}, since the step from $Z_0$ is Gaussian given $\xi_0$,
    \begin{equation*}
        Z_1 \sim \mathcal N\Bigl( Z_0 - h\, \nabla \hat U(Z_0;\xi_0) ,\, 2h I_d \Bigr) .
    \end{equation*}
    Let $h_0$ also satisfy $4 M h_0 \Le 1$, let $|Z_0|^2 \Le R_W$, and let
    \begin{equation*}
        E_0 := \Bigl\{ \xi : |\nabla \hat U(Z_0;\xi)| \Le 4M \bigl( 1 + |Z_0| \bigr) \Bigr\} ,
    \end{equation*}
    so that Markov's inequality and the moment bound of Assumption~\ref{ass: sg} give $\mathbb P_\xi(E_0^c) \Le 4^{-12} \Le \frac12$. For $\xi_0 \in E_0$ and $x$ in the closed unit ball $\mathbb B$, since $4Mh \Le 1$,
    \begin{equation*}
        \bigl| x - Z_0 + h\, \nabla \hat U(Z_0;\xi_0) \bigr| \Le 2 \bigl( 1 + R_W^{\frac12} \bigr) ,
    \end{equation*}
    and the Gaussian density of $Z_1$ given $\xi_0$, with $\nu$ the uniform distribution on $\mathbb B$ and $\mathbb P_\xi(E_0) \Ge \frac12$, gives for every Borel set $A$
    \begin{equation*}
        \mathcal P_h(Z_0, A) \Ge \frac12\, (4\pi h)^{-\frac d2} \exp\Bigl( -\frac{(1 + R_W^{\frac12})^2}{h} \Bigr) \mathrm{Leb}(\mathbb B)\, \nu(A) , \qquad |Z_0|^2 \Le R_W ,
    \end{equation*}
    which is the minorization, with $\alpha_h$ the minimum of $\frac12$ and the factor of $\nu(A)$.

    Harris's ergodic theorem then gives the unique invariant distribution $\bar\pi_h$ and \eqref{eq: harris} with $\bar\pi_h$ in place of $\mu_h$. The moment bound follows from the truncations $\min\bigl\{ |x|^{12}, \ell \bigr\}$, and the almost sure limit from the Harris recurrence and the strong law of large numbers \citep[Theorem~17.0.1(i)]{meyn1993markov} applied to $g = f$, which is Lipschitz by Assumption~\ref{ass: f}. Both arguments are those of Appendix~\ref{sec: invariant proof}, with the bound on $\mathbb E |Z_k|^{12}$ above, at $Z_0 = 0$, in place of Lemma~\ref{lem: moments}. The passage from \eqref{eq: sglm markov} to the samples of \eqref{eq: sglm} is not needed, since $(Z_k)_{k \Ge 0}$ is itself the chain.
\end{proof}

\subsection{Proof of Theorem~\ref{thm: sgld bias}: numerical bias of SGLD}
\label{sec: sgld bias proof}

The proof uses the local weak error of \eqref{eq: sgld} against \eqref{eq: langevin}, expanded to the order $h^2$ with explicit coefficients. The following lemma is the counterpart of Lemma~\ref{lem: sg weak} for \eqref{eq: sglm markov}, with \eqref{eq: langevin} in place of the modified dynamics \eqref{eq: modified langevin}, and its proof follows the same steps.

\begin{lemma}
\label{lem: sgld weak}
    Let Assumptions~\ref{ass: U} and \ref{ass: sg} hold and let $G \in C^4(\mathbb R^d)$ have bounded derivatives of orders $1$ to $4$. Then, for every $h \in (0,1]$, the chain \eqref{eq: sgld} satisfies, for every $k \Ge 0$,
    \begin{equation}
    \begin{aligned}
        & \biggl| \hat{\mathbb E}_k\bigl[ G(Z_{k+1}) \bigr] - (e^{h\mathcal L} G)(Z_k) - \frac{h^2}{2}\, \nabla^2 G(Z_k) : \bigl( 2\nabla^2 U(Z_k) + \mathcal C(Z_k) \bigr) \\
        & + \frac{h^2}{2}\, \nabla G(Z_k) \cdot \bigl( \nabla^2 U(Z_k)\, \nabla U(Z_k) - \nabla \Delta U(Z_k) \bigr) \biggr| \Le C h^3 \sum_{j=1}^{4} \norm{\nabla^j G}_\infty \bigl( 1 + |Z_k| \bigr)^{6} .
    \end{aligned}
        \label{eq: sgld weak}
    \end{equation}
\end{lemma}

\begin{proof}
    By time homogeneity, it suffices to take $k = 0$. The proof expands $G(Z_1)$ about $Z_0 + D_0$ to the order $2$ in $h\eta_0$, replaces $\nabla^2 G(Z_0 + D_0)$ by $\nabla^2 G(Z_0)$ in the term of order $h^2$, and compares $G(Z_0 + D_0)$ with $G(X_h)$ through the moments of $D_0$ and $X_h - Z_0$. With $\hat{\mathcal H}_0 = \sigma\bigl( \hat{\mathcal F}_0, (B_s)_{s \Le h} \bigr)$ as in the proof of Lemma~\ref{lem: sg weak}, $Z_0$, the increment $D_0$ of \eqref{eq: sgld step} and the solution $X_h$ are $\hat{\mathcal H}_0$-measurable, and $\xi_0$ is independent of $\hat{\mathcal H}_0$, so Assumption~\ref{ass: sg} and \eqref{eq: covariance} give
    \begin{equation*}
        \mathbb E\bigl[ \eta_0 \,\big|\, \hat{\mathcal H}_0 \bigr] = 0 , \qquad \mathbb E\bigl[ \eta_0 \otimes \eta_0 \,\big|\, \hat{\mathcal H}_0 \bigr] = \mathcal C(Z_0) .
    \end{equation*}
    Taylor's formula for $G$ at $Z_0 + D_0$, to the order $2$ in $h\eta_0$, as in the proof of Lemma~\ref{lem: sg weak}, and the bound $\hat{\mathbb E}_0 |\eta_0|^3 \Le C \bigl( 1 + |Z_0| \bigr)^3$ from Assumption~\ref{ass: sg} give
    \begin{equation}
        \Bigl| \hat{\mathbb E}_0\bigl[ G(Z_1) - G(Z_0 + D_0) \bigr] - \frac{h^2}{2}\, \hat{\mathbb E}_0\bigl[ \nabla^2 G(Z_0 + D_0) \bigr] : \mathcal C(Z_0) \Bigr| \Le C h^3 \norm{\nabla^3 G}_\infty \bigl( 1 + |Z_0| \bigr)^3 .
        \label{eq: sgld noise}
    \end{equation}
    Next, $\nabla^2 G(Z_0 + D_0)$ is replaced by $\nabla^2 G(Z_0)$ in \eqref{eq: sgld noise} at the cost $C h^3$. By Taylor's formula,
    \begin{equation*}
        \bigl\| \nabla^2 G(Z_0 + D_0) - \nabla^2 G(Z_0) - \nabla^3 G(Z_0)[D_0] \bigr\| \Le \frac12 \norm{\nabla^4 G}_\infty |D_0|^2 ,
    \end{equation*}
    and, since $\nabla^3 G(Z_0)$ and $\mathcal C(Z_0)$ are $\hat{\mathcal F}_0$-measurable while $B_h$ has mean $0$ and is independent of $\hat{\mathcal F}_0$,
    \begin{equation*}
        \hat{\mathbb E}_0\bigl[ \nabla^3 G(Z_0)[D_0] : \mathcal C(Z_0) \bigr] = -h\, \nabla^3 G(Z_0)[\nabla U(Z_0)] : \mathcal C(Z_0) .
    \end{equation*}
    With $|\nabla U(Z_0)| \Le M \bigl( 1 + |Z_0| \bigr)$ and $\hat{\mathbb E}_0 |D_0|^2 \Le C h \bigl( 1 + |Z_0| \bigr)^2$, \eqref{eq: trace bound} and $\Sigma \Le M^2$ give
    \begin{align}
        & \frac{h^2}{2} \Bigl| \hat{\mathbb E}_0\bigl[ \nabla^2 G(Z_0 + D_0) - \nabla^2 G(Z_0) \bigr] : \mathcal C(Z_0) \Bigr| \Le \frac{\Sigma h^2}{2} \Bigl( h \norm{\nabla^3 G}_\infty |\nabla U(Z_0)| + \frac12 \norm{\nabla^4 G}_\infty\, \hat{\mathbb E}_0 |D_0|^2 \Bigr) \notag \\
        & \quad \Le C h^3 \bigl( \norm{\nabla^3 G}_\infty + \norm{\nabla^4 G}_\infty \bigr) \bigl( 1 + |Z_0| \bigr)^2 .
        \label{eq: sgld hessian replacement}
    \end{align}

    The state $Z_0 + D_0$ is compared with $X_h$ as the state $Y_1$ after one step of \eqref{eq: lm markov} is compared with $\tilde Y_0(h)$ in the proof of Lemma~\ref{lem: weak}. Since $X_h = Z_0 + \sqrt2\, B_h - \int_0^h \nabla U(X_s)\, \D s$, the local error is
    \begin{equation*}
        Z_0 + D_0 - X_h = \int_0^h \bigl( \nabla U(X_s) - \nabla U(Z_0) \bigr) \D s ,
    \end{equation*}
    and Assumption~\ref{ass: U} with the bound on $X_s - Z_0$ in \eqref{eq: exact sizes} for \eqref{eq: langevin} give, with $\pnorm{\cdot}_4$ taken under $\hat{\mathbb E}_0$,
    \begin{equation}
        \pnorm{Z_0 + D_0 - X_h}_4 \Le C h^{\frac32} \bigl( 1 + |Z_0| \bigr) , \qquad \pnorm{D_0}_4 + \pnorm{X_h - Z_0}_4 \Le C \sqrt h \bigl( 1 + |Z_0| \bigr) .
        \label{eq: sgld sizes}
    \end{equation}
    With $Z_0$, $Z_0 + D_0$, $X_h$ and $X_h - Z_0$ in place of $Y_0$, $Y_1$, $\tilde Y_0(h)$ and $\tilde W_h$, the expansion \eqref{eq: G expansion full} reads
    \begin{equation}
        G(Z_0 + D_0) - G(X_h) = \sum_{j=1}^{3} \frac{1}{j!}\, \nabla^j G(Z_0) : \Bigl( D_0^{\otimes j} - (X_h - Z_0)^{\otimes j} \Bigr) + \varrho_4^{D_0} ,
        \label{eq: sgld G expansion}
    \end{equation}
    with a remainder $\varrho_4^{D_0}$ that satisfies
    \begin{equation}
        |\varrho_4^{D_0}| \Le \frac16\, \norm{\nabla^4 G}_\infty\, \bigl| Z_0 + D_0 - X_h \bigr|\, \max\bigl\{ |D_0| ,\, |X_h - Z_0| \bigr\}^3 .
        \label{eq: sgld remainder}
    \end{equation}
    Taking $\hat{\mathbb E}_0$ in \eqref{eq: sgld G expansion}, H\"older's inequality with \eqref{eq: sgld remainder} and \eqref{eq: sgld sizes} gives $\hat{\mathbb E}_0 |\varrho_4^{D_0}| \Le C h^3 \norm{\nabla^4 G}_\infty \bigl( 1 + |Z_0| \bigr)^4$, so
    \begin{equation}
        \biggl| \hat{\mathbb E}_0\bigl[ G(Z_0 + D_0) - G(X_h) \bigr] - \sum_{j=1}^{3} \frac{1}{j!}\, \nabla^j G(Z_0) : \hat{\mathbb E}_0\Bigl[ D_0^{\otimes j} - (X_h - Z_0)^{\otimes j} \Bigr] \biggr| \Le C h^3 \norm{\nabla^4 G}_\infty \bigl( 1 + |Z_0| \bigr)^4 .
        \label{eq: sgld em expansion}
    \end{equation}

    The moments of $D_0$ in \eqref{eq: sgld em expansion} are explicit. With all derivatives of $U$ evaluated at $Z_0$ and $1 \Le i,j,l \Le d$, since $\sqrt2\, B_h \sim \mathcal N(0, 2h I_d)$ and $\nabla U(Z_0)$ is fixed,
    \begin{equation*}
    \begin{gathered}
        \hat{\mathbb E}_0[D_0] = -h\, \nabla U , \qquad \hat{\mathbb E}_0\bigl[ D_0^{\otimes 2} \bigr] = 2h\, I_d + h^2\, \nabla U^{\otimes 2} , \\
        \Bigl( \hat{\mathbb E}_0\bigl[ D_0^{\otimes 3} \bigr] \Bigr)_{ijl} = -2h^2 \bigl( \delta_{ij}\, \partial_l U + \delta_{il}\, \partial_j U + \delta_{jl}\, \partial_i U \bigr) - h^3\, \partial_i U\, \partial_j U\, \partial_l U .
    \end{gathered}
    \end{equation*}
    The moments of $X_h - Z_0$ are expanded in \eqref{eq: exact m1}--\eqref{eq: exact m3} of Lemma~\ref{lem: exact moments} for \eqref{eq: langevin}, with $U$, $Z_0$, $X_h$ and $\hat{\mathbb E}_0$ in place of $\tilde U_h$, $Y_0$, $\tilde Y_0(h)$ and $\mathbb E_0$. Next we subtract them from the moments of $D_0$, for $j = 1, 2, 3$.
    \begin{enumerate}
    \item For $j = 1$, the terms $-h\, \nabla U$ cancel, and the term of order $h^2$ in \eqref{eq: exact m1} remains with the opposite sign:
    \begin{equation}
        \Bigl| \hat{\mathbb E}_0\bigl[ D_0 - (X_h - Z_0) \bigr] + \frac{h^2}{2} \bigl( \nabla^2 U\, \nabla U - \nabla \Delta U \bigr) \Bigr| \Le C h^3 \bigl( 1 + |Z_0| \bigr)^3 .
        \label{eq: sgld m1}
    \end{equation}
    \item For $j = 2$, the terms $2h\, I_d$ and $h^2\, \nabla U^{\otimes 2}$ cancel, and the term $-2h^2\, \nabla^2 U$ in \eqref{eq: exact m2} remains with the opposite sign:
    \begin{equation}
        \Bigl\| \hat{\mathbb E}_0\bigl[ D_0^{\otimes 2} - (X_h - Z_0)^{\otimes 2} \bigr] - 2h^2\, \nabla^2 U \Bigr\| \Le C h^3 \bigl( 1 + |Z_0| \bigr)^4 .
        \label{eq: sgld m2}
    \end{equation}
    \item For $j = 3$, the terms of order $h^2$ cancel, and the term $-h^3\, \nabla U^{\otimes 3}$ of $\hat{\mathbb E}_0\bigl[ D_0^{\otimes 3} \bigr]$ has norm $h^3 |\nabla U(Z_0)|^3 \Le C h^3 \bigl( 1 + |Z_0| \bigr)^3$:
    \begin{equation}
        \Bigl\| \hat{\mathbb E}_0\bigl[ D_0^{\otimes 3} - (X_h - Z_0)^{\otimes 3} \bigr] \Bigr\| \Le C h^3 \bigl( 1 + |Z_0| \bigr)^6 .
        \label{eq: sgld m3}
    \end{equation}
    \end{enumerate}
    Inserting \eqref{eq: sgld m1}--\eqref{eq: sgld m3} in \eqref{eq: sgld em expansion} gives
    \begin{align}
        & \Bigl| \hat{\mathbb E}_0\bigl[ G(Z_0 + D_0) - G(X_h) \bigr] - h^2\, \nabla^2 G(Z_0) : \nabla^2 U(Z_0) + \frac{h^2}{2}\, \nabla G(Z_0) \cdot \bigl( \nabla^2 U(Z_0)\, \nabla U(Z_0) - \nabla \Delta U(Z_0) \bigr) \Bigr| \notag \\
        & \quad \Le C h^3 \sum_{j=1}^{4} \norm{\nabla^j G}_\infty \bigl( 1 + |Z_0| \bigr)^6 .
        \label{eq: sgld em weak}
    \end{align}
    Since $(e^{h\mathcal L} G)(Z_0) = \hat{\mathbb E}_0[G(X_h)]$, the sum of \eqref{eq: sgld noise}, \eqref{eq: sgld hessian replacement} and \eqref{eq: sgld em weak} proves \eqref{eq: sgld weak}.
\end{proof}

\begin{proof}[Proof of Theorem~\ref{thm: sgld bias}]
    Here $\phi_h$ is the solution \eqref{eq: discrete poisson solution} of \eqref{eq: discrete poisson}. Let $h_0 \Le \lambda^{-1}$, with $\lambda$ the rate of Theorem~\ref{thm: u decay}, so that \eqref{eq: poisson bound} holds. Let $Z_0 \sim \bar\pi_h$, independent of $(B_t)_{t \Ge 0}$ and of $(\xi_k)_{k \Ge 0}$, so that $Z_1 \sim \bar\pi_h$, and add $Z_0$ to the generators of $\hat{\mathcal F}_k$ in \eqref{eq: sg filtration}, as in the proof of Theorem~\ref{thm: numerical bias}, so that $\hat{\mathcal F}_0 = \sigma(Z_0)$. The discrete Poisson equation \eqref{eq: discrete poisson} at $Z_0$, with $(e^{h\mathcal L} \phi_h)(Z_0) = \hat{\mathbb E}_0[\phi_h(X_h)]$, gives after taking expectations
    \begin{equation}
        \bar\pi_h(f) - \pi(f) = \frac1h\, \mathbb E\bigl[ \phi_h(Z_0) - \phi_h(X_h) \bigr] = \frac1h\, \mathbb E\bigl[ \phi_h(Z_1) - \phi_h(X_h) \bigr] ,
        \label{eq: sgld bias identity}
    \end{equation}
    where $\mathbb E\, \phi_h(Z_0) = \mathbb E\, \phi_h(Z_1)$ since $Z_1 \sim \bar\pi_h$. The expectations are finite, as in the proof of Theorem~\ref{thm: numerical bias}, since $\phi_h$ is Lipschitz with constant $C L$ by \eqref{eq: poisson bound}, $\bar\pi_h(|x|^{12}) \Le C$ by Lemma~\ref{lem: sgld invariant}, and $\mathbb E|X_h - Z_0| \Le C$ by the bound on $X_s - Z_0$ in \eqref{eq: exact sizes} for \eqref{eq: langevin}, averaged over $Z_0 \sim \bar\pi_h$. By \eqref{eq: poisson bound}, $G = \phi_h$ has derivatives of orders $1$ to $4$ bounded by $C L$, so \eqref{eq: sgld weak} of Lemma~\ref{lem: sgld weak} at $k = 0$, averaged over $Z_0 \sim \bar\pi_h$, gives
    \begin{equation*}
    \begin{aligned}
        & \biggl| \mathbb E\bigl[ \phi_h(Z_1) - \phi_h(X_h) \bigr] - \frac{h^2}{2}\, \mathbb E\Bigl[ \nabla^2 \phi_h(Z_0) : \bigl( 2\nabla^2 U(Z_0) + \mathcal C(Z_0) \bigr) \Bigr] \\
        & + \frac{h^2}{2}\, \mathbb E\Bigl[ \nabla \phi_h(Z_0) \cdot \bigl( \nabla^2 U(Z_0)\, \nabla U(Z_0) - \nabla \Delta U(Z_0) \bigr) \Bigr] \biggr| \Le C L h^3\, \mathbb E \bigl( 1 + |Z_0| \bigr)^6 \Le C L h^3 .
    \end{aligned}
    \end{equation*}
    Dividing by $h$ and inserting \eqref{eq: sgld bias identity} gives \eqref{eq: sgld bias}. Finally, \eqref{eq: trace bound}, $\Sigma \Le M^2$, \eqref{eq: poisson bound} and Assumption~\ref{ass: U} give
    \begin{equation*}
    \begin{aligned}
        & \biggl| \mathbb E\Bigl[ \nabla^2 \phi_h(Z_0) : \bigl( 2\nabla^2 U(Z_0) + \mathcal C(Z_0) \bigr) \Bigr] \biggr| \\
        & + \biggl| \mathbb E\Bigl[ \nabla \phi_h(Z_0) \cdot \bigl( \nabla^2 U(Z_0)\, \nabla U(Z_0) - \nabla \Delta U(Z_0) \bigr) \Bigr] \biggr| \Le C L\, \mathbb E \bigl( 1 + |Z_0| \bigr) \Le C L ,
    \end{aligned}
    \end{equation*}
    which with \eqref{eq: sgld bias} gives $|\bar\pi_h(f) - \pi(f)| \Le C L h$.
\end{proof}

\subsection{Proof of Corollary~\ref{cor: sgld coefficient}: coefficient of the numerical bias of SGLD}
\label{sec: sgld coefficient proof}

\begin{proof}[Proof of Corollary~\ref{cor: sgld coefficient}]
    Let
    \begin{equation}
    \begin{aligned}
        \bar A_h(x) & := \nabla^2 \phi_h(x) : \bigl( 2\nabla^2 U(x) + \mathcal C(x) \bigr) - \nabla \phi_h(x) \cdot \bigl( \nabla^2 U(x)\, \nabla U(x) - \nabla \Delta U(x) \bigr) , \\
        \bar A_0(x) & := \nabla^2 \phi_0(x) : \bigl( 2\nabla^2 U(x) + \mathcal C(x) \bigr) - \nabla \phi_0(x) \cdot \bigl( \nabla^2 U(x)\, \nabla U(x) - \nabla \Delta U(x) \bigr) ,
    \end{aligned}
        \label{eq: sgld coefficient integrand}
    \end{equation}
    so that dividing \eqref{eq: sgld bias} by $h$ gives
    \begin{equation*}
        \Bigl| \frac{\bar\pi_h(f) - \pi(f)}{h} - \frac12\, \bar\pi_h(\bar A_h) \Bigr| \Le C L h .
    \end{equation*}
    It suffices to prove $\bar\pi_h(\bar A_h) \to \pi(\bar A_0)$ and to compute $\pi(\bar A_0)$.

    The measures $\bar\pi_h$ converge weakly to $\pi$ as $h \to 0$. For $\varphi \in C_c^\infty(\mathbb R^d)$, which satisfies Assumption~\ref{ass: f} with some constant $L_\varphi$, Theorem~\ref{thm: sgld bias} gives $|\bar\pi_h(\varphi) - \pi(\varphi)| \Le C L_\varphi h$, and Lemma~\ref{lem: sgld invariant} gives, uniformly in $h$,
    \begin{equation}
        \bar\pi_h\Bigl( \bigl( 1 + |x| \bigr)\, \mathbf 1_{\{|x| > \rho\}} \Bigr) \Le \rho^{-11}\, \bar\pi_h\Bigl( \bigl( 1 + |x| \bigr)^{12} \Bigr) \Le C \rho^{-11} , \qquad \rho > 0 .
        \label{eq: sgld tail}
    \end{equation}
    The approximation of a bounded continuous $b$ by $\varphi \in C_c^\infty(\mathbb R^d)$ in Appendix~\ref{sec: bias coefficient proof}, with $\bar\pi_h$ and \eqref{eq: sgld tail} in place of $\mu_h$ and \eqref{eq: mu tail}, then gives $\bar\pi_h(b) \to \pi(b)$ for every bounded continuous $b$.

    The gradient and the Hessian of $\phi_h$ converge to those of $\phi_0$ uniformly on compact sets. The argument for the convergence of $\nabla^2 \tilde\phi_h$ in Appendix~\ref{sec: bias coefficient proof}, through \eqref{eq: modified dynkin} and \eqref{eq: g consistency}, applies with $U$, $\mathcal L$, $f$, $\pi(f)$, $\phi_h$, the solution $X_s$ of \eqref{eq: langevin} and its initial value $X_0$ in place of $\tilde U_h$, $\tilde{\mathcal L}_h$, $\tilde f_h$, $\tilde\pi_h(\tilde f_h)$, $\tilde\phi_h$, $\tilde Y_0(s)$ and $Y_0$. It then uses \eqref{eq: discrete poisson}, \eqref{eq: poisson bound} and Assumption~\ref{ass: U} in place of \eqref{eq: modified poisson}, \eqref{eq: modified poisson bound} and Lemma~\ref{lem: modified class}, and the differences $\nabla \tilde U_h - \nabla U$, $\tilde f_h - f$ and $\tilde\pi_h(\tilde f_h) - \pi(f)$ estimated there are $0$. It gives, for every sequence $h_n \to 0$, a subsequence along which the derivatives of $\phi_{h_n}$ of orders $1$ to $3$ converge uniformly on compact sets to those of a function that differs from $\phi_0$ by a constant, hence
    \begin{equation*}
        \sup_{|x| \Le \rho} \Bigl( |\nabla \phi_h(x) - \nabla \phi_0(x)| + \norm{\nabla^2 \phi_h(x) - \nabla^2 \phi_0(x)} \Bigr) \to 0 , \qquad h \to 0 , \quad \rho > 0 .
    \end{equation*}
    By \eqref{eq: sgld coefficient integrand}, \eqref{eq: trace bound}, $\Sigma \Le M^2$ and Assumption~\ref{ass: U},
    \begin{equation*}
        |\bar A_h(x) - \bar A_0(x)| \Le C \norm{\nabla^2 \phi_h(x) - \nabla^2 \phi_0(x)} + C \bigl( 1 + |x| \bigr) |\nabla \phi_h(x) - \nabla \phi_0(x)| ,
    \end{equation*}
    so $\bar A_h \to \bar A_0$ uniformly on compact sets, and \eqref{eq: poisson bound} and \eqref{eq: continuous poisson bound} give
    \begin{equation*}
        |\bar A_h(x)| + |\bar A_0(x)| \Le C L \bigl( 1 + |x| \bigr) , \qquad x \in \mathbb R^d .
    \end{equation*}
    The function $\bar A_0$ is continuous, since $\nabla \phi_0$, $\nabla^2 \phi_0$ and the derivatives of $U$ are continuous and $\mathcal C(x)$ is Lipschitz by \eqref{eq: covariance lipschitz}. For every $\rho > 0$,
    \begin{equation*}
        \bigl| \bar\pi_h(\bar A_h) - \pi(\bar A_0) \bigr| \Le \sup_{|x| \Le \rho} |\bar A_h(x) - \bar A_0(x)| + C L\, \bar\pi_h\Bigl( \bigl( 1 + |x| \bigr)\, \mathbf 1_{\{|x| > \rho\}} \Bigr) + \bigl| \bar\pi_h(\bar A_0) - \pi(\bar A_0) \bigr| .
    \end{equation*}
    The bounded continuous function $\max\bigl\{ -C L (1+\rho), \min\{ \bar A_0, C L (1+\rho) \} \bigr\}$ equals $\bar A_0$ on $\{|x| \Le \rho\}$ and differs from $\bar A_0$ by at most $C L (1 + |x|)\, \mathbf 1_{\{|x| > \rho\}}$, so the weak convergence of $\bar\pi_h$ to $\pi$, \eqref{eq: sgld tail} and the moments of $\pi$ in Lemma~\ref{lem: dissipative} give $\limsup_{h \to 0} |\bar\pi_h(\bar A_0) - \pi(\bar A_0)| \Le C L \rho^{-11}$. With \eqref{eq: sgld tail} for $C L\, \bar\pi_h\bigl( (1 + |x|)\, \mathbf 1_{\{|x| > \rho\}} \bigr)$,
    \begin{equation*}
        \limsup_{h \to 0}\, \bigl| \bar\pi_h(\bar A_h) - \pi(\bar A_0) \bigr| \Le C L \rho^{-11} , \qquad \rho > 0 ,
    \end{equation*}
    and letting $\rho \to \infty$ gives $\bar\pi_h(\bar A_h) \to \pi(\bar A_0)$.

    It remains to compute $\pi(\bar A_0)$. Since $\phi_0$ has bounded derivatives of orders $1$ to $3$ by \eqref{eq: continuous poisson bound}, \eqref{eq: sgld integration by parts} with $\psi = \phi_0$ gives
    \begin{equation*}
        \pi(\bar A_0) = \pi\Bigl( \nabla^2 \phi_0(x) : \bigl( 2\nabla^2 U(x) + \mathcal C(x) \bigr) \Bigr) - \pi\bigl( \nabla^2 \phi_0(x) : \nabla^2 U(x) \bigr) = \pi\Bigl( \nabla^2 \phi_0(x) : \bigl( \nabla^2 U(x) + \mathcal C(x) \bigr) \Bigr) ,
    \end{equation*}
    which with $\bar\pi_h(\bar A_h) \to \pi(\bar A_0)$ proves \eqref{eq: sgld coefficient}.
\end{proof}

\section{Proofs for the explicit coefficient in Wasserstein-1 error}
\label{app: w1}

Throughout this appendix, the constant $C$ depends only on $m$, $M$, $R$ and $d$. Let
\begin{equation}
    \mathcal D := \bigl\{ \omega \in C^4(\mathbb R^d) : \norm{\nabla^j \omega}_\infty < \infty ,\ j = 1, 2, 3, 4 \bigr\} , \qquad L^2_0(\pi) := \bigl\{ g \in L^2(\pi) : \pi(g) = 0 \bigr\} ,
    \label{eq: test class}
\end{equation}
so $\mathcal D$ is the set of the test functions $\omega$ of \eqref{eq: riesz}, and $L^2_0(\pi)$ is a Hilbert space with the inner product of $L^2(\pi)$.

\subsection{Proof of Theorem~\ref{thm: riesz}: well-posedness of the Riesz representor}
\label{sec: riesz well-posedness}

Lemma~\ref{lem: poincare} follows from \citep{bakry2008simple}. By \eqref{eq: dissipative bounds}, $U$ is bounded from below, and $\langle x, \nabla U(x) \rangle \Ge \frac m4 |x|^2 - C \Ge |x|$ outside a ball whose radius depends only on $m$, $M$, $R$ and $d$. Hence $U$ satisfies the drift condition of \citep[Corollary~1.6]{bakry2008simple}, and \citep[Theorem~1.4]{bakry2008simple} gives \eqref{eq: poincare} for smooth $g$. In the proof of \citep[Corollary~1.6]{bakry2008simple}, the Lyapunov function is an exponential of $|x|$ outside a larger ball, whose radius depends only on the previous radius and on $d$, and the constants of the inequality \citep[(1.3)]{bakry2008simple} that it satisfies depend on $U$ only through $\nabla U$ on that ball. By \citep[Theorem~1.4]{bakry2008simple}, $C_P$ depends only on these constants and on the Poincar\'e constant of $\pi$ restricted to the same ball, which \citep[Section~1]{bakry2008simple} bounds through the oscillation of $U$ on the ball. Since $|\nabla U(x)| \Le M(1 + |x|)$, $C_P$ depends only on $m$, $M$, $R$ and $d$. Mollification and dominated convergence, with the moments of Lemma~\ref{lem: dissipative}, extend \eqref{eq: poincare} to every $g \in C^1(\mathbb R^d)$ with a bounded gradient. The next lemma is the identity of Bochner under $\pi$, used in the proofs of Theorems~\ref{thm: riesz} and \ref{thm: batch threshold}.

\begin{lemma}
\label{lem: bochner}
    Let Assumption~\ref{ass: U} hold. For every $\omega \in \mathcal D$,
    \begin{equation}
        \pi\bigl( (\mathcal L \omega)^2 \bigr) = \pi\bigl( \nabla^2 \omega : \nabla^2 \omega \bigr) + \pi\bigl( \nabla \omega \cdot \nabla^2 U\, \nabla \omega \bigr) .
        \label{eq: bochner}
    \end{equation}
\end{lemma}

\begin{proof}
    Differentiating \eqref{eq: generator} gives $\partial_i \mathcal L \omega = \mathcal L \partial_i \omega - \sum_{j=1}^d \partial_i \partial_j U\, \partial_j \omega$ for each component $i$. The functions $\omega$, $\mathcal L \omega$ and $\partial_i \omega$ and their derivatives of orders $1$ and $2$ grow at most linearly, so \eqref{eq: generator ibp} with $g_1 = \mathcal L \omega$, $g_2 = \omega$, and then with $g_1 = g_2 = \partial_i \omega$, gives
    \begin{equation*}
    \begin{aligned}
        \pi\bigl( (\mathcal L \omega)^2 \bigr) & = -\sum_{i=1}^d \pi\bigl( \partial_i \omega\, \partial_i \mathcal L \omega \bigr) = -\sum_{i=1}^d \pi\bigl( \partial_i \omega\, \mathcal L \partial_i \omega \bigr) + \pi\bigl( \nabla \omega \cdot \nabla^2 U\, \nabla \omega \bigr) \\
        & = \sum_{i=1}^d \pi\bigl( |\nabla \partial_i \omega|^2 \bigr) + \pi\bigl( \nabla \omega \cdot \nabla^2 U\, \nabla \omega \bigr) ,
    \end{aligned}
    \end{equation*}
    which is \eqref{eq: bochner}.
\end{proof}

\begin{proof}[Proof of Theorem~\ref{thm: riesz}]
    The proof bounds $|\pi(\nabla^2 \omega : \mathcal A)|$ by the norm of $\mathcal L \omega$ in $L^2(\pi)$, through \eqref{eq: poincare} and \eqref{eq: bochner}, and applies the Riesz representation theorem on the dense set of the functions $-\mathcal L \omega$, $\omega \in \mathcal D$, in $L^2_0(\pi)$. Let $\omega \in \mathcal D$. Since $\omega$ grows at most linearly, \eqref{eq: generator ibp} and \eqref{eq: poincare} give
    \begin{equation*}
        \pi\bigl( |\nabla \omega|^2 \bigr) = \pi\bigl( (\omega - \pi(\omega)) (-\mathcal L \omega) \bigr) \Le C_P^{\frac12}\, \pi\bigl( |\nabla \omega|^2 \bigr)^{\frac12}\, \pi\bigl( (\mathcal L \omega)^2 \bigr)^{\frac12} ,
    \end{equation*}
    so $\pi(|\nabla \omega|^2) \Le C_P\, \pi((\mathcal L \omega)^2)$, and \eqref{eq: bochner} with $\nabla^2 U \succeq -M I_d$ gives
    \begin{equation}
        \pi\bigl( \nabla^2 \omega : \nabla^2 \omega \bigr) \Le \pi\bigl( (\mathcal L \omega)^2 \bigr) + M\, \pi\bigl( |\nabla \omega|^2 \bigr) \Le (1 + M C_P)\, \pi\bigl( (\mathcal L \omega)^2 \bigr) .
        \label{eq: hessian l2}
    \end{equation}
    Since $\mathcal A(x) : \mathcal A(x) \Le d\, \norm{\mathcal A(x)}^2$, the Cauchy--Schwarz inequality under $\pi$ gives
    \begin{equation*}
        \bigl| \pi\bigl( \nabla^2 \omega : \mathcal A \bigr) \bigr| \Le \sqrt d\, \norm{\mathcal A}_\infty\, \pi\Bigl( \bigl( \nabla^2 \omega : \nabla^2 \omega \bigr)^{\frac12} \Bigr) \Le \sqrt d\, \norm{\mathcal A}_\infty\, \pi\bigl( \nabla^2 \omega : \nabla^2 \omega \bigr)^{\frac12} ,
    \end{equation*}
    and inserting the bound on $\pi(\nabla^2 \omega : \nabla^2 \omega)$ in \eqref{eq: hessian l2} gives
    \begin{equation}
        \bigl| \pi\bigl( \nabla^2 \omega : \mathcal A \bigr) \bigr| \Le \sqrt{d\, (1 + M C_P)}\; \norm{\mathcal A}_\infty\; \pi\bigl( (\mathcal L \omega)^2 \bigr)^{\frac12} , \qquad \omega \in \mathcal D .
        \label{eq: riesz functional bound}
    \end{equation}

    The functions $-\mathcal L \omega$, $\omega \in \mathcal D$, lie in $L^2_0(\pi)$ by their linear growth and \eqref{eq: generator ibp} with $g_1 = 1$. Every $g \in C_c^\infty(\mathbb R^d)$ satisfies Assumption~\ref{ass: f}, so the solution $\phi_0$ of the continuous Poisson equation \eqref{eq: continuous poisson} with $g$ in place of $f$ belongs to $\mathcal D$ by \eqref{eq: continuous poisson bound}, and $g - \pi(g) = -\mathcal L \phi_0$. Since the functions $g - \pi(g)$, $g \in C_c^\infty(\mathbb R^d)$, are dense in $L^2_0(\pi)$, so are the functions $-\mathcal L \omega$, $\omega \in \mathcal D$. By \eqref{eq: riesz functional bound}, $-\mathcal L \omega \mapsto \pi(\nabla^2 \omega : \mathcal A)$ is a well-defined linear functional on the functions $-\mathcal L \omega$, $\omega \in \mathcal D$, bounded in the norm of $L^2_0(\pi)$, and linear in $\mathcal A$. It extends to $L^2_0(\pi)$, and the Riesz representation theorem gives the unique $\psi_{\mathcal A} \in L^2_0(\pi)$ that satisfies \eqref{eq: riesz} for every $\omega \in \mathcal D$. By this uniqueness and the linearity of the functional in $\mathcal A$, the map $\mathcal A \mapsto \psi_{\mathcal A}$ is linear, and \eqref{eq: riesz functional bound} gives \eqref{eq: riesz bound}.
\end{proof}

\subsection{Auxiliary estimates for Theorem~\ref{thm: kr expansion}}
\label{sec: kr auxiliary}

Throughout this subsection and the next, $h \Le h_0$ with $h_0 \Le \lambda^{-1}$, $\lambda$ the rate of Lemma~\ref{lem: smoothing}. For a signed measure $\nu$ on $\mathbb R^d$ with $\nu(\mathbb R^d) = 0$ and $\int (1 + |x|)\, |\nu|(\D x) < \infty$, such as $\psi\pi$ with $\psi \in L^2_0(\pi)$, the Kantorovich--Rubinstein norm of Section~\ref{subsec: notation} satisfies
\begin{equation}
    \norm{\nu}_{\mathrm{KR}} = \sup\biggl\{ \int_{\mathbb R^d} f\, \D\nu : f \text{ satisfies Assumption~\ref{ass: f}},\ \norm{\nabla f}_\infty \Le 1 \biggr\} .
    \label{eq: kr test functions}
\end{equation}
The supremum on the right side of \eqref{eq: kr test functions} is at most $\norm{\nu}_{\mathrm{KR}}$, since its functions $f$ are Lipschitz with $\norm{\nabla f}_\infty \Le 1$, and it is at least $\norm{\nu}_{\mathrm{KR}}$, since the mollifications of a Lipschitz function $g$ with $\norm{\nabla g}_\infty \Le 1$ satisfy Assumption~\ref{ass: f} with the same bound on the gradient and converge to $g$ uniformly.

In the proofs of Theorems~\ref{thm: kr expansion} and \ref{thm: batch threshold}, the test function $f$ satisfies Assumption~\ref{ass: f} and $\norm{\nabla f}_\infty \Le 1$, and the estimates do not depend on the constant $L$ of Assumption~\ref{ass: f}. The modified test function \eqref{eq: heat shift} is $\tilde f_h = f * \mathcal N(0, \frac h2 I_d)$, where $\mathcal N(0, \frac h2 I_d)$ also denotes the density of this Gaussian distribution, so
\begin{equation*}
    \nabla^j \tilde f_h = (\nabla f) * \nabla^{j-1}\mathcal N\Bigl( 0, \frac h2 I_d \Bigr) , \qquad \int_{\mathbb R^d} \Bigl\| \nabla^{j-1}\mathcal N\Bigl( 0, \frac h2 I_d \Bigr)(y) \Bigr\| \,\D y \Le C h^{-\frac{j-1}{2}} , \qquad j = 1, 2, 3, 4 ,
\end{equation*}
and
\begin{equation}
    \norm{\nabla \tilde f_h}_\infty \Le 1 , \qquad
    \norm{\nabla^j \tilde f_h}_\infty \Le C h^{-\frac{j-1}{2}} , \qquad j = 2, 3, 4 .
    \label{eq: heat derivatives}
\end{equation}

The semigroup estimates of Appendix~\ref{sec: semigroup} and Lemma~\ref{lem: fourth from lipschitz} combine into one bound for the exact and the modified dynamics. The proof technique is essentially that of \citep[Lemma~10]{ye2026mse}, which splits the times into $0 < t \Le 1$ and $t > 1$.

\begin{lemma}
\label{lem: smoothing}
    Let Assumption~\ref{ass: U} hold. There is $\lambda > 0$, depending only on $m$, $M$, $R$ and $d$, such that, for every $h \in (0,h_0]$, every $g \in \mathcal D$ of \eqref{eq: test class} and every $t > 0$,
    \begin{equation}
        \norm{\nabla^j (e^{t\mathcal L} g)}_\infty + \norm{\nabla^j (e^{t\tilde{\mathcal L}_h} g)}_\infty \Le C \min(1,t)^{-\frac{j-1}{2}}\, e^{-\lambda t}\, \norm{\nabla g}_\infty , \qquad j = 1, 2, 3, 4 .
        \label{eq: smoothing}
    \end{equation}
\end{lemma}

\begin{proof}
    For $j = 1$, the bound on $e^{t\mathcal L} g$ is \eqref{eq: grad decay}. For $j = 2, 3, 4$ and $0 < t \Le 1$, \eqref{eq: smoothing 23} and \eqref{eq: fourth from lipschitz} give
    \begin{equation*}
        \norm{\nabla^j (e^{t\mathcal L} g)}_\infty \Le C\, t^{-\frac{j-1}{2}}\, \norm{\nabla g}_\infty \Le C e^{\lambda} \min(1,t)^{-\frac{j-1}{2}}\, e^{-\lambda t}\, \norm{\nabla g}_\infty ,
    \end{equation*}
    since $\min(1,t) = t$ and $e^{-\lambda t} \Ge e^{-\lambda}$, which is \eqref{eq: smoothing} for $e^{t\mathcal L}$ with $C e^{\lambda}$ in place of $C$. For $j = 2, 3, 4$ and $t > 1$, the semigroup property gives $e^{t\mathcal L} g = e^{\mathcal L} \bigl( e^{(t-1)\mathcal L} g \bigr)$, and $e^{(t-1)\mathcal L} g \in \mathcal D$, since \eqref{eq: short time} at time $t - 1$ bounds its derivatives of orders $1$ to $4$. Hence \eqref{eq: smoothing 23} and \eqref{eq: fourth from lipschitz} at time $1$, with $e^{(t-1)\mathcal L} g$ in place of $g$, bound $\norm{\nabla^j (e^{t\mathcal L} g)}_\infty$ by $C \norm{\nabla (e^{(t-1)\mathcal L} g)}_\infty$, and \eqref{eq: grad decay} at time $t - 1$ bounds this gradient, so that
    \begin{equation*}
        \norm{\nabla^j (e^{t\mathcal L} g)}_\infty \Le C \bigl\| \nabla (e^{(t-1)\mathcal L} g) \bigr\|_\infty \Le C e^{-\lambda (t-1)} \norm{\nabla g}_\infty = C e^{\lambda} e^{-\lambda t} \norm{\nabla g}_\infty .
    \end{equation*}
    Since $\min(1,t) = 1$ for $t > 1$, this is \eqref{eq: smoothing} for $e^{t\mathcal L}$ with $C e^{\lambda}$ in place of $C$.

    The same argument applies to $e^{t\tilde{\mathcal L}_h}$. Lemmas~\ref{lem: semigroup}, \ref{lem: short time} and \ref{lem: fourth from lipschitz} use only the bounds of Assumption~\ref{ass: U} up to the order $5$ and the convexity outside the ball of radius $R$, and by Lemma~\ref{lem: modified class}, for every $h \in (0, h_0]$, $\tilde U_h$ satisfies them with $2M$ and $\frac m2$ in place of $M$ and $m$. These lemmas therefore hold for $e^{t\tilde{\mathcal L}_h}$ with a rate and a constant that depend only on $m$, $M$, $R$ and $d$, as in Corollary~\ref{cor: modified poisson}. Since the left side of \eqref{eq: smoothing} is the sum of the norms for $e^{t\mathcal L}$ and for $e^{t\tilde{\mathcal L}_h}$, \eqref{eq: smoothing} holds with $\lambda$ the smaller of the rates for $U$ and for $\tilde U_h$, and $C$ twice the larger of the constants.
\end{proof}

In the proofs of Theorem~\ref{thm: kr expansion}, the bounds \eqref{eq: smoothing} at the times $t = ih$ are summed over $i \Ge 0$ through the following elementary bounds.

\begin{lemma}
\label{lem: time sums}
    For every $h \in (0,h_0]$,
    \begin{equation}
        h \sum_{i=0}^{\infty} \min\bigl( 1, (i+1)h \bigr)^{-\frac{j-1}{2}} e^{-\lambda i h} \Le
        \begin{cases}
            C , & j = 1, 2 , \\
            C \bigl( 1 + \log\frac1h \bigr) , & j = 3 , \\
            C h^{-\frac12} , & j = 4 .
        \end{cases}
        \label{eq: time sums}
    \end{equation}
\end{lemma}

\begin{proof}
    Split the sum at $(i+1) h = 1$. The terms with $(i+1)h > 1$, that is, $i \Ge \lfloor 1/h \rfloor$, have the minimum $1$, and since $\lambda h \Le 1$,
    \begin{equation*}
        h \sum_{i = \lfloor 1/h \rfloor}^{\infty} e^{-\lambda i h} \Le \frac{h}{1 - e^{-\lambda h}} \Le \frac{2}{\lambda} .
    \end{equation*}
    For the terms with $(i+1)h \Le 1$, put $n = i+1$ and bound $e^{-\lambda ih}$ by $1$:
    \begin{equation*}
        h \sum_{n=1}^{\lfloor 1/h \rfloor} (nh)^{-\frac{j-1}{2}} = h^{\frac{3-j}{2}} \sum_{n=1}^{\lfloor 1/h \rfloor} n^{-\frac{j-1}{2}} \Le
        \begin{cases}
            1 , & j = 1 , \\
            2 , & j = 2 , \\
            1 + \log\frac1h , & j = 3 , \\
            3 h^{-\frac12} , & j = 4 .
        \end{cases}
    \end{equation*}
    Adding the bound $2/\lambda$ gives \eqref{eq: time sums}.
\end{proof}

The gradient noise enters a step of \eqref{eq: sgld} with the factor $h$, but in the next lemma it contributes only at the order $h^{\frac32}$: the Gaussian increment of the step smooths the test function, and the noise has mean $0$.

\begin{lemma}
\label{lem: sgld lipschitz step}
    Let Assumptions~\ref{ass: U} and \ref{ass: sg} hold. Then, for every Lipschitz function $f$ with $\norm{\nabla f}_\infty \Le 1$ and every $h \in (0,1]$, the chain \eqref{eq: sgld} satisfies
    \begin{equation}
        \bigl| \hat{\mathbb E}_0\bigl[ f(Z_1) \bigr] - (e^{h\mathcal L} f)(Z_0) \bigr| \Le C h^{\frac32} \bigl( 1 + |Z_0| \bigr) .
        \label{eq: sgld lipschitz step}
    \end{equation}
\end{lemma}

\begin{proof}
    Write the step from $Z_0$ as in \eqref{eq: sgld step}, with $D_0$, $\eta_0$ and the solution $X_s$ of \eqref{eq: langevin} from $X_0 = Z_0$ as in Appendix~\ref{app: sgld}. Then
    \begin{equation}
        \hat{\mathbb E}_0\bigl[ f(Z_1) \bigr] - (e^{h\mathcal L} f)(Z_0) = \hat{\mathbb E}_0\bigl[ f(Z_0 + D_0 - h\eta_0) - f(Z_0 + D_0) \bigr] + \hat{\mathbb E}_0\bigl[ f(Z_0 + D_0) - f(X_h) \bigr] .
        \label{eq: sgld step split}
    \end{equation}

    Since $\norm{\nabla f}_\infty \Le 1$, the bound on $Z_0 + D_0 - X_h$ in \eqref{eq: sgld sizes} gives
    \begin{equation*}
        \Bigl| \hat{\mathbb E}_0\bigl[ f(Z_0 + D_0) - f(X_h) \bigr] \Bigr| \Le \hat{\mathbb E}_0 |Z_0 + D_0 - X_h| \Le C h^{\frac32} \bigl( 1 + |Z_0| \bigr) .
    \end{equation*}
    For the expectation of $f(Z_0 + D_0 - h\eta_0) - f(Z_0 + D_0)$ in \eqref{eq: sgld step split}, condition on $\hat{\mathcal F}_0$ and $\eta_0$. Since $\sqrt2\, B_h \sim \mathcal N(0, 2h I_d)$ is independent of $\sigma(\hat{\mathcal F}_0, \eta_0)$,
    \begin{equation*}
        \mathbb E\bigl[ f(Z_0 + D_0 - h\eta_0) - f(Z_0 + D_0) \,\big|\, \hat{\mathcal F}_0, \eta_0 \bigr] = (e^{h\Delta} f)\bigl( Z_0 - h\, \nabla U(Z_0) - h\eta_0 \bigr) - (e^{h\Delta} f)\bigl( Z_0 - h\, \nabla U(Z_0) \bigr) ,
    \end{equation*}
    with $e^{h\Delta} f = f * \mathcal N(0, 2h I_d)$, for which $\norm{\nabla^2 e^{h\Delta} f}_\infty \Le C h^{-\frac12}$ as in \eqref{eq: heat derivatives}. Taylor's formula in $h\eta_0$ gives
    \begin{equation*}
        \Bigl| \mathbb E\bigl[ f(Z_0 + D_0 - h\eta_0) - f(Z_0 + D_0) \,\big|\, \hat{\mathcal F}_0, \eta_0 \bigr] + h\, \nabla (e^{h\Delta} f)\bigl( Z_0 - h\, \nabla U(Z_0) \bigr) \cdot \eta_0 \Bigr| \Le C h^{\frac32}\, |\eta_0|^2 .
    \end{equation*}
    Since $\nabla (e^{h\Delta} f)\bigl( Z_0 - h\, \nabla U(Z_0) \bigr)$ is $\hat{\mathcal F}_0$-measurable and $\hat{\mathbb E}_0[\eta_0] = 0$, taking $\hat{\mathbb E}_0$ removes the term linear in $\eta_0$, and $\hat{\mathbb E}_0 |\eta_0|^2 \Le \Sigma \Le M^2$ by Assumption~\ref{ass: sg} gives
    \begin{equation*}
        \Bigl| \hat{\mathbb E}_0\bigl[ f(Z_0 + D_0 - h\eta_0) - f(Z_0 + D_0) \bigr] \Bigr| \Le C h^{\frac32}\, \hat{\mathbb E}_0 |\eta_0|^2 \Le C h^{\frac32} .
    \end{equation*}
    With the bound $C h^{\frac32} (1 + |Z_0|)$ on the expectation of $f(Z_0 + D_0) - f(X_h)$, \eqref{eq: sgld step split} then gives \eqref{eq: sgld lipschitz step}.
\end{proof}

In the proof of \eqref{eq: kr sgld}, the integrand of the sum in \eqref{eq: sgld leading} carries $\nabla^2 U\, \nabla U - \nabla \Delta U$, whose gradient grows linearly, so the integrand satisfies only the bound \eqref{eq: local lipschitz} below, and $\bar\pi_h$ is replaced by $\pi$ in the coefficient of $\frac{h^2}{2}$ in \eqref{eq: sgld leading} through the following truncation.

\begin{lemma}
\label{lem: local lipschitz}
    Let $\mu$ and $\nu$ be probability measures on $\mathbb R^d$ with $\mu(|x|^{12}) + \nu(|x|^{12}) < \infty$ and $\mathcal W_1(\mu, \nu) \Le 1$, and let $F : \mathbb R^d \to \mathbb R$ satisfy, for some $K \Ge 0$,
    \begin{equation}
        |F(x) - F(x')| \Le K \bigl( 1 + |x| + |x'| \bigr) |x - x'| , \qquad x, x' \in \mathbb R^d .
        \label{eq: local lipschitz}
    \end{equation}
    Then $|\mu(F) - \nu(F)| \Le 3 K \bigl( 1 + \mu(|x|^{12}) + \nu(|x|^{12}) \bigr)\, \mathcal W_1(\mu, \nu)^{\frac{10}{11}}$.
\end{lemma}

\begin{proof}
    If $\mathcal W_1(\mu,\nu) = 0$, then $\mu = \nu$. Otherwise let $\rho := \mathcal W_1(\mu,\nu)^{-\frac1{11}} \Ge 1$, and let $T_\rho(x) := x$ for $|x| \Le \rho$ and $T_\rho(x) := \rho x / |x|$ for $|x| > \rho$, the retraction onto the closed ball of radius $\rho$, which is $1$-Lipschitz. By \eqref{eq: local lipschitz}, $F \circ T_\rho$ is Lipschitz with constant $K(1 + 2\rho) \Le 3K\rho$, so the Kantorovich--Rubinstein duality gives
    \begin{equation*}
        \bigl| \mu(F \circ T_\rho) - \nu(F \circ T_\rho) \bigr| \Le 3 K \rho\, \mathcal W_1(\mu,\nu) .
    \end{equation*}
    For $|x| > \rho$, $|T_\rho(x)| = \rho < |x|$ and $|x - T_\rho(x)| \Le |x|$, so \eqref{eq: local lipschitz} and $|x| > \rho \Ge 1$ give
    \begin{equation*}
        |F(x) - F(T_\rho(x))| \Le K \bigl( 1 + 2|x| \bigr) |x| \Le 3 K |x|^2 \Le 3 K \rho^{-10} |x|^{12} ,
    \end{equation*}
    while $F = F \circ T_\rho$ on $\{|x| \Le \rho\}$. Hence
    \begin{equation*}
        |\mu(F) - \nu(F)| \Le 3 K \rho\, \mathcal W_1(\mu,\nu) + 3 K \rho^{-10} \bigl( \mu(|x|^{12}) + \nu(|x|^{12}) \bigr) = 3 K \bigl( 1 + \mu(|x|^{12}) + \nu(|x|^{12}) \bigr)\, \mathcal W_1(\mu,\nu)^{\frac{10}{11}} ,
    \end{equation*}
    with $F$ integrable under $\mu$ and $\nu$ since $|F(x)| \Le |F(0)| + K(1 + |x|)|x|$.
\end{proof}

\subsection{Proof of Theorem~\ref{thm: kr expansion}: expansions in the Kantorovich--Rubinstein norm}
\label{sec: kr expansion proof}

Since $\mathcal C(x)$ is positive semidefinite, $\norm{\mathcal C}_\infty \Le \Sigma \Le M^2$ by \eqref{eq: covariance}, and $\norm{\nabla^2 U}_\infty \Le M$. The bound \eqref{eq: riesz bound} of Theorem~\ref{thm: riesz}, the Cauchy--Schwarz inequality and the moments of Lemma~\ref{lem: dissipative} then give
\begin{equation}
    \pi\bigl( |\psi_{\mathcal C}|\, (1 + |x|)^2 \bigr) + \pi\bigl( |\psi_{\nabla^2 U} + \psi_{\mathcal C}|\, (1 + |x|)^2 \bigr) \Le C .
    \label{eq: psi moments}
\end{equation}

\begin{proof}[Proof of \eqref{eq: kr sglm} in Theorem~\ref{thm: kr expansion}]
    Recall that $\mu_h$ is the invariant distribution of \eqref{eq: sglm markov} and $\tilde\pi_h \propto e^{-\tilde U_h}$ is the modified target, the invariant distribution of \eqref{eq: modified langevin}. The proof writes $\mu_h(\tilde f_h) - \tilde\pi_h(\tilde f_h)$, by the telescoping formula \eqref{eq: sglm telescope}, as a series of local errors of Lemma~\ref{lem: sg weak} on the functions
    \begin{equation}
        \tilde{\mathcal V}_i := e^{ih\tilde{\mathcal L}_h} \tilde f_h , \qquad i = 0, 1, 2, \dots ,
        \label{eq: propagated test functions}
    \end{equation}
    moves the leading term of each to an integral against $\pi$, and sums the series through \eqref{eq: riesz} into $\frac h2\, \pi(f \psi_{\mathcal C})$.
    
    Let $f$ satisfy Assumption~\ref{ass: f} and $\norm{\nabla f}_\infty \Le 1$.
    Let $Y_0 \sim \mu_h$, independent of $(B_t)_{t \Ge 0}$ and $(\xi_k)_{k \Ge 0}$, so that $Y_1 \sim \mu_h$. As in the proof of Theorem~\ref{thm: numerical bias}, add $Y_0$ to the generators of $\hat{\mathcal F}_k$ in \eqref{eq: sg filtration}, so that $\hat{\mathcal F}_0 = \sigma(Y_0)$, while $\chi_0$ and $\xi_0$ remain independent of $\hat{\mathcal F}_0$. Then $Z_0 = Y_0 + \sqrt{h/2}\, \chi_0 \sim \pi_h$. By \eqref{eq: heat shift invariant}, $\pi_h(f) = \mu_h(\tilde f_h)$, so
    \begin{equation}
        \pi_h(f) - \pi(f) = \bigl( \mu_h(\tilde f_h) - \tilde\pi_h(\tilde f_h) \bigr) + \bigl( \tilde\pi_h(\tilde f_h) - \pi(f) \bigr) ,
        \label{eq: sglm split}
    \end{equation}
    and \eqref{eq: shifted bias} bounds the bracket $\tilde\pi_h(\tilde f_h) - \pi(f)$ of \eqref{eq: sglm split} by $C h^2$.

    For $i = 0$, $\tilde{\mathcal V}_0 = \tilde f_h$ and \eqref{eq: heat derivatives} bounds its derivatives, and for $i \Ge 1$, \eqref{eq: smoothing} at $t = ih$ with $\min(1, ih) \Ge \frac12 \min(1, (i+1)h)$ bounds those of $\tilde{\mathcal V}_i$, so that, for $i \Ge 0$ and $j = 1, 2, 3, 4$,
    \begin{equation}
        \norm{\nabla^j \tilde{\mathcal V}_i}_\infty \Le C \min\bigl( 1, (i+1)h \bigr)^{-\frac{j-1}{2}} e^{-\lambda i h} ,
        \label{eq: propagated derivatives}
    \end{equation}
    and $\tilde{\mathcal V}_i \in \mathcal D$. Since $\tilde Y_0(h)$ solves the modified dynamics \eqref{eq: modified langevin} from $Y_0$, driven by Brownian increments independent of $\hat{\mathcal F}_0$, $\hat{\mathbb E}_0[\tilde{\mathcal V}_i(\tilde Y_0(h))] = \tilde{\mathcal V}_{i+1}(Y_0)$, and since $Y_0$ and $Y_1$ have the law $\mu_h$,
    \begin{equation*}
        \mathbb E\bigl[ \tilde{\mathcal V}_i(Y_1) - \tilde{\mathcal V}_i(\tilde Y_0(h)) \bigr] = \mu_h(\tilde{\mathcal V}_i) - \mu_h(\tilde{\mathcal V}_{i+1}) .
    \end{equation*}
    The modified target $\tilde\pi_h$ is invariant under $e^{t\tilde{\mathcal L}_h}$ and has a moment of order $1$ bounded uniformly in $h$, as in the proof of Lemma~\ref{lem: shifted bias}, so \eqref{eq: propagated derivatives} at $j = 1$ gives
    \begin{equation}
        \bigl| \tilde{\mathcal V}_N(y) - \tilde\pi_h(\tilde f_h) \bigr| \Le \norm{\nabla \tilde{\mathcal V}_N}_\infty \int_{\mathbb R^d} |y - y'|\, \tilde\pi_h(\D y') \Le C e^{-\lambda N h} \bigl( 1 + |y| \bigr) ,
        \label{eq: sglm relaxation}
    \end{equation}
    and $\mu_h(\tilde{\mathcal V}_N) \to \tilde\pi_h(\tilde f_h)$ as $N \to \infty$. Summing over $i$,
    \begin{equation}
        \mu_h(\tilde f_h) - \tilde\pi_h(\tilde f_h) = \sum_{i=0}^{\infty} \mathbb E\bigl[ \tilde{\mathcal V}_i(Y_1) - \tilde{\mathcal V}_i(\tilde Y_0(h)) \bigr] .
        \label{eq: sglm telescope}
    \end{equation}

    By \eqref{eq: sg weak} at $k = 0$ with $G = \tilde{\mathcal V}_i$, the tower property for the $\hat{\mathcal F}_0$-measurable factor $\nabla^2 \tilde{\mathcal V}_i(Y_0)$, and the moment bound $\mu_h(|y|^{12}) \Le C$ of Lemma~\ref{lem: invariant},
    \begin{equation}
        \Bigl| \mathbb E\bigl[ \tilde{\mathcal V}_i(Y_1) - \tilde{\mathcal V}_i(\tilde Y_0(h)) \bigr] - \frac{h^2}{2}\, \mathbb E\bigl[ \nabla^2 \tilde{\mathcal V}_i(Y_0) : \mathcal C(Z_0) \bigr] \Bigr| \Le C h^3 \sum_{j=1}^{4} \norm{\nabla^j \tilde{\mathcal V}_i}_\infty .
        \label{eq: sglm local}
    \end{equation}
    By \eqref{eq: propagated derivatives} and \eqref{eq: time sums}, the remainders of \eqref{eq: sglm local} satisfy
    \begin{equation*}
        C h^3 \sum_{i=0}^{\infty} \sum_{j=1}^{4} \norm{\nabla^j \tilde{\mathcal V}_i}_\infty \Le C h^2 \Bigl( 1 + \log\frac1h + h^{-\frac12} \Bigr) \Le C h^{\frac32} .
    \end{equation*}
    With \eqref{eq: sglm split}, \eqref{eq: shifted bias} and \eqref{eq: sglm telescope},
    \begin{equation}
        \biggl| \pi_h(f) - \pi(f) - \frac{h^2}{2} \sum_{i=0}^{\infty} \mathbb E\bigl[ \nabla^2 \tilde{\mathcal V}_i(Y_0) : \mathcal C(Z_0) \bigr] \biggr| \Le C h^{\frac32} .
        \label{eq: sglm leading}
    \end{equation}
    By \eqref{eq: trace bound}, \eqref{eq: propagated derivatives} and \eqref{eq: time sums} at $j = 2$, the sum in \eqref{eq: sglm leading} satisfies
    \begin{equation*}
        \frac{h^2}{2} \sum_{i=0}^{\infty} \Bigl| \mathbb E\bigl[ \nabla^2 \tilde{\mathcal V}_i(Y_0) : \mathcal C(Z_0) \bigr] \Bigr| \Le \frac{\Sigma h^2}{2} \sum_{i=0}^{\infty} \norm{\nabla^2 \tilde{\mathcal V}_i}_\infty \Le C \Sigma h .
    \end{equation*}
    Hence $|\pi_h(f) - \pi(f)| \Le C(\Sigma h + h^{\frac32})$ for every $f$ satisfying Assumption~\ref{ass: f} and $\norm{\nabla f}_\infty \Le 1$, and \eqref{eq: kr test functions} gives
    \begin{equation}
        \mathcal W_1(\pi_h, \pi) \Le C \bigl( \Sigma h + h^{\frac32} \bigr) \Le C h .
        \label{eq: w1 rough}
    \end{equation}

    In the coefficient of $\frac{h^2}{2}$ in \eqref{eq: sglm leading}, the expectation can be replaced by an integral against $\pi$. Since $Z_0 - Y_0 = \sqrt{h/2}\, \chi_0$, \eqref{eq: trace bound} gives
    \begin{equation}
        \Bigl| \mathbb E\bigl[ \nabla^2 \tilde{\mathcal V}_i(Y_0) : \mathcal C(Z_0) \bigr] - \pi_h\bigl( \nabla^2 \tilde{\mathcal V}_i : \mathcal C \bigr) \Bigr| = \Bigl| \mathbb E\bigl[ \bigl( \nabla^2 \tilde{\mathcal V}_i(Y_0) - \nabla^2 \tilde{\mathcal V}_i(Z_0) \bigr) : \mathcal C(Z_0) \bigr] \Bigr| \Le C \sqrt h\, \norm{\nabla^3 \tilde{\mathcal V}_i}_\infty .
        \label{eq: sglm point replacement}
    \end{equation}
    By \eqref{eq: trace bound}, \eqref{eq: covariance lipschitz} and $|A : B| \Le d \norm{A} \norm{B}$ for $A, B \in \mathbb R^{d \times d}$, the function $\nabla^2 \tilde{\mathcal V}_i : \mathcal C$ is Lipschitz,
    \begin{equation*}
        \bigl| \nabla^2 \tilde{\mathcal V}_i(x) : \mathcal C(x) - \nabla^2 \tilde{\mathcal V}_i(x') : \mathcal C(x') \bigr| \Le \bigl( \Sigma \norm{\nabla^3 \tilde{\mathcal V}_i}_\infty + 4 d M^2 \norm{\nabla^2 \tilde{\mathcal V}_i}_\infty \bigr) |x - x'| ,
    \end{equation*}
    so the Kantorovich--Rubinstein duality and \eqref{eq: w1 rough} give
    \begin{equation}
        \bigl| \pi_h(\nabla^2 \tilde{\mathcal V}_i : \mathcal C) - \pi(\nabla^2 \tilde{\mathcal V}_i : \mathcal C) \bigr| \Le C h \bigl( \norm{\nabla^2 \tilde{\mathcal V}_i}_\infty + \norm{\nabla^3 \tilde{\mathcal V}_i}_\infty \bigr) .
        \label{eq: sglm measure replacement}
    \end{equation}
    By \eqref{eq: sglm point replacement}, \eqref{eq: sglm measure replacement}, \eqref{eq: propagated derivatives} and \eqref{eq: time sums} at $j = 2, 3$, the replacements of $Y_0$ by $Z_0$ and of $\pi_h$ by $\pi$ in the coefficient of $\frac{h^2}{2}$ in \eqref{eq: sglm leading} introduce an additional error of at most
    \begin{equation*}
    \begin{aligned}
        & \frac{h^2}{2} \sum_{i=0}^{\infty} \Bigl( C \sqrt h\, \norm{\nabla^3 \tilde{\mathcal V}_i}_\infty + C h \bigl( \norm{\nabla^2 \tilde{\mathcal V}_i}_\infty + \norm{\nabla^3 \tilde{\mathcal V}_i}_\infty \bigr) \Bigr) \\
        & \quad \Le C h^{\frac32} \Bigl( 1 + \log\frac1h \Bigr) + C h^2 \Bigl( 1 + \log\frac1h \Bigr) \Le C h^{\frac32} \Bigl( 1 + \log\frac1h \Bigr) ,
    \end{aligned}
    \end{equation*}
    and \eqref{eq: sglm leading} gives
    \begin{equation}
        \biggl| \pi_h(f) - \pi(f) - \frac{h^2}{2} \sum_{i=0}^{\infty} \pi\bigl( \nabla^2 \tilde{\mathcal V}_i : \mathcal C \bigr) \biggr| \Le C h^{\frac32} \Bigl( 1 + \log\frac1h \Bigr) .
        \label{eq: sglm under pi}
    \end{equation}

    Since $\tilde{\mathcal V}_i \in \mathcal D$, \eqref{eq: riesz} with $\mathcal A = \mathcal C$ and $\omega = \tilde{\mathcal V}_i$ gives
    \begin{equation}
        \pi\bigl( \nabla^2 \tilde{\mathcal V}_i : \mathcal C \bigr) = \pi\bigl( \psi_{\mathcal C}\, (-\mathcal L \tilde{\mathcal V}_i) \bigr) .
        \label{eq: sglm riesz}
    \end{equation}
    The generators \eqref{eq: generator} and \eqref{eq: modified generator} differ only in the drift, so
    \begin{equation}
        -\mathcal L \tilde{\mathcal V}_i = -\tilde{\mathcal L}_h \tilde{\mathcal V}_i - \bigl( \nabla \tilde U_h - \nabla U \bigr) \cdot \nabla \tilde{\mathcal V}_i ,
        \label{eq: generator difference}
    \end{equation}
    and \eqref{eq: modified class} bounds the drift term of \eqref{eq: generator difference}:
    \begin{equation*}
        \bigl| \bigl( \nabla \tilde U_h(x) - \nabla U(x) \bigr) \cdot \nabla \tilde{\mathcal V}_i(x) \bigr| \Le C h \bigl( 1 + |x| \bigr) \norm{\nabla \tilde{\mathcal V}_i}_\infty , \qquad x \in \mathbb R^d .
    \end{equation*}
    Inserting \eqref{eq: generator difference} into \eqref{eq: sglm riesz} and using the bound on $\pi\bigl( |\psi_{\mathcal C}|\, (1 + |x|)^2 \bigr)$ in \eqref{eq: psi moments} gives
    \begin{equation*}
        \Bigl| \pi\bigl( \nabla^2 \tilde{\mathcal V}_i : \mathcal C \bigr) - \pi\bigl( \psi_{\mathcal C}\, (-\tilde{\mathcal L}_h \tilde{\mathcal V}_i) \bigr) \Bigr| \Le C h\, \pi\bigl( |\psi_{\mathcal C}|\, (1 + |x|) \bigr) \norm{\nabla \tilde{\mathcal V}_i}_\infty \Le C h \norm{\nabla \tilde{\mathcal V}_i}_\infty .
    \end{equation*}
    By \eqref{eq: propagated derivatives} and \eqref{eq: time sums} at $j = 1$, $h \sum_{i=0}^{\infty} \norm{\nabla \tilde{\mathcal V}_i}_\infty \Le C$, so
    \begin{equation*}
        \frac{h^2}{2} \sum_{i=0}^{\infty} \Bigl| \pi\bigl( \nabla^2 \tilde{\mathcal V}_i : \mathcal C \bigr) - \pi\bigl( \psi_{\mathcal C}\, (-\tilde{\mathcal L}_h \tilde{\mathcal V}_i) \bigr) \Bigr| \Le C h^3 \sum_{i=0}^{\infty} \norm{\nabla \tilde{\mathcal V}_i}_\infty \Le C h^2 ,
    \end{equation*}
    and \eqref{eq: sglm under pi} gives
    \begin{equation}
        \biggl| \pi_h(f) - \pi(f) - \frac{h^2}{2} \sum_{i=0}^{\infty} \pi\bigl( \psi_{\mathcal C}\, (-\tilde{\mathcal L}_h \tilde{\mathcal V}_i) \bigr) \biggr| \Le C h^{\frac32} \Bigl( 1 + \log\frac1h \Bigr) .
        \label{eq: sglm modified generator}
    \end{equation}

    It remains to sum $h (-\tilde{\mathcal L}_h \tilde{\mathcal V}_i)$ over $i$ in \eqref{eq: sglm modified generator}. For this sum, $Y_0 \in \mathbb R^d$ is a deterministic initial point of $\tilde Y_0(s)$, in place of $Y_0 \sim \mu_h$. By \eqref{eq: modified dynkin}, with $\tilde{\mathcal V}_i$ and $\tilde{\mathcal L}_h \tilde{\mathcal V}_i$ in place of $\tilde\phi_h$ and $\tilde g_h$, which applies since $\nabla \tilde{\mathcal V}_i$ is bounded, and $e^{h\tilde{\mathcal L}_h} \tilde{\mathcal V}_i = \tilde{\mathcal V}_{i+1}$,
    \begin{equation*}
        \tilde{\mathcal V}_i(Y_0) - \tilde{\mathcal V}_{i+1}(Y_0) = \int_0^h \mathbb E\bigl[ (-\tilde{\mathcal L}_h \tilde{\mathcal V}_i)(\tilde Y_0(s)) \bigr] \D s ,
    \end{equation*}
    hence
    \begin{equation}
        h\, (-\tilde{\mathcal L}_h \tilde{\mathcal V}_i)(Y_0) - \bigl( \tilde{\mathcal V}_i(Y_0) - \tilde{\mathcal V}_{i+1}(Y_0) \bigr) = \int_0^h \mathbb E\bigl[ (\tilde{\mathcal L}_h \tilde{\mathcal V}_i)(\tilde Y_0(s)) - (\tilde{\mathcal L}_h \tilde{\mathcal V}_i)(Y_0) \bigr] \D s .
        \label{eq: dynkin error}
    \end{equation}
    Differentiating \eqref{eq: modified generator} gives
    \begin{equation*}
        \nabla (\tilde{\mathcal L}_h \tilde{\mathcal V}_i) = \nabla \Delta \tilde{\mathcal V}_i - \nabla^2 \tilde U_h\, \nabla \tilde{\mathcal V}_i - \nabla^2 \tilde{\mathcal V}_i\, \nabla \tilde U_h ,
    \end{equation*}
    and \eqref{eq: modified class bounds} gives
    \begin{equation*}
        \norm{\nabla^2 \tilde U_h(z)} \Le 2M , \qquad |\nabla \tilde U_h(z)| \Le 2M (1 + |z|) , \qquad z \in \mathbb R^d .
    \end{equation*}
    Hence, for $z \in \mathbb R^d$ and $0 \Le s \Le h$,
    \begin{equation*}
    \begin{gathered}
        \bigl| \nabla (\tilde{\mathcal L}_h \tilde{\mathcal V}_i)(z) \bigr| \Le C (1 + |z|) \sum_{j=1}^{3} \norm{\nabla^j \tilde{\mathcal V}_i}_\infty , \\[-8pt]
        \bigl| (\tilde{\mathcal L}_h \tilde{\mathcal V}_i)(\tilde Y_0(s)) - (\tilde{\mathcal L}_h \tilde{\mathcal V}_i)(Y_0) \bigr| \Le C \bigl( 1 + |Y_0| + |\tilde Y_0(s) - Y_0| \bigr) |\tilde Y_0(s) - Y_0| \sum_{j=1}^{3} \norm{\nabla^j \tilde{\mathcal V}_i}_\infty .
    \end{gathered}
    \end{equation*}
    The bound on $\tilde Y_0(s) - Y_0$ in \eqref{eq: exact sizes} with $p = 2$ gives
    \begin{equation*}
        \mathbb E\bigl[ (1 + |Y_0| + |\tilde Y_0(s) - Y_0|)\, |\tilde Y_0(s) - Y_0| \bigr] \Le C \sqrt h\, (1+|Y_0|)^2 , \qquad 0 \Le s \Le h ,
    \end{equation*}
    so integrating over $s \in [0,h]$ in \eqref{eq: dynkin error} gives, for every $Y_0 \in \mathbb R^d$,
    \begin{equation*}
        \Bigl| h\, (-\tilde{\mathcal L}_h \tilde{\mathcal V}_i)(Y_0) - \bigl( \tilde{\mathcal V}_i(Y_0) - \tilde{\mathcal V}_{i+1}(Y_0) \bigr) \Bigr| \Le C h^{\frac32} \bigl( 1 + |Y_0| \bigr)^2 \sum_{j=1}^{3} \norm{\nabla^j \tilde{\mathcal V}_i}_\infty .
    \end{equation*}
    Integrating this bound in $Y_0$ against $|\psi_{\mathcal C}|\, \pi$, with \eqref{eq: psi moments}, \eqref{eq: propagated derivatives} and \eqref{eq: time sums} at $j = 1, 2, 3$,
    \begin{equation*}
        \frac h2 \sum_{i=0}^{\infty} \Bigl| \pi\Bigl( \psi_{\mathcal C}\, \bigl( h\, (-\tilde{\mathcal L}_h \tilde{\mathcal V}_i) - \tilde{\mathcal V}_i + \tilde{\mathcal V}_{i+1} \bigr) \Bigr) \Bigr| \Le C h^{\frac32} \Bigl( 1 + \log\frac1h \Bigr) .
    \end{equation*}
    The partial sums of the telescoping series are
    \begin{equation*}
        \sum_{i=0}^{N-1} \pi\bigl( \psi_{\mathcal C}\, (\tilde{\mathcal V}_i - \tilde{\mathcal V}_{i+1}) \bigr) = \pi\bigl( \psi_{\mathcal C} \tilde f_h \bigr) - \pi\bigl( \psi_{\mathcal C} \tilde{\mathcal V}_N \bigr) ,
    \end{equation*}
    and $\pi(\psi_{\mathcal C}) = 0$, \eqref{eq: sglm relaxation} and \eqref{eq: psi moments} give $\pi(\psi_{\mathcal C} \tilde{\mathcal V}_N) \to 0$ as $N \to \infty$. The series of $\pi(\psi_{\mathcal C}\, (-\tilde{\mathcal L}_h \tilde{\mathcal V}_i))$ converges absolutely: by $|\tilde{\mathcal L}_h \tilde{\mathcal V}_i(x)| \Le C (1+|x|) (\norm{\nabla \tilde{\mathcal V}_i}_\infty + \norm{\nabla^2 \tilde{\mathcal V}_i}_\infty)$, \eqref{eq: psi moments} and \eqref{eq: propagated derivatives},
    \begin{equation*}
        \sum_{i=0}^{\infty} \Bigl| \pi\bigl( \psi_{\mathcal C}\, (-\tilde{\mathcal L}_h \tilde{\mathcal V}_i) \bigr) \Bigr| \Le C \sum_{i=0}^{\infty} \bigl( \norm{\nabla \tilde{\mathcal V}_i}_\infty + \norm{\nabla^2 \tilde{\mathcal V}_i}_\infty \bigr) \Le C h^{-\frac12} \sum_{i=0}^{\infty} e^{-\lambda i h} < \infty .
    \end{equation*}
    Therefore
    \begin{equation}
        \biggl| \frac{h^2}{2} \sum_{i=0}^{\infty} \pi\bigl( \psi_{\mathcal C}\, (-\tilde{\mathcal L}_h \tilde{\mathcal V}_i) \bigr) - \frac h2\, \pi\bigl( \psi_{\mathcal C} \tilde f_h \bigr) \biggr| \Le C h^{\frac32} \Bigl( 1 + \log\frac1h \Bigr) .
        \label{eq: sglm telescope psi}
    \end{equation}
    Finally, \eqref{eq: heat shift} and $\norm{\nabla f}_\infty \Le 1$ give $|\tilde f_h - f| \Le \sqrt{h/2}\; \mathbb E|\chi| \Le \sqrt{d h / 2}$, so $\frac h2\, |\pi(\psi_{\mathcal C}\, (\tilde f_h - f))| \Le C h^{\frac32}$ by \eqref{eq: psi moments}. With \eqref{eq: sglm modified generator} and \eqref{eq: sglm telescope psi},
    \begin{equation*}
        \Bigl| \pi_h(f) - \pi(f) - \frac h2\, \pi( f \psi_{\mathcal C} ) \Bigr| \Le C h^{\frac32} \Bigl( 1 + \log\frac1h \Bigr) .
    \end{equation*}
    Since $\pi_h(f) - \pi(f) - \frac h2\, \pi(f \psi_{\mathcal C}) = \int f\, \D(\pi_h - \pi - \frac h2\, \psi_{\mathcal C}\, \pi)$, the supremum over the functions $f$ satisfying Assumption~\ref{ass: f} and $\norm{\nabla f}_\infty \Le 1$, with \eqref{eq: kr test functions}, proves \eqref{eq: kr sglm}.
\end{proof}

\begin{proof}[Proof of \eqref{eq: kr sgld} in Theorem~\ref{thm: kr expansion}]
    The proof follows that of \eqref{eq: kr sglm}, with the local error \eqref{eq: sgld weak} in place of \eqref{eq: sg weak}, Lemma~\ref{lem: sgld lipschitz step} for the term $i = 0$, and Lemma~\ref{lem: local lipschitz} for the move from $\bar\pi_h$ to $\pi$. Let $f$ satisfy Assumption~\ref{ass: f} and $\norm{\nabla f}_\infty \Le 1$, and let
    \begin{equation}
        \mathcal V_i := e^{ih\mathcal L} f , \qquad i = 0, 1, 2, \dots ,
        \label{eq: sgld test functions}
    \end{equation}
    so that $\mathcal V_0 = f$ and $\mathcal V_i \in \mathcal D$ by \eqref{eq: short time}. By \eqref{eq: smoothing} at $t = ih$, as for \eqref{eq: propagated derivatives}, for $i \Ge 1$ and $j = 1, 2, 3, 4$,
    \begin{equation}
        \norm{\nabla^j \mathcal V_i}_\infty \Le C \min\bigl( 1, (i+1)h \bigr)^{-\frac{j-1}{2}} e^{-\lambda i h} .
        \label{eq: sgld propagated}
    \end{equation}
    With the transition kernel $\mathcal P_h$ of \eqref{eq: sgld}, as in Appendix~\ref{app: sgld}, $(\mathcal P_h g)(Z_0) = \hat{\mathbb E}_0\bigl[ g(Z_1) \bigr]$ for a given $Z_0$. Since $\bar\pi_h$ is invariant under $\mathcal P_h$ and $e^{h\mathcal L} \mathcal V_i = \mathcal V_{i+1}$,
    \begin{equation*}
        \bar\pi_h(\mathcal V_i) - \bar\pi_h(\mathcal V_{i+1}) = \bar\pi_h\bigl( \mathcal P_h \mathcal V_i - e^{h\mathcal L} \mathcal V_i \bigr) ,
    \end{equation*}
    and $|\mathcal V_N(x) - \pi(f)| \Le C e^{-\lambda N h} (1 + |x|)$ by \eqref{eq: smoothing} at $j = 1$ and the invariance of $\pi$ under $e^{t\mathcal L}$, as in \eqref{eq: sglm relaxation}. Summing over $i$,
    \begin{equation}
        \bar\pi_h(f) - \pi(f) = \sum_{i=0}^{\infty} \bar\pi_h\bigl( \mathcal P_h \mathcal V_i - e^{h\mathcal L} \mathcal V_i \bigr) .
        \label{eq: sgld telescope}
    \end{equation}
    For the term $i = 0$ of \eqref{eq: sgld telescope}, averaging \eqref{eq: sgld lipschitz step} over $Z_0 \sim \bar\pi_h$, with $\bar\pi_h(|x|) \Le C$ from Lemma~\ref{lem: sgld invariant}, gives
    \begin{equation*}
        \bigl| \bar\pi_h\bigl( \mathcal P_h f - e^{h\mathcal L} f \bigr) \bigr| \Le C h^{\frac32} \bigl( 1 + \bar\pi_h(|x|) \bigr) \Le C h^{\frac32} .
    \end{equation*}
    For $i \Ge 1$, averaging \eqref{eq: sgld weak} at $k = 0$, with $G = \mathcal V_i$, over $Z_0 \sim \bar\pi_h$, with $\bar\pi_h(|x|^{12}) \Le C$, gives
    \begin{equation*}
        \biggl| \bar\pi_h\bigl( \mathcal P_h \mathcal V_i - e^{h\mathcal L} \mathcal V_i \bigr) - \frac{h^2}{2}\, \bar\pi_h\Bigl( \nabla^2 \mathcal V_i : \bigl( 2 \nabla^2 U + \mathcal C \bigr) - \nabla \mathcal V_i \cdot \bigl( \nabla^2 U\, \nabla U - \nabla \Delta U \bigr) \Bigr) \biggr| \Le C h^3 \sum_{j=1}^{4} \norm{\nabla^j \mathcal V_i}_\infty ,
    \end{equation*}
    and, by \eqref{eq: sgld propagated} and \eqref{eq: time sums},
    \begin{equation*}
        C h^3 \sum_{i=1}^{\infty} \sum_{j=1}^{4} \norm{\nabla^j \mathcal V_i}_\infty \Le C h^2 \Bigl( 1 + \log\frac1h + h^{-\frac12} \Bigr) \Le C h^{\frac32} .
    \end{equation*}
    Hence
    \begin{equation}
        \biggl| \bar\pi_h(f) - \pi(f) - \frac{h^2}{2} \sum_{i=1}^{\infty} \bar\pi_h\Bigl( \nabla^2 \mathcal V_i : \bigl( 2 \nabla^2 U + \mathcal C \bigr) - \nabla \mathcal V_i \cdot \bigl( \nabla^2 U\, \nabla U - \nabla \Delta U \bigr) \Bigr) \biggr| \Le C h^{\frac32} .
        \label{eq: sgld leading}
    \end{equation}
    By Assumption~\ref{ass: U} and \eqref{eq: trace bound}, for every $x \in \mathbb R^d$ and every symmetric matrix $A$,
    \begin{equation*}
    \begin{gathered}
        \bigl| \nabla^2 U(x)\, \nabla U(x) - \nabla \Delta U(x) \bigr| \Le M^2 (1 + |x|) + d M \Le C (1 + |x|) , \\
        \bigl| A : \bigl( 2 \nabla^2 U(x) + \mathcal C(x) \bigr) \bigr| \Le (2 d M + \Sigma) \norm{A} .
    \end{gathered}
    \end{equation*}
    Hence the integrand in \eqref{eq: sgld leading} satisfies
    \begin{equation*}
        \Bigl| \nabla^2 \mathcal V_i : \bigl( 2 \nabla^2 U + \mathcal C \bigr) - \nabla \mathcal V_i \cdot \bigl( \nabla^2 U\, \nabla U - \nabla \Delta U \bigr) \Bigr| \Le C \bigl( 1 + |x| \bigr) \bigl( \norm{\nabla \mathcal V_i}_\infty + \norm{\nabla^2 \mathcal V_i}_\infty \bigr) ,
    \end{equation*}
    and the moment bound $\bar\pi_h(|x|) \Le C$ of Lemma~\ref{lem: sgld invariant}, \eqref{eq: sgld propagated} and \eqref{eq: time sums} at $j = 1, 2$ give
    \begin{equation*}
        \frac{h^2}{2} \sum_{i=1}^{\infty} \biggl| \bar\pi_h\Bigl( \nabla^2 \mathcal V_i : \bigl( 2 \nabla^2 U + \mathcal C \bigr) - \nabla \mathcal V_i \cdot \bigl( \nabla^2 U\, \nabla U - \nabla \Delta U \bigr) \Bigr) \biggr| \Le C h^2 \sum_{i=1}^{\infty} \bigl( \norm{\nabla \mathcal V_i}_\infty + \norm{\nabla^2 \mathcal V_i}_\infty \bigr) \Le C h .
    \end{equation*}
    With \eqref{eq: sgld leading}, $|\bar\pi_h(f) - \pi(f)| \Le C (h + h^{\frac32}) \Le C h$ for every $f$ satisfying Assumption~\ref{ass: f} and $\norm{\nabla f}_\infty \Le 1$, and \eqref{eq: kr test functions} gives
    \begin{equation}
        \mathcal W_1(\bar\pi_h, \pi) \Le C h .
        \label{eq: w1 rough sgld}
    \end{equation}

    In the coefficient of $\frac{h^2}{2}$ in \eqref{eq: sgld leading}, the measure $\bar\pi_h$ can be replaced by $\pi$. By Assumption~\ref{ass: U}, the derivatives of $U$ of orders $2$ to $4$ are bounded, so $\nabla^2 U$ and $\nabla \Delta U$ are bounded and Lipschitz, and
    \begin{equation}
        \bigl| \bigl( \nabla^2 U\, \nabla U - \nabla \Delta U \bigr)(x) - \bigl( \nabla^2 U\, \nabla U - \nabla \Delta U \bigr)(x') \bigr| \Le C \bigl( 1 + |x| + |x'| \bigr) |x - x'| , \qquad x, x' \in \mathbb R^d .
        \label{eq: drift correction lipschitz}
    \end{equation}
    With \eqref{eq: drift correction lipschitz}, \eqref{eq: trace bound}, \eqref{eq: covariance lipschitz} and $|A : B| \Le d \norm{A} \norm{B}$, for each $i$, the integrand in the $\bar\pi_h$ bracket of \eqref{eq: sgld leading} satisfies \eqref{eq: local lipschitz} with $K = C \sum_{j=1}^{3} \norm{\nabla^j \mathcal V_i}_\infty$. Lemma~\ref{lem: local lipschitz}, with $\mu = \bar\pi_h$, $\nu = \pi$, the moments of Lemmas~\ref{lem: sgld invariant} and \ref{lem: dissipative}, and $\mathcal W_1(\bar\pi_h, \pi) \Le C h \Le 1$ by \eqref{eq: w1 rough sgld} with $C h_0 \Le 1$, gives, for $i \Ge 1$,
    \begin{equation}
    \begin{aligned}
        & \Bigl| (\bar\pi_h - \pi)\Bigl( \nabla^2 \mathcal V_i : \bigl( 2 \nabla^2 U + \mathcal C \bigr) - \nabla \mathcal V_i \cdot \bigl( \nabla^2 U\, \nabla U - \nabla \Delta U \bigr) \Bigr) \Bigr| \\
        & \quad \Le C\, \mathcal W_1(\bar\pi_h, \pi)^{\frac{10}{11}} \sum_{j=1}^{3} \norm{\nabla^j \mathcal V_i}_\infty \Le C h^{\frac{10}{11}} \sum_{j=1}^{3} \norm{\nabla^j \mathcal V_i}_\infty .
    \end{aligned}
    \label{eq: sgld measure replacement}
    \end{equation}
    By \eqref{eq: sgld propagated} and \eqref{eq: time sums} at $j = 1, 2, 3$, $h \sum_{i=1}^{\infty} \sum_{j=1}^{3} \norm{\nabla^j \mathcal V_i}_\infty \Le C (1 + \log\frac1h)$, so, by \eqref{eq: sgld measure replacement}, replacing $\bar\pi_h$ by $\pi$ in the coefficient of $\frac{h^2}{2}$ in \eqref{eq: sgld leading} introduces an additional error of at most
    \begin{equation*}
        \frac{h^2}{2} \sum_{i=1}^{\infty} C h^{\frac{10}{11}} \sum_{j=1}^{3} \norm{\nabla^j \mathcal V_i}_\infty \Le C h^{1 + \frac{10}{11}} \Bigl( 1 + \log\frac1h \Bigr) \Le C h^{\frac32} .
    \end{equation*}
    Under $\pi$, \eqref{eq: sgld integration by parts} with $\mathcal V_i$ in place of $\psi$ and \eqref{eq: riesz} give
    \begin{equation*}
        \pi\Bigl( \nabla^2 \mathcal V_i : \bigl( 2 \nabla^2 U + \mathcal C \bigr) - \nabla \mathcal V_i \cdot \bigl( \nabla^2 U\, \nabla U - \nabla \Delta U \bigr) \Bigr) = \pi\bigl( \nabla^2 \mathcal V_i : ( \nabla^2 U + \mathcal C ) \bigr) = \pi\bigl( ( \psi_{\nabla^2 U} + \psi_{\mathcal C} )\, (-\mathcal L \mathcal V_i) \bigr) ,
    \end{equation*}
    so we obtain
    \begin{equation}
        \biggl| \bar\pi_h(f) - \pi(f) - \frac{h^2}{2} \sum_{i=1}^{\infty} \pi\bigl( ( \psi_{\nabla^2 U} + \psi_{\mathcal C} )\, (-\mathcal L \mathcal V_i) \bigr) \biggr| \Le C h^{\frac32} .
        \label{eq: sgld under pi}
    \end{equation}

    The sum in \eqref{eq: sgld under pi} telescopes as in the proof of \eqref{eq: kr sglm}, with $\mathcal L$, $\mathcal V_i$ of \eqref{eq: sgld test functions}, the solution $X_s$ of \eqref{eq: langevin} and its initial value $X_0$ in place of $\tilde{\mathcal L}_h$, $\tilde{\mathcal V}_i$ of \eqref{eq: propagated test functions}, $\tilde Y_0(s)$ and $Y_0$, with $\psi_{\nabla^2 U} + \psi_{\mathcal C}$ in place of $\psi_{\mathcal C}$, and with \eqref{eq: sgld propagated} and Assumption~\ref{ass: U} in place of \eqref{eq: propagated derivatives} and \eqref{eq: modified class bounds}. The identity \eqref{eq: dynkin error} holds in this form, and so does the bound on its left side by $C h^{\frac32} (1 + |X_0|)^2 \sum_{j=1}^{3} \norm{\nabla^j \mathcal V_i}_\infty$. Since the sum starts at $i = 1$, the partial sums of the telescoping series are $\pi((\psi_{\nabla^2 U} + \psi_{\mathcal C})\, \mathcal V_1) - \pi((\psi_{\nabla^2 U} + \psi_{\mathcal C})\, \mathcal V_N)$. By $\pi(\psi_{\nabla^2 U} + \psi_{\mathcal C}) = 0$, the bound $|\mathcal V_N(x) - \pi(f)| \Le C e^{-\lambda N h} (1 + |x|)$ and \eqref{eq: psi moments}, $\pi((\psi_{\nabla^2 U} + \psi_{\mathcal C})\, \mathcal V_N) \to 0$ as $N \to \infty$. Hence
    \begin{equation*}
        \biggl| \frac{h^2}{2} \sum_{i=1}^{\infty} \pi\bigl( ( \psi_{\nabla^2 U} + \psi_{\mathcal C} )\, (-\mathcal L \mathcal V_i) \bigr) - \frac h2\, \pi\bigl( ( \psi_{\nabla^2 U} + \psi_{\mathcal C} )\, \mathcal V_1 \bigr) \biggr| \Le C h^{\frac32} \Bigl( 1 + \log\frac1h \Bigr) .
    \end{equation*}
    Finally, the bound on $X_s - X_0$ in \eqref{eq: exact sizes} for \eqref{eq: langevin} gives $|\mathcal V_1(x) - f(x)| = |\mathbb E^x[f(X_h) - f(X_0)]| \Le \mathbb E^x |X_h - X_0| \Le C \sqrt h\, (1 + |x|)$, so $\frac h2\, |\pi((\psi_{\nabla^2 U} + \psi_{\mathcal C})\, (\mathcal V_1 - f))| \Le C h^{\frac32}$ by \eqref{eq: psi moments}. With \eqref{eq: sgld under pi},
    \begin{equation*}
        \Bigl| \bar\pi_h(f) - \pi(f) - \frac h2\, \pi\bigl( f\, ( \psi_{\nabla^2 U} + \psi_{\mathcal C} ) \bigr) \Bigr| \Le C h^{\frac32} \Bigl( 1 + \log\frac1h \Bigr) ,
    \end{equation*}
    and the supremum over the functions $f$ satisfying Assumption~\ref{ass: f} and $\norm{\nabla f}_\infty \Le 1$, with \eqref{eq: kr test functions}, proves \eqref{eq: kr sgld}.
\end{proof}

\subsection{Proof of Corollary~\ref{cor: w1 limits}: limits of the Wasserstein-1 errors}
\label{sec: w1 limits proof}

\begin{proof}[Proof of Corollary~\ref{cor: w1 limits}]
    Both limits of Corollary~\ref{cor: w1 limits} follow from \eqref{eq: w1 limits} as in Section~\ref{subsec: w1}, since $\pi_h$, $\bar\pi_h$ and $\pi$ have finite first moments by Lemmas~\ref{lem: invariant}, \ref{lem: sgld invariant} and \ref{lem: dissipative}.

    For the positivity, a function $\psi \in L^2_0(\pi)$ with $\norm{\psi\, \pi}_{\mathrm{KR}} = 0$ has $\pi(g \psi) = 0$ for every $g \in C_c^\infty(\mathbb R^d)$, hence $\psi = 0$ $\pi$-almost everywhere. By \eqref{eq: psi hessian}, $\psi_{\nabla^2 U}$ is continuous and unbounded, since $|\Delta U| \Le d M$ while $|\nabla U(x)| \Ge \langle x, \nabla U(x) \rangle / |x| \to \infty$ as $|x| \to \infty$ by \eqref{eq: dissipative bounds}. Hence $\psi_{\nabla^2 U}$ is nonzero on a nonempty open set, which has positive $\pi$-measure since $\pi(x) > 0$ for every $x \in \mathbb R^d$, and $\norm{\psi_{\nabla^2 U}\, \pi}_{\mathrm{KR}} > 0$. For $\psi_{\mathcal C}$, let $r > 0$ and
    \begin{equation}
        \varphi_r(x) := \bigl( r^2 + |x|^2 \bigr)^{\frac12} , \qquad x \in \mathbb R^d ,
        \label{eq: phi r}
    \end{equation}
    so that $\varphi_r \in \mathcal D$ and
    \begin{equation}
    \begin{gathered}
        |\nabla \varphi_r| \Le 1 , \qquad \nabla^2 \varphi_r = \frac{\varphi_r^2 I_d - x \otimes x}{\varphi_r^3} \succeq \frac{r^2 I_d}{\varphi_r^3} , \qquad \norm{\nabla^2 \varphi_r} \Le \frac{1}{\varphi_r} , \\
        \Delta \varphi_r = \frac{d-1}{\varphi_r} + \frac{r^2}{\varphi_r^3} , \qquad |\nabla \Delta \varphi_r| = \frac{(d-1)\, |x|}{\varphi_r^3} + \frac{3 r^2 |x|}{\varphi_r^5} \Le \frac{d+1}{2 r^2} .
    \end{gathered}
        \label{eq: phi r bounds}
    \end{equation}
    Since $|\nabla U(x)| \Le M(1 + |x|) \Le \sqrt2\, M \varphi_1(x)$, the gradient of $\mathcal L \varphi_1 = \Delta \varphi_1 - \nabla U \cdot \nabla \varphi_1$ satisfies, by \eqref{eq: phi r bounds}, $|\nabla (\mathcal L \varphi_1)| \Le \frac{d+1}{2} + M + \sqrt2\, M \Le d + 3M$. With \eqref{eq: riesz} for $\omega = \varphi_1$ and $\mathcal C(x)$ positive semidefinite,
    \begin{equation*}
        \norm{\psi_{\mathcal C}\, \pi}_{\mathrm{KR}} \Ge \frac{\pi\bigl( \psi_{\mathcal C}\, (-\mathcal L \varphi_1) \bigr)}{d + 3M} = \frac{\pi\bigl( \nabla^2 \varphi_1 : \mathcal C \bigr)}{d + 3M} \Ge \frac{1}{d + 3M}\, \pi\biggl( \frac{\mathrm{Tr}\bigl[ \mathcal C(x) \bigr]}{(1 + |x|^2)^{\frac32}} \biggr) .
    \end{equation*}
    The trace $\mathrm{Tr}\bigl[ \mathcal C(x) \bigr]$ is nonnegative and continuous by \eqref{eq: covariance lipschitz}; if $\mathcal C(x)$ does not vanish identically, the trace is positive on a nonempty open set, and $\pi\Bigl( \mathrm{Tr}\bigl[ \mathcal C(x) \bigr]\, (1 + |x|^2)^{-\frac32} \Bigr) > 0$, so $\norm{\psi_{\mathcal C}\, \pi}_{\mathrm{KR}} > 0$.
\end{proof}

\subsection{Proof of Lemma~\ref{lem: w1 counterexample}: failure of the comparison}
\label{sec: w1 counterexample proof}

\begin{proof}[Proof of Lemma~\ref{lem: w1 counterexample}]
    The proof checks the assumptions, writes the Kantorovich--Rubinstein norm in dimension $1$ as the integral \eqref{eq: kr one dimension}, computes $\psi_{\mathcal C}$ explicitly, and evaluates both norms of Lemma~\ref{lem: w1 counterexample} by quadrature. Assumption~\ref{ass: U} holds with $m = 1$, $R = 0$ and every $M \Ge 1$. The estimator is unbiased, $|\nabla \hat U(x;\xi)| \Le |x| + 10$, its Hessian $1 + 10\, \xi\, a'(x)$ is bounded by $1 + 10 \norm{a'}_\infty$, and its variance is $100\, a(x)^2 \Le 100$, so Assumption~\ref{ass: sg} holds for $M$ large enough.

    Here $\pi$ is the standard Gaussian distribution, with $\pi'(x) = -x\, \pi(x)$. For $\psi \in L^2_0(\pi)$ and $s \in \mathbb R$, the Cauchy--Schwarz inequality on $(-\infty, s]$ for $s \Le 0$, and on $[s, \infty)$ for $s > 0$ since $\pi(\psi) = 0$, gives
    \begin{equation*}
        \biggl| \int_{-\infty}^s \psi(y)\, \pi(y)\, \D y \biggr| \Le \pi(\psi^2)^{\frac12} e^{-\frac14 s^2} ,
    \end{equation*}
    so for $f$ Lipschitz with $|f'| \Le 1$ an integration by parts gives
    \begin{equation*}
        \pi(f \psi) = -\int_{\mathbb R} f'(s) \int_{-\infty}^s \psi(y)\, \pi(y)\, \D y \,\D s \Le \int_{\mathbb R} \biggl| \int_{-\infty}^s \psi(y)\, \pi(y)\, \D y \biggr| \,\D s ,
    \end{equation*}
    with equality for $f(s) = -\int_0^s \mathrm{sign}\bigl( \int_{-\infty}^{s'} \psi(y)\, \pi(y)\, \D y \bigr) \D s'$. Hence
    \begin{equation}
        \norm{\psi\, \pi}_{\mathrm{KR}} = \int_{\mathbb R} \biggl| \int_{-\infty}^s \psi(y)\, \pi(y)\, \D y \biggr| \,\D s .
        \label{eq: kr one dimension}
    \end{equation}

    By \eqref{eq: psi hessian}, $\psi_{\nabla^2 U}(x) = \frac12 (x^2 - 1)$, and
    \begin{equation}
        \int_{-\infty}^s \psi_{\nabla^2 U}(y)\, \pi(y)\, \D y = -\frac12\, s\, \pi(s) , \qquad s \in \mathbb R .
        \label{eq: G hessian}
    \end{equation}
    For $\mathcal A = a^2$, let $J(x) := \int_{-\infty}^x y\, a(y)^2\, \D y$ and $\zeta := \pi(a^2) - \pi(J)$. The function $J - a^2 + \zeta$ is smooth and bounded, satisfies $\pi(J - a^2 + \zeta) = 0$, and has the derivative $x\, a^2 - (a^2)'$, so for $\omega \in \mathcal D$, \eqref{eq: generator ibp} gives
    \begin{equation*}
        \pi\bigl( (J - a^2 + \zeta) (-\mathcal L \omega) \bigr) = \pi\bigl( (x\, a^2 - (a^2)')\, \omega' \bigr) = -\int_{\mathbb R} \omega'(x)\, \bigl( a(x)^2 \pi(x) \bigr)'\, \D x = \pi\bigl( \omega''\, a^2 \bigr) ,
    \end{equation*}
    which is \eqref{eq: riesz} for $\mathcal A = a^2$. By the uniqueness in Theorem~\ref{thm: riesz},
    \begin{equation}
        \psi_{a^2} = J - a^2 + \zeta .
        \label{eq: psi shells}
    \end{equation}
    Since $\mathcal C = 100\, a^2$, the linearity in Theorem~\ref{thm: riesz} gives $\psi_{\mathcal C} = 100\, \psi_{a^2}$. Quadrature in \eqref{eq: kr one dimension}, with \eqref{eq: G hessian} and \eqref{eq: psi shells}, then gives the values in Lemma~\ref{lem: w1 counterexample}, so \eqref{eq: w1 comparison} fails.
\end{proof}

\subsection{Proof of Theorem~\ref{thm: batch threshold}: threshold on the batch size}
\label{sec: batch threshold proof}

\begin{proof}[Proof of Theorem~\ref{thm: batch threshold}]
    The proof bounds $\norm{( \psi_{\nabla^2 U} + \psi_{\mathcal C} )\, \pi}_{\mathrm{KR}}$ from below by testing \eqref{eq: riesz} with a function of the distance to the minimizer of $U$, bounds $\norm{\psi_{\mathcal C}\, \pi}_{\mathrm{KR}}$ from above through \eqref{eq: bochner}, and compares the $2$ bounds. Let $x_*$ be the minimizer of $U$, and let $\varphi_r$ be the function of \eqref{eq: phi r} with $r^2 = d/m$ and $x - x_*$ in place of $x$. The gradient of $\mathcal L \varphi_r = \Delta \varphi_r - \nabla U \cdot \nabla \varphi_r$ satisfies
    \begin{equation}
        |\nabla (\mathcal L \varphi_r)| \Le 2\, \norm{\nabla^2 U}_\infty + \frac{d+1}{2r^2} \Le 3\, \norm{\nabla^2 U}_\infty ,
        \label{eq: lipschitz phi R}
    \end{equation}
    since \eqref{eq: phi r bounds} bounds $|\nabla \Delta \varphi_r|$ by $\frac{d+1}{2r^2} = \frac{(d+1)\, m}{2d} \Le m \Le \norm{\nabla^2 U}_\infty$, and, by $|\nabla \varphi_r| \Le 1$, $\norm{\nabla^2 \varphi_r} \Le \varphi_r^{-1}$ and $|\nabla U(x)| = |\nabla U(x) - \nabla U(x_*)| \Le \norm{\nabla^2 U}_\infty\, \varphi_r(x)$,
    \begin{equation*}
        |\nabla(\nabla U \cdot \nabla \varphi_r)| \Le \norm{\nabla^2 U}\, |\nabla \varphi_r| + \norm{\nabla^2 \varphi_r}\, |\nabla U| \Le 2 \norm{\nabla^2 U}_\infty .
    \end{equation*}

    By \eqref{eq: lipschitz phi R}, the function $-\mathcal L \varphi_r / (3\, \norm{\nabla^2 U}_\infty)$ is $1$-Lipschitz, and \eqref{eq: riesz} with $\omega = \varphi_r$ gives $\pi((-\mathcal L \varphi_r)\, \psi_{\mathcal A}) = \pi(\nabla^2 \varphi_r : \mathcal A)$ for $\mathcal A = \nabla^2 U$ and $\mathcal A = \mathcal C$. As $\nabla^2 \varphi_r$ and $\mathcal C(x)$ are positive semidefinite, $\pi(\nabla^2 \varphi_r : \mathcal C) \Ge 0$, and as $\nabla^2 U \succeq m I_d$, $\nabla^2 \varphi_r : \nabla^2 U \Ge m\, \Delta \varphi_r$. Hence
    \begin{equation}
        \norm{( \psi_{\nabla^2 U} + \psi_{\mathcal C} )\, \pi}_{\mathrm{KR}} \Ge \frac{\pi\bigl( (-\mathcal L \varphi_r)\, ( \psi_{\nabla^2 U} + \psi_{\mathcal C} ) \bigr)}{3\, \norm{\nabla^2 U}_\infty} \Ge \frac{m\, \pi( \Delta \varphi_r )}{3\, \norm{\nabla^2 U}_\infty} .
        \label{eq: threshold lower}
    \end{equation}
    The identity \eqref{eq: generator ibp} with $g_1 = 1$ and $g_2 = \frac12 |x - x_*|^2$ gives $\pi(\nabla U \cdot (x - x_*)) = d$, and $\nabla U(x) \cdot (x - x_*) \Ge m |x - x_*|^2$, so $\pi(|x - x_*|^2) \Le d/m$. As functions of $|x - x_*|^2$, $(r^2 + |x - x_*|^2)^{-\frac12}$ and $(r^2 + |x - x_*|^2)^{-\frac32}$ are convex and decreasing, so Jensen's inequality and $r^2 = d/m$ give
    \begin{equation*}
        \pi( \Delta \varphi_r ) \Ge \frac{d-1}{(2d/m)^{\frac12}} + \frac{d/m}{(2d/m)^{\frac32}} = \Bigl( \frac{m}{2d} \Bigr)^{\frac12} \Bigl( d - \frac12 \Bigr) \Ge \Bigl( \frac{m d}{8} \Bigr)^{\frac12} ,
    \end{equation*}
    and \eqref{eq: threshold lower} gives
    \begin{equation}
        \norm{( \psi_{\nabla^2 U} + \psi_{\mathcal C} )\, \pi}_{\mathrm{KR}} \Ge \frac{m^{\frac32} \sqrt d}{6 \sqrt2\, \norm{\nabla^2 U}_\infty} \Ge \frac{m^{\frac32} \sqrt d}{9\, \norm{\nabla^2 U}_\infty} .
        \label{eq: threshold lower explicit}
    \end{equation}

    For $\norm{\psi_{\mathcal C}\, \pi}_{\mathrm{KR}}$ in \eqref{eq: w1 comparison}, let $f$ satisfy Assumption~\ref{ass: f} and $\norm{\nabla f}_\infty \Le 1$, and let $\phi_0$ solve \eqref{eq: continuous poisson}. By \eqref{eq: bochner} with $\omega = \phi_0 \in \mathcal D$ and $\nabla^2 U \succeq m I_d$, and \eqref{eq: generator ibp},
    \begin{equation}
        \pi\bigl( \nabla^2 \phi_0 : \nabla^2 \phi_0 \bigr) + m\, \pi\bigl( |\nabla \phi_0|^2 \bigr) \Le \pi\bigl( (f - \pi(f))^2 \bigr) = \pi\bigl( \nabla f \cdot \nabla \phi_0 \bigr) \Le \pi\bigl( |\nabla \phi_0|^2 \bigr)^{\frac12} .
        \label{eq: bochner phi0}
    \end{equation}
    By \eqref{eq: bochner phi0}, $\pi(|\nabla \phi_0|^2)^{\frac12} \Le m^{-1}$ and then $\pi(\nabla^2 \phi_0 : \nabla^2 \phi_0) \Le m^{-1}$. As for \eqref{eq: riesz functional bound}, \eqref{eq: riesz} with $\omega = \phi_0$ and the Cauchy--Schwarz inequality give $|\pi(f \psi_{\mathcal C})| = |\pi(\nabla^2 \phi_0 : \mathcal C)| \Le \sqrt d\, \norm{\mathcal C}_\infty\, m^{-\frac12}$, and the supremum over $f$ with \eqref{eq: kr test functions} gives
    \begin{equation}
        \norm{\psi_{\mathcal C}\, \pi}_{\mathrm{KR}} \Le \frac{\sqrt d\, \norm{\mathcal C}_\infty}{m^{\frac12}} .
        \label{eq: threshold upper}
    \end{equation}
    Hence, whenever $\norm{\mathcal C}_\infty \Le m^2 / (9\, \norm{\nabla^2 U}_\infty)$, \eqref{eq: threshold upper} and \eqref{eq: threshold lower explicit} give
    \begin{equation*}
        \norm{\psi_{\mathcal C}\, \pi}_{\mathrm{KR}} \Le \frac{\sqrt d\, \norm{\mathcal C}_\infty}{m^{\frac12}} \Le \frac{m^{\frac32} \sqrt d}{9\, \norm{\nabla^2 U}_\infty} \Le \norm{( \psi_{\nabla^2 U} + \psi_{\mathcal C} )\, \pi}_{\mathrm{KR}} ,
    \end{equation*}
    which is \eqref{eq: w1 comparison}. For SGLM and SGLD with a batch of $P$ independent stochastic gradients, the averaged estimator satisfies Assumption~\ref{ass: sg} with $\Sigma/P$ in place of $\Sigma$ (Section~\ref{subsec: mse analysis}), so Corollary~\ref{cor: w1 limits} applies with its covariance $\mathcal C(x)/P$, and $\norm{\mathcal C}_\infty / P \Le m^2 / (9\, \norm{\nabla^2 U}_\infty)$ exactly when $P$ satisfies \eqref{eq: batch threshold}. The representor of $\mathcal C(x)/P$ is $\psi_{\mathcal C}/P$ by the linearity in Theorem~\ref{thm: riesz}, and \eqref{eq: w1 comparison} with $\psi_{\mathcal C}/P$ in place of $\psi_{\mathcal C}$ is \eqref{eq: w1 comparison batch}.
\end{proof}

\end{document}